\documentclass{amsbook}

\usepackage{amssymb,amsmath, cite} 
\usepackage{latexsym}
\usepackage[usenames,dvipsnames]{color}
\usepackage{tikz}
\usepackage[all]{xy}

\usepackage{amsmath, amsfonts, amsthm, amssymb, cite}
\usepackage[all]{xy}

\newtheorem{thm}{Theorem}[chapter]
\newtheorem{lem}[thm]{Lemma}
\newtheorem{prop}[thm]{Proposition}
\newtheorem{cor}[thm]{Corollary}

\theoremstyle{definition}
\newtheorem{dfn}[thm]{Definition}
\newtheorem{ex}[thm]{Example}

\theoremstyle{remark}
\newtheorem{rmk}[thm]{Remark}

\numberwithin{section}{chapter}
\numberwithin{equation}{chapter}

\def\cyclic{\mathop{\kern0.9ex{{+}
\kern-2.2ex\raise-.28ex\hbox{\Large\hbox
{$\circlearrowright$}}}}\limits}

\def\buildrel#1_#2^#3{\mathrel{\mathop{\kern 0pt#1}\limits_{#2}^{#3}}}

\newcommand{\Aut}{\mbox{$\mathtt{Aut}$}}

\renewcommand{\mid}{\mbox{$\mathtt{mid}$}}
\newcommand{\End}{\mbox{$\mathtt{End}$}}

\newcommand{\Span}{\mbox{$\mathtt{span}$}}

\newcommand{\Id}{\mbox{$\mathtt{Id}$}}
\newcommand{\id}{\mbox{$\mathtt{id}$}}

\newcommand{\Ad}{\mbox{$\mathtt{Ad}$}}
\newcommand{\ad}{\mbox{$\mathtt{ad}$}}

\newcommand{\C}{\mathbb C} 
 
\newcommand{\Y}{\mathbb Y} 
\renewcommand{\L}{\mathbb L} 
\newcommand{\bC}{{\bf C}} 
\newcommand{\bB}{{\bf B}}

\newcommand{\bS}{{\bf S}}

\newcommand{\bE}{{\bf E}}
\newcommand{\bD}{{\bf D}}

\newcommand{\bR}{{\bf R}}
\newcommand{\bK}{{\bf K}}

\newcommand{\bm}{{\bf m}}

\newcommand{\A}{\mathbb A} 

\newcommand{\B}{\mathbb B} 
\renewcommand{\S}{\mathbb S}

\newcommand{\R}{\mathbb R}

\newcommand{\N}{\mathbb N}

\newcommand{\sech}{\mbox{\rm sech}} 
 
 \DeclareMathOperator{\arcsinh}{arcsinh}
 
\newcommand{\Tr}{\mbox{\rm Tr}} 

\newcommand{\g}{{\mathfrak{g}}{}} 
 
\renewcommand{\u}{{\mathfrak{u}}{}} 
\renewcommand{\l}{{\mathfrak{l}}{}}
\renewcommand{\k}{{\mathfrak{k}}{}} 
\newcommand{\p}{{\mathfrak{p}}{}} 
\newcommand{\fB}{{\mathfrak{B}}{}} 
\newcommand{\fD}{{\mathfrak{D}}{}} 
\newcommand{\fL}{{\mathfrak{L}}{}} 
 
\newcommand{\fV}{{\mathfrak{V}}{}} 
\newcommand{\fY}{{\mathfrak{Y}}{}}

\newcommand{\fsp}{{\mathfrak{sp}}{}} 
\newcommand{\q}{{\mathfrak{q}}{}} 
\newcommand{\n}{{\mathfrak{n}}{}} 
\newcommand{\s}{{\mathfrak{s}}{}}

\newcommand{\fd}{{\mathfrak{d}}{}} 
\newcommand{\h}{{\mathfrak{h}}{}} 
\renewcommand{\a}{{\mathfrak{a}}{}} 
\renewcommand{\b}{{\mathfrak{b}}{}} 
\newcommand{\z}{{\mathfrak{z}}{}}

\renewcommand{\k}{{\mathfrak{k}}{}}

\newcommand{\CS}{\mathcal S}
\newcommand{\CE}{\mathcal E}
\newcommand{\CK}{\mathcal K}
\newcommand{\CR}{\mathcal R}

\newcommand{\CD}{\mathcal D}
\newcommand{\CJ}{\mathcal J}
\newcommand{\CI}{\mathcal I}
\newcommand{\CH}{\mathcal H}
\newcommand{\CB}{\mathcal B}
\newcommand{\CU}{\mathcal U}
\newcommand{\CF}{\mathcal F}

\newcommand{\CA}{\mathcal A}
\newcommand{\CC}{\mathcal C}

\newcommand{\vf}{\varphi}

\newcommand{\Jac}{\rm Jac}

\begin{document}

\frontmatter

\title{Deformation  Quantization for Actions of K\"ahlerian  Lie Groups}



\author{Pierre Bieliavsky}
\address{University of Louvain, Belgium}
\email{pierre.bieliavsky@gmail.com}
\thanks{}

\author{Victor Gayral}
\address{University of Reims, France}
\email{victor.gayral@univ-reims.fr}
\thanks{Work
supported by the Belgian Interuniversity Attraction Pole (IAP) within the framework ``Dynamics, Geometry and Statistical Physics'' (DYGEST)}
\date{ August 24, 2012, and, in revised form, November 4, 2013.}



\begin{abstract}
Let $\B$ be a Lie group admitting a left-invariant negatively curved K\"ahlerian
structure. Consider a strongly continuous action $\alpha$ of  $\B$ on a Fr\'echet algebra 
$\CA$. Denote by $\CA^\infty$ the associated   Fr\'echet algebra 
of smooth vectors for this action. In the Abelian case $\B=\R^{2n}$ and 
$\alpha$ isometric, Marc Rieffel proved in  \cite{Ri}
that Weyl's operator symbol composition formula (the so called Moyal product)
yields a deformation  through Fr\'echet  algebra structures
$\{\star_{\theta}^\alpha\}_{\theta\in\R}$ on $\CA^\infty$. When $\CA$ is a $C^*$-algebra, 
every deformed Fr\'echet algebra $(\CA^\infty,\star^\alpha_\theta)$ admits a compatible
pre-$C^*$-structure, hence yielding a deformation theory at the level
of $C^*$-algebras too. In this memoir, we prove both analogous 
statements for general negatively
curved K\"ahlerian groups.  
The construction relies on the one hand on combining a 
non-Abelian version of oscillatory integral
on tempered Lie groups with geometrical objects coming from invariant WKB-quantization of 
solvable symplectic symmetric spaces, and, on the second hand, in establishing a non-Abelian 
version of the Calder\'on-Vaillancourt Theorem. In particular, we give an oscillating kernel formula 
for  WKB-star products on symplectic symmetric spaces
that fiber over an exponential Lie group. 
\end{abstract}

\maketitle

\chapter*{Warning}

The published version of the present article contains a mistake which invalidates the proof of the invariance of the $K$-theory under 
the deformation (Theorem 8.50 in the printed version and Theorem 7.50 in the former arXiv version). 
The mistake is located in Proposition 8.47 (Proposition 7.47 in the former arXiv version), where it is claimed that the map $\R^2\to\R^2$,
$(a_1,a_2)\mapsto (e^{-2a_2}\sinh(2a_1), -e^{-2a_1}\sinh(2a_2))$ is diffeomorphism. In fact, this map only defines a diffeomorphism
from $\R^2$ onto a proper open subspace of $\R^2$. As a consequence, the deformed $C^*$-algebra $A_{\theta,\bm}$ needs not to be 
Strongly Morita equivalent to the reduced crossed-product $\mathbb B\ltimes (\CK(\CH_\chi)\otimes A)$ and the $K$-theory 
needs not to be an invariant of the deformation.  

In this updated version, we have removed the whole contents of  Section 7.7 (Section 8.7 of the published version). 

\tableofcontents

\mainmatter

\chapter*{Introduction}
The general idea of deforming a given theory by use of  its symmetries goes back  to Drinfel'd. 
One paradigm being that
the data of a \emph{Drinfel'd twist} based on a bi-algebra acting on an associative algebra $\A$,
produces an associative deformation of $\A$. 
In the context of Lie theory, one considers for instance the category of module-algebras
over the universal enveloping algebra $\CU(\g)$ of the Lie algebra $\g$ of a given Lie group $G$. 
In that situation, the notion 
of Drinfel'd twist is in a one to one  correspondence with the one of  left-invariant formal 
star-product  $\star_\nu$ on the space of formal power series $C^\infty(G)[[\nu]]$, see  \cite{Drinfeld}. 
Disposing of such a twist, every $\CU(\g)$-module-algebra $\A$ may then be formally deformed
into an associative algebra $\A[[\nu]]$.
 It is important to observe that, within this situation, the symplectic leaf $\B$ 
through the unit element $e$ of $G$ in the characteristic
foliation of the (left-invariant) Poisson structure directing the star-product $\star_\nu$, 
always consists of an immersed Lie subgroup
of $G$.  The Lie group $\B$ therefore carries a left-invariant \emph{symplectic} structure. This stresses 
the importance of \emph{symplectic Lie groups} (i.e$.$ connected Lie groups endowed with
invariant 
symplectic forms) as 
semi-classical approximations of Drinfel'd twists attached to Lie algebras.

\quad

 In the present memoir, we address the question of designing \emph{non-formal} Drinfel'd twists for 
 actions of symplectic Lie groups $\B$ that underly negatively curved \emph{K\"ahlerian  Lie 
 groups},
 i.e$.$ Lie groups that admit a left-invariant K\"ahlerian structure of negative curvature.
These groups exactly correspond to the normal $\bf j$-algebras defined by Pyatetskii-Shapiro in his work on automorphic 
forms \cite{PS}. In particular, this class of groups contains all Iwasawa factors $AN$ of Hermitian type simple Lie groups $G=KAN$.

Roughly speaking, one looks for a smooth one-parameter family of complex valued smooth 
two-point functions on the group, $\{K_\theta\}_{\theta\in\R}\subset C^\infty(\B\times\B,\C)$, 
with the property 
that, for every strongly continuous and isometric
 action $\alpha$ of $\B$ on a $C^*$-algebra 
$A$, the  following formula
\begin{equation}\label{UDFK}
a\star_\theta^\alpha b
:=\int_{\B\times \B}\, K_\theta(x,y)\,\alpha_x(a)\,\alpha_y(b)\,{\rm d}x\,{\rm d}y\;,
\end{equation}
defines a one-parameter deformation of  the  $C^*$-algebra $A$.

 The above program was realized by Marc Rieffel in the particular case of the Abelian Lie group 
 $\B=\R^{2n}$ in \cite{Ri}. More precisely,
Rieffel proved that for \emph{any} strongly continuous and isometric action of $\R^{2n}$ on  \emph{any} 
Fr\'echet algebra $\CA$, the 
associated Fr\'echet subalgebra $\CA^\infty$ of
smooth vectors for this action, is deformed by the rule (\ref{UDFK}), where the two-point kernel
there, consists of the Weyl symbol composition kernel: 
$$K_\theta(x,y):=\theta^{-2n}\,\exp\big\{\tfrac{i}{\theta}\omega^0(x,y)\big\}\;,$$
associated to a translation invariant symplectic structure  $\omega^0$ on $\R^{2n}$. 
The associated star-product therefore corresponds here to Moyal's product. 
In the special case where the Fr\'echet algebra $\CA$ is a $C^*$-algebra, Rieffel also constructed
a deformed $C^*$-structure, so that $(\CA^\infty,\star^\alpha_\theta)$ becomes a pre-$C^*$-algebra,
which in turn yields  a deformation theory at the level of $C^*$-algebras too. Many 
further results  have
been proved then (for example  continuity of the field of deformed $C^*$-algebras 
 \cite{Ri}, invariance of the $K$-theory 
  \cite{Ri2}...), and many 
 applications have been found (for instance in  locally compact quantum groups
 \cite{Ri3},  quantum fields theory \cite{biz,biz2},  spectral triples  \cite{GGIVS}...).

 In this  memoir, we  investigate the deformation theory of $C^*$-algebras endowed
 with an isometric action of a  negatively curved  K\"ahlerian  Lie group. Most of the results we present
 here are of a pure analytical nature. Indeed, once a family $\{K_\theta\}_{\theta\in\R}$ 
 of associative   (i.e$.$ such that the associated deformed product \eqref{UDFK}
 is at least formally associative) two-point functions has been found, in order
  to give a precise meaning to   the associated multiplication rule, there is no doubt that
  the integrals in \eqref{UDFK} need to be interpreted in a suitable (here oscillatory) sense.
  Indeed, there is no reason to expect the two-point function $K_\theta$ to be integrable:
  it is typically  not even bounded in the non-Abelian case!
  Thus, already in the case of an isometric action on a $C^*$-algebra, we have to face a serious
  analytical difficulty.
 We stress that contrarily to the case of $\R^{2n}$, in the situation of a non-Abelian group 
 action, this is an highly non-trivial feature of our deformation theory. 
  
  \quad
  
The memoir is organized as follows.

  \quad

In chapter 1, we start by introducing non-Abelian and unbounded versions of 
Fr\'echet-valued symbol 
spaces on 
a Lie group $G$, with Lie algebra $\g$:
$$
\CB^{\underline\mu}
(G,\CE):=\left\{f\in C^\infty(G,\CE)\,:\,\forall X\in\CU(\g),\,\forall j\in\N,\,\exists C>0\,:\,\|
\widetilde X f\|_j\leq C\,\mu_j\right\}\;,
$$
where $\CE$ is a Fr\'echet space,  $\underline\mu:=\{\mu_j\}_{j\in\N}$ is a countable family of 
specific positive functions on $G$, called \emph{weights} (see Definition \ref{weight-def})
affiliated to a countable set of semi-norms  $\{\|.\|_j\}_{j\in\N}$ defining the Fr\'echet topology on 
$\CE$ and where $\widetilde X$ is the left-invariant differential operator on $G$ associated 
to an element $X\in\CU(\g)$.
For example, $\CB^{1}(G,\C)$ consists of the smooth vectors for the right regular representation of 
$G$ on the space of 
bounded right-uniformly continuous functions on $G$ (the uniform structures on $G$ are
 generally not balanced in our non-Abelian situation) and it coincides with Laurent Schwartz's 
 space $\CB$ when $G=\R^n$. We shall also mention that  function spaces on Lie groups
  of a similar type
 are considered in \cite{LS}
and in \cite{LW} in the context of actions of $\R^d$ on locally convex 
algebras. \\
We then define a notion of oscillatory integrals on Lie groups $G$ that are endowed with a specific 
type of smooth function
$S\in C^\infty(G,\R)$ (see Definitions \ref{temp-grp}, \ref{TEMPPAIR} and  \ref{TEMPADM}). We call 
such a pair $(G,S)$ an \emph{admissible tempered pair}. The main
result of this chapter is that associated to an admissible tempered pair $(G,S)$, and given a 
growth-controlled function $\bm$, the oscillatory integral 
$$
\CD(G,\CE)\to\CE\,,\quad F\mapsto\int_G\bm\,e^{iS}\,F\;,
$$
canonically extends from $\CD(G,\CE)$, the space of smooth compactly supported functions,  to our symbol space $\CB^{\underline\mu}(G,\CE)$.
This construction is explained in Definition \ref{OI}, which  turns out to apply in our situation as a direct consequence of 
Proposition \ref{PROPIP2}, the main technical result of this chapter.

\quad

In chapter 2, we consider an arbitrary \emph{normal $\bf j$-group} $\B$ (i.e$.$ a connected
simply connected Lie group whose Lie algebra is a 
normal $\bf j$-algebra---see Definition \ref{normal-j-alg}). The main result of this chapter, Theorem
\ref{TASP},
shows that its square $\B\times\B$ canonically underlies an admissible tempered pair 
$(\B\times\B,S_{\rm can}^\B)$.  When elementary, every normal $\bf j$-group has a canonical simply 
transitive action on a specific
solvable symplectic symmetric space. The two-point function $S_{\rm can}^\B$ we consider here 
comes from an 
earlier work of one of us. It  consists of the sum of the phases $S_{\rm can}^{\S_j}$
of the oscillatory kernels associated to invariant star-products on  solvable symplectic 
symmetric space \cite{BiMas,Bi07}, 
 in the Pyatetskii-Shapiro decomposition \cite{PS} of a normal $\bf j$-group $\B$  into a sequence of split extensions of elementary normal
  $\bf j$-factors: $\B=(\dots((\S_1\ltimes\S_2)\ltimes\S_3)\ltimes\dots)\ltimes\S_N$.
The two-point phase function $S_{\rm can}^\S$ in the elementary
 case, then consists of the symplectic area of 
the unique geodesic triangle in $\S$ (viewed as a  solvable symplectic symmetric space), whose geodesic edges admit $e, x$ and $y$ as midpoints ($e$ denotes the 
unit element of the group $\S$):
$$
S_{\rm can}^\S(x_1,x_2):={\rm Area}\left(\Phi_\S^{-1}(e,x_1,x_2)\right)\;,
$$
with
$$
\Phi_\S:\S^3\to\S^3\,,\quad (x_1,x_2,x_3)\mapsto \big(\mid(x_1,x_2),\mid(x_2,x_3),
\mid(x_3,x_1)\big)\;,
$$
where $\mid(x,y)$ denotes the geodesic midpoint between $x$ and $y$  in $\S$ 
(again uniquely defined in our situation).

\quad

 In chapter 3, we consider  an arbitrary
normal $\bf j$-group $\B$, and define the above-mentioned oscillatory kernels $K_{\theta}$ simply by 
tensorizing oscillating kernels found in \cite{BiMas} on elementary $\bf j$-factors. The resulting 
kernel has the form 
$$
K_\theta={\theta^{-\dim\B}}\,\bm_{\rm can}^\B\,\exp\big\{\tfrac{i}{\theta}S_{\rm can}^\B\big\}\;,
$$
where $S_{\rm can}^\B$ is the two-point phase mentioned in the description of chapter 2 above, 
and $\bm_{\rm can}^\B=\bm_{\rm can}^{\S_1}\otimes\dots\otimes\bm_{\rm can}^{\S_N}$, where $\bm_{\rm can}^{\S_j}=\Jac_{\Phi^{-1}_{\S_j}}^{1/2}$
denotes  the square root of the Jacobian of the ``medial triangle" map $\Phi_\S^{-1}$. In particular, it defines an oscillatory integral on every symbol 
space of the type
$\CB^{\underline\mu}(\B\times\B,\CB^{\underline\nu}(\B,\CE))$.
When valued in a Fr\'echet \emph{algebra} $\CA$, this yields a non-perturbative and associative star-product 
$\star_{\theta}$ on the union of 
all symbol spaces $\CB^{\underline\mu}(\B,\CA)$.

\quad

In chapter 4, we consider any \emph{tempered} action of a normal $\bf j$-group $\B$ on a Fr\'echet 
algebra $\CA$. By tempered action
we mean a strongly continuous action $\alpha$ of $\B$ by automorphisms on $\CA$, such that for 
every semi-norm $\|.\|_j$ there is a  weight (``tempered" for a suitable notion of temperedness)  $\mu_j^\alpha$ such that 
$\|\alpha_g(a)\|_j\,\leq\,\mu_j^\alpha(g)\,\|a\|_j
$ for all $a\in\CA$ and $g\in\B$. In that case, the space of smooth vectors $\CA^\infty$ for $\alpha$
naturally identifies with a subspace of $\CB^{\underline{\hat\mu}}(\B,\CA^\infty)$ (where $\underline{\hat\mu}$ 
is affiliated 
to $\underline\mu^\alpha=\{\mu_j^\alpha\}_{j\in\N}$) through the injective map:
$$
\alpha:\CA^\infty\to\CB^{\underline{\hat\mu}}(\B,\CA^\infty)\,,\quad a\mapsto[g\mapsto\alpha_g(a)]\;.
$$
We stress that even in the case of an isometric action, and contrarily to the Abelian 
situation,  the map $\alpha$ always takes values in a symbol space $\CB^{\underline\mu}$, 
for a non-trivial sequence of weights
$\underline\mu$. This explains why the non-Abelian framework forces to consider such symbol spaces.
 Applying the results of chapter 3 to this situation, we get a new associative product on 
 $\CA^\infty$ defined by the 
formula 
$$
a\star^\alpha_{\theta} b:=\big(\alpha(a)\star_{\theta}\alpha(b)\big)(e)\;.
$$
Then, the main result of this chapter, stated as Theorem \ref{UDF} in the text, is the following
fact:

\quad

\noindent
{\bf Universal Deformation Formula for Actions of K\"ahlerian Lie Groups on Fr\'echet Algebras:} \\
{\it Let  $(\CA,\alpha,\B)$ be a Fr\'echet algebra endowed with 
a tempered action  of a normal $\bf j$-group.  
Then,  $(\CA^\infty,\star_{\theta}^\alpha)$ is an associative Fr\'echet algebra, (abusively) called the 
Fr\'echet deformation of $\CA$.}

\quad

Following the terminology introduced in \cite{GGS}, the word ``universal" refers to the fact that our deformation procedure applies to \emph{any} Fr\'echet algebra  the K\"ahlerian Lie group acts on.
The word ``universal" does not refer to the possibility that our construction might
yield all such deformation procedures valid for a given K\"ahlerian Lie group.
However, regarding the last sentence, the following remark can nevertheless be made. In order to get rid of technicalities let us consider, for a short moment, the purely formal 
framework of formal Drinfel'd twists based on the bi-algebra underlying the enveloping algebra $\CU(\b)$ of the Lie group $\B$. Given such a Drinfel'd twist 
$F\in\CU(\b)\otimes\CU(\b)[[\nu]]$, one defines the \emph{internal symmetry} of the twist as the group $G(F)$,  consisting of  diffeomorphisms that preserve the twist:
$$
G(F):=\{\varphi:\B\to\B\;|\;\varphi_\star\tilde{F}=\tilde{F}\}\;,
$$
where $\tilde{F}$ denotes the left-invariant (formal sequence of) operator(s) associated to the twist $F$. Within the present work
we treat (all) the situations where, in the elementary normal case, the internal symmetry equals the automorphism group of a 
symmetric symplectic space structure on $\B$ whose underlying affine connection consists in the canonical torsion free invariant connection on the group $\B$
(see \cite{Bi07}).

\quad

The rest of the memoir is devoted to the construction of a pre-$C^*$-structure on 
$(\CA^\infty,\star_{\theta}^\alpha)$, in the case $\CA$ is a $C^*$-algebra.
The  
method we use is a generalization of Unterberger's Fuchs calculus \cite{Un84}
and fits within the general framework of ``Moyal quantizer'' as defined in
\cite{CGV} (see also  \cite[Section 3.5]{GBFV}).

\quad 

 In chapter 5, we define a special class of symplectic symmetric spaces which naturally 
give rise to explicit WKB-quantizations (i.e$.$ invariant star-products representable through oscillatory 
kernels)
that underlie pseudo-differential
operator calculus. Roughly speaking, an \emph{elementary symplectic symmetric space}
is a symplectic symmetric space that consists of the total space of a fibration in flat fibers over a Lie 
group $Q$ of exponential type.
In that case, a  variant of Kirillov's orbit method yields  a unitary and self-adjoint representation 
on an Hilbert-space $\CH$,
of the symmetric space $M$:
$$
\Omega:M\to \mathcal U_{\rm sa}(\CH)\;,
$$
with associated ``quantization rule'':
$$
\Omega: L^1(M)\to \CB(\CH)\,,\quad F\mapsto\Omega(F):=\int_MF(x)\Omega(x)\,{\rm d}x\;.
$$
Weighting the above mapping by the multiplication by a (growth controlled) function $\bm$ defined on the base $Q$ yields a pair of adjoint maps:
$$
\Omega_\bm:L^2(M)\to{\mathcal L}^2(\CH)\quad\mbox{\rm and}\quad\sigma_\bm:{\mathcal L}^2(\CH)\to L^2(M)\;,
$$
where ${\mathcal L}^2(\CH)$ denotes the Hilbert space  of Hilbert-Schmidt operators on $\CH$. Both of the above maps are equivariant under
the \emph{whole} automorphism group of $M$. Note that this last feature very much contrasts with the usual notion of coherent-state quantization
for groups (as opposed to symmetric spaces). The corresponding ``Berezin transform" $B_\bm:=\sigma_\bm\circ\Omega_\bm$ is explicitly controlled. In particular, when invertible,
the associated star-product $F_1\star F_2:=B_\bm^{-1}\sigma_\bm(\Omega_\bm(F_1)\Omega_\bm(F_2))$ is of oscillatory (WKB) type and its associated
kernel is explicitly determined. Note that, because entirely explicit, this chapter yields a proof of Weinstein's conjectural form for star-product WKB-kernels on symmetric spaces \cite{We} in the situation considered here. The chapter ends with considerations on extending the construction to semi-direct products.

\vspace{3mm}

Chapter 6 is entirely devoted to applying the construction of chapter 5 to the particular case of K\"ahlerian Lie groups with negative curvature.
Such a Lie group is always a normal ${\bf j}$-group is the sense of Pyatetskii-Shapiro. Each of its elementary factors admits the structure
of an elementary symplectic symmetric space. Accordingly to \cite{Bi07},  the obtained non-formal star products coincide with the one described in chapters 2-4.

\vspace{3mm}

Chapter 7 deals with the deformation theory for $C^*$-algebras. We eventually prove the following statement:

\quad

\noindent
{\bf Universal Deformation Formula for Actions of K\"ahlerian Lie Groups on $C^*$-Algebras:} \\
{\it Let  $(A,\alpha,\B)$ be a $C^*$-algebra endowed with 
a strongly continuous and isometric action  of a normal $\bf j$-group.  
Then,  there exists a canonical  $C^*$-norm on the   Fr\'echet algebra 
$(A^\infty,\star_{\theta}^\alpha)$. Its  $C^*$-closure  is  (abusively) called the 
$C^*$-deformation of $A$.  }

\quad

 The above statement follows from a non-Abelian generalization of the 
Calder\'on-Vaillancourt
Theorem in the context of the usual Weyl pseudo-differential calculus on $\R^{2n}$
(see Theorem \ref{C-def}). Our result 
asserts that the element $\Omega_\bm(F)$
associated to a function $F$ in $\CB^1(\B,A)$ naturally consists of an
element of the spatial tensor product of $A$ by $\CB(\CH)$.
 Its proof relies on a combination of a resolution of the identity 
obtained from wavelet analysis considerations
(see section \ref{WA})
 and further properties of our oscillatory integral defined in chapter \ref{OSIL}.
 We finally  prove that the $K$-theory is an invariant of the deformation.

\qquad 

{\bf Acknowledgments}
We warmly thank the referee for his positive criticism. His comments, remarks and suggestions
  have greatly improved this work.

\chapter*{Notations and conventions}

Given a 
Lie group $G$, with Lie algebra $\g$, we denote by ${\rm d}_G(g)$ a left invariant Haar 
measure. In the non-unimodular case, we consider  the modular function $\Delta_G$, defined by   
the relation:
$$
{\rm d}_G(g)\Delta_G(g):={\rm d}_G(g^{-1})\;.
$$  
Unless otherwise specified, $L^p(G)$, $p\in[1,\infty]$, will always denote the Lebesgue  space 
associated with the choice of a  left-invariant Haar measure made above. We also denote by
$\CD(G)$ the space of smooth compactly supported functions on $G$ and by $\CD'(G)$ the dual 
space of distributions.

We use the notations $L^\star$ and $R^\star$, for the left and right regular actions:
\begin{equation*}
L^\star_gf(g'):=f(g^{-1}g')\,,\qquad R^\star_gf(g'):=f(g'g)\;.
\end{equation*}
By $\widetilde X$ and $\underline X$, we mean the left-invariant and right-invariant 
vector fields on $G$ associated to  the elements $X$ and $-X$ of $\g$:
\begin{equation*}
\widetilde X:=\frac d{dt}\Big|_{t=0}\,R^\star_{e^{tX}}\,,\qquad  \underline X:=\frac d{dt}\Big|_{t=0}\,L^
\star_{e^{tX}}\;.
\end{equation*}
Given a    element $X$
of the universal enveloping algebra $\CU(\g)$ of $\g$, 
we adopt the same notations $\widetilde{X}$ and $\underline X$ for  the associated left- and
 right-invariant differential 
operators on $G$.
More generally, if $\alpha$ is an action of $G$ on a  topological 
vector space $\CE$,  we consider the 
infinitesimal 
form of the action, given for $X\in\g$ by:
\begin{equation*}
X^\alpha(a):=\frac d{dt}\Big|_{t=0}\,\alpha_{e^{tX}}(a)\,,\qquad a\in \CE^\infty\;,
\end{equation*}
and extended  to the whole universal enveloping algebra $\CU(\g)$, by declaring that the map
$\CU(\b)\to \End(\CA^\infty)$, $X\mapsto X^\alpha$ is an algebra homomorphism.
Here, $\CE^\infty$ denotes the set of smooth vectors for the action:
$$
\CE^\infty:=\big\{a\in\CE\;:\; [g\mapsto\alpha_g(a)]\in C^\infty(G,\CE)\big\}\;.
$$

 Let $\Delta_{\CU}$ be the ordinary co-product of $\CU(\g)$. 
We make use of the Sweedler notation:
$$
\Delta_{\CU}(X)=\sum_{(X)} X_{(1)}\otimes X_{(2)}\in\CU(\g)\otimes\CU(\g)\,,\qquad X\in\CU(\g)\;,
$$
and accordingly, for $f_1,f_2\in C^\infty(G)$ and $X\in\CU(\g)$, we write
\begin{equation}
\label{Sweedler}
\widetilde X(f_1f_2)= \sum_{(X)} \big(\widetilde X_{(1)}\,f_1\big)\,\big(\widetilde X_{(2)}\,f_2\big)\,,
\quad \underline X(f_1f_2)= \sum_{(X)} \big(\underline X_{(1)}\,f_1\big)\,\big(\underline X_{(2)}
\,f_2\big)\;.
\end{equation}
More generally, we use the notation
\begin{align}
\label{itS}
(\Delta_{\CU}\otimes\Id)\circ\Delta_{\CU}(X) &=\sum_{(X)}\sum_{(X_{(1)})} (X_{(1)})_{(1)}
\otimes (X_{(1)})_{(2)}\otimes X_{(2)}\\
&=:\sum_{(X)} X_{(11)}\otimes X_{(12)}\otimes X_{(2)}\nonumber\;,
\end{align}
and obvious generalization of it.

To a fixed  ordered basis $\{X_1,\dots,X_m\}$ of the Lie algebra $\g$, we associate a
PBW basis of $\CU(\g)$:
\begin{equation}
\label{PBW}
\{ X^\beta,\,\beta\in\N^{m}\}\,,\qquad X^\beta:=X_1^{\beta_1}X_2^{\beta_2}\dots X_m^{\beta_m}\;.
\end{equation}
This induces a filtration 
$$
\CU(\g)=\bigcup_{k\in\N}\CU_k(\g)\,,\qquad\CU_k(\g)\subset \CU_l(\g)\,,\quad
k\leq l \;,
$$
 in terms of the subsets
\begin{equation}
\label{filter}
\CU_k(\g):=\Big\{\sum_{|\beta|\leq k} C_\beta\,X^\beta\,,\;C_\beta\in\R\Big\}\,,\quad k\in\N\;,
\end{equation}
where $|\beta|:=\beta_1+\dots+\beta_m$.
For $\beta,\beta_1,\beta_2\in\N^m$, we define the `structure constants' 
$\omega_\beta^{\beta_1,\beta_2}\in\R$ of $\CU(\g)$, by
\begin{equation}
\label{structure-constants}
X^{\beta_1} X^{\beta_2}=\sum_{|\beta|\leq|\beta_1|+|\beta_2|} 
\omega_\beta^{\beta_1,\beta_2}\,X^\beta\in \CU_{|\beta_1|+|\beta_2|}(\g)\;.
\end{equation}
We endow the finite dimensional vector space $\CU_k(\g)$, with the $\ell^1$-norm  $|.|_k$
within the basis $\{X^\beta,|\beta|\leq k\}$:
\begin{equation}
\label{norm-CUk}
|X|_k:=\sum_{|\beta|\leq k} |C_\beta|\qquad \mbox{if}\qquad X=\sum_{|\beta|\leq k} C_\beta\,
 X^\beta\in\CU_k(\g)\;.
\end{equation}
We  observe that the  family of norms $\{|.|_k\}_{k\in\N}$ is compatible with the 
filtered 
structure of $\CU(\g)$, in the sense
that if $X\in\CU_k(\g)$, then $|X|_k=|X|_l$ whenever $l\geq k$.
Considering a  subspace $V\subset\g$, we also denote by $\CU(V)$ the \emph{unital} 
subalgebra of $\CU(\g)$ generated by $V$:
\begin{align}
\label{CUVE}
\CU(V)={\Span}\Big\{X_1X_2\dots X_n\,:\, X_j\in V\,,n\in\mathbb N\Big\}\;,
\end{align}
that we may filtrate using the induced filtration of $\CU(\g)$. We also observe that the co-product 
preserves  the latter subalgebras,  in the sense that $\Delta_\CU\big(\CU(V)\big)\subset\CU(V)\otimes\CU(V)$.

Regarding the uniform structures on a locally compact group $G$, we say that a function
$f:G\to \mathbb C$ is right (respectively left) uniformly continuous if for all $\varepsilon>0$, there exists $U$, an open 
neighborhood of the neutral element $e$, such that for all $(g,h)\in G\times G$ we have
$$
|f(g)-f(h)|\leq\varepsilon\,,\quad\mbox{whenever} \quad g^{-1}h\in U\quad\big(\mbox{respectively}\;
\;hg^{-1}\in U\big)\;.
$$
We call a Lie group $G$ (with Lie algebra $\g$) exponential, if the exponential map
$\exp:\g\to G$ is a global diffeomorphism.

Let $f_1,f_2$ be two real valued functions on  $G$.  We say that $f_1$ and $f_2$
have the same behavior, that we write 
$f_1\asymp f_2$,
when
there exist $0<c\leq C$ such that for all $x\in G$, we have $cf_1(x)\leq f_2(x)\leq Cf_1(x)$.

\chapter{Oscillatory integrals}
\label{OSIL}

This chapter is the most technical part of this memoir. To outline  its content, we shall first 
explain the situation on which the general theory is designed. Consider $\B$ a simply connected
Lie group endowed with a left invariant K\"ahlerian structure with negative curvature. Let also
$\alpha$ be a strongly continuous and isometric action of $\B$ on a $C^*$-algebra $A$. In this
context and from previous works of one of us, we have at our disposal
 two functions $\bm_{\rm can}\in C^\infty(\B
\times\B,\R_+^*)$ and $S_{\rm can}\in C^\infty(\B\times\B,\R)$ such that setting for  $\theta\in\R^*$
$$
K_\theta={\theta^{-\dim\B}}\,\bm_{\rm can}\,\exp\big\{\tfrac{i}{\theta}S_{\rm can}\big\}
\in C^\infty(\B\times\B,\C)\;,
$$
the following formula
\begin{equation}
\label{sous}
a\star_\theta^\alpha b:=\int_{\B\times\B} K_\theta(x,y) \,\alpha_x(a)\,\alpha_y(b)\,{\rm d}x\,{\rm d}y\,,
\qquad a,b\in A\;,
\end{equation}
formally defines a one-parameter family of associative algebra structures on $A$. Since the
function $\bm_{\rm can}$ is unbounded and since 
the map $\alpha(a):=[x\mapsto\alpha_x(a)]$ is constant
in norm, the only hope to go beyond the formal level is to define the integral sign in \eqref{sous}
in an oscillatory way. What we are precisely looking for, is a pair 
 $(A_0, \underline\bD)$, where 
\begin{itemize}
\item  $A_0$ is a dense and $\alpha$-stable Fr\'echet subalgebra
of $A$ with a topology (finer than the uniform one)  
determined by a countable set of semi-norms $\{\|.\|_j\}_{j\in\N}$ 
\item $\underline\bD:=\{\bD_j\}_{j\in\N}$ is a 
countable family of differential operators on $\B\times \B$, such that 
for all $a,b\in A_0$ and all $j\in\N$, the image
under $\bD_j$  of the map 
$[(x,y)\in\B\times\B\mapsto \bm_{\rm can}(x,y) \,\alpha_x(a)\,\alpha_y(b)\in A_0]$ belongs to 
$L^1(\B\times\B, (A_0,\|.\|_j))$ and, denoting by $\bD_j^*$ the formal adjoint of $\bD_j$, such that
$$
\bD_j^* \exp\big\{\tfrac{i}{\theta}S_{\rm can}\big\}=\exp\big\{\tfrac{i}{\theta}S_{\rm can}\big\}\;.
$$
\end{itemize}

Once such a pair is found, there is a clear  way to give a meaning of $\star_\theta^\alpha$
on $A_0$, namely by the $\|.\|_j$-absolutely convergent $A$-valued integral:
$$
a\star_\theta^\alpha b:={\theta^{-\dim\B}}\int_{\B\times\B} \exp\big\{\tfrac{i}{\theta}S_{\rm can}(x,y)\big\} \,
\bD_{j;x,y}\big(\bm_{\rm can}(x,y) \,\alpha_x(a)\,\alpha_y(b)\big)\,{\rm d}x\,{\rm d}y\;.
$$
There is an obvious candidate for $A_0$, which is $A^\infty$,  the Fr\'echet subalgebra of $A$ consisting of  smooth
vectors for the action $\alpha$:
$$
A^\infty:=\big\{a\in A\;:\; \alpha(a)=[x\mapsto\alpha_x(a)]\in C^\infty(\B,A)\big\}\;.
$$
Recall that by the strong continuity assumption of the action, $A^\infty$ is dense in $A$.
Then, the crucial observation, proved in Lemma \ref{embed}, is that we have
a continuous embedding:
$$
\alpha\otimes\alpha: A^\infty\times A^\infty\to\CB^{\underline\mu}(\B\times\B,A^\infty)\,,\quad (a,b)\mapsto[(g,h)\mapsto \alpha_g(a)\alpha_h(b)]\;,
$$
where the symbol space $\CB^{\underline\mu}(\B\times\B,A^\infty)$ is defined in \eqref{Bfam}
and where $\underline\mu=\{\mu_j\}_{j\in\N}$ is a  family of unbounded functions
on $\B\times\B$, called weights (see Definition \ref{weight-def} below). We should already
stress that
it is because the group $\B$ is non Abelian (and noncompact) that we are forced to consider
such symbol spaces associated to unbounded weights.
 Therefore, the problem becomes the construction of
a family of differential operators $\underline\bD=\{\bD_j\}_{j\in\N}$, such that
for all $j\in\N$, we have
$$
\bD_j:\CB^{\mu_j}(\B\times\B,A^\infty)\to L^1\big(\B\times\B,(A^\infty,\|.\|_j)\big)\;,
$$
and
$$
\bD_j^* \exp\big\{\tfrac{i}{\theta}S_{\rm can}\big\}=\exp\big\{\tfrac{i}{\theta}S_{\rm can}\big\}\;.
$$

The construction of the family of differential operators $\underline
\bD$, and therefore of the associated notion of
oscillatory integral, is  rather involved and this  is the whole subject of this chapter. Since in this
memoir we eventually need an oscillatory integral for different groups (for $\B\times\B$ in chapters
\ref{NFSP} and \ref{DFA} and for $\B$ in chapter \ref{SCSTAR}) and since we intend  to apply
the whole ideas of the present construction to  situations involving other groups and other
kernels, we have decided to formulate
the results of this chapter in the greatest possible generality.

\section{Symbol spaces}\label{SS}
In this preliminary section,  
we consider a non-Abelian, weighted and Fr\'echet space valued version of the Laurent 
Schwartz space $\CB$ of smooth functions that, together with all of their derivatives, are bounded. 
For reasons that will become clear later
(cf$.$ chapter \ref{QPSSS}), we refer to such function spaces as {\em symbol spaces}.
They are constructed out of a family of specific functions on a Lie group $G$, that we call
{\em weights}. The prototype of a weight for a non-Abelian Lie group is given in  
Definition \ref{MWS}.
The key properties of these symbol spaces are established in Lemmas \ref{symbols} and
\ref{SmoothFamily}. In Lemma \ref{OUTPUTLEM}, we show on a fundamental example, how such spaces
naturally appear in the context of non-Abelian Lie group actions.

\begin{dfn} 
\label{weight-def}
Consider   a connected  real Lie group $G$ with Lie algebra $\g$.
An element $\mu\in C^\infty(G,\R_+^*)$ is called a {\bf weight}  if it satisfies the following properties:

\noindent (i) For every element $X\in\CU(\g)$, there exist $C_L,C_R>0$ such that:
$$
|\widetilde{X}\mu|\leq C_L\,\mu\quad\mbox{and}\quad |\underline{X}\,\mu|\leq C_R\,\mu\;.
$$
\noindent (ii) There exist  $C,L,R\in\R_+^*$  such that for all $g,h\in G$ we have:
$$
\mu(gh)\leq C\,\mu(g)^L\,\mu(h)^R\;.
$$
\noindent A pair $(L,R)$ as in item (ii) is called a {\bf sub-multiplicative degree} of the weight $\mu$. A weight
with sub-multiplicative degree $(1,1)$ is called a {\bf sub-multiplicative weight}.
\end{dfn}
\begin{rmk}
\label{vee}
For $\mu\in C^\infty(G)$, we set $\mu^{\!\vee}\!(g):=\mu(g^{-1})$. Then, from the relation
 $\widetilde X \,\mu^{\!\vee}=(\underline X\,\mu)^{\!\vee}$  for all $X\in\CU(\g)$, we see that $\mu$ is a 
 weight of sub-multiplicative degree $(L,R)$ if and only if
 $\mu^{\!\vee}$ is a weight of sub-multiplicative degree $(R,L)$. Moreover, a product of two weights 
is a weight and a (positive) power of a weight is a weight. Also, the tensor product of two weights
is a weight on the direct product group.
\end{rmk}

 In the following, we construct a canonical and non-trivial weight for  non-Abelian Lie 
groups.  This specific weight is an important object as it will naturally and repeatedly appear 
 in all our analysis (see for instance Lemmas  \ref{OUTPUTLEM}, \ref{genout},  \ref{belle},
\ref{cinqcinq} and  \ref{embed}).
 
 \begin{dfn}
 \label{MWS}
 Choosing an Euclidean norm $|.|$ on $\g$, for $x\in G$, we let $|\Ad_{x}|$ be the operator norm
of the adjoint action of $G$ on $\g$.
The function
$$
\mathfrak{d}_G:G\to\R^*_+\,,\quad x\mapsto
\sqrt{1\,+\,\left|\Ad_{x}\right|^2\,+\,\left|\Ad_{x^{-1}}\right|^2}\;,
$$
is called the {\bf modular weight} of $G$.
 \end{dfn}
 \begin{lem}
The modular weight $\fd_G$ is  a sub-multiplicative weight on $G$. \end{lem}
\begin{proof}
To prove that $\fd_G$ is a weight,
we start   from the relations for $X\in\g$ and $x\in G$:
\begin{align*}
\widetilde X|\Ad_x|^2&=2\sup_{Y\in\g,\,|Y|=1}\langle\Ad_x\circ\ad_X(Y),\Ad_x(Y)\rangle\;,\\
\widetilde X|\Ad_{x^{-1}}|^2&=-2\sup_{Y\in\g,\,|Y|=1}\langle\ad_X\circ\Ad_{x^{-1}}(Y),\Ad_{x^{-1}}(Y)
\rangle\;,\\
\underline X|\Ad_x|^2&=-2\sup_{Y\in\g,\,|Y|=1}\langle\ad_X\circ\Ad_x(Y),\Ad_x(Y)\rangle\;,\\
\underline X|\Ad_{x^{-1}}|^2&=2\sup_{Y\in\g,\,|Y|=1}\langle\Ad_{x^{-1}}\circ\ad_X(Y),\Ad_{x^{-1}}(Y)
\rangle\;.
\end{align*}
Then, we get by induction and  for every $X\in\CU(\g)$ of strictly positive
homogeneous degree:
$$
\big|\widetilde{X}\,\fd_G(x)\big|\,,\big|\underline{X}\,\fd_G(x)\big|\leq \left|\ad_X\right|\frac{\left|\Ad_{x}\right|^2
+\left|\Ad_{x^{-1}}\right|^2}{\sqrt{1\,+\,\left|\Ad_{x}\right|^2\,+\,\left|\Ad_{x^{-1}}\right|^2}}\leq  |\ad_X|\,\fd_G(x)\;,
$$
where, for $X\in \CU(\g)$, we denote by $|\ad_{X}|$  the operator norm
of the adjoint action of $\CU(\g)$ on $\g$. This implies condition (i) of Definition \ref{weight-def}.
Last, the inequality
$$|\Ad_{gh}|=|\Ad_g\circ\Ad_h|\leq |\Ad_{g}||\Ad_{h}|\,,\quad g,h\in G\;,$$
implies the sub-multiplicativity  of $\fd_G$.
\end{proof}

We next give further informations regarding the modular weight of a semi-direct product of groups.
These results will be used in Lemmas \ref{fd-explicit}, \ref{mwdi}, Proposition \ref{defofdef} and
 Lemma \ref{F-eta}.

\begin{lem}\label{fd-semi}
Let $G_j$, $j=1,2$, be two connected Lie groups with  Lie algebras $\g_j$, let
$\bR:G_1\to\mbox{\rm Aut}(G_2)$ be an extension
homomorphism and consider the associated semi-direct product  
$G:=G_1\ltimes_\bR G_2$ whose Lie algebra
is denoted by $\g=\g_1\ltimes\g_2$.
Consider the adjoint actions $\Ad^j:G_j\to GL(\g_j)$ and $\Ad:G\to GL(\g)$
and define the mapping 
$$
\Phi:G_2\to {\rm Hom}(\g_1,\g_2)\,,\quad
g_2\mapsto\big[X_1\mapsto \Ad_{g_2}(X_1)-X_1\big]\;.
$$
Then, in the parametrization $G_2\times G_1\to G$, $(g_2,g_1)
\mapsto g_1g_2$,  we have the following behavior of the modular weight of the semi-direct product:
\begin{align*}
\fd_{G}\asymp\Big[(g_2,g_1)\mapsto
\fd_{G_1}(g_1)+\big(1&+
|\Phi(g_2)\circ\Ad^1_{g_1}|^2+| \Ad^2_{g_2}\circ\bR_{g_1}|^2\\&+
|\bR_{g_1^{-1}}\circ\Phi(g_2^{-1})|^2+| \bR_{g_1^{-1}}\circ\Ad^2_{g_2^{-1}}|^2
\big)^{1/2}\Big]\;,
\end{align*}
where we use the same notation for the extension homomorphism and its derivative:
$$
\g_2\to \g_2\,,\quad X\mapsto\frac d{dt}\Big|_{t=0} \bR_{g_1}(e^{tX})\,,\quad g_1\in G_1\;.
$$
In particular, we get the following behaviors of  the restrictions of $\fd_{G}$ to the subgroups 
$G_1$ and $G_2$:
\begin{align*}
\fd_{G}\big|_{G_1}&\asymp\fd_{G_1}+\Big[g_1\mapsto \big(1+|\bR_{g_1}|^2
+|\bR_{g_1^{-1}}|^2\big)^{1/2}\Big]
\\
\fd_{G}\big|_{G_2}&\asymp\fd_{G_2}+\Big[g_2\mapsto\big(1+|\Phi(g_2)|^2
+|\Phi(g_2^{-1})|^2\big)^{1/2}\Big]\;.
\end{align*}
Moreover, in the case of a direct product, we have the behavior
$$
\fd_{G_1\times G_2}\asymp \fd_{G_1}\otimes 1+1\otimes\fd_{G_2}\;.
$$
\end{lem}
\begin{proof}
Fix  Euclidean structures on $\g_1$ and on $\g_2$ and induce  one on $\g_1\oplus\g_2$ by 
declaring that $\g_1$ is orthogonal to $\g_2$. A direct computation shows that
in the vector decomposition $\g=\g_1\oplus\g_2$, the operator $\Ad_{g_2g_1}$ takes the following matrix form:
$$
\Ad_{g_2g_1}=\left(\begin{array}{cc} \Ad^1_{g_1} & 0\\ \Phi(g_2)\circ\Ad^1_{g_1} &
 \Ad^2_{g_2}\circ
\bR_{g_1}\end{array}\right)\;,
$$
with inverse given by
$$
\Ad_{(g_2g_1)^{-1}}=\left(\begin{array}{cc} \Ad^1_{g_1^{-1}} & 0\\ \bR_{g_1^{-1}}
\circ\Phi(g_2^{-1}) & 
\bR_{g_1^{-1}}\circ\Ad^2_{g_2^{-1}}\end{array}\right)\;.
$$
But on the finite dimensional vector space $\End(\g)$,  the operator norm of $\Ad_{g_2g_1}$
is equivalent to its Hilbert-Schmidt norm, and the latter reads
$$
\big(|\Ad^1_{g_1}|^2+| \Phi(g_2)\circ\Ad^1_{g_1}|^2+| \Ad^2_{g_2}
\circ\bR_{g_1}|^2\big)^{1/2}\;.
$$
Similarly, the operator norm of $\Ad_{(g_2g_1)^{-1}}$ is equivalent to 
$$
\big(|\Ad^1_{g_1^{-1}}|^2+ | \bR_{g_1^{-1}}\circ\Phi(g_2^{-1})|+ 
|\bR_{g_1^{-1}}\circ\Ad^2_{g_2^{-1}}|^2\big)^{1/2}\;,
$$
proving the first claim. The remaining statements follow immediately.
 \end{proof}

We should  mention that the modular function, $\Delta_G$, is also  a sub-multiplicative weight
on $G$. 
Indeed the multiplicativity property implies that for every
 $X\in\CU(\g)$ and $x\in G$: 
 $$
 \big(\widetilde X\,\Delta_G\big)(x)=\big(\widetilde X\,\Delta_G\big)(e) \,\Delta_G(x)\,,\qquad
  \big(\underline X\,\Delta_G\big)(x)=\big(\underline X\,\Delta_G\big)(e) \,\Delta_G(x)\;.
 $$
 However, this weight will not be of much interest in what follows.

 The next notion will play a key role to establish density results for our symbol spaces.
We assume from now on the Lie group $G$ to be non-compact.
\begin{dfn}
\label{weight}
Given two weights $\mu$ and $\hat\mu$, we  say that $\hat\mu$ {\bf dominates} $\mu$, which we denote by $\mu\prec\hat\mu$, if 
$$
\lim_{g\to\infty}\frac{\mu(g)}{\hat\mu(g)}=0\;.
$$
\end{dfn}
\begin{rmk}
 For negatively curved K\"ahlerian 
Lie groups, the modular weight  has the crucial property to dominate the constant weight $1$
(see Corollary \ref{MWP} and Lemma \ref{mwdi}).
\end{rmk}

We now  let $\CE$ be a complex Fr\'echet space with topology underlying a countable family of  semi-norms $\{\|.\|_j\}_{j\in\N}$.
Given a weight  $\mu$, 
we first consider  the following space of $\CE$-valued functions on $G$:
\begin{align}
\label{Bun}
&\CB^\mu(G,\CE):=\\
&\qquad\left\{F\in C^\infty(G,\CE)\,:\,\forall X\in\CU(\g),\,\forall j\in\N,\,\exists C>0\,:\,\|\widetilde X F\|_j\leq C\,\mu\right\}\;.\nonumber
\end{align}
When $\CE=\C$ (respectively when  $\mu=1$, respectively when $\CE=\C$ and $\mu=1$), we denote $\CB^\mu(G,\CE)$ by $\CB^\mu(G)$ 
 (respectively  by $\CB(G,\CE)$, respectively by $\CB(G)$).
We endow the space $\CB^\mu(G,\CE)$  with the natural topology associated to the following semi-norms:
\begin{equation}
\label{norms}
 \|F\|_{j,k,\mu}:=\sup_{X\in\,\CU_k(\g)}\,\sup_{g\in G}\Big\{\frac{\|\widetilde{X}F(g)\|_j}
 {\mu(g)\,|X|_k}\Big\}\,,\qquad \,\,j,k\in\N\;,
 \end{equation}
 where $\CU(\g)=\cup_{k\in\N}\CU_k(\g)$ is the filtration described in \eqref{filter} and 
 $|.|_k$ is the norm on $\CU_k(\g)$ defined in \eqref{norm-CUk}. Note that for
 $X=\sum_{|\beta|\leq k} C_\beta\,X^\beta\in\CU_k(\g)$, we have
 $$
 \frac{\|\widetilde{X}F(g)\|_j}{|X|_k}\,\leq\,\frac{\sum_{|\beta|\leq k} |C_\beta| \|\widetilde{X^\beta}F(g)\|_j}{
 \sum_{|\beta|\leq k} |C_\beta|}\,\leq\, \max_{|\beta|\leq k} \|\widetilde{X^\beta}F(g)\|_j\;,
 $$
 and hence
 \begin{align}
 \label{eqeq}
  \|F\|_{j,k,\mu}\,\leq \,\max_{|\beta|\leq k}\,\sup_{g\in G}\frac{ \|\widetilde{X^\beta}F(g)\|_j}{\mu(g)}
  =\max_{|\beta|\leq k}\, \|\widetilde{X^\beta} F\|_{j,0,\mu}\;,
 \end{align}
 which shows  that the semi-norms \eqref{norms} are well defined on $\CB^\mu(G,\CE)$.
 When $\CE=\C$ (respectively when  $\mu=1$, respectively when $\CE=\C$ and $\mu=1$), we 
 denote the semi-norms
 \eqref{eqeq} by $\|.\|_{k,\mu}$, $k\in\N$ (respectively by $\|.\|_{j,k}$, $j,k\in\N$, respectively 
 by $\|.\|_{k}$, $k\in\N$).

 The basic properties of the spaces $\CB^\mu(G,\CE)$ are 
established in the next lemma. They are essentially standard but are used all over the text.
In particular, in the last item, we prove that $\CD(G,\CE)$ is a dense subset of $\CB^\mu(G,\CE)$
for the induced topology of $\CB^{\hat\mu}(G,\CE)$, for $\hat\mu$ an arbitrary weight which
dominates $\mu$. This fact will be used in a crucial way  for the construction of the oscillatory integral
given in Definition \ref{OI}.
 
\quad

Let  $C_b(G,\CE)$ be the Fr\'echet space of $\CE$-valued continuous and bounded   functions 
on $G$. The topology we consider on the latter is the one associated to  the semi-norms 
$\|F\|_{j}:=\sup_{g\in G}\|F(g)\|_j$, $j\in\N$. This space carries an action of $G$ by 
right-translations. This action is of course isometric but not necessarily strongly continuous.  Consider 
therefore its closed subspace $C_{ru}(G, \CE)$ constituted by the right-uniformly continuous 
functions.  

\begin{lem}\label{symbols}
Let $(G,\CE)$ as above and let $\mu$, $\nu$ and $\hat\mu$ be three  weights on $G$.
\begin{enumerate}
\item[(i)] The right regular action $R^\star$ of $G$ on $C_{ru}(G, \CE)$ is isometric and strongly 
continuous. 
\item[(ii)] Let $C_{ru}(G, \CE)^\infty$ be the subspace of $C_{ru}(G, \CE)$ of smooth vectors for the 
right regular action. Then $C_{ru}(G, \CE)^\infty$ identifies with $\CB(G,\CE)$ as topological vector 
spaces. In particular, $\CB(G,\CE)$ is Fr\'echet.
\item[(iii)] The left regular action $L^\star$ of $G$ on $\CB(G,\CE)$ is isometric.
\item[(iv)] The bilinear map:
\begin{equation*}
\CB^{\mu}(G)\times\CB^{\nu}(G,\CE)\to\CB^{\mu\nu}(G,\CE)\,,\quad (u,F)\mapsto[g\in G\mapsto u(g)
\,F(g)\in\CE]\;,
\end{equation*}
is continuous.
\item[(v)] The map
$$
\CB^\mu(G,\CE)\to \CB(G,\CE)\,,\quad F\mapsto \mu^{-1} F\;,
$$
is an homeomorphism. In particular, the space $\CB^\mu(G,\CE)$ is Fr\'echet as well.
\item[(vi)] For every $X\in\CU(\g)$, the associated left invariant differential operator $\widetilde X$ 
acts continuously on $\CB^\mu(G,\CE)$.
\item[(vii)]  If there exists $C>0$ such that $\mu \leq C\hat\mu$, then $\CB^{\mu}(G,\CE)\subset 
\CB^{\hat\mu}(G,\CE)$, continuously.
\item[(viii)] Assume that $\mu\prec\hat\mu$. Then the closure of $\CD(G,\CE)$
in $\CB^{\hat\mu}(G,\CE)$ contains $\CB^{\mu}(G,\CE)$. In particular, the space $\CD(G,\CE)$ is a 
dense subset of 
$\CB^{\mu}(G,\CE)$ for the induced topology of $\CB^{\hat\mu}(G,\CE)$.
\end{enumerate}
\end{lem}
\begin{proof}
(i) Recall that $G$ being locally compact and countable at infinity, the space $C_b(G,\CE)$ is 
Fr\'echet (by the same argument as in the proof of \cite[Proposition 44.1 and Corollary 1]{Tr}). The 
subspace $C_{ru}(G, \CE)$ is then closed as a uniform limit of right-uniformly continuous functions is 
right-uniformly continuous. 
Thus $C_{ru}(G, \CE)$ endowed with the induced topology is a Fr\'echet space as well.\\
Being isometric on $C_b(G,\CE)$, the right action is consequently isometric on $C_{ru}(G, \CE)$ too.
Moreover, for any converging sequence $\{g_n\}\subset G$, with limit $g\in G$, and any 
$F\in C_{ru}(G, \CE)$, we have 
$\|(R^\star_{g_n}-R^\star_g)F\|_{j}=\sup_{g_0\in G}\|F(g_0g_n)-F(g_0g)\|_j$ which tends to zero 
due to the right-uniform continuity of $F$. Hence the right regular action $R^\star$ is strongly 
continuous on $C_{ru}(G, \CE)$.

(ii) Note that an element $F\in C_{ru}(G, \CE)^\infty$ is such that the map
$[g\mapsto R^\star_gF]$ 
is smooth as a $C_{ru}(G, \CE)$-valued function on $G$. In particular, for every $X\in\CU(\g)$, 
$\widetilde{X}F$ is bounded and smooth. This clearly  gives  the inclusion $C_{ru}(G, \CE)^\infty
\subset\CB(G,\CE)$.
Reciprocally, $G$ acts on $\CB(G,\CE)$ via the right regular representation. Indeed,  for all 
$g\in G$ and $X\in\CU(\g)$, we have 
$$
\widetilde{X}\circ R^\star_g=R^\star_g\circ\widetilde{(\Ad_{g^{-1}}X)}\;,
$$
and hence for $j,k\in\N$ and $F\in\CB(G,\CE)$, we deduce
\begin{align*}
\|R^\star_gF\|_{j,k}&=\sup_{X\in\,\CU_k(\g)}\,\sup_{g'\in G}
\Big\{\frac{\|\widetilde{(\Ad_{g^{-1}} X)} F(g'g)\|_j}
 {|X|_k}\Big\}\\
& =\sup_{X\in\,\CU_k(\g)}\,\sup_{g'\in G}
\Big\{\frac{\|\widetilde{(\Ad_{g^{-1}} X)} F(g')\|_j}
 {|X|_k}\Big\}\leq |\Ad_{g^{-1}}|_k \|F\|_{j,k}\;,
\end{align*}
where $|\Ad_{g}|_k$ denotes the operator norm of the adjoint action of $G$ on the
(finite dimensional) Banach space $(\CU_k(\g),|.|_k)$.
 Now we have the inclusion
$\CB(G,\CE)\subset C_{ru}(G,\CE)$. Indeed,  for $F\in\CB(G)$, $g\in G$ and for fixed $X\in\g$, one observes that 
\begin{align*}
\big|F(g\exp(tX))-F(g)\big|&=\Big|\int_0^t\, \frac{\rm d}{{\rm d}\tau}(F(g\exp(\tau X)))\,{\rm d}\tau\Big|
\\&=\Big|\int_0^t\, \widetilde XF(g\exp(\tau X))\,{\rm d}\tau\Big|\leq |X|_1\,\|F\|_{1} \,|t|\;,
\end{align*}
where $\|F\|_1=\sup_{X\in\CU_1(\g)}\sup_{g\in G}|\widetilde X F|/|X|_1$.
Hence, we obtain  the right-uniform continuity of $F$. To show that $F\in\CB(G)$ is a differentiable vector for the right-action, we observe that
\begin{align*}
\Big|\frac{1}{t}\big(F(g\exp(t X))-F(g)\big)-\big(\widetilde{X}F\big)(g)\Big|&\leq\int_0^1\Big|\big(\widetilde{X}F\big)(g\exp(t\tau X))-\big(\widetilde{X}F\big)(g)\Big|\, {\rm d}\tau\\
&\leq \int_0^1\int_0^{t\tau}\Big|\big(\widetilde{X}^2F\big)(g\exp(\tau' X))\Big|\, {\rm d}\tau'\, {\rm d}\tau\\
&\leq
|t|\,\sup_{g\in G}\left\{\big|\widetilde{X}^2F\big|(g)\right\}\leq |X^2|_2\,\|F\|_{2}\,|t|\;,
\end{align*}
which tends to zero together with $t$. (Here $\|F\|_2=\sup_{X\in\CU_2(\g)}
\sup_{g\in G}|\widetilde X F|/|X|_2$.)
 This yields differentiability at the unit element. One gets it everywhere else by observing that
\begin{equation}
\label{tech1}
\widetilde{X}_g(R^\star_g F)=R^\star_g(\widetilde{X}F)\,,\qquad \forall\,X\in\CU(\g)\,,\quad
\forall\,g\in G\,,\quad
\forall\,F\in\CB(G)\;.
\end{equation}
An induction on the order of derivation implies $\CB(G)\subset C_{ru}(G)^\infty$.
The $\CE$-valued case is entirely similar. \\
The assertion concerning the topology follows from the definition of the topology on smooth vectors \cite{War}. 

(iii) The fact that $G$ acts isometrically on $\CB(G,\CE)$ via the left regular representation, follows from
\begin{align*}
\|L^\star_gF\|_{j,k}&=\sup_{X\in\,\CU_k(\g)}\sup_{g'\in G}\frac{\|\widetilde{X}\big(L^\star_g 
F\big)(g')\|_j}{|X|_k}=\sup_{X\in\,\CU_k(\g)}\sup_{g'\in G}\frac{\|\big(L^\star_g\widetilde{X} 
F\big)(g')\|_j}{|X|_k}\\
&=\sup_{X\in\,\CU_k(\g)}\sup_{g'\in G}\frac{\|\widetilde{X} F(g^{-1}g')\|_j}{|X|_k} =
\sup_{X\in\,\CU_k(\g)}\sup_{g'\in G}\frac{\|\widetilde{X} F(g')\|_j}{|X|_k}=\|F\|_{j,k}\;.
\end{align*}

 (iv) Let  $u\in\CB^{\mu}(G)$ and $F\in\CB^{\nu}(G,\CE)$. Using  Sweedler's notation 
 \eqref{Sweedler}, we have for $j,k\in\N$:
 \begin{align*}
 \|uF\|_{j,k,\mu\nu}&=\sup_{X\in\,\CU_k(\g)}\,
 \sup_{g\in G} \frac{\|\widetilde X(uF)(g)\|_j}{\mu(g)\nu(g)|X|_k}\\
& \leq \sup_{X\in\,\CU_k(\g)}\,
 \sup_{g\in G} \sum_{(X)} \frac{\big|\big(\widetilde X_{(1)}u\big)(g)\big|\,
 \big\|\big(\widetilde X_{(2)}F\big)(g)\big\|_j}
 {\mu(g)\nu(g)|X|_k}\\
& \leq\Big( \sup_{X\in\,\CU_k(\g)}\sum_{(X)} \frac{|X_{(1)}|_k\,|X_{(2)}|_k}{|X|_k}\Big)\, \|u\|_{k,\mu}
  \|F\|_{j,k,\nu}\;.
 \end{align*}
 Now, for $X=\sum_{|\beta|\leq k} C_\beta \,X^\beta\in\CU_k(\g)$, expanded in the PBW basis \eqref{PBW},
  we have 
 \begin{align*}
 \Delta_\CU(X)&=\sum_{|\beta|\leq k} C_\beta \sum_{\gamma\leq\beta}
 \binom {\beta}{ \gamma}X^\gamma\otimes X^{\beta-\gamma}
\;,
 \end{align*}
 which, with $m$ the dimension of $\g$, implies that
 \begin{equation}
 \label{compat}
 \sum_{(X)} |X_{(1)}|_k\,|X_{(2)}|_k\leq \sum_{|\beta|\leq k} |C_\beta| \sum_{\gamma\leq\beta}
\binom{\beta}{ \gamma}
 \leq 2^{mk} \sum_{|\beta|\leq k} |C_\beta|=2^{mk}|X|_k\;.
 \end{equation}
 Hence we get
 $$
  \|uF\|_{j,k,\mu\nu}\leq  2^{mk}\|u\|_{k,\mu}
  \|F\|_{j,k,\nu}\;,
  $$
 proving  the continuity.

(v) Note first that since $\mu$ is a weight,
  for every $X\in\CU(\g)$ there is a constant $C>0$ such that 
$|\widetilde X(\mu^{-1})|\leq C \mu^{-1}$. This easily  follows by induction once one has noticed that 
for $X\in\g$, we have $\widetilde X( \mu^{-1})=-(\widetilde X\mu) \mu^{-2}$. This entails that
 if   $\mu\in\CB^\mu(G)$ then 
$\mu^{-1}\in\CB^{\mu^{-1}}(G)$, even if $\mu^{-1}$ may not be a weight in the sense of Definition 
\ref{weight-def}, as it may not have any sub-multiplicative degree. Then, (v) 
follows from (iv), where the existence of a sub-multiplicative degree is not used in that proof.

(vi) This assertion follows immediately  from equation \eqref{eqeq} and  (vii) follows easily from equation
 \eqref{norms}.

(viii) Choose an increasing sequence $\{C_n\}_{n\in\N}$ of relatively compact open subsets in $G$, such that $\lim_n C_n=G$. Pick $0\leq\psi\in\CD(G)$ of $L^1(G)$-norm one and define 
\begin{equation}
\label{en}
e_n:=\int_G \psi(g)\,R^\star_g(\chi_n) \,{\rm d}_G(g)\;,
\end{equation}
where $\chi_n$ denotes the characteristic function of $C_n$. It is clear that $e_n$ is an increasing 
family of smooth compactly supported functions, which by Lebesgue dominated convergence,
   converges point-wise to the unit function. Moreover, for all $F\in\CB^{\mu}(G,\CE)$ and
   $\hat\mu$ a weight dominating $\mu$, we have
\begin{align*}
\|(1-e_n)F\|_{j,0,\hat\mu}&=\sup_{g\in G}\Big\{\frac1{\hat\mu(g)}\big(1-e_n(g)\big)\|F(g)\|_j\Big\}
\\&\leq 
\|F\|_{j,0,\mu}\;
\sup_{g\in G}\Big\{\frac{\mu(g)}{\hat\mu(g)}\big(1-e_n(g)\big)\Big\}\;,
\end{align*}
which converges  to zero when $n$ goes to infinity, since $\mu/\hat\mu\to 0$ when $g\to\infty$ and for fixed $g\in G$, 
$1-e_n(g)$ decreases to zero when $n\to\infty$.
We need to show that the same property holds true for all the semi-norms
$\|.\|_{j,k,\hat\mu}$, $k\geq 1$. We use an induction.
First note that if $ X\in\g$, then we have
\begin{align*}
\widetilde X e_n=\frac d{dt}\Big|_{t=0}R^\star_{e^{tX}}(e_n)&=\frac d{dt}\Big|_{t=0}\int_G \psi(g)\,R^\star_{e^{tX}}R^\star_g(\chi_n) \,{\rm d}_G(g)\\&
=
\frac d{dt}\Big|_{t=0}\int_G \psi(e^{-tX}g)\,R^\star_g(\chi_n) \,{\rm d}_G( g)\\
&=\int_G \big(\underline{X}\psi\big)(g)\,
R^\star_g(\chi_n) \,{\rm d}_G( g)\;.
\end{align*}
A routine inductive argument then gives
\begin{equation}
\label{routine}
\widetilde X e_n=\int_G \big(\underline{X}\psi\big)(g)\,
R^\star_g(\chi_n) \,{\rm d}_G( g)\,,\quad\forall\,X\in\CU(\g)\;,
\end{equation}
which entails
$$
\|\widetilde X e_n\|_\infty\leq \|\underline{X}\psi\|_1<\infty,\quad \forall X\in\CU(\g)\;.
$$
(This means that the sequence $\{e_n\}_{n\in\N}$ belongs to $\CB(G)$, uniformly in $n$.)
Now, assume that $\|(1-e_n)F\|_{j,k,\hat\mu}\to 0$, $n\to\infty$, for a given $k\in\N$,
 for all $F\in \CB^{\mu}(G,\CE)$ and all  $j\in\N$. From the same reasoning as those leading to \eqref{eqeq}
 and with $X^\beta$ the element of the PBW basis of $\CU(\g)$ defined in \eqref{PBW}, we 
 see that
 $$
 \|(1-e_n)F\|_{j,k+1,\hat\mu}\leq \|(1-e_n)F\|_{j,k,\hat\mu}
 +\max_{|\beta|=k+1}\,\|\widetilde {X^\beta}\big((1-e_n)F\big)
 \|_{j,0,\hat\mu}\;.
 $$
 We only need to show that the second term in the inequality above goes to zero when $n\to\infty$, as
 the first does by induction hypothesis.
 Writing $X^\beta=X^\gamma X$, with $|\gamma|=k$ and $X\in\g$,
by virtue of the Leibniz rule, we get
$$
\widetilde {X^\gamma} \widetilde X\big((1-e_n)F\big)=-\widetilde {X^\gamma}\big((\widetilde X e_n)F\big)+\widetilde 
{X^\gamma}\big((1-e_n)\,\widetilde X F\big)\;.
$$
Note that
$$
\|\widetilde {X^\gamma}\big((1-e_n)\widetilde X F\big)\|_{j,0,\hat\mu}\leq
\|(1-e_n)\widetilde X F\|_{j,k,\hat\mu}\;,
$$
which converges to zero when $n\to\infty$ by induction hypothesis, since
 $\widetilde X F\in \CB^{\mu}(G,\CE)$ and $|\gamma|=k$.
 Regarding the first term, we have using Sweedler's notations \eqref{Sweedler} and for a finite sum:
$$
\widetilde {X^\gamma}\big((\widetilde X e_n)F\big)=\sum_{({X^\gamma})} \big(\widetilde {X^\gamma_{(1)}} \widetilde X e_n\big)\big(\widetilde {X^\gamma_{(2)}} F\big)\;.
$$
Note that  $\int_G \underline{P}\underline{X}\,\psi\, {\rm d}_G=0$ for any $P\in\CU(\g)$, $X\in\g$ any $\psi\in\CD(G)$. Indeed, this follows from an inductive argument starting with 
$$
\int_G \underline X\psi (g)\,{\rm d}_G( g)=\frac d{dt}\Big|_{t=0}\int_G L^\star_{e^{tX}}\big(\psi\big)(g)\, {\rm d}_G( g)=\frac d{dt}\Big|_{t=0}\int_G \psi(g)\, {\rm d}_G( g)=0\;,
$$
for all $X\in\g$.
Using \eqref{routine}, we arrive at
$$
\widetilde {X^\gamma}\big((\widetilde X e_n)F\big)=
\sum_{({X^\gamma})}\Big(
\int_G \big(\underline {X^\gamma}_{(1)} \underline X  \psi\big)(g)\big(R^\star_g(\chi_n)-1\big)\big)
{\rm d}_G(g)\Big)\widetilde {X^\gamma}_{(2)} F\;,
$$
which converges  to zero in the norms $\|.\|_{j,0,\hat\mu}$, $j\in\N$, since it is a finite sum of terms of 
the form $(1-e_n)F$ 
(with possibly re-defined $F$'s in $\CB^{\mu}(G,\CE)$ and  $\psi$'s in $\CD(G)$).
 \end{proof}

\begin{rmk}
\label{stress}
We stress that as $\mu\prec\hat\mu$ implies that $\mu\leq C\hat\mu$ for some $C>0$ (since
$G$ is locally compact), in Lemma \ref{symbols}, item (viii) strengthens item (vii). In particular,
when $\mu\prec\hat\mu$, $\CB^{\mu}(G,\CE)$ is continuously contained in $\CB^{\hat\mu}(G,\CE)$.
This fact is of upmost importance for the continuity of the oscillatory integral, given in
 Definition 
\ref{OI}.
\end{rmk}

\begin{rmk}
 On $\CB(G,\CE)$,  the left regular action is generally not strongly continuous and  the right 
regular action is never isometric
unless $G$ is Abelian.
\end{rmk}

We now generalize  the spaces $\CB^\mu(G,\CE)$, by allowing a 
certain 
behavior at infinity of the $\CE$-valued functions on $G$, which is not necessarily uniform  with 
respect to the semi-norm index. We first introduce  some more notations.

\begin{dfn}
Let  $J$ be a countable set and let $\underline\mu:=\{\mu_j\}_{j\in J}$ be 
an associated family of weights  on a  Lie group $G$. We  
denote by $(\underline L,\underline R):=\{(L_j,R_j)\}_{j\in J}$ the associated family 
of sub-multiplicative degrees. Given two families of weights $\underline\mu,\underline{\hat\mu}$,
we say that $\underline{\hat\mu}$ dominates $\underline\mu$,  denoted by $\underline \mu\prec  
\underline{\hat\mu}$, if $\hat\mu_j$ dominates $\mu_j$  for all $j\in J$. The term by term 
product  (respectively tensor product) of two families of weights $\underline\mu$ and $\underline\nu$,
is denoted by $\underline\mu.\underline\nu$ (respectively by $\underline\mu\otimes\underline\nu$).
\end{dfn}

Consider a Fr\'echet space $\CE$  and 
$\underline\mu$  a 
countable family of weights  on $G$. We then define
\begin{align}
\label{Bfam}
&\CB^{\underline\mu}(G,\CE):=\\
&\qquad\left\{F\in C^\infty(G,\CE)\,:\,\forall X\in\CU(\g),\,\forall j\in\N,\,\exists C>0\,:\,\|
\widetilde X F\|_j\leq C\,\mu_j\right\}\;.\nonumber
\end{align}
We endow the latter space with the following set of  the semi-norms:
\begin{equation}
\label{norms2}
 \|F\|_{j,k,\underline\mu}:=\sup_{X\in\,\CU_k(\g)}\,\sup_{g\in G}\Big\{\frac{\|\widetilde{X}F(g)\|_j}
 {\mu_j(g)\,|X|_k}\Big\}\,,\qquad \,\,j,k\in\N\;,
 \end{equation}
 As expected, the space $\CB^{\underline\mu}(G,\CE)$ is  Fr\'echet  for the topology 
induced by the
semi-norms \eqref{norms2} and most of the properties of Lemma \ref{symbols} remain true.

\begin{lem} 
\label{SmoothFamily}
Let $\big(G,\CE\big)$ as above and let $\underline\mu$,
$\underline\nu$ and 
$\underline{\hat\mu}$ be three families of weights  on $G$.
\begin{enumerate}
\item[(i)] The space $\CB^{\underline\mu}(G,\CE)$ is Fr\'echet.
\item[(ii)] Let $(\underline L,\underline R)$ be the sub-multiplicative degree  of $\underline\mu$. Then, for every $g\in G$ the 
left-translation $L^\star_g$ defines a continuous map from $\CB^{\underline\mu}(G,\CE)$ to 
$\CB^{\underline\lambda}(G,\CE)$, where $\underline\lambda:=\{\lambda_j\}_{j\in\N}$ 
with $\lambda_j:=\mu_j^{R_j}$.
\item[(iii)] The bilinear map:
\begin{equation*}
\CB^{\underline\mu}(G)\times\CB^{\underline{\nu}}(G,\CE)
\to\CB^{\underline\mu.\underline{\nu}}(G,\CE)\,,\quad (u,F)\mapsto[g\in G\mapsto u(g)\,F(g)\in\CE]\;,
\end{equation*}
is  continuous.
\item[(iv)] For every $X\in\CU(\g)$, the  left invariant differential operator $\widetilde X$ acts continuously on $\CB^{\underline\mu}(G,\CE)$.
\item[(v)]  If for every $j\in\N$, there exists $C_j>0$ such that $\mu_j \leq C_j\hat\mu_j$, then 
$\CB^{\underline\mu}(G,\CE)\subset \CB^{\underline{\hat\mu}}(G,\CE)$.
\item[(vi)] Assume that $\underline\mu\prec\underline{\hat\mu}$. Then, the closure of  $\CD(G,\CE)$
in $\CB^{\underline{\hat\mu}}(G,\CE)$ contains $\CB^{\underline\mu}(G,\CE)$. In particular, 
$\CD(G,\CE)$ is a dense subset of $\CB^{\underline\mu}(G,\CE)$ for the induced topology of 
$\CB^{\underline{\hat\mu}}(G,\CE)$.
\end{enumerate}
\end{lem}
\begin{proof}
(i) For each $j\in\N$, define $\|.\|^\sim_j:=\sum_{k=0}^j\|.\|_k$. Clearly,  the topologies on
$\CE$ associated with the families of semi-norms 
$\{\|.\|_j\}_{j\in\N}$ and $\{\|.\|_j^\sim\}_{j\in\N}$ are equivalent.
Thus, we may assume without loss of generality that the family of semi-norms $\{\|.\|_j\}_{j\in\N}$ is increasing.
We start by recalling the standard realization of the Fr\'echet space $(\CE,\{\|.\|_j\}_{j\in\N})$ as a projective limit. One considers the null spaces
$V_j:=\{v\in\CE\;|\;\|v\|_j=0\}$ and form the normed quotient spaces $\dot{\CE}_j:=\CE/V_j$. Denoting by $\CE_j$ the Banach completion of the latter,
the family of semi-norms being increasing, one gets, for every pair of indices $i\leq j$, a natural continuous linear mapping
$g_{ji}:\CE_j\to\CE_i$. The Fr\'echet space $\CE$ is then isomorphic to the subspace $\tilde{\CE}$ of the product space $\prod_j\CE_j$ constituted 
by the elements $(x)\in\prod_j\CE_j$ such that $x_i=g_{ji}(x_j)$. Within this setting, the subspace $\tilde{\CE}$ is endowed with the projective topology associated with the family of maps
$\{f_j:\tilde{\CE}\to\CE_j:(x)\mapsto x_j\}$ (i.e$.$ the coarsest topology that renders continuous each of the $f_j$'s--- see e.g$.$ \cite[pp. 50-52]{Schaeffer}).
 Within this context, we then observe that the topology on $\CB^{\underline\mu}(G,\CE)\;
\simeq\;\CB^{\underline\mu}(G,\tilde{\CE})$ induced by the semi-norms
(\ref{norms2}) consists of the projective topology associated with the mappings $\phi_j:\CB^{\underline\mu}(G,\tilde{\CE})\to\CB^{\mu_j}(G,\CE_j):F\mapsto f_j\circ  F$.
Next we consider a Cauchy sequence $\{F_n\}_{n\in\N}$ in $\CB^{\underline\mu}(G,\tilde{\CE})$. Since every space $\CB^{\mu_j}(G,\CE_j)$ is Fr\'echet,
each sequence $\{f_j\circ F_n\}_{n\in\N}$ converges in $\CB^{\mu_j}(G,\CE_j)$ to an element denoted by $F^j$. Moreover, for every $g\in G$, one has
\begin{align*}
\|g_{ji}\circ F^j(g)-F^i(g)\|_i&=\|g_{ji}\circ F^j(g)-f_i\circ F_n(g)+f_i\circ F_n(g)-F^i(g)\|_i\\
&\leq\|g_{ji}\circ \big(F^j-f_j\circ F_n\big)(g)\|_i+\|f_i\circ F_n(g)-F^i(g)\|_i\;,
\end{align*}
which can be rendered as small as we want since every $g_{ji}$ is continuous. Hence
$g_{ji}\circ F^j=F^i$ which amounts to say that $\CB^{\underline\mu}(G,\tilde{\CE})$ is complete.

(ii) Let $F\in \CB^{\underline\mu}(G,\CE)$  and $g\in G$ and set $\lambda_j:=\mu_j^{R_j}$.
 We have for $j,k\in\N$:
\begin{align*}
\|L^\star_g \,F\|_{j,k,\underline\lambda}&=\sup_{X\in\,\CU_k(\g)}
\sup_{g'\in G}\frac{\|\widetilde X\big(L^\star_g \,F\big)(g')\|_j}
{\mu_j(g')^{R_j}|X|_k}=\sup_{X\in\,\CU_k(\g)}\sup_{g'\in G}\frac{\|\big(L^\star_g \,\widetilde XF\big)(g')\|_j}
{\mu_j(g')^{R_j}|X|_k}\\&=\sup_{X\in\,\CU_k(\g)}\sup_{g'\in G}\frac{\|\big( \widetilde XF\big)(g^{-1}g')\|_j}
{\mu_j(g')^{R_j}|X|_k}\leq \mu_j(g^{-1})^{L_j}\,\|F\|_{j,k,\underline\mu}\;.
\end{align*}

 Items  (iii), (iv), (v) and (vi) are proved in the same way as their counterparts in Lemma \ref{symbols}.
\end{proof}

\begin{rmk}
\label{stress2}
In the same fashion as in Remark \ref{stress}, we see that when $\underline\mu\prec
\underline{\hat\mu}$, then 
$\CB^{\underline\mu}(G,\CE)\subset \CB^{\underline{\hat\mu}}(G,\CE)$ continuously.
\end{rmk}

In Lemma \ref{OUTPUTLEM}, we show on a first example,
how the notion of $\CB$-spaces for families of weights, 
naturally appears in the context of non-Abelian Lie group actions.
We start by a preliminary result.
We  fix an Euclidean structure on $\g$, such that the basis 
$\{X_1,\dots,X_n\}$ (from which we have constructed the PBW basis \eqref{PBW})  is orthonormal.
\begin{lem}
\label{lemtwo}
For $g\in G$ and $k\in\N$, denote by $|\Ad_g|_k$ the operator norm of the adjoint action $\Ad$
of $G$ on the finite dimensional  Banach space $\big(\CU_k(\g),|.|_k\big)$. Then, for each $k\in\N$, there exists 
 a constant $C_k>0$, such that for all $g\in G$
  $$
 |\Ad_g|_k\leq C_k\,\fd_G(g)^k\;.
 $$
 where $\fd_G\in C^\infty(G)$ is the modular weight (given in Definition \ref{MWS}).
 \end{lem}
\begin{proof}
Note first  that for all $k\in\N$, there
exists a constant $\omega_k>0$ such that for all $X\in\CU_{k_1}(\g)$ and 
$Y\in\CU_{k_2}(\g)$,  we have
\begin{align}
\label{fou}
|X\,Y|_{k_1+k_2}\leq \omega_{k_1+k_2}\,|X|_{k_1}\,|Y|_{k_2}\;.
\end{align}
Indeed, 
observe that if 
$$
X=\sum_{|\beta|\leq k_1} C^1_\beta\, X^\beta\in\CU_{k_1}(\g)\quad\mbox{and}\quad 
Y=\sum_{|\beta|\leq k_2} C^2_\beta\, X^\beta\in\CU_{k_2}(\g)\;,
$$
 we have
$$
X\,Y=\sum_{|\beta_1|\leq k_1,|\beta_2|\leq k_2} C^1_{\beta_1}\,C^2_{\beta_2}
\sum_{|\beta|\leq|\beta_1+\beta_2|} \omega_\beta^{\beta_1,\beta_2} \, X^\beta\;,
$$
where the constants $ \omega_\beta^{\beta_1,\beta_2}$ are defined in \eqref{structure-constants}.
The sub-additivity of the norm $|.|_{k_1+k_2}$ then  entails that
\begin{align*}
|X\,Y|_{k_1+k_2}&
\leq\sum_{|\beta_1|\leq k_1,|\beta_2|\leq k_2} |C^1_{\beta_1}|\,|C^2_{\beta_2}|
\,\sum_{|\beta|\leq|\beta_1+\beta_2|} |\omega_\beta^{\beta_1,\beta_2}| \;.
\end{align*}
Thus, it leads to defining
$$
\omega_{k_1+k_2}:=\sup_{|\beta_1+\beta_2|\leq k_1+k_2}\; 
\sum_{|\beta|\leq|\beta_1+\beta_2|} |\omega_\beta^{\beta_1,\beta_2}|\;,
$$
and the inequality \eqref{fou} is proved.
Next, for
 $$
 X=\sum_{|\beta|\leq k} C_\beta\, X_1^{\beta_1}\dots X_m^{\beta_m}\in\CU_{k}(\g)\;,
 $$
  we have
 $$
 \Ad_g(X)=\sum_{|\beta|\leq k} C_\beta\, \big(\Ad_g(X_1)\big)^{\beta_1}\dots 
 \big(\Ad_g(X_m)\big)^{\beta_m}\in\CU_{k}(\g)\;,
 $$
 and thus by the previous considerations, we deduce
 \begin{align*}
 \big| \Ad_g(X)\big|_k&\leq \sum_{|\beta|\leq k} |C_\beta|\, \Big(\prod_{j=2}^{|\beta|}\omega_j\Big)
 \,|\Ad_g(X_1)|_1^{\beta_1}\,
 |\Ad_g(X_2)|_1^{\beta_2}\dots 
|\Ad_g(X_m)|_1^{\beta_m}\;.
 \end{align*}
 As the restriction of the norm $|.|_1$ from $\CU_1(\g)$ to $\g$ coincides with the $\ell^1$-norm
 of $\g$ within the basis $\{X_1,...,X_m\}$, we deduce for $j=1,\dots,m$ and with $|\Ad_g|$ the operator norm
 of $\Ad_g$ with respect to the Euclidean structure of $\g$ chosen:
 $$
 |\Ad_g(X_j)|_1\leq\sqrt m\, |\Ad_g(X_j)|\leq  \sqrt m\, |\Ad_g|\,|X_j|_\g=\sqrt m\, |\Ad_g|\;,
 $$
 as $X_j\in\g$ belongs to the unit sphere of $\g$ for the Euclidean  norm $|.|_\g$. This implies
 $$
 \big| \Ad_g(X)\big|_k\leq  m^{k/2}\,  \Big(\sup_{|\beta|\leq k}\prod_{j=2}^{|\beta|} \omega_j\Big)\,|\Ad_g|^k\;,
 $$
and the result follows from Definition \ref{MWS}.
 \end{proof}

\begin{lem}\label{OUTPUTLEM}
Let $\underline\mu$ be a family of weights on $G$, with 
 sub-multiplicative degree $(\underline L,\underline R)$. Then 
the linear map
$$
{\mathcal R}:=\Big[F\in C^\infty(G,\CE)\mapsto\big[g\mapsto R^\star_gF\big]\in C^\infty\big(G,C^\infty(G,\CE)\big)\Big]\;,
$$
is continuous from $\CB^{\underline\mu}(G,\CE)$ to 
$
\CB^{\underline\nu}\big(G,\CB^{\underline\lambda}(G,\CE)\big)
$, where $\underline\nu:=\{\nu_{j,k}\}_{j,k\in\N}$ and $\underline\lambda:=\{\lambda_j\}_{j\in\N}$ with
$$
\nu_{j,k}:=\mu_j^{R_j}\,\fd_G^k\qquad\mbox{and}\qquad \lambda_j:=\mu_j^{L_j}\;.
$$
 More precisely, labeling by $(j,k)\in\N^2$
the semi-norm $\|.\|_{j,k,\underline\lambda}$ of $\CB^{\underline\lambda}(G,\CE)$, 
for each $(j,k,k')\in\N^3$,
there exists a constant $C>0$, such that for all $F\in \CB^{\underline\mu}(G,\CE)$, we have
$$
\|{\mathcal R}(F)\|_{(j,k),k',\underline\nu}\leq C\|F\|_{j,k+k',\underline\mu}\;.
$$
\end{lem}
\begin{proof}
Using the relation \eqref{tech1}, we obtain for $X\in \CU_{k'}(\g)$, $F\in \CB^{\underline\mu}(G,\CE)$ 
and $g\in G$:
$$
\|\widetilde X_g R^\star_g(F)\|_{j,k,\underline\lambda}
=\| R^\star_g(\widetilde XF)\|_{j,k,\underline\lambda}
=\sup_{Y\in\,\CU_{k}(\g)}\sup_{x\in G}\frac{\|\widetilde Y_x R^\star_g
\big(\widetilde XF\big)(x)\|_j}{\mu_j^{L_j}(x)|Y|_{k}}\;.
$$
Moreover, since for any $Y\in\CU(\g)$ and $g\in G$, we have
$
R^\star_{g^{-1}}\widetilde Y R^\star_g=\widetilde {\Ad_{g^{-1} }Y}
$
and since  $F\in\CB^{\underline\mu}(G,\CE)$ and  $\underline\mu$ has sub-multiplicative 
with degree $(\underline L,\underline R)$, we get
\begin{align*}
\|\widetilde X_g R^\star_g(F)\|_{j,k,\underline\lambda}&=\sup_{Y\in\,\CU_{k}(\g)} \sup_{x\in G}\frac{\|\big(\widetilde 
{\Ad_{g^{-1} }Y}\widetilde XF\big)(xg)\|_j}{\mu_j^{L_j}(x)|Y|_{k}}
\\&\leq C\|F\|_{j,k+k',\underline\mu} 
|\Ad_{g^{-1}}|_{ k} |X|_{k'}\sup_{x\in G}\frac{
\mu_j(xg)}{\mu_j^{L_j}(x)}\\
&\leq C\|F\|_{j,k+k',\underline\mu}\, |\Ad_{g^{-1}}|_{ k} |X|_{k'}\,\mu_j^{R_j}(g)\;,
\end{align*}
and one concludes using Lemma \ref{lemtwo}.
\end{proof}
\section{Tempered pairs}
\label{Tp}
In this section, we establish the main technical result of  the first part of this memoir  (Proposition
 \ref{PROPIP2}), on which the construction of our oscillatory integral (and thus of our universal 
 deformation formula for Fr\'echet algebras) essentially relies. To this aim,
 we start by introducing the class of {\em tempered Lie groups} (Definition \ref{temp-grp})
 and the sub-class of {\em tempered pairs} (Definition \ref{TEMPPAIR}). In Lemma \ref{fd-tempered},
 we give simple but important consequences for the modular weight, modular
function and Haar measure, when the group is tempered. The rest of this section
is devoted to the proof of Proposition \ref{PROPIP2}.

\begin{lem}
\label{lem:equiv-tempered-left-inv-or-not}
Let $G$ be a connected real Lie group and  $\psi:\R^{m}\to G$ be a global diffeomorphism. Then 
the multiplication and inverse operations seen through $\psi$ are tempered functions\footnote{By tempered function, we mean a smooth function
 whose every derivative  is bounded by a polynomial function. These functions are sometimes called ``slowly increasing".} (in the 
ordinary sense of $\R^m$) if and only if for 
every element $A\in\CU(\g)$ their derivatives along $\widetilde{A}$ are bounded by
a  function which is polynomial within the chart $\psi$.
\end{lem}
\begin{proof}
Denote $m(x,y)=m_x(y)=\psi^{-1}(\psi(x)\cdot \psi(y))$ and $\iota(x)=\psi^{-1}(\psi(x)^{-1})$ the 
multiplication and inverse of $G$ seen through $\psi\in{\rm Diff}(\R^m,G)$.
Denote by $x_e$ the transportation of the neutral element $e$ of $G$: $\psi(x_e):=e$.
Lastly,  identifying naturally $T_{x_e}(\R^m)$ with $\R^m$, for  every $X\in\R^{m}$ denote by 
$$
\widetilde X^
\psi_x:={{m_x}_\star}_{x_e}(X)={{\psi^{-1}}_\star} \widetilde{\left({\psi_\star}_{x_e}X\right)}_{\psi(x)}
\;,
$$
 the left invariant 
vector field corresponding to ${\psi_\star}_{x_e}X\in\g$.

 Assume $m$ and $\iota$ are tempered in the usual sense. Then for $X\in\R^{m}$, by definition
 \[ \widetilde X^\psi_x =  \frac{d}{dt} m(x,x_e+tX) \Big|_{t=0} \;, \]
which is a linear combination of partial derivatives of $m$ all of them being bounded by some 
polynomials in $x$ since $m$ is tempered. In the same way, the derivatives of left-invariant vector 
fields are linear combinations of higher partial derivatives of compositions of $m$ with itself in the 
second variable, which are also bounded by some polynomials. Hence the left-invariant vector 
fields are tempered, and consequently so are the left-invariant derivatives of $m$ and $\iota$.

Conversely, assume $m$ and $\iota$ are tempered in the sense of left-invariant vector fields.  We 
will see that the constant vector fields on $\R^{m}$ are linear combinations of left-invariant vector 
fields, the coefficients  being tempered functions.
Indeed, we have $X=\left({{m_x}_\star}_{x_e}\right)^{-1}\big(\widetilde X^\psi_x\big)$ and the matrix 
elements of that inverse matrix are finite sums and products of the matrix elements of the original 
one, which are tempered, divided by its determinant. Thus all we have to check is that the inverse 
of the determinant is a tempered function. But 
$1/\det(m_{x_{\star x_e}}) = \det(m_{\iota(x)_{\star x}})$  is tempered since $m$ and $\iota$ are.
\end{proof}

The preceding observation yields us to introduce the following notion:
\begin{dfn}
\label{temp-grp}
A Lie group $G$ is called {\bf tempered} if there exists a global coordinate system $\psi:\R^{m}\to G$ 
where the multiplication and inverse operations are tempered functions. A smooth function  $f$
on a tempered Lie group 
is  called a {\bf tempered function} if $f\circ\psi$ is tempered.
\end{dfn}
\begin{rmk}
 Every tempered Lie group, being diffeomorphic to an Euclidean space, is connected and simply 
 connected.
Moreover, by arguments similar to those of Lemma \ref{lem:equiv-tempered-left-inv-or-not}, a smooth function  $f$
on a tempered Lie group 
is   tempered  if and only if for any $X\in\CU(\g)$,  its derivative along $\widetilde X$ is 
bounded by a polynomial function  within the global chart $\R^{m}\to G$.
\end{rmk}

\begin{ex}
For any  (simply connected) nilpotent Lie group, the exponential coordinates,
$\g\to G:X\mapsto\exp(X)$, provides $G$ with a structure of a tempered Lie group. Indeed, in the case of a nilpotent Lie group, the Baker-Campbell-Hausdorff series is finite.
\end{ex}

\begin{rmk}
\label{LRORNOT}
 Observe also that we can replace in Lemma 
\ref{lem:equiv-tempered-left-inv-or-not} left-invariant differential operators by right-invariant
one. In fact, for a tempered group, left-invariant  vector fields
are linear combinations of right-invariant one
with tempered coefficients and vice versa.
\end{rmk}

\begin{lem}
\label{fd-tempered}
Let $G$ be a tempered Lie group. Then the modular weight $\fd_G$ 
and the modular function $\Delta_G$ are tempered.
Moreover, in the transported coordinates,
 every Haar measure on $G$ is a  multiple of a Lebesgue measure  on $\R^m$
by a tempered density. 
\end{lem}
\begin{proof}
The conjugate action $\bC:G\times G\to G:(g,x)\mapsto gxg^{-1}$ is a tempered map when read in the global coordinate system.
Therefore, the evaluation of the restriction of its tangent mapping to the first factor $G\times\{e\}$ on the constant section $0\oplus X$,
$X\in\g$, of $T(G\times G)$ consists of a tempered mapping:
$$
G\to\g\,,\quad g\mapsto\bC_{\star(g,e)}(0_g\oplus X)\;.
$$
The latter coincides with $g\mapsto\Ad_g(X)$. Varying $X$ in the finite dimensional
vector space $\g$, yields the tempered map $\Ad:G\to\End(\g)$. This shows that $\fd_G$
is tempered.

Transporting the group structure of $G$ to $\R^n$ by means of  the global coordinates,
it is clear that any Haar measure on $\R^n$ (for the transported group law) is absolutely continuous
with respect to the Lebesgue measure. Let ${\rm d}_G(\xi)$ be a left invariant Haar measure on $G$ transported to $\R^n$. Let also $\rho:\R^n\to \R$ be the Radon-Nikodym derivative of ${\rm d}_G(\xi)$ with respect to ${\rm d}\xi$, the Lebesgue measure on $\R^n$. Let $\xi_e\in\R^n$ be the transported neutral element of $G$. By left-invariance of the Haar measure ${\rm d}_G(\xi)$, we get
$$
\rho\big(m(\xi',\xi)\big)=\rho(\xi) |{\rm Jac}_{L_{\xi'}^\star}|(\xi)\,,\quad\forall \xi,\xi'\in \R^n\;,
$$
where $m(.,.)$ denotes the transported multiplication law on $\R^n$ and $L_\xi^\star$ stands for the associated left translation operator on $\R^n$. Letting $\xi\to\xi_e$, we deduce
$$
\rho(\xi)=\rho(\xi_e) |{\rm Jac}_{L^\star_{\xi}}|(\xi_e)\,,\quad\forall \xi\in \R^n\;,
$$
and we conclude by Lemma \ref{lem:equiv-tempered-left-inv-or-not} using the fact 
that the multiplication law is tempered. 

 Next, we let $\iota$ the inversion map of $G$ transported to $\R^n$. We have in the 
transported coordinates:
$$
\Delta_G(\xi)=\frac{ {\rm d}_G\big(\iota(\xi)\big)}{{\rm d}_G(\xi)}=\frac{ {\rm d}_G\big(\iota(\xi)\big)}{{\rm d}\xi} \frac{{\rm d}\xi}{{\rm d}_G(\xi)}
={\rm Jac}_{\iota}(\xi)\,\frac{ {\rm d}_G\big(\iota(\xi)\big)}{{\rm d}\big(\iota(\xi)\big)} \frac{{\rm d}\xi}{{\rm d}_G(\xi)}={\rm Jac}_{\iota}(\xi)\,\frac{\rho\big(\iota(\xi)\big)}{\rho(\xi)}\;,
$$
and we conclude using what precedes and the temperedness of the inversion map on $G$.
\end{proof}

 We now consider the data of a pair $(G,S)$ where 
$G$ is a connected real Lie group with real  Lie algebra $\g$ and
  $S$ is a real-valued smooth function on $G$.

\begin{dfn}
\label{TEMPPAIR}
The pair $(G,S)$ is called {\bf tempered}  if the following two properties are satisfied:
\begin{enumerate}
\item[(i)] The map
\begin{equation}\label{TEMPCOORD}
\phi:G\,\to\,\g^\star\,,\quad x\,\mapsto\,\Big[\,\g\to\R,\;X\mapsto\big(\widetilde{X}.\,S\big)(x)\,\Big]\;,
\end{equation}
is a global diffeomorphism.
\item[(ii)] The inverse map $\phi^{-1}:\g^\star\simeq\R^m\to G$ endows $G$ with the structure
of a tempered Lie group.
\end{enumerate}
\end{dfn}
\begin{rmk}
\label{S-temp}
Within the above situation, the function $S$ is itself automatically tempered. Indeed, in the proof of Lemma \ref{lem:equiv-tempered-left-inv-or-not}, we have seen 
that directional derivatives in the coordinate system $\phi$ are expressed as linear combinations
of left-invariant vector fields with tempered coefficients. Hence, within a basis $\{X_j\}_{j=1,...,N}$ of $\g$, denoting
$x_j\:=\big(\widetilde{X_j}S\big)( x)$, we have
$\partial_{x_j}S(x)=\sum_km_j^k(x)x_k$ where the matrix $(m_j^k)_{j,k}$ have tempered
entries. This implies that the partial derivatives
(of every strictly 
 positive order) of $S$ are tempered. In spherical  coordinates $(r,\theta)$  (associated to the 
$x_j$'s) one  observes that
$\partial_rS(r,\theta)=R(\theta)_1^k\partial_{x_k}S$ where $\theta$ belongs to the unit sphere $S^{N-1}$ and where
$R(\theta)$ is a (rotation) matrix that smoothly depends on $\theta$. Hence: 
$$
|S(x)|=\Big|C+\int_{r_0}^r\partial_\rho S(\rho,\theta){\rm d}\rho\Big|
\leq |C|+\int_{r_0}^r\left|R(\theta)_1^k\right|\,\rho^{n_k}{\rm d}\rho\;,
$$
which for large $x$ is smaller than a multiple of some positive power of $r$.  Therefore the function 
$S$ has polynomial growth too.
\end{rmk}

Given a tempered pair $(G,S)$, with $\g$ the Lie algebra of $G$, we now consider a vector space decomposition:
\begin{equation}\label{DECOMP}
\g=\bigoplus_{n=0}^NV_n\;,
\end{equation}
and for every $n=0,\dots, N$,
an ordered basis $\{e^n_{j}\}_{j=1,\dots,\dim (V_n)}$ of $V_n$. We  get global coordinates on $G$:
\begin{equation}\label{COORDINATES}
x_n^{j}:=\big(\widetilde{e_{j}^n}. S\big)(x)\,,\qquad n=0,\dots,N\,,\quad j=1,\dots,\dim (V_n)\;.
\end{equation}
We choose a scalar product on each $V_n$ and let $|.|_n$ be the associated Euclidean norm.
 Given an element $A\in\CU(\g)$,  
 we let  ${\widetilde A}^*$ be
the formal adjoint of the left-invariant differential operator $\widetilde A$, with respect to the inner 
product of $L^2(G)$. We make the obvious observation that ${\widetilde A}^*$  is still left-invariant. 
Indeed, for $\psi,\vf\in C^\infty_c(G)$ and $g\in G$, since $L_g^\star$ is unitary on $L^2(G)$,
we have
$$
\langle L^\star_g {\widetilde A}^*\psi,\vf\rangle=\langle  \psi,{\widetilde A}L^\star_{g^{-1}}\vf\rangle=
\langle  \psi,L^\star_{g^{-1}}{\widetilde A}\vf\rangle=\langle  {\widetilde A}^*L^\star_g\psi,\vf\rangle\;.
$$
Moreover, we make the following requirement of compatibility of the adjoint map on 
$L^2(G)$ with respect to the ordered decomposition \eqref{DECOMP}.  
Namely, denoting for every $n\in\{0,...,N\}$:
\begin{align}
\label{Vn}
V^{(n)}:=\bigoplus_{k=0}^nV_k\;,
\end{align}
we assume:
\begin{equation}
\label{adjoint}
\forall\, n=0,\dots,N,\quad\forall \,A\in
\CU(V_n)\,,\quad\exists\, B\in \CU(V^{(n)})\quad\mbox{such that}\quad
\widetilde A^*=\widetilde B\;,
\end{equation}
where the space $\CU(V^{(n)})$ is  the subalgebra of $\CU(\g)$ generated by $V^{(n)}$,
as defined in \eqref{CUVE}. We now pass to regularity assumptions regarding the function $S$. 

\begin{dfn}\label{TEMPADM} Set 
\begin{equation}
\label{bE}
\bE:=\exp\{iS\}\;.
\end{equation}
 A tempered pair $(G,S)$ is called {\bf admissible}, if there exists a decomposition (\ref{DECOMP})
with associated coordinate system (\ref{COORDINATES}), such that
for every $n=0,\dots,N$, there
exists an element $X_n\in\CU(V_n)\subset\CU(\g)$ whose associated multiplier $\alpha_n$, defined as 
\begin{equation*}
\widetilde{X}_n\bE=:\alpha_n\bE\;,
\end{equation*}
satisfies the following properties:
\begin{enumerate}
\item[(i)] There exist $C_n>0$ and $\rho_n>0$ such that:
$$
\left|\alpha_n\right|\;\geq\;C_n\big(1\,+\,\left|x_n\right|_n^{\rho_n}\big)\;,
$$
where $x_n:=(x_n^{j})_{j=1,\dots,{\rm dim}(V_n)}$.
\item[(ii)] For all $n=0,\dots, N$, there exists a tempered function $\mu_n\in C^\infty(G,\R^*_+)$ such that:

(ii.1) For every $A\in
\CU(V^{(n)})\subset
\CU(\g)$ there exists $C_A>0$ such that:
\begin{equation*}
\big|\widetilde{A}\,{\alpha_n}\big|\leq C_A\,\left|\alpha_n\right|\,\mu_n\;.
\end{equation*}

(ii.2) n $\mu_n$ is independent of the variables $\{x^{j}_r\}_{j=1,\dots,\dim(V_r)}$, for all $r\leq n$:
\begin{equation*}
\frac{\partial \mu_n}{\partial{x^{j}_r}}=0\,,\quad\forall r\leq n\,,\,\,\,\forall j=1,\dots,\dim(V_r)\;.
\end{equation*}
\end{enumerate}
\end{dfn}

We now give a sequence of Lemmas allowing the construction of  a sequence of
differential operator
$\{\bD_j\}_{j\in\N}$ such that,
for $\underline\mu$ a family of tempered weights, $\bD_j$ sends continuously $\CB^{\mu_j}(G,\CE_j)$
to $L^1(G,\CE_j)$, where $\CE_j$ is the semi-normed space $(\CE,\|.\|_j)$ and which is such that
$\bD_j^*\bE=\bE$. Thus is all that we need to construct our oscillatory integral.
We start with a preliminary result, which gives an upper bound for powers of derivatives of the 
inverse of a multiplier,
in the context of admissible tempered pairs. 
\begin{lem}
\label{POWER} 
Fix $n=0,\dots,N$.
Let $\alpha\in C^\infty(G)$ be non-vanishing and  $1\leq \mu\in C^\infty(G)$ such that for every $A\in
\CU(V^{(n)})$ there exists $C>0$ with $|\widetilde A \alpha|\leq C \,\mu\,|\alpha|$. Fixing 
$X\in\CU(V^{(n)})$, a monomial of homogeneous degree $M\in\N$, we consider the 
differential operator
$$
D_{X,\alpha}: C^\infty(G)\to C^\infty(G)\,,\quad \Phi\mapsto\widetilde{X}\Big(\frac{\Phi}{\alpha}
\Big)\;.
$$
Then, for every $r\in\N$, there exist an element $X'\in\CU(\g)$ of maximal homogeneous degree 
bounded by $rM$ and a constant $C>0$  such that for every $\Phi\in C^\infty(G)$ we have:
$$
\left|D_{X,\alpha}^r\Phi\right|\leq C\,\frac{\mu^{r^2M}}{|\alpha|^r}\,\big|\widetilde{X}'\Phi\big|\;.
$$
\end{lem}
\begin{proof}
We start by recalling Fa\`a di Bruno's formula:
$$
\frac{{\rm d}^r}{{\rm d}t^r}\Big(\frac{1}{\bf f}\Big)=\frac{1}{\bf f}\,\sum_{\vec{M}}\,C^r_{\vec{M}}\,
\prod_{j=1}^r\Big(\frac{{\bf f}^{(j)}}{{\bf f}}\Big)^{M_j}\,,\quad{\bf f}\in C^\infty(\R)\;,
$$
where $\vec{M}=(M_1,\dots,M_r)$ runs along partitions of $r$  (i.e$.$ $r=\sum_{j=1}^rjM_j$) and 
where $C^r_{\vec{M}}$
is some combinatorial coefficient.
 Within Sweedler's notations \eqref{Sweedler}, Fa\`a  di Bruno formula then yields for $\Phi\in  C^\infty(G)$:
$$
D_{X,\alpha}\Phi=\sum_{(X)}\Big(\widetilde{X}_{(1)}\frac{1}{\alpha}\Big)\;\Big(\widetilde{X}_{(2)}
\Phi\Big)=\sum_{(X)}\frac{1}{\alpha}\sum\prod\Big(\frac{\widetilde{X}_j\alpha}{\alpha}\Big)^{M_j}\;
\Big(\widetilde X_{(2)}\Phi\Big)\;,
$$
where the second sum and product run over partitions of $M_{(1)}:=\deg( X_{(1)})\leq M$ and 
where the element $X_j$ is of homogeneous degree
$j=1,\dots,r$  and contains the  combinatorial
coefficients of Fa\`a di Bruno's formula. Of course, we also have that $X_{(1)}$, $X_{(2)}$ and $X_j$ all belong to $
\CU(V^{(n)})$. Thus, $|\widetilde{X}_j\alpha|\leq C(X) \,\mu\,|\alpha|$, the estimation is 
satisfied for $r=1$. 
 For $r=2$, we observe:
\begin{align*}
D_{X,\alpha}^2\Phi&=\sum_{(X)}\Big(\widetilde{X}_{(1)}\frac{1}{\alpha}\Big)\;\widetilde{X}_{(2)}
\Big(\sum_{(X)}\widetilde{X}_{(1)}\frac{1}{\alpha}\;\widetilde{X}_{(2)}\Phi\Big)\\
&=\sum_{(X),(X_{(2)})}\Big(\widetilde{X}_{(1)}\frac{1}{\alpha}\Big)\;\sum_{(X)}
\Big(\widetilde{X}_{(21)}\widetilde{X}_{(1)}\frac{1}{\alpha}\Big)\;\Big(\widetilde{X}_{(22)}
\widetilde{X}_{(2)}\Phi\Big)\;.
\end{align*}
Fa\`a di Bruno' s formula for $\frac{1}{\alpha}$ then yields the assertion for $r=2$. 
Iterating this procedure, we get that
\begin{equation}\label{ITER}
D_{X,\alpha}^r\Phi=\sum_{(X)}\prod_{j=1}^r\Big(\widetilde{X}^{(j)}\frac{1}{\alpha}\Big)\;
\Big(\widetilde{X}' \Phi\Big)\;,
\end{equation}
for some elements $X^{(j)},X'\in\CU(V^{(n)})$ where the maximal homogeneous 
degree of $X^{(j)}$ is bounded by $jM$. (We have absorbed the combinatorial
coefficients of Fa\`a di Bruno' s formula in $X^{(j)}$ and $X'$.)
Therefore,  Fa\`a di Bruno's formula yields for every $j=1,\dots,r$:
\begin{align*}
\Big|\widetilde{X}^{(j)}\frac{1}{\alpha}\Big|=\Big|\frac{1}{\alpha}\sum\prod
\Big(\frac{\widetilde{X}^{(j)}_k.\alpha}{\alpha}
\Big)^{\deg(X^{(j)}_k )}\Big|&\leq C\,\frac{\mu^{\sum_k\deg( X^{(j)} _k)}}{|\alpha|}\\
&\leq
C\,\frac{\mu^{\deg( X^{(j)} )}}{|\alpha|}
\leq C\,\frac{\mu^{jM}}{|\alpha|}\;.
\end{align*}
Therefore since $\frac{r(r+1)}{2}\leq r^2$ and $\mu\geq 1$, we get the (rough) estimation:
$$
\left|D_{X,\alpha}^r\Phi\right|\leq C\sum_{(X)}\prod_{j=1}^r\frac{1}{|\alpha|}\,\mu^{jM}\;\big|
\widetilde{X}'\,\Phi\big|
\leq C'\sum_{(X)}\frac{1}{|\alpha|^r}\,\mu^{r^2M}\;\big|\widetilde{X}'\,\Phi\big|\;,
$$
which delivers the proof.
\end{proof}

 We now fix an admissible tempered pair $(G,S)$ and for all $n=0,\dots, N$,
 we let $X_n\in\CU(V_n)$ as given in Definition \ref{TEMPADM} 
and we  let $\alpha_n,\mu_n\in C^\infty(G)$ be the associated multiplier and tempered function. 
Accordingly to the previous notations, we introduce the differential operators:
\begin{align}
\label{Dn}
D_n:=D_{X_n^*,\alpha_n}:C^\infty(G)\to C^\infty(G)\,,\quad\Phi\mapsto \widetilde{X}_n^*\,\Big(\frac{\Phi}{\alpha_n}
\Big)\;.
\end{align}
Recall that by assumption, there exists $Y_n\in \CU(V^{(n)})$ such that 
$\widetilde{X}_n^*=\widetilde{Y}_n$ and thus, we can apply  Lemma \ref{POWER} to these
operators.
For every $r_n\in\N$, accordingly to the expression (\ref{ITER}), we write for $\Phi\in C^\infty(G)$:
$$
D^{r_n}_n\Phi=\sum_{(X_n)}\prod_{j=1}^{r_n}\Big(\widetilde{X}_n^{(j)}\frac{1}{\alpha_n}\Big)\;
\Big(\widetilde{X}_n'\Phi\Big)\;,
$$
where ${X}_n^{(j)}\in \CU(V^{(n)})$ and its  homogeneous degree is bounded by $jM_n$, with $M_n$ the maximal homogeneous degree of $X_n$ and where the one of $X_n'$ is bounded by $r_nM_n$.
Setting 
\begin{equation}
\label{Psin}
\Psi_n=\prod_{j=1}^{r_n}\Big(\widetilde{X}_n^{(j)}\frac{1}{\alpha_n}\Big)\;,
\end{equation}
we then write (abusively since in fact it is a finite sum of such terms):
\begin{equation}
\label{PX}
D^{r_n}_n=:\Psi_n\,\widetilde{X}_n' \;.
\end{equation}
Given a $N+1$-tuple of  integers $\vec  r=(r_0,\dots,r_N)$, we will be led to consider the operator 
\begin{equation}
\bD_{\vec r}:=\bD_{r_0,\dots,r_N}:=D_0^{r_0}\,D_1^{r_1}\,\dots\,D_N^{r_N}\;.
\label{bD}
\end{equation}
Using  the iterated Sweedler's notation  \eqref{itS}, a recursive use of \eqref{PX} yields
\begin{align*}
\bD_{\vec r}=&
\sum_{(X_0')}\sum_{( X_1')}
\dots\sum_{(X_{N-2}')}\sum_{(X_{N-1}')}\\
&\big[\Psi_0\big]\big[(\widetilde X_0')_{(1)}\Psi_1\big]
\big[(\widetilde X_0')_{(21)}(\widetilde X_1')_{(1)}\Psi_2\big]
\big[(\widetilde X_0')_{(221)}
(\widetilde X_1')_{(21)}(\widetilde X_2')_{(1)}\Psi_3\big]\dots\\
&\big[(\widetilde X_0')_{(22\dots21)}(\widetilde X_1')_{(2\dots21)}
\dots(\widetilde X_{N-3}')_{(21)}(\widetilde X_{N-2}')_{(1)}\Psi_{N-1}\big]\\
&
\big[(\widetilde X_0')_{(222\dots21)}(\widetilde X_1')_{(22\dots21)}
\dots (\widetilde X_{N-3}')_{(221)}(\widetilde X_{N-2}')_{(21)}(\widetilde X_{N-1}')_{(1)}
\Psi_{N}\big]\\
&(\widetilde X_0')_{(222\dots2)}\,(\widetilde X_1')_{(22\dots2)}
\dots(\widetilde X_{N-3}')_{(222)}\,(\widetilde X_{N-2}')_{(22)}\,(\widetilde X_{N-1}')_{(2)}
\,\widetilde X_N'
\;.
\end{align*}
Set then, 
\begin{align}
\label{Psi}
\Psi_{1,0}&:=(\widetilde X_0')_{(1)}\Psi_1\\
\Psi_{2,1,0}&:=(\widetilde X_0')_{(21)}(\widetilde X_1')_{(1)}\Psi_2\nonumber\\
\Psi_{3,2,1,0}&:=(\widetilde X_0')_{(221)}
(\widetilde X_1')_{(21)}(\widetilde X_2')_{(1)}\Psi_3\nonumber\\
&\;\;\vdots\nonumber\\
\Psi_{N-1,\dots,1,0}&:=(\widetilde X_0')_{(22\dots21)}(\widetilde X_1')_{(2\dots21)}
\dots(\widetilde X_{N-3}')_{(21)}(\widetilde X_{N-2}')_{(1)}\Psi_{N-1}\nonumber\\
\Psi_{N,N-1,\dots,1,0}&:=(\widetilde X_0')_{(222\dots21)}(\widetilde X_1')_{(22\dots21)}
\dots(\widetilde X_{N-2}')_{(21)}(\widetilde X_{N-1}')_{(1)}
\Psi_{N}\nonumber\;,
\end{align}
and
$$
{X}'_{N,\dots,0}:=(X_0')_{(222\dots2)}\,( X_1')_{(22\dots2)}
\dots( X_{N-3}')_{(222)}\,(X_{N-2}')_{(22)}\,( X_{N-1}')_{(2)}\, X_N'\;,
$$
in terms of which we have (with the same abuse of notations as in \eqref{bD} above):
\begin{equation}
\label{bDPhi}
\bD_{\vec r}=\Psi_0\,\Psi_{1,0}\,\Psi_{2,1,0}\,\Psi_{3,2,1,0}\dots
\Psi_{N-1,\dots,1,0}\,\Psi_{N,\dots,0}\,\widetilde{X}'_{N,N-1,\dots,1,0}\;.
\end{equation}

\begin{lem}
\label{POWER2}
Fix $n=0,\dots,N$ and
let $\alpha\in C^\infty(G)$ and $\mu\in C^\infty(G,\R^*_+)$ satisfying the hypothesis of 
 Lemma \ref{POWER}. For $j=1,\dots,r$ and $r\in\N^*$, fix also $X^{(j)}\in\CU(V^{(n)})$ and define
$$
\Psi:=\prod_{j=1}^r\Big(\widetilde{X}^{(j)}\frac{1}{\alpha}\Big)\;,
$$
where $\deg({X}^{(j)})\leq jM$, for a given $M\in\N^*$. 
Consider a monomial $Y\in\CU(V^{(n)})$, then we have 
 $$ \widetilde{Y}\,\Psi=\sum_{(Y)}\prod_{j=1}^r\Big(\widetilde{Y}^{(j)}\frac{1}{\alpha}\Big)\quad\mbox{\rm with}\quad \deg({Y}^{(j)})\leq jM+\deg(Y)\;,$$
and moreover there exists $C>0$ such that
$$ \big|\widetilde{Y}\,\Psi\big|\leq C\,\frac{\mu^{r^2M+r\deg(Y)}}{|\alpha|^r}\;.$$
\end{lem}
\begin{proof}
The equality is immediate. Regarding the inequality, we first note that by virtue of 
Fa\`a di Bruno's formula, we have
for a finite sum:
$$
\widetilde{Y}^{(j)}\frac{1}{\alpha}=\frac{1}{\alpha}\sum\prod_{k=1}^{\deg(Y^{(j)})}\Big(\frac{\widetilde{Y}^{(j)}_k\alpha}{\alpha}\Big)^{\deg(Y^{(j)}_k)}\;.
$$
Hence
$$
\Big|\widetilde{Y}^{(j)}\frac{1}{\alpha}\Big|\leq C\,\frac{\mu^{\sum_k\deg(Y^{(j)}_k)}}{|\alpha|}\leq C\,\frac{\mu^{\deg(Y^{(j)})}}{|\alpha|}\leq C\,\frac{\mu^{jM+\deg(Y)}}{|\alpha|}\;.
$$
We then conclude as in the proof of Lemma \ref{POWER}.
\end{proof}
 From the lemmas above, we deduce an estimate for the `coefficient functions' appearing in the expression of the differential operator $\bD_{\vec r}$ in \eqref{bDPhi}.
\begin{cor}
\label{ESTIMATION}
Let $(G,S)$ be an admissible tempered pair with decomposition $\g=\bigoplus_{n=0}^NV_n$ 
and accordingly 
to Definition \ref{TEMPADM}, for $n=0,\dots,N$,  we let 
$(X_n,\alpha_n,\mu_n)\in\CU(V_n)\times C^\infty(G)\times C^\infty(G)$ be the associated differential operator, multiplier and 
tempered function. Then, for $k=0,\dots, N$ and $r_k\in\N^*$,  with $\Psi_{k,\dots,0}\in C^\infty(G)$ defined in \eqref{Psi}, we have
$$
\left|\Psi_{k,\dots,0}\right|\leq C_k\,\frac{\mu_k^{r_k^2M_k+r_k\sum_{j=0}^{k-1}r_jM_j}}{|\alpha_k|^{r_k}}\;,
$$
for some finite non-negative constant $C_k$ and where $M_n:=\deg(X_n)$, $n=0,\dots,N$.
\end{cor}
\begin{proof}
Observe that
$$
\Psi_{k,\dots,0}=\prod_{j=0}^{k-1}(\widetilde{X'_j})_{(2\dots21)}\Psi_k\;,
$$
where $\Psi_k$ is defined in \eqref{Psin}. Since  $({X}'_j)_{(2\dots21)}\in\CU(V^{(k)})$
with homogeneous degree of is bounded by $r_jM_j$ for every $j=0,\dots,k-1$, the estimate 
 we need follows from Lemma \ref{POWER2}.
 \end{proof}
 
 We can now state the main technical results of this chapter.

\begin{prop}\label{PROPIP}
Let $(G,S)$ be an admissible tempered pair and let $\mu$ be a tempered weight. Then,  there exists
 $\vec r=(r_0,\dots,r_N)\in\N^{N+1}$ such that  for every element $F\in \CB^\mu(G)$, the function 
 $\bD_{\vec r} F$ belongs to $ L^1(G)$.
More precisely, there exist a finite constant $C>0$ and $K\in\N$ with $K\leq \sum_{k=0}^Nr_kM_k$ 
and $M_k=\deg(X_k)$ (with $X_k\in\CU_{M_k}(\g)$ as given in Definition \ref{TEMPADM}), such that 
for all $F\in \CB^\mu(G)$, we have:
\begin{equation*}
\|\bD_{\vec r} F\|_1\leq C\,\sup_{X\in\,\CU_K(\g)}\,\sup_{g\in G}\Big\{\frac{\big| \widetilde X\,F(x)\big|}
{\mu(g)\,|X|_K}\Big\}=
C\,\|F\|_{K,\mu}\;.
\end{equation*}
\end{prop}
\begin{proof}
By Lemma \ref{fd-tempered},  in the coordinates \eqref{COORDINATES}, the  Radon-Nikodym 
derivative of the left Haar measure on $G$ with respect to the Lebesgue measure on $\g^\star$, is 
bounded by a polynomial in  $\{x_n^{j}\,,\,\,j=1,\dots,{\rm dim}(V_n)\,,\,\,n=0,\dots,N\}$. By the assumption 
of temperedness of the weight $\mu$, the latter is also bounded by a polynomial in the same 
coordinates. Now, observe
from \eqref{bDPhi}, that we have for any $\vec r=(r_1,\dots, r_N)$ and for 
$K=\deg(X'_{N,\dots,0})\leq \sum_{k=0}^Nr_kM_k$:
\begin{align}\label{LDOMO}
|\bD_{\vec r} \,F|
&\leq|\Psi_0|\,|\Psi_{1,0}|\,|\Psi_{2,1,0}|\,\dots\,|\Psi_{N,\dots,0}|\,\big|\widetilde{X}'_{N,\dots,0}\,F\big|
\nonumber\\
&\leq \,|\Psi_0|\,|\Psi_{1,0}|\,|\Psi_{2,1,0}|\,\dots\,|\Psi_{N,\dots,0}|\,\mu\;|X'_{N,\dots,0}|_K\,
\|F\|_{K,\mu}\;.
\end{align}
This will gives the estimate, if we prove that the function in front of $\|F\|_{K,\mu}$ 
in \eqref{LDOMO} is integrable for a suitable choice of $\vec r\in\N^{N+1}$.
We  prove  a stronger result, namely that given $\vec R=(R_0,\dots,R_N)\in\N^{N+1}$, there exists $\vec r=(r_0,\dots,r_N)\in\N^{N+1}$ such that the associated  functions $\Psi_{k,\dots,0}$ (which depend on $\vec r$) satisfy:
$$
|\Psi_0(x)|\,|\Psi_{1,0}(x)|\,|\Psi_{2,1,0}(x)|\,\dots\,|\Psi_{N,\dots,0}(x)|\leq \frac{C}{(1+|x_0|)^{R_0}\dots(1+|x_N|)^{R_N}}\;.
$$
From Corollary \ref{ESTIMATION} and writing $r_k^2M_k+r_k\sum_{j=0}^{k-1}r_jM_j=r_k\sum_{j=0}^{k}r_jM_j$, we obtain the following estimation:
\begin{equation*}
|\Psi_0|\,|\Psi_{1,0}|\,|\Psi_{2,1,0}|\,\dots\,|\Psi_{N,\dots,0}|\leq 
C\,\prod_{k=0}^N\frac{\mu_k^{r_k\sum_{j=0}^{k}r_jM_j}}{\alpha_k^{r_k}}\;.
\end{equation*}
Moreover, by assumption of temperedness, see Definition \ref{TEMPADM} (ii.1), there exist $\rho_0,\dots,\rho_N>0$ such that
\begin{equation*}
|\Psi_0(x)|\,|\Psi_{1,0}(x)|\,|\Psi_{2,1,0}(x)|\,\dots\,|\Psi_{N,\dots,0}(x)|\leq 
C\,\prod_{k=0}^N\frac{\mu_k(x)^{r_k\sum_{j=0}^{k}r_jM_j}}{(1+|x_k|)^{\rho_k r_k}}\;.
\end{equation*}
From the hypothesis of Definition \ref{TEMPADM} (ii.2), we deduce that the element $\mu_N$ is constant. Indicating the variable dependence into parentheses, one also has
\begin{align*}
\mu_{N-1}&=\mu_{N-1}(x_N)\,,\quad
\mu_{N-2}=\mu_{N-2}(x_{N-1},x_N)\,,\quad\dots\\
\mu_{1}&=\mu_{1}(x_2,\dots,x_N)\,,\quad
\mu_{0}=\mu_{0}(x_1,x_2,\dots,x_N)\;.
\end{align*}
 Denoting by $m_n$, $n=0,\dots,N$, the degree of a polynomial function that, in the variables \eqref{COORDINATES}, dominates the tempered function $\mu_n$, we obtain
the sufficient conditions:
$$
r_0\rho_0\;\geq\;R_0\quad\mbox{\rm and}\quad 
\rho_n \,r_n\,-\,\sum_{k=0}^{n-1}\Big(\,m_k\,r_k\,\sum_{j=0}^kr_jM_j\Big)\;\geq\;R_n\,,\quad 
n=1,\dots,N\;,
$$
which are always achievable. 
\end{proof}
Let now $\CE$ be a complex Fr\'echet space, with topology  associated with a countable
 family of semi-norms $\{\|.\|_j\}_{j\in\N}$. An immediate  modification of its proof,  lead us to the 
 following  version of Proposition \ref{PROPIP} (the only difference with the former is that
 now the index $K\in\N$ may depend on $j$ via the order of the tempered weight $\mu_j$).

\begin{prop}\label{PROPIP2}
Let $(G,S)$ be an admissible tempered pair, $\CE$ be a complex Fr\'echet space  and let 
$\underline\mu$ be a
family of tempered weights. Then for all $j\in\N$, there exist
 $\vec r_j\in\N^{N+1}$, $C_j>0$ and $k_j\in\N$, such that  for every element 
 $F\in \CB^{\underline\mu}(G,\CE)$, we have
\begin{equation*}
\int_G\|\bD_{\vec r_j} F(g)\|_j\,{\rm d}_G(g)\leq C_j\,\sup_{X\in\,\CU_{K_j}(\g)}\sup_{g\in G}\Big\{
\frac{\| \widetilde X\,F(x)\|_j}{\mu_j(g)\,|X|_{k_j}}\Big\}=:C_j\,\|F\|_{j,k_j,\underline\mu}\;.
\end{equation*}
\end{prop}

\section{An oscillatory integral for admissible  tempered pairs}
\label{OITP}

Our notion of oscillatory integral
 is an immediate  consequence of  Proposition \ref{PROPIP2} together with the following  
 immediate but essential  
 observation:
\begin{prop}
Let $(G,S)$ be an admissible tempered pair, let $\bE$ be the function on $G$ defined in \eqref{bE} 
and, for $\vec r\in\N^{N+1}$, let  $\bD_{\vec r}$
be  the differential operator  
given in \eqref{bD}. Then, with $\bD_{\vec r}^*$ the formal adjoint of $\bD_{\vec r}$
with respect to the inner product 
of $L^2(G)$,  we have:
$$
\bD_{\vec r}^*\,\bE=\bE\,,\qquad\forall \vec r\in\N^{N+1}\;,
$$
\end{prop}
\begin{proof}
By construction, we have
$$
\bD_{\vec r}=D_0^{r_0}\,D_1^{r_1}\,\dots\,D_N^{r_N}\;,
$$
where the operators $D_n$, $n=0,\dots,N$, are given in \eqref{Dn}.
Hence it suffices to prove the identity above for $\vec r=(0,\dots,0,1,0,\dots,0)$.
But this immediately follows by construction since
$$
D_n^*= \frac{1}{\alpha_n}\widetilde{X}_n
\qquad\mbox{and}\qquad
\widetilde{X}_n\,\bE=\alpha_n\,\bE\;.
$$
This completes the proof.
\end{proof}

\begin{dfn}
\label{OI}
Let $(G,S)$ be an admissible tempered pair,  $\mu$ a tempered weight, $\bm$ an element of 
$\CB^\mu(G)$ and 
$\underline\mu,\underline{\hat\mu}$ two families of tempered weights such that
$\underline\mu\prec\underline{\hat\mu}$ (hence $\underline\mu.\mu\prec\underline{\hat\mu}.\mu$).  
Associated to the family of  weights $\underline{\hat\mu}.\mu$,
let $\vec r_j\in\N^{N+1}$, $j\in\N$,
  as given in Proposition \ref{PROPIP2} and 
let $\bD_{\vec r_j}$ be the differential operators given in \eqref{bD}.
Performing integrations by parts,  the Dunford-Pettis theorem \cite{Gr} 
 yields a $\CB^{\underline{\hat\mu}}(G,\CE)$-continuous mapping
$$
\CD(G,\CE)\to\CE\;,\quad F\mapsto\int_G\bm\bE\,F=\int_G\bE\,\bD_{\vec r_j}\big(\bm\,F\big)\;.
$$
Then by  Lemma \ref{SmoothFamily} (v)
(and Remark \ref{stress2}),  the latter  extends to the following 
continuous  linear mapping:
$$
\widetilde{\int_G \bm\,\bE}:\CB^{\underline{\mu}}(G,\CE)\to \CE\;,
$$
that we refer to as an {\bf oscillatory integral}.
\end{dfn}

Our next aim is to prove that the oscillatory integral on $\CB^{\underline{\mu}}(G,\CE)$,
does not depend on the choices made.  
So let $\mu$, $\underline{\mu}$ and  $\underline{\hat\mu}$ as in Definition \ref{OI}.  Fix 
$F\in \CB^{\underline{\mu}}(G,\CE)$ and chose a sequence $\{F_n\}_{n\in\N}$ of elements of 
$\CD(G,\CE)$ converging to $F$ for the topology of $\CB^{\underline{\hat\mu}}(G,\CE)$. By definition
of the oscillatory integral and undoing the integrations by parts at the level of smooth
compactly supported $\CE$-valued functions, we first observe that by continuity, we have:
\begin{align*}
\widetilde{\int_G \bm\,\bE}\big(F\big)&=\widetilde{\int_G \bm\,\bE}\big(\lim_{n\to\infty}F_n\big)
\\&=\lim_{n\to\infty}\widetilde{\int_G \bm\,\bE}\big(F_n\big)=
\lim_{n\to\infty}\int_G \bm(g)\,\bE(g)\,F_n(g)\,{\rm d}_G(g)\;,
\end{align*}
where the first limit is in $\CB^{\underline{\hat\mu}}(G,\CE)$  and the last two  are in $\CE$.
Then,  the estimate of Proposition \ref{PROPIP2} immediately implies that the limit
above is independent of the approximation sequence $\{F_n\}_{n\in\N}$ chosen. 
This shows that the oscillatory integral 
does not depend on the differential operators in $\bD_{\vec r_j}$ used to define
the extension (in the topology of $\CB^{\underline{\hat\mu}}(G,\CE)$) of the oscillatory integral from $\CD(G,\CE)$ to $\CB^{\underline{\mu}}(G,\CE)$.
Last,  to see that 
the oscillatory integral mapping is also independent of the choice of the family of dominant 
weights $\underline{\hat\mu}$ chosen, it suffices to  
remark that the approximation sequence constructed
in the proof of Lemma \ref{symbols} (viii) can be used for 
any family $\underline{\hat\mu}$ such that  $\mu_j\prec\hat\mu_j$. 
Of course this will hold provided that we can always find dominant weights. This is certainly the
case if there exists a weight dominating the constant weight $1$.
Thus we have proved:
\begin{prop}
\label{uniqueness}
Let $(G,S)$ be an admissible tempered pair, $\CE$ be a complex Fr\'echet space,
$\mu$ be tempered weights, $\bm\in\CB^\mu(G)$ and $\underline\mu$ be a family of tempered
weights. Assuming that there exists
a tempered weight $\mu_c$ which dominates the constant weight $1$,  then  the oscillatory
integral mapping 
$$
\widetilde{\int_G \bm\,\bE}:\CB^{\underline{\mu}}(G,\CE)\to \CE\;,
$$
does not depend on the choice of the integers $\vec r_j\in\N^{N+1}$ and on the family
of dominant weights
$\underline{\hat\mu}$
 given in Definition \ref{OI}. Moreover, given $F\in \CB^{\underline{\mu}}(G,\CE)$, we have
\begin{equation*}
\widetilde{\int_G \bm\,\bE}\big(F\big)=
\lim_{n\to\infty}\int_G \bm(g)\,\bE(g)\,F_n(g)\,{\rm d}_G(g)\;,
\end{equation*}
where $\{F_n\}_{n\in\N}$ is an arbitrary sequence in $\CD(G,\CE)$, converging to $F$ in the 
topology of $\CB^{\underline{\hat\mu}}(G,\CE)$, for an arbitrary sequence of weights 
$\underline{\hat\mu}$, which 
dominates $\underline{\mu}$.
\end{prop}
\begin{rmk}
 Note that Proposition \ref{uniqueness} does not assert   that the oscillatory integral on 
$\CB^{\underline{\mu}}(G,\CE)$ is the unique continuous extension of its restriction to 
$\CD(G,\CE)$.  
\end{rmk}
We observe that the existence of a tempered weight that dominates the constant weight
$1$ implies that every weight is dominated, which
is crucial for the construction of the oscillatory integral, as observed above. 
This leads us to introduce the notion of tameness below.
We will see that this property holds
for negatively curved K\"ahlerian groups, where we can simply use  the modular weight
(see Corollary \ref{MWP} and Lemma \ref{mwdi}). 

\begin{dfn}
\label{tame}
A tempered Lie group $G$, with associated diffeomorphism $\phi:G\to\g^\star$ is called {\bf tame} 
 if there exist an Euclidean norm $|.|$ on $\g^\star$, a tempered weight $\mu_\phi$ and two positive 
 constants $C,\rho$ such that 
$$
C\left(1+|\phi|^2\right)^{\frac{\rho}{2}}\leq\mu_\phi\;.
$$
\end{dfn}
\begin{rmk}
Observe that when a tempered Lie group $G$ is tame, then there exists a tempered weight that
dominates the constant weight $1$. Indeed, $1\prec\mu_\phi$. This is the main reason for 
introducing this notion.
\end{rmk}

In the constant family case, i.e$.$ for $\CB^\mu(G,\CE)$, we can express the oscillatory integral as 
an absolutely convergent one for each semi-norm $\|.\|_j$:
$$
\widetilde{\int_G \bm\,\bE}\big(F\big)=\int_G\bE\,\bD_{\vec r}\big(\bm\,F\big)\,,
\qquad \forall F\in \CB^\mu(G,\CE)\;,
$$
where the label $\vec r\in\N^{N+1}$ of the differential operator $\bD_{\vec r}$ is given by 
Proposition \ref{PROPIP} and does not depend on the index $j\in\N$ which labels the semi-norms
defining the topology of $\CE$ nor on the element $F\in\CB^\mu(G,\CE)$.
However, we cannot  access to such a formula  in the case of $\CB^{\underline{\mu}}(G,\CE)$,
since in this case the label $\vec r\in\N^{N+1}$  depends on the index $j\in\N$.
However, and from Proposition \ref{PROPIP2}, we have the following weaker statement:

\begin{prop}
\label{weak-conv}
Let $(G,S)$ an admissible and tame tempered pair, $\CE$ a  Fr\'echet space,
$\mu,\underline{\mu}$ be tempered weights and $\bm\in\CB^\mu(G)$. Then for every $j\in\N$, there exists
$\vec r_j\in\N^{N+1}$, such that for all $F\in \CB^{\underline\mu}(G,\CE)$, we have
$$
\widetilde{\int_G \bm\,\bE}\big(F\big)=\int_G\bE\,\bD_{\vec r_j}\big(\bm\,F\big)\;,
$$
where the right hand side is an absolutely convergent integral for the semi-norm $\|.\|_j$.
\end{prop}

 We close this section with a natural result  on the compatibility of the
oscillatory integral with continuous linear maps between Fr\'echet spaces.
\begin{lem}
\label{continuous-maps}
Let $(G,S)$ be an admissible and tame tempered pair, $(\CE,\{\|.\|_j\})$ and $(\CF,\{\|.\|_j'\})$ two Fr\'echet spaces and 
 $T:\CE\to\CF$ a continuous linear map.
Define the map $\hat T$ from 
$C(G,\CE)$ to $C(G,\CF)$, by setting
$$
\big(\hat T F\big)(g):= T\big(F(g)\big)\;.
$$
Then, for any family $\underline\mu$ of tempered weights, there exists another family of tempered
weights $\underline\nu$, such that  $\hat T$ 
is continuous from $\CB^{\underline\mu}(G,\CE)$
 to $\CB^{\underline\nu}(G,\CF)$.
Moreover for any $\bm\in\CB^\mu(G)$, with 
$\mu$ another tempered weight, we have
\begin{align}
\label{mapT}
T\,\Big(\widetilde{\int_G \bm\,\bE}\big(F\big)\Big)=\widetilde{\int_G \bm\,\bE}\big(\hat T F\big)\;.
\end{align}
\end{lem}
\begin{proof}
  By continuity of $T$, 
 for all $j\in\N$ there exist $l(j)\in\N$ 
and $C_j>0$, such that for all $a\in\CE$, we have $\|T(a)\|_j'\leq C_j\|a\|_{l(j)}$. 
This immediately implies the continuity of $\hat T$. Indeed, setting $\nu_j:=\mu_{l(j)}$, then for 
$F\in \CB^{\underline\nu}(G,\CE)$ and $j,k\in\N$, we have
\begin{align*}
\|\hat T \,F\|_{j,k,\underline\nu}'&=\sup_{X\in\,\CU_k(\g)}\sup_{g\in G}\frac{\|\widetilde X T\big(F(g)\big)\|_j'}
{\mu_{l(j)}|X|_k}
=\sup_{X\in\,\CU_k(\g)}\sup_{g\in G}\frac{\| T\big(\widetilde XF(g)\big)\|_j'}{\mu_{l(j)}|X|_k}\\
&\leq C_j
\sup_{X\in\,\CU_k(\g)}\sup_{g\in G}\frac{\| \big(\widetilde XF(g)\big)\|_{l(j)}}{\mu_{l(j)}|X|_k}=C_j\|F\|_{l(j),k,\underline\mu}\;.
\end{align*}
Repeating the arguments  for $ \CB^{\underline\mu.\mu_\phi}(G,\CE)$ ($\mu_\phi$ is the tempered
weight associated to tameness)   instead of $ \CB^{\underline\mu}(G,\CE)$, we see that both sides of 
\eqref{mapT} define continuous linear maps from $ \CB^{\underline\mu.\mu_\phi}(G,\CE)$ to $\CF$. 
Moreover, it is easy to see that they coincide on $\CD(G,\CE)$ and thus they coincide on the closure 
of $\CD(G,\CE)$
inside $ \CB^{\underline\mu.\mu_\phi}(G,\CE)$, which contains  $ \CB^{\underline\mu}(G,\CE)$
 by Lemma \ref{SmoothFamily} (vi), as 
$\underline\mu.\mu_\phi$ dominates $\underline\mu$.
\end{proof}

\section{A Fubini Theorem for semi-direct products}\label{FTSDP}

 The aim of this section is to prove a Fubini type result for  the oscillatory integral 
on a semi-direct product of tempered pairs. We start with following observation:
\begin{lem}
\label{Fubini-basis}
Let $\CE$  be a  Fr\'echet space,
$G_1$, $G_2$ be two Lie groups with Lie algebras $\g_1$, $\g_2$
 and $\bR\in{\rm Hom}\big(G_1,\Aut(G_2)\big)$ be an 
extension homomorphism. Consider $\underline\mu$, a family of weights on the semi-direct product
$G_1\ltimes_\bR G_2$, with  sub-multiplicative degree $(\underline L,\underline R)$. 
Set also $\underline\mu_{1}$ and
$\underline\mu_{2}$ for its restrictions  to the subgroups  $G_1$ and $G_2$
and let $\fd_1$ be the restriction of the modular weight (cf$.$ Definition \ref{MWS}) of 
$G_1\ltimes_\bR G_2$ to $G_1$. Then the map  
\begin{align}
\label{hat}
C^\infty(G_1\ltimes_\bR G_2,\CE)&\to C^\infty\big(G_1,C^\infty( G_2,\CE)\big)\;,\nonumber\\
F&\mapsto\hat F:=\big[g_1\in G_1
\mapsto [g_2\in G_2\mapsto F(g_2g_1)]\big]
\;,
\end{align}
sends continuously $\CB^{\underline\mu}(G_1\ltimes_\bR G_2,\CE)$ to 
$\CB^{\underline\nu}\big(G_1,\CB^{\underline\lambda}(G_2,\CE)\big)$, where
$\underline\nu:=\{\nu_{j,k}\}_{j,k\in\N}$
and 
$\underline\lambda:=\{\lambda_j\}_{j\in\N}$ with 
$$
\nu_{j,k}:=\mu_{1,j}^{R_j}\,\fd_1^k\,,\qquad \lambda_j:=\mu_{2,j}^{L_j}\;.
$$
\end{lem} 
\begin{proof}
First, observe that
for $g\in G_1\ltimes_\bR G_2$ with $g=g_2g_1$, $g_1\in G_1$, $g_2\in G_2$, 
$F\in C^\infty(G_1\ltimes_\bR G_2)$ and for $X^1\in \g_1$, $X^2\in\g_2$, we have
\begin{align}
\label{truff}
\widetilde{X^1}_{.g_1}\,\hat F(g_1,g_2)=\widetilde{X^1}_{.g}\,F(g)\,,\quad 
\widetilde{X^2}_{.g_2}\,\hat F(g_1,g_2)=\widetilde{\bR_{g_1^{-1}} (X^2)}_{.g}\,F(g)\;,
\end{align}
where we use the same notation for the extension homomorphism and its derivative:
$$
\g_2\to \g_2\,,\quad X\mapsto\frac d{dt}\Big|_{t=0} \bR_{g_1}(e^{tX})\,,\quad g_1\in G_1\;.
$$
From this, it follows that the restriction of a weight on $G_1\ltimes_\bR G_2$ to $G_1$ or $G_2$
is still a weight on $G_1$ or $G_2$. Indeed, given $\mu$ a weight on $G_1\ltimes_\bR G_2$, call
$\mu^ i$, $ i=1,2$, its restriction to the subgroup $G_ i$ and given $X\in\CU(\g_ i)$
call $X_ i$ its image in $\CU(\g_1\ltimes \g_2)$. 
Then, Equation \eqref{truff} yields $\widetilde X\mu^ i
=\big(\widetilde X_ i \mu\big)^ i$, $ i=1,2$, which together with $(\mu^{\!\vee})^ i=(\mu^ i)^{\!\vee}$ (where 
$\mu^{\!\vee}\!(g):=\mu(g^{-1})$)
 implies that the first condition of 
Definition \ref{weight-def} is satisfied. Sub-multiplicativity at the level of each subgroups $G_ i$,
$ i=1,2$, follows from sub-multiplicativity at the level of $G_1\ltimes_\bR G_2$  (with the same
sub-multiplicativity degree).

 Moreover, \eqref{truff} also implies that for 
$F\in\CB^{\underline\mu}(G_1\ltimes_\bR G_2,\CE)$, $X^1\in \CU_{k_1}(\g_1)$, 
$X^2\in\CU_{k_2}(\g_2)$ and 
$k_1,k_2,j\in\N$, we have for $g_2g_1\in G_1\ltimes_\bR G_2$:
\begin{align*}
&\|\widetilde{X^1}_{.g_1}\,\widetilde{X^2}_{.g_2}\,\hat F(g_1,g_2)\|_j
=\|\big(\widetilde{X^1}\widetilde{\bR_{g_1^{-1}} (X^2)}\,
F\big)(g_2g_1)\|_j\\
&\leq C(k_1,k_2)\, |X^1|_{k_1}\,|X^2|_{k_2}\,|\bR_{g_1^{-1}}|_{k_2}\,\sup_{Y\in\,\CU_{k_1+k_2}
(\g_1\ltimes \g_2)}
\frac{\|\widetilde Y F(g_2g_1)\|_j}{|Y|_{k_1+k_2}}\\
&\leq C(k_1,k_2)\,|X^1|_{k_1}\,|X^2|_{k_2}\,|\bR_{g_1^{-1}}|_{k_2}\,\mu_j(g_2g_1)\,
\|F\|_{j,k_1+k_2,\underline\mu}\\
&\leq C'(k_1,k_2)\, |X^1|_{k_1}\,|X^2|_{k_2}\,\fd_1(g_1)^{k_2}\,\mu_j(g_1)^{R_j}\,\mu_j(g_2)^{L_j}\,
\|F\|_{j,k_1+k_2,\underline\mu}\;,
\end{align*}
by Lemma \ref{lemtwo}, since for $g_1\in G_1$, $\bR_{g_1}$ coincides with the restriction of 
$\Ad_{g_1}$ to $\g_2$. Thus, labeling  by $(j,k_2)\in\N^2$ the semi-norms 
$\|.\|_{j,k_2,\underline\lambda}$ of $\CB^{\underline\lambda}(G_2,\CE)$, we finally get:
$$
\|\hat F\|_{(j,k_2),k_1,\underline\nu}\leq C'(k_1,k_2)\, \|F\|_{j,k_1+k_2,\underline\mu}\;,
$$
which completes the proof.
\end{proof}

Now, assume that the groups $G_1$ and $G_2$
come from admissible tempered pairs $(G_1,S_1)$ and $(G_2,S_2)$. Parametrizing $g=g_2g_1\in
G_1\ltimes_\bR G_2$ with $g_i\in G_i$, $i=1,2$, we can then set
\begin{equation}
\label{double-phase}
S:G_1\ltimes_\bR G_2\to \R\,,\quad g_2g_1\mapsto S_1(g_1)+S_2(g_2)\;,
\end{equation}
and using the notation \eqref{bE},  we set accordingly
$$
\bE(g_2g_1):=\exp\{iS(g_2g_1)\}=\bE_1(g_1)\,\bE_2(g_2)\;.
$$
For $g_1\in G_1$, we let $|\bR_{g_1}|$ be the operator norm of $\bR_{g_1}\in\End(\g_2)$.
Assume further that the map $G_1\to\R^+$, $g_1\mapsto|\bR_{g_1}|$ is tempered.
Then, by Lemma \ref{fd-semi}  and Lemma \ref{fd-tempered}, we deduce that
   $\fd_1$, the restriction of the modular weight on $G_1\ltimes_\bR G_2$ to 
$G_1$, is  tempered too. 
Thus for $\bm\in\CB^\mu(G_1\ltimes_\bR G_2)$ with $\mu$ a tempered weight,
and  with 
 $\hat \bm$ the associated function on $G_1\times G_2$ as constructed in \eqref{hat},  Lemma 
 \ref{Fubini-basis} shows that the map
$$
\CB^{\underline\mu}(G_1\ltimes_\bR G_2,\CE)\to\CE\;,\quad F\mapsto \widetilde{\int_{G_2}\bE_2}\,
\Big(\widetilde{ \int_{G_1}\bE_1}\,\big(\hat\bm\,\hat F\big)\Big)\;,
$$
is well defined as a continuous linear map. Thus under these circumstances, this map could be used as a definition for the oscillatory integral on the semi-direct product $G_1\ltimes_\bR G_2$.
Moreover, when the pair $\big(G_1\ltimes_\bR G_2,S\big)$ is also tempered and admissible
and when the extension homomorphism preserves the Haar measure ${\rm d}_{G_2}$, then 
the map above coincides with the oscillatory integral on  $G_1\ltimes_\bR G_2$, 
as given in Definition \ref{OI}.
This is our Fubini-type result in the context of semi-direct product of tempered pairs:

\begin{prop}
\label{Fubini}
Within the context of Lemma \ref{Fubini-basis}, assume further that the groups $G_1$ and $G_2$
come from admissible and tame
 tempered pairs $(G_1,S_1)$ and $(G_2,S_2)$ and, with $S$  defined in
\eqref{double-phase}, that $\big(G_1\ltimes_\bR G_2,S\big)$ is admissible, tame  and tempered too.
 Assume last 
 that the extension homomorphism  
 $$
 \bR\in{\rm Hom}\big(G_1,\Aut(G_2)\big)\;,
 $$ 
 preserves the Haar measure ${\rm d}_{G_2}$ and be such that  the map 
 $[g_1\mapsto|\bR_{g_1}|]$ is tempered on $(G_1,S_1)$.
 Let also $\mu,\underline\mu$ 
 be tempered weights on the semi-direct product $G_1\ltimes_\bR G_2$. Then, for
 $F\in\CB^{\underline\mu}(G_1\ltimes_\bR G_2,\CE)$, $\bm\in\CB^\mu(G_1\ltimes_\bR G_2)$, with 
 $\hat F$ and $\hat \bm$ the associated functions on $G_1\times G_2$ as in \eqref{hat},
 we have
 \begin{equation}
 \label{eq:Fubini}
 \widetilde{\int_{G_1\ltimes_\bR G_2}\bE\bm}\,\big(F\big)=\widetilde{\int_{G_2}\bE_2}\,\Big(\widetilde{
 \int_{G_1}\bE_1}\,\big(\hat\bm\,\hat F\big)\Big)\;.
 \end{equation}
\end{prop}
\begin{proof}
Since the map $[g_1\mapsto|\bR_{g_1}|]$ is tempered on $(G_1,S_1)$, 
by Lemma \ref{fd-semi}  and Lemma \ref{fd-tempered} one deduces
that  $\fd_1$, the restriction of the modular weight on $G_1\ltimes_\bR G_2$ to 
$G_1$ is  tempered on $(G_1,S_1)$.
  Thus by Lemma \ref{Fubini-basis}, the right hand-side
of \eqref{eq:Fubini} is well defined as a continuous linear map from 
$\CB^{\underline\mu.\mu_\phi}(G_1\ltimes_\bR G_2,\CE)$ to 
$\CE$
($\mu_\phi$ is the tempered weight on $G_1\ltimes_\bR G_2$ associated with tameness). 
 Note also that, by our assumptions, the pair $(G_1\ltimes_\bR G_2,S)$ is tempered and
 admissible, the left hand side of \eqref{eq:Fubini} is also well 
defined as a continuous linear map from $\CB^{\underline\mu.\mu_\phi}(G_1\ltimes_\bR G_2,\CE)$ to 
$\CE$, too.
Now, take $F\in\CD(G_1\ltimes_\bR G_2,\CE)$ and associate to it 
$\hat F\in\CD\big(G_1,\CD( G_2,\CE)\big)$
as in \eqref{hat}. By construction, we have
\begin{align*}
 \widetilde{\int_{G_1\ltimes_\bR G_2}\bE\bm}\,( F)=\int_{G_1\ltimes_\bR G_2}
 \bE(g)\,\bm(g)\,F(g)\,{\rm d}_{G_1\ltimes_\bR G_2}(g)\;.
\end{align*}
Since the extension homomorphism $\bR$ preserves ${\rm d}_{G_2}$, we have for $g_1\in G_1$,
$g_2\in G_2$:
$$
  {\rm d}_{G_1\ltimes_\bR G_2}(g_2g_1)={\rm d}_{G_1}(g_1)\,{\rm d}_{G_2}(g_2)\;,
$$
which, by the ordinary Fubini Theorem,  implies that
\begin{align*}
\widetilde{\int_{G_1\ltimes_\bR G_2}\bE\bm}\,( F)&=\int_{ G_2} \bE_2(g_2)\,
\Big(\int_{G_1} \bE_1(g_1)
 \,\bm(g_2g_1)\,F(g_2g_1)\,{\rm d}_{G_1}(g_1)\Big)\,{\rm d}_{G_2}(g_2)\\
&=\widetilde{\int_{G_2}\bE_2}\,\Big(\widetilde{
 \int_{G_1}\bE_1}\,(\hat\bm\,\hat F)\Big)\;.
\end{align*}
Thus, both sides
of \eqref{eq:Fubini} are continuous linear maps from 
$\CB^{\underline\mu.\mu_\phi}(G_1\ltimes_\bR G_2,\CE)$ 
to $\CE$  and coincide on $\CD(G_1\ltimes_\bR G_2,\CE)$. Therefore, these maps coincide on the 
closure of
$\CD(G_1\ltimes_\bR G_2,\CE)$ inside $\CB^{\underline\mu.\mu_\phi}(G_1\ltimes_\bR G_2,\CE)$. 
One concludes using
Lemma \ref{SmoothFamily} (vi), which shows that the latter closure contains
$\CB^{\underline\mu}(G_1\ltimes_\bR G_2,\CE)$.
\end{proof}

\section{A Schwartz space for tempered pairs}
\label{sec:Schwartz}

In  this section, we introduce
a Schwartz type functions space, out of an admissible tempered pair $(G,S)$ and prove that it is 
Fr\'echet and nuclear. Our notion of Schwartz space is of course closely related, if not in many cases 
equivalent, to 
 other notions of Schwartz space on Lie groups, but the point here is that it is formulated in terms of 
 the phase function
$S$ only. This is this formulation that allows to immediately implement the compatibility with our 
notion of oscillatory integral.

\begin{dfn}
\label{SCHH}
Let $(G,S)$ be a tempered pair. For all $X\in \CU(\g)$, we let 
$\alpha_X:=\bE^{-1}\,\widetilde X\,\bE\in C^\infty(G)$, where $\bE$ is defined in \eqref{bE}. 
Then we set
$$
\CS^S(G):=\big\{f\in C^\infty(G)\;:\; \forall\,X,Y\in\CU(\g)\,,\;\forall\,n\in\N\,,\;\sup_{x\in G}
\big|\alpha_X^n(x)\,\big(\widetilde Y\,f\big)(x)\big|<\infty\big\}\;. 
$$
\end{dfn}
We first prove that this space is isomorphic to the ordinary Schwartz space of the Euclidean space 
$\g^\star$.
\begin{lem}
Let $\phi:G\to \g^\star$ be the diffeomorphism underlying Definition \ref{TEMPPAIR}, associated to  
an admissible tempered pair $(G,S)$. Fixing an Euclidean structure
on $\g^\star$,  denote by $\CS(\g^\star)$ the ordinary Schwartz space of $\g^\star$. 
Then, $\CS^S(G)$ coincides
with
$$
\CS^\phi(G):=\big\{f\in C^\infty (G)\;:\;f\circ\phi^{-1}\in\CS(\g^\star)\big\}\;.
$$
 In particular, endowed with the transported topology, 
$\CS^S(G)$ is a nuclear Fr\'echet space. If moreover the pair $(G,S)$ is tame,
then the transported topology is equivalent to the one associated with the semi-norms
\begin{align*}
\|.\|_{k,n}\;:\;f\in\CS^S(G,\CE)\mapsto\sup_{X\in\,\CU_k(\g)}\sup_{x\in G}\Big\{\frac
{\mu_\phi(x)^n\big|\widetilde X\,f(x)
\big|}{|X|_k}\Big\}\;,\qquad k,n\in\N\;,
\end{align*}
\label{SCH}
\end{lem}
\begin{proof}
Recall that $f\in \CS^\phi(G)$ if and only if for all $\alpha,\beta\in\N^{\dim (G)}$, we have
\begin{equation}
\label{e1}
\sup_{\xi\in\g^\star}\big|\xi^\alpha\,\partial^\beta(f\circ\phi^{-1})(\xi)\big|<\infty\;,
\end{equation}
while $f\in \CS^S(G)$ if and only if for all $X, Y\in\CU(\g)$ and all $n\in\N$
\begin{equation}
\label{e2}
\sup_{x\in G}\big|\alpha_X^n(x)\,\big(\widetilde Y\,f\big)(x)\big|<\infty\;.
\end{equation}
Fix $\{X_j\}_{j=1}^{\dim(G)}$ a basis of $\g$ and let $\{\xi_j\}_{j=1}^{\dim(G)}$ be the dual basis on $\g^\star$. From
the same methods as in Lemma \ref{fd-tempered}, one can construct an invertible matrix $M(\xi)$ which is tempered with tempered inverse and which is such that in the $\phi$-coordinates
$$
\widetilde X_j =\sum_{i=1}^{\dim(G)} M(\xi)_{j,i}\partial_{\xi_i}\;.
$$
Since by Remark \ref{S-temp}
 $S$ is tempered,  for all $X\in\CU(\g)$, the associated multiplier $\alpha_X=
 \bE^{-1}\,\widetilde{X}\bE$ in $\phi$-coordinates is bounded by a polynomial function on 
 $\g^\star$. Last, since the pair $(G,S)$ is admissible,
associated to the vector space decomposition $\g=\bigoplus_{k=0}^N V_k$, there exist  elements $X_k\in\CU(V_k)$ and constants $C_k,\rho_k>0$ such that with $\alpha_k=\bE^{-1}\,
\widetilde{X_k}\bE$ and with $|.|_k$  the Euclidean norm on $V_k$, we have
$$
|x_n|_n\leq \big(C_n^{-1}|\alpha_n|-1\big)^{1/\rho_n}.
$$
Summing up over $k=1,\dots, N$ gives the existence of $C,\rho>0$, such that
$$
|\xi|\leq C\Big(1+\sum_{k=0}^N\big|\alpha_k\big(\phi^{-1}(\xi)\big)\big|\Big)^\rho\;.
$$
Putting these three facts together gives the equality between the two sets of functions on $G$ and the
equivalence of the topologies associated with the semi-norms \eqref{e1} and \eqref{e2}.

  Last, 
assume that the pair $(G,S)$ is also tame, with associated
weight $\mu_\phi$. Then there exist
constants $C_1,C_2,\rho_1,\rho_2$ such that, with $\phi:G\to \g^\star$ 
the diffeomorphism underlying Definition \ref{TEMPPAIR}, we have
$$
C_1(1+|\phi|^2)^{\rho_1/2}\leq\mu_\phi\leq C_2(1+|\phi|^2)^{\rho_2/2}\;.
$$
which is enough to prove the last claim.
\end{proof}

More generally, when $\CE$ is a complex Fr\'echet space with topology underlying a countable set 
of 
semi-norms $\{\|.\|_j\}_{j\in\N}$, we define the $\CE$-valued Schwartz space associated to an 
admissible tempered and tame pair $(G,S)$
as the Fr\'echet space completion of $\CD(G,\CE)$ for the topology underlying the family
of semi-norms:
\begin{align*}
\|f\|_{j,k,n}:=\sup_{X\in\,\CU_k(\g)}\sup_{x\in G}\Big\{\frac
{\mu_\phi(x)^n\big\|\widetilde X\,f(x)
\big\|_j}{|X|_k}\Big\}\;,\qquad j,k,n\in\N\;,
\end{align*}
Note that by nuclearity of  $\CS^S(G)$, we have $\CS^S(G,\CE)=\CS^S
(G)\hat\otimes\CE$ (for any
completed tensor product).

\section{Bilinear mappings from the oscillatory integral}\label{BMOI}

We now present several results which establish most of the
analytical properties we will need  to construct our universal deformation 
formula for actions of K\"ahlerian groups on Fr\'echet algebras. In all that   follows,
when considering  a Fr\'echet algebra $\CA$ with topology underlying a countable set of 
semi-norms $\{\|.\|_j\}_{j\in\N}$, we will always assume the latter to be  sub-multiplicative, 
i.e$.$ 
$$
\|ab\|_j\leq\|a\|_j\|b\|_j\,,\qquad\forall\,a,b\in\CA\,,\quad \forall\,j\in\N\;.
$$
We start with  a crucial fact. Its proof being very similar to those of 
Lemma \ref{OUTPUTLEM}, we  omit it. 

\begin{lem}\label{genout}
Let $\CA$ be a Fr\'echet algebra 
 and let $\underline\mu_1,\underline\mu_2$ be two
families of weights on $G$ with sub-multiplicative degree respectively denoted by
 $(\underline L_{1},\underline R_{1})$ 
and $( \underline L_{2}, \underline R_{2})$.
Then the bilinear mapping
\begin{align*}
&{\mathcal R}\otimes {\mathcal R}:C^\infty(G,\CA)\times C^\infty(G,\CA)\to
C^\infty\big(G\times G,C^\infty(G,\CA)\big)\;,\\
&(F_1, F_2) \mapsto
\Big[(x,y)\in G\times G\mapsto(R^\star_xF_1)(R^\star_yF_2)
:=\big[g\in G\mapsto F_1(gx)F_2(gy)\big]\in\CA\Big]\;,
\end{align*}
is  continuous from $\CB^{\underline\mu_1}(G,\CA)\times\CB^{\underline\mu_2}(G,\CA)$
to $\CB^{\underline\nu}
\big(G\times G\,,\,\CB^{\underline\lambda}(G,\CA)\big)$, 
where the families
$\underline\nu:=\{\nu_{j,k}\}_{j,k\in\N}$ 
and $\underline\lambda:=\{\lambda_j\}_{j\in\N}$ are given by
$$
\nu_{j,k}:=\big( \mu_{1,j}\otimes\mu_{2,j}\big)^{\max(R_{1,j},R_{2,j})}\;\fd_{G\times G}^k\,,\qquad
\lambda_j:=\mu_{1,j}^{L_{1,j}}\,\mu_{2,j}^{ L_{2,j}}\;.
$$
More precisely, labeling by $(j,k)\in
\N^2$ the semi-norm $\|.\|_{j,k,\underline\lambda}$ of 
$\CB^{\underline\lambda}(G,\CA)$, for all $(j,k,k')$ in $\N^3$, 
there exists $C>0$ such that for all $F_1\in \CB^{\underline\mu_1}(G,\CA)$,
$F_2\in \CB^{\underline{\mu}_2}(G,\CA)$, we have
$$
\big\|{\mathcal R}\otimes{\mathcal R}(F_1,F_2)\big\|_{(j,k),k',\underline\nu}
\leq C\,\|F_1\|_{j,k+k',\underline\mu_1}\,\|F_2\|_{j,k+k',\underline{\mu}_2}\;.
$$
\end{lem}

\begin{thm}\label{OSCKERN}    
Let $(G\times G,S)$ be an admissible and tame tempered pair. Let also 
$\bm\in\CB^\mu(G\times G)$ for some tempered weight $\mu$ on $G\times G$
and  let $\underline\mu_1,\underline{\mu}_2$ be two
families of weights on $G$ with sub-multiplicative degree respectively denoted by 
 $(\underline L_{1},\underline R_{1})$ 
and $( \underline L_{2}, \underline R_{2})$, such that the family of  weights 
$\underline\mu_{1}\otimes \underline\mu_{2}$ is  tempered 
on $G\times G$.
Then, for any Fr\'echet algebra $\CA$, the oscillatory integral
\begin{align}
\label{starS}
\star_S:=\Big[(F_1,F_2)\mapsto\widetilde{\int_{G\times G}\bm\bE}\;\circ\;{\mathcal R}\otimes {\mathcal R} \,
(F_1,F_2)\Big]\;,
\end{align}
defines a continuous bilinear map  from 
$\CB^{\underline\mu_1}(G,\CA)\times\CB^{\underline{\mu}_2}(G,\CA)$ 
to $\CB^{\underline\lambda}(G,\CA)$, where the family $\underline\lambda:=\{\lambda_j\}_{j\in\N}$, 
is given by
$$
\lambda_j:=\mu_{1,j}^{L_{1,j}}\,\mu_{2,j}^{ L_{2,j}}\;.
$$ 
More precisely, for any $(j,k)\in\N^2$ there exist $C>0$ and $l\in\N$ such that for any 
$F_1\in\CB^{\underline\mu_1}(G,\CA)$
and $F_2\in\CB^{\underline{\mu}_2}(G,\CA)$, we have
$$
\|F_1\star_S F_2\|_{j,k,\underline\lambda}\leq C\,\|F_1\|_{j,l,\underline\mu_1}\,
\|F_2\|_{j,l,\underline{\mu}_2}\;.
$$
In particular, one has a  continuous bilinear product (not necessarily associative!):
$$
\star_S:\CB(G,\CA)\times
\CB(G,\CA)\to\CB(G,\CA)\;.
$$
\end{thm}
\begin{proof} By Lemma \ref{genout}, the map
$${\mathcal R}\otimes {\mathcal R} :\CB^{\underline\mu_1}(G,\CA)\times\CB^{\underline{\mu}_2}(G,\CA)\to
\CB^{\underline\nu}
\big(G\times G\,,\,\CB^{\underline\lambda}(G,\CA)\big)\;,$$
where $\nu_{j,k}=\big( \mu_{1,j}\otimes\mu_{2,j}\big)^{\max(R_{1,j},R_{2,j})}
\;\fd_{G\times G}^k$ and 
$\lambda_j=\mu_{1,j}^{L_{1,j}}\mu_{2,j}^{ L_{2,j}}$,
is a  continuous bilinear mapping.
By tameness,  the family of tempered weights
$\mu.\underline\nu$ 
is dominated.
Hence the oscillatory integral composed with the map ${\mathcal R}\otimes {\mathcal R} $ is well defined
 as a   continuous bilinear mapping
 from $\CB^{\underline\mu_1}(G,\CA)\times\CB^{\underline{\mu}_2}(G,\CA)$ to 
 $\CB^{\underline\lambda}(G,\CA)$.
  The precise estimate follows by putting 
 together Proposition \ref{PROPIP2}, Lemma \ref{genout}
  and the 
 continuity relation $\|.\|_{j,k,\underline{\hat\mu}}\leq C \|.\|_{j,k,\underline\mu}$, 
 for any family of weights $\underline{\hat\mu}$
 that dominates another family  $\underline\mu$.
\end{proof}
 We now discuss some issues regarding associativity of the  bilinear mapping $\star_S$.
To this aim, we need to show how to compute the product $F_1\star_S F_2$  as the limit  of a double 
sequence of products of smooth compactly supported functions.
\begin{lem}
\label{specific-limits}
Within the context of Theorem \ref{OSCKERN}, for $F_1\in \CB^{\underline\mu_1}(G,\CA)$ and 
$F_2\in \CB^{\underline{\mu}_2}(G,\CA)$, we let $\{F_{1,n}\}$, $\{F_{2,n}\}$ be two sequences in 
$\CD(G,\CA)$ 
converging respectively to $F_1$ and $F_2$ for the topologies of 
$\CB^{\underline{\hat\mu}_1}(G,\CA)$ and 
$\CB^{\underline{\hat\mu}_2}(G,\CA)$ where $\underline{\hat\mu}_1\succ  \underline{\mu}_1$,
$\underline{\hat\mu}_2\succ  \underline{\mu}_2$ are such that $\underline{\hat\mu}_1\otimes
\underline{\hat\mu}_2$ is tempered on $G\times G$. Setting then
 $\lambda_j:=\mu_{1,j}^{L_{1,j}}\mu_{2,j}^{ L_{2,j}}$, we have the equalities in 
$\CB^{\underline\lambda}(G,\CA)$:
$$
F_1\star_S F_2=\lim_{n_1\to\infty}\lim_{n_2\to\infty} F_{1,n_1}\star_S F_{2,n_2}=
\lim_{n_2\to\infty}\lim_{n_1\to\infty} F_{1,n_1}
\star_S F_{2,n_2}\;.
$$
\end{lem}
\begin{proof}
 Note that the family 
$\underline{\hat\nu}:=\{({\hat\mu}_{1,j}\otimes {\hat \mu}_{2,j})^{\max(R_{1,j},R_{2,j})}
\fd_{G\times G}^k\}_{(j,k)}$ is tempered (by assumption) and dominates
the (tempered by assumption)
 family $\underline\nu:=\{({\mu}_{1,j}\otimes { \mu}_{2,j})^{\max(R_{1,j},R_{2,j})}\,
\fd_{G\times G}^k\}_{(j,k)}$ and consequently (see Remark \ref{stress2}),
 we may 
view ${\mathcal R}\otimes {\mathcal R} (F_1,F_2)$ as an element of 
$\CB^{\underline{\hat\nu}}
\big(G\times G\,,\,\CB^{\underline\lambda}(G,\CA)\big)$,
with $\underline\lambda:=\{\mu_{1,j}^{L_{1,j}}\mu_{2,j}^{ L_{2,j}}\}_{j\in\N}$.
By the estimate of Theorem \ref{OSCKERN}, we know that for all $j,k\in\N$, 
there exists $l\in\N$
such that
$$
\|F_1\star_SF_2\|_{j,k,\underline\lambda}\leq C(k,j)\,\|F_1\|_{j,l,\underline{\hat\mu}_1}\,
\|F_2\|_{j,l,\underline{\hat\mu}_2}\;.
$$
We then write
\begin{align*}
&\|F_1\star_S F_2- F_{1,n_1}\star_S F_{2,n_2}\|_{j,k,\underline\lambda}\\&\leq
\|(F_1-F_{1,n_1})\star_S F_2\|_{j,k,\underline\lambda}+
\|F_{1,n_1}\star_S(F_2-F_{2,n_2})\|_{j,k,\underline\lambda}\\
&\leq C(k,j)\big( \|F_1-F_{1,n_1}\|_{j,l,\underline{\hat\mu}_1}\| F_2\|_{j,l,\underline{\hat\mu}_2}+
\|F_{1,n_1}\|_{j,l,\underline{\hat\mu}_1}\|F_2-F_{2,n_2}\|_{j,l,\underline{\hat\mu}_2}\big)\\
&\leq C'(k,j)\big( \|F_1-F_{1,n_1}\|_{j,l,\underline{\hat\mu}_1}\| F_2\|_{j,l,\underline{\mu}_2}+
\|F_{1,n_1}\|_{j,l,\underline{\hat\mu}_1}\|F_2-F_{2,n_2}\|_{j,l,\underline{\hat\mu}_2}\big)\;,
\end{align*}
 where in the last  inequality we used Remark \ref{stress2} which shows that
the semi-norm $\|.\|_{j,k,\underline{\hat\mu}_2}$ is
dominated by $ \|.\|_{j,k,\underline{\mu}_2} $. The latter remark also shows that the
numerical  sequence $\{\|F_{1,n_1}\|_{j,l,\underline{\hat\mu}_1}\}_{n_1\in\N}$
is bounded since
$$
\|F_{1,n_1}\|_{j,l,\underline{\hat\mu}_1}\leq \|F_{1,n_1}-F_1\|_{j,l,\underline{\hat\mu}_1}
+\|F_1\|_{j,l,\underline{\hat\mu}_1}
\leq \|F_{1,n_1}-F_1\|_{j,l,\underline{\hat\mu}_1}
+C\|F_1\|_{j,l,\underline{\mu}_1}\;.
$$
This completes the proof. 
\end{proof}
\begin{rmk}
\label{cor-explicit}
In other words, within the setting of  Lemma \ref{specific-limits}, we have  in 
$\CB^{\underline\lambda}(G,\CA)$:
$$
F\star_SF'=\lim_{m,n\to\infty}\int_{G\times G} \bE(x,x')\,\bm(x,x')\,R^\star_{x}(F_m)\,
R^\star_{x'}(F'_n)\,{\rm d}_G(x)\,{\rm d}_G(x')\;,
$$
for suitable approximation sequences $\{F_n\},\{F'_n\}\subset \CD(G,\CA)$.
\end{rmk}

\begin{dfn}
\label{wass}
Within the context of Theorem \ref{OSCKERN}, we say that the product $\star_S$, given in \eqref{starS}, is {\bf weakly associative} when for all $\psi_1,\psi_2,\psi_3\in\CD(G,\CA)$, one has
$(\psi_1\star_S\psi_2)\star_S\psi_3=\psi_1\star_S(\psi_2\star_S\psi_3)$ in $\CB(G,\CA)$.
\end{dfn}

\begin{prop}\label{WASS}
Within the context of Theorem \ref{OSCKERN}, weak associativity implies strong associativity in 
the sense that, when weakly associative, 
 for every further 
 family of   weights $\underline\mu_3$ on $G$ with sub-multiplicative degree denoted by 
 $(\underline L_{3},\underline R_{3})$ such that $\underline\mu_2\otimes\underline\mu_3$
 is tempered on $G\times G$.
 Then, for every element $(F_1,F_2,F_3)\in\CB^{\underline\mu_1}(G,\CA)\times
\CB^{\underline\mu_2}(G,\CA)
\times\CB^{\underline\mu_3}(G,\CA)$,
one has the equality $(F_1\star_SF_2)\star_S F_3=F_1\star_S (F_2\star_S F_3)$ in 
$\CB^{\underline\lambda}(G,\CA)$ for 
$\underline\lambda=\{\mu_{1,j}^{L_{1,j}^2}{\mu_{2,j}}^{L_{2,j}^2}{\mu_{3,j}}^{ L_{3,j}^2}\}_{j\in\N}$.
\end{prop}
\begin{proof}
Let $\mu_\phi$ be a tempered weight on $G\times G$ which dominates the constant weight $1$ (it exists 
by assumption of tameness.)
Consider the element $\nu_\phi\in C^\infty(G)$ defined by $\nu_\phi(g):=\mu_\phi(g,e)$. 
The latter is then (by a direct application of Lemma \ref{fd-semi}) a
 weight on $G$ that dominates $1$. Moreover,  it is easy to see that $\nu_\phi\otimes\nu_\phi$
 is tempered on $G\times G$.
Hence, all the family of weights $\underline\mu_1$, $\underline\mu_2$ and $\underline\mu_3$ are 
dominated e.g$.$ by
 $\underline{\hat\mu}_1:=\nu_\phi.\underline\mu_1$, $\underline{\hat\mu}_2:=
 \nu_\phi.\underline\mu_2$
  and $\underline{\hat\mu}_3:=\nu_\phi.\underline\mu_3$
and $\underline{\hat\mu}_1\otimes\underline{\hat\mu}_2$ and $\underline{\hat\mu}_2\otimes
\underline{\hat\mu}_3$ are tempered on $G\times G$. 
The assumptions of Lemma \ref{specific-limits} are therefore satisfied.
Let us  consider sequences of smooth compactly supported elements 
$\{\Phi_{1,n}\}_{n\in\N}\,,\{\Phi_{2,n}\}_{n\in\N}$ and $\{\Phi_{3,n}\}_{n\in\N}$
that converge to the elements $F_1$, $F_2$ and $F_3$ respectively in 
$\CB^{\underline{\hat\mu}_1}(G,\CA)\,,\,\CB^{\underline{\hat\mu}_2}(G,\CA)$ and 
$\CB^{\underline{\hat\mu}_3}(G,\CA)$.
Using separate continuity of $\star_S$ and Lemma \ref{specific-limits}, we observe the 
following equality:
$$
\lim_{n_1\to\infty}\Big(\lim_{n_2\to\infty}\Big(\lim_{n_3\to\infty}(\Phi_{1,n_1}\star_S
\Phi_{2,n_2})\star_S\Phi_{3,n_3}\Big)\Big)=(F_1\star_S F_2)\star_S F_3\;,
$$
in 
$\CB^{\underline\lambda}(G,\CA)$ for 
$\underline\lambda=\{\mu_{1,j}^{L_{1,j}^2}{\mu_{2,j}}^{L_{2,j}^2}{\mu_{3,j}}^{ L_{3,j}^2}\}$. 
One then concludes
using  weak associativity and the commutativity of the limits, as shown in Lemma \ref{specific-limits}.
\end{proof}

In section \ref{sec:Schwartz}, we have seen how to associate in a canonical way a Schwartz type 
functions space to a tempered, admissible and tame pair. Hence, starting with such a pair 
$(G\times G,S)$, we get a Schwartz space on $G\times G$. But we can also define a one-variable 
Schwartz space using the continuity of the partial evaluation maps:
\begin{dfn}
\label{S1P}
Let $(G\times G,S)$ be a tempered admissible
and tame pair and
$\CA$ be a Fr\'echet algebra. We define the $\CA$-valued Schwartz space on $G$ associated to 
$S$ by
$$
\CS^{S}(G,\CA):=\Big\{\big[g\in G\mapsto f(g,e)\big]\,,\,\,\, f\in\CS^{S}(G \times G,\CA)\Big\}\;.
$$  
We endow  the latter with the topology induced by the semi-norms:
\begin{align}
\label{SNS2}
\|f\|_{j,k,n}:=\sup_{X\in\,\CU_k(\g)}\sup_{x\in G}\Big\{\frac
{\mu_{\phi,1}(x)^n\big\|\widetilde X\,f(x)
\big\|_j}{|X|_k}\Big\}\;,\qquad j,k,n\in\N\;,
\end{align}
with $\mu_{\phi,1}(x):=\mu_\phi(x,e)$ and  $\mu_\phi$ is the tempered  weight on $G\times G$
associated with the tameness (Definition \ref{tame}).
\end{dfn}

 The next Lemma shows that the right action on the space of $\CA$-valued Schwartz 
functions, leads us  to a $\CB$-type space for family of weights too.
\begin{lem}
\label{belle}
Let $(G\times G,S)$ be a tame and admissible tempered pair, $\CA$ 
be a Fr\'echet algebra and
 $\underline\mu$  be a
family of  weights on $G$ with sub-multiplicative degree  $(\underline L,\underline R)$,
which is bounded by a power of $\mu_{\phi,1}$ (see Definition \ref{S1P}).
Then, there exists a  sequence $\{M_n\}_{n\in\N}$ of integers, such that
for  all elements $F\in\CB^{\underline\mu}(G,\CA)$ and $\varphi\in\CS^S(G,
\CA)$, the element $({\mathcal R}\otimes {\mathcal R} )(F,\varphi)$ (defined in Lemma \ref{genout}) belongs 
to $\CB^{\underline\nu}(G\times G,\CS^S(G,\CA))$ where $\underline\nu:=\{\nu_{j,k,n}\}_{j,k,n\in\N}$
with 
$$
\nu_{j,k,n}:=\big(\mu_j^{R_j}\otimes\mu_{\phi,1}^{{\!\vee} M_n }\big)
\fd_{G\times G}^{2k}\,,\qquad j,k,n\in\N\;.
$$
\end{lem}
\begin{proof} 
Denote by $\|.\|_{j,k,n}$, $(j,k,n)\in\N^3$, the semi-norms \eqref{SNS2} of $\CS^S(G,\CA)$.
Then, the semi-norms of  an element 
$$
\Phi=\big[(x,y)\in G\times G\mapsto[g\in G\mapsto \Phi(x,y;g)\in\CA]\big]\in 
\CB^{\underline\nu}\big(G\times G,\CS^S(G,\CA)\big)\;,
$$
are given by
\begin{align*}
&\|\Phi\|_{(j,k,n),(k_1,k_2),\underline\nu}=\\
&\quad\sup_{x,y,g\in G}
\sup_{X\in\,
\CU_k(\g)}\sup_{(Y_1,Y_2)\in\,
\CU_{k_1}(\g)\times\CU_{k_2}(\g)} \frac{\mu_{\phi,1}(g)^n\|
\widetilde{X}_g(\widetilde{Y}_1\otimes\widetilde{Y}_2)_{(x,y)}\,\Phi(x,y;g)\|_j}{\nu_{j,k,n}(x,y)}\;.
\end{align*}
Using Sweedler's notation \eqref{Sweedler}, we have for $X\in\CU_k(\g)$, $Y_1\in\CU_{k_1}(\g)$,
$Y_2\in\CU_{k_2}(\g)$:
\begin{align*}
&\widetilde{X}_g.\left((\widetilde{Y}_1\otimes\widetilde{Y}_2)_{(x,y)}
\left(R_x^\star F(g)\, R_y^\star\varphi(g)\right)\right)=\\
&\qquad\qquad
\sum_{(X)}\big(\widetilde{(\Ad_{x^{-1}}X_1)} \widetilde Y_1 F\big)(gx)\,
\big(\widetilde{(\Ad_{y^{-1}}X_2)} \widetilde Y_2 \vf\big)(gy)\;,
\end{align*}
which yields the following estimation  for  arbitrary $N\in\N$:
\begin{align*}
&\|\widetilde{X}_g.\left((\widetilde{Y}_1\otimes\widetilde{Y}_2)_{(x,y)}.
\left(R_x^\star F(g)\, R_y^\star\varphi(g)\right)\right)\|_j\\
&\quad\leq \sum_{(X)}|X_{(1)}|_k|X_{(2)}|_k |\Ad_{x^{-1}}|_k\,|\Ad_{y^{-1}}|_k\,|Y_1|_{k_1}|Y_2|_{k_2}\\
&\quad\qquad\qquad\times\sup_{Z_1\in\,\CU_{k+k_1}(\g)}\frac{\|\widetilde Z_1
F(gx)\|_j}{|Z_1|_{k+k_1}}\,
\sup_{Z_2\in\,\CU_{k+k_2}(\g)}\frac{\|\widetilde Z_2\vf(gy)\|_j}{|Z_2|_{k+k_2}}
\\
&\quad\leq \sum_{(X)}|X_{(1)}|_k|X_{(2)}|_k |\Ad_{x^{-1}}|_k\,|\Ad_{y^{-1}}|_k\,|Y_1|_{k_1}|Y_2|_{k_2}\\
&\quad\qquad\qquad\times\mu_j(gx)\,\mu_{\phi,1}^{-N}(gy)\,
\|F\|_{j,k+k_1,\underline\mu}\,\|\vf\|_{j,k+k_2,N}\;,
\end{align*}
which by Lemma \ref{fd-semi}, Lemma \ref{lemtwo} and the estimate \eqref{compat} is bounded by a constant times
\begin{equation}
\label{try}
|X|_k|Y_1|_{k_1}|Y_2|_{k_2}\,\fd_{G\times G}(x,y)^{2k}\,
\mu_j(gx)\,\mu_{\phi,1}^{-N}(gy)\,
\|F\|_{j,k+k_1,\underline\mu}\,\|\vf\|_{j,k+k_2,N}\;.
\end{equation}
Setting $(L,R)$ for the sub-multiplicative degree of $\mu_\phi$ and
using
 $$
 \mu_{\phi,1}^{-1}(gy)\leq C^{1/L}\,\mu_{\phi,1}^{-1/L}(g)
\mu_{\phi,1}^{R/L}(y^{-1})\,,\qquad y,g\in G\;,
$$
we see that \eqref{try} is (up to a constant)  bounded by
\begin{align*}
&|X|_k|Y_1|_{k_1}|Y_2|_{k_2}\,\fd_{G\times G}(x,y)^{2k}\,
\mu_j^{L_j}(g)\mu_j^{R_j}(x)\,\mu_{\phi,1}^{-N/L}(g)\mu_{\phi,1}^{N R/L}(y^{-1})\\
&\qquad\qquad\qquad\qquad\qquad\qquad\qquad\qquad\qquad\qquad
\times\|F\|_{j,k+k_1,\underline\mu}\,\|\vf\|_{j,k+k_2,N}\;.
\end{align*}
Now, given $n\in\N$, chose $N_n,M_n\in\N$  such that 
$$
\mu_j^{L_j}\mu_{\phi,1}^{-N_n/L}
\leq \mu_{\phi,1}^{-n}\;,
\qquad\mbox{and}\qquad   \mu_{\phi,1}^{N_nR/L}\leq \mu_{\phi,1}^{M_n}\;.
$$
 Then, set
$$
\nu_{j,k,n}(x,y):=\mu_j^{R_j}(x)\mu_{\phi,1}^{{\!\vee} M_n }(y)
\fd_{G\times G}^{2k}(x,y)\,,\qquad j,k,n\in\N\;.
$$
This entails that \eqref{try} is bounded by
$$
|X|_k|Y_1|_{k_1}|Y_2|_{k_2}\,\mu_{\phi,1}^{-n}(g) \mu_{j,k,n}(x,y)\,
\|F\|_{j,k+k_1,\underline\mu}\,\|\vf\|_{j,k+k_2,N_n}\;,
$$
which, for a finite constant $C(j,k,n,k_1,k_2)>0$, finally gives
$$
\|({\mathcal R}\otimes {\mathcal R} )(F,\varphi)\|_{(j,k,n),(k_1,k_2),\underline\nu}\leq C(j,k,n,k_1,k_2)
\|F\|_{j,k+k_1,\underline\mu}\,\|\vf\|_{j,k+k_2,N_n}\;,
$$
proving the claim.
\end{proof}

Observe that $1\otimes\mu_{\phi,1}$ is tempered on $G\times G$. By Remark \ref{vee} and since
the inversion map on a tempered group is a tempered map, we deduce that 
$1\otimes\mu_{\phi,1}^\vee$ is tempered  on $G\times G$ too. Hence, when $\underline \mu\otimes
1$ is also tempered on $G\times G$, so is the family of weights $\underline\nu$ given in  the 
Lemma \ref{belle}.
We then deduce the following important consequence of the latter:
\begin{prop}
\label{action-schwartz}
Let $(G\times G,S)$ be a tame and admissible tempered pair, $\CA$ 
 a Fr\'echet algebra,
 $\underline\mu$   a
family of tempered weights on $G$, such that the family of
weights $\underline\mu\otimes 1$ is tempered on $G\times G$
and  $\bm\in\CB^\mu(G\times G)$ for some tempered weight $\mu$ on $G\times G$.
Then the bilinear map 
$\star_S$, defined in \eqref{starS},  is   continuous   on $\CS^S(G,\CA)$ and
one has the  continuous bilinear map:
$$
\star_S:\CB^{\underline\mu}(G,\CA)\times \CS^S(G,\CA)\to \CS^S(G,\CA)\,,\quad(F,\varphi)\mapsto L_{\star_{S}}(F):\vf\mapsto F\star_S \vf\;.
$$
\end{prop}

\begin{rmk}
\label{conv-expr}
In the context of the proposition above, observe that the restriction to 
 $ \CS^S(G,\CA)\times \CS^S(G,\CA)$
 of the  bilinear product $\star_S$ \eqref{starS}, we  have the (point-wise and semi-norm-wise)
   absolutely convergent expression: 
$$
\vf_1\star_S\vf_2=\int_{G\times G} \bm(x_1,x_2)\,\bE(x_1,x_2)\,R_{x'_1}^\star(\vf_1)\,
R_{x'_2}^\star(\vf_2)\,
{\rm d}_G(x_1)\,{\rm d}_G(x_1)\;.
$$
\end{rmk}

\chapter{Tempered pairs for K\"ahlerian Lie groups}
\label{GeometrySSS}
The aim of this chapter is to endow each negatively curved  K\"ahlerian Lie group with the structure of 
a tempered, tame and admissible pair.
Recall that a
 Lie group $G$ is called a {\bf K\"ahlerian Lie group} when it is endowed with an invariant K\"ahler 
structure, i.e$.$ a left-invariant complex structure ${\bf J}$ together 
with a left-invariant Riemannian metric ${\bf g}$ such that the triple $(G,{\bf J},{\bf g})$ constitutes a 
K\"ahler manifold. Within the present memoir, we will
be concerned with K\"ahlerian Lie groups whose sectional curvature is negative. We call them 
{\bf negatively curved}.

In section \ref{PStheory}, we briefly review the theory of normal $\bf j$-algebras and associated
normal $\bf j$-groups, which in turn gives a classification result for negatively curved 
K\"ahlerian Lie groups. In section \ref{H-group}, we explain how an elementary normal 
$\bf j$-group is naturally endowed with a structure of a symplectic symmetric space. Such 
structure is  the core  of the construction of chapter \ref{QPSSS}. It is also in this context that
we have a clear geometric construction for the phase  and amplitude of the kernel underlying
our deformation formula. It is then in sections \ref{2P} and \ref{2pencore} that we prove
that the phase mentioned above  endows a negatively curved K\"ahlerian Lie group with the
structure of an admissible and tempered pair.

\section{Pyatetskii-Shapiro's  theory }
\label{PStheory}
The following definition, due to Pyatetskii-Shapiro \cite{PS}, describes the infinitesimal structure of 
negatively curved K\"ahlerian Lie groups.
\begin{dfn}
\label{normal-j-alg}
A {\bf normal $\bf j$-algebra} is a triple  $(\b,\alpha,{\bf j})$ where 
\begin{enumerate}
\item[(i)]\label{axiomI}
$\b$ is a solvable Lie algebra which is split over the reals,  
i.e$.$ $\ad_X$ has only real eigenvalues for 
all $X\in\b$,

\item[(ii)]\label{axiomII}
$\bf j$ is an endomorphism of $\b$ such that ${\bf j}^2=-1$ and
\[ [X,Y] + {\bf j}[{\bf j}X,Y] + {\bf j}[X,{\bf j}Y] - [{\bf j}X,{\bf j}Y] = 0\,,\quad X,Y\in\b \;,\]

\item[(iii)]\label{axiomIII}
$\alpha$ is a linear form on $\b$ such that
\[ \alpha([{\bf j}X,X]) > 0\;\text{ if }\;X\neq 0 \text{\quad and \quad} \alpha([{\bf j}X,{\bf j}Y])=\alpha([X,Y]) 
\,,\quad X,Y\in\b \;.\]
\end{enumerate}
\end{dfn}
We quote the following structure result from \cite{PS}.
\begin{prop}
The Lie algebra of a negatively curved K\"ahlerian Lie group always carries a structure of normal 
${\bf j}$-algebra.
\end{prop}
 If $\b'$ is a  subalgebra of $\b$ which is invariant by $\bf j$, then 
$(\b',\alpha|_{\b'},{\bf j}|_{\b'})$ is again a normal $\bf j$-algebra, called a {\bf $\bf j$-subalgebra} of 
$(\b,\alpha,{\bf j})$. A $\bf j$-subalgebra whose underlying Lie algebra $\b'$ is an ideal of $\b$ is 
called a {$\bf j$-ideal}.
\begin{ex}\label{SU1n}{\rm 
Every Iwasawa factor $AN$ of the simple Lie group $SU(1,n)$ is naturally a negatively curved 
K\"ahlerian Lie group. Indeed, denoting by $K\simeq U(n)$ a maximal compact subgroup of 
$SU(1,n)$,
one knows that the associated symmetric space $G/K$ is a negatively curved K\"ahlerian 
$SU(1,n)$-manifold. The associated Iwasawa decomposition $SU(1,n)=ANK$ then yields
a global diffeomorphism between $G/K$ and $AN$. Transporting to $AN$ the K\"ahler structure of 
$G/K$ under the latter diffeomorphism, then endows $AN$ with a negatively curved K\"ahlerian
Lie group structure, called \emph{elementary} after Pyatetskii-Shapiro. 
}
\end{ex}
The infinitesimal structure underlying an elementary normal $\bf j$-group (cf$.$ the above Example \ref{SU1n}) may be precisely described as follows.
Let $(V,\omega^0)$ be a symplectic vector space   of real dimension $2d$.  We consider  the associated Heisenberg Lie algebra  $\h:=V\oplus\R E$. That is, $\h$ is the central extension of the Abelian Lie algebra $V$, with brackets given by
$$
[v_1,v_2]:=\omega^0(v_1,v_2)\,E\,, \quad v_1,v_2\in V\,,\qquad[E,X]:=0\,,\quad X\in\h\;.
$$ 
\begin{dfn}
\label{ENLA}
 Let  $\a$ be a one-dimensional real Lie algebra, with generator  $H$.
We consider the split extension of Lie algebras:
$$
0\,\to\,\h\,\to\,\s:=\a\ltimes_{\rho_\h}\h\,\to\,\a\,\to\,0\;,
$$ 
with extension homomorphism $\rho_\h\,:\,\a\,\to\,{\rm Der}(\h)$ given by
\begin{equation}
\label{split}
\rho_\h(H)\big(v+ t\,E\big):=[H,v+ t\,E]:=v+ 2t\,E\,,\qquad v\in V\,,\quad t\in\R\;.
\end{equation}
The Lie algebra $\s$ is called {\bf elementary normal}.
Last, we denote by $\S$ the connected simply connected Lie group whose Lie algebra is  $\s$ and we call
the latter an {\bf elementary normal $\bf j$-group}.
\end{dfn}
Note that  $\S$ is a solvable group of real dimension $2d+2$ and if $V=\{0\}$, $\S$ is isomorphic to the affine
group of the real line.
 It turns out that every negatively curved K\"ahlerian Lie group can be decomposed into elementary pieces: at the infinitesimal level, one has the following result, due to Pyatetskii-Shapiro \cite{PS}.
\begin{prop}\label{P-S}
 Let $(\b,\alpha,\bf j)$ be a normal $\bf j$-algebra. Then, there exist  $\z$, a one-dimensional ideal of $\b$ and $V$,  a vector subspace  of $\b$, such that setting $\a:={\bf j}\z$, the algebra $\s:=\a\oplus V\oplus \z$ underlies an elementary normal $\bf j$-ideal of $\b$. Moreover,
the associated extension sequence
$$
0\longrightarrow\s\longrightarrow\b\longrightarrow\b'\longrightarrow0\;,
$$
is split as a sequence of normal $\bf j$-algebras and such that:
\begin{equation}
\label{P-S_a} 
[\b',\a\oplus\z]=0\quad\mbox{and}\quad
[\b', V]\subset V\;.
\end{equation}
\end{prop}
In particular, every normal $\bf j$-algebra $\b$ admits a decomposition as a sequence of split extensions of elementary normal $\bf j$-algebras $\s_i$, $i=1,\dots,N$, of real  dimension $2d_i+2$, $d_i\in\N$:
\begin{equation*}
\big(\dots\big(\s_N\ltimes\s_{N-1}\big)\ltimes\dots\ltimes\s_2\big)\ltimes\s_1\;,
\end{equation*}
such that for all $i=1,\dots,N-1$
$$
\big[\big(\s_N\ltimes\dots\big)\ltimes\s_{i+1},\a_{i}\oplus\z_{i}\big]=0\quad\mbox{and}\quad
\big[\big(\s_N\ltimes\dots\big)\ltimes\s_{i+1}, V_{i}\big]\subset V_{i}\;.
$$
\begin{dfn}
\label{norma-j-group}
A {\bf normal $\bf j$-group} $\B$, consists of a connected simply connected Lie group that admits a normal $\bf j$-algebra as Lie algebra, i.e$.$ $\B=\exp\{\b\}$, where $\b$ is a normal $\bf j$-algebra.
\end{dfn}

At the group level, for $i=1,\dots,N-1$, call $\bR^i$ the extension homomorphism at each step:
\begin{align}
\label{AES}
\bR^i\in{\rm Hom}\big((\S_N\ltimes\dots)\ltimes\S_{i+1},\Aut(\S_i)\big)\;.
\end{align}
The conditions given in \eqref{P-S_a} implies that $\bR^i$ takes values in ${\rm Sp}(V_i,\omega^0_i)$, where
$(V_i,\omega^0_i)$ denotes the symplectic vector space attached  to $\S_i$. 

\section{Geometric structures on elementary normal $\bf j$-groups}
\label{H-group}
In this section, we review the properties of a symplectic symmetric space structure every elementary normal $\bf j$-group is naturally endowed with.
The phase function with respect to which an admissible tempered pair will be associated to later on, was defined in \cite{Bi07} in terms 
of this symplectic symmetric space structure. We start with the definition of a symplectic symmetric space as in \cite{Bi95} which is 
an adaptation to the symplectic case of the notion of symmetric space as introduced by O. Loos \cite{Lo}.

\begin{dfn}
\label{def:SSS}
A {\bf symplectic symmetric space} is a triple $(M,s,\omega)$ where 
\begin{enumerate}
\item[(i)] $M$ is a connected smooth manifold,

\item[(ii)]  $s$ is a smooth map
$$s:M\times M\to M\,,\quad(x,y)\mapsto s_x(y):=s(x,y)\;,$$
 such that:

\noindent (ii.1) For every $x\in M$, the partial map $s_x:M\to M$ is an involutive diffeomorphism 
admitting $x$ as isolated fixed point. The diffeomorphism $s_x$ is called the {\bf symmetry} 
at point $x$.

\noindent (ii.2) For all points $x$ and $y$ in $M$, the following relation holds:
$$
s_x\circ s_y \circ s_x=s_{s_x(y)}\;.
$$
\item[(iii)] $\omega$ is a closed\footnote{The closedness condition is, in fact, redundant,
see   \cite{Bi95}.} and non-degenerate  two-form on $M$ that is invariant 
under the symmetries:
$$
s_x^\star\,\omega=\omega\,,\quad\forall x\in M\;.
$$
\end{enumerate}
A {\bf morphism} between two  symplectic symmetric spaces is defined as a
symplectomorphism  that intertwines the symmetries.
\end{dfn}

Symplectic symmetric spaces always carry a preferred Lie group of transformations \cite{Bi95}:

\begin{dfn}
\label{autSSS}
The {\bf automorphism group} of a symplectic symmetric space $(M,s,\omega)$
is constituted by the symplectomorphisms
$\vf\in\mbox{\rm Symp}(M,\omega)$ which are covariant under the 
symmetries:
$$
\vf\circ s_x=s_{\vf(x)}\circ\vf,\quad\forall x\in M\;.
$$
It is a Lie subgroup of $\mbox{\rm Symp}(M,\omega)$ that 
acts transitively on $M$ and it is denoted by 
$\Aut(M,s,\omega)$.  Its Lie algebra is called
the {\bf derivation algebra} of $(M,s,\omega)$ and is denoted by $\mathfrak{aut}(M,s,\omega)$.
\end{dfn}

We now pass to the particular case of a given $2d+2$-dimensional elementary normal 
$\bf j$-group $\S$ with associated 
symplectic form $\omega^\S$.
Let $a,t\in\R$ and $v\in V\simeq\R^{2d}$.  The following identification will always be understood:
$$
\R^{2d+2}\to\s\;,\quad  x:=(a,v,t)\mapsto aH+v+tE\;.
$$
The following result is extracted from \cite{BiMas,BCSV,Bi07}:
\begin{prop}
\label{proprietes}
Let $\S$ be an elementary normal $\bf j$-group.
\begin{enumerate}
\item[(i)] The map
\begin{equation}
\label{chartS}
\s\to\S\,,\,\quad(a,v,t)\mapsto \exp(aH)\exp(v+tE)=\exp(aH)\exp(v)\exp(tE)\;,
\end{equation}
is a global Darboux chart on $(\S,\omega^\S)$ in which the symplectic structure  reads:
$$
\omega^\S:=2{\rm d}a\wedge {\rm d}t\,+\,\omega^0\;.
$$
\item[(ii)] Setting furthermore
\begin{align}
\label{SSs}
& s_{(a,v,t)}(a',v',t'):=\\
\big(2a-a',2v\cosh(a-a')-v',&2t\cosh(2a-2a')-t'+\omega^0(v,v')\sinh(a-a')\big)
\nonumber\;,
\end{align}
defines a symplectic symmetric space structure $(\S,s,\omega^\S)$ on the elementary normal 
$\bf j$-group $\S$.
\item[(iii)] The left action $L_x:\S\to\S$, $x'\mapsto x.x'$, defines a injective Lie group homomorphism
$$
L:\S\to\Aut(\S,s,\omega^\S)\;.
$$
In the coordinates \eqref{chartS}, we have
$$
x.x'=(a,v,t).(a',v',t')=\big(a+a',e^{-a'}v+v',e^{-2a'}t+t'+\tfrac12e^{-a'}\omega^0(v,v')\big)\;.
$$
and
$$
x^{-1}=(a,v,t)^{-1}=(-a,-e^a v,-e^{2a}t)\;.
$$
\item[(iv)] The action $\bR\;:\;\mbox{\rm Sp}(V,\omega^0)\times \S\to\S$, 
$(A,(a,v,t))\mapsto\bR_A(a,v,t):=(a,Av,t)$ by automorphisms of the normal $\bf j$-group $\S$
induces an injective Lie group homomorphism:
$$
\bR:\mbox{\rm Sp}(V,\omega^0)\to\Aut(\S,s,\omega^\S)\,,\quad A\mapsto\bR_A\;.
$$
  In fact,
 $\mbox{\rm Sp}(V,\omega^0)\simeq\Aut(\S)\cap \Aut(\S,s,\omega^\S)$.
 \end{enumerate}
\end{prop}
Note that  in the coordinates \eqref{chartS}, the modular function of $\S$, $\Delta_{\S}$, reads $e^{(2d+2)a}$.

We now pass to the definition of the three-point phase on $\S$. For this we need 
the notion of ``double geodesic triangle" as introduced by A. Weinstein \cite{We} and Z. Qian
\cite{Qi}.
\begin{dfn} 
\label{def-midpoint}
Let $(M,s)$ be a symmetric space.
A {\bf midpoint map} on $M$ is a smooth map
$$
M\times M\to M\,,\quad(x,y)\mapsto\,\mid(x,y)\;,
$$
such that, for all points $x,y$ in $M$:
$$
s_{\mid(x,y)}(x)=y\;.
$$
\end{dfn}
\begin{rmk}
When it exists, such a midpoint map on a symmetric space $(M,s)$ is necessarily unique (see 
Lemma 2.1.6 of  \cite{YVthese}).
\end{rmk}

\begin{rmk}
\label{MDPT}
Observe that in the case where the partial maps
$s^y:M\to M$, $x\mapsto s_x(y)$
are global diffeomorphisms of $M$, a midpoint map exists and is given by:
$$
\mid(x,y):=\left(\,s^x\,\right)^{-1}(y)\;.
$$
Note that in this case, every $\vf\in\Aut(M,s)$ intertwines the midpoints. Indeed, since for all 
$x,y\in M$ we have $\vf(s_y(x))=s_{\vf(y)}\big(\vf(x)\big)$, we get
$$
\vf\big(\mid(x,y)\big)=\mid\big(\vf(x)\,,\,\vf(y)\big)\;.
$$
\end{rmk}
An immediate computation shows that a midpoint map always exists on the symplectic symmetric
space attached to an elementary normal $\bf j$-group:
\begin{lem}
For the symmetric space $(\S,s)$  underlying an elementary normal $\bf j$-group, the associated 
partial maps are global  diffeomorphisms. In the coordinates \eqref{chartS}, we have:
\begin{align*}
&\big(s^{(a_0,v_0,t_0)}\big)^{-1}(a,v,t)=\\
&\Big(\frac{a+a_0}2,\frac{v+v_0}{2\cosh(\tfrac{a-a_0}2)},
\frac{t+t_0}{2\cosh({a-a_0})}-\omega^0(v,v_0)\tfrac{\sinh(\tfrac{a-a_0}2)}{4\cosh({a-a_0})
\cosh(\tfrac{a-a_0}2)} \Big)\;.
\end{align*}
\end{lem}
The following statement is proved in \cite{Bi07}.
\begin{prop}
\label{LB}
Let $\S$ be an elementary normal $\bf j$-group.
\begin{enumerate}
\item[(i)] The K\"ahler   manifold $\S$ is strictly geodesically complete:  two points determine 
a unique geodesic arc.
\item[(ii)] The ``medial triangle" three-point function
$$
\Phi:\S^3\to\S^3\,,\quad (x_1,x_2,x_3)\mapsto \big(\mid(x_1,x_2),\mid(x_2,x_3),
\mid(x_3,x_1)\big)\;,
$$
is a $\S$-equivariant (under the left regular action) global diffeomorphism. 
\end{enumerate}
\end{prop}
Since our space $\S$ has trivial de Rham 
cohomology in degree two,  any three points $(x,y,z)\in\S^3$ 
define an oriented
geodesic triangle $T(x,y,z)$ whose symplectic area is well-defined by integrating the two-form 
$\omega^\S$
on any surface admitting $T(x,y,z)$ as boundary. With a slight abuse of notation, we set
$$
{\rm Area}(x,y,z):=\int_{T(x,y,z)}\omega^\S\;.
$$
\begin{dfn}
\label{df-Scan}
The {\bf canonical two-point phase} associated to an elementary normal $\bf j$-group  is defined by
$$
S_{\rm can}^\S(x_1,x_2):={\rm Area}\left(\Phi^{-1}(e,x_1,x_2)\right)\;\in\;C^\infty(\S^2,\R)\;,
$$
where $e:=(0,0,0)$ denotes the unit element in $\S$.
In the coordinates \eqref{chartS}, one has the explicit expression:
\begin{equation}\label{SEXPL}
S_{\mathrm{can}}^\S(x_1,x_2)=t_2\sinh2a_1-t_1\sinh2a_2+\omega^0(v_1,v_2)\cosh a_1 
\cosh a_2\,\;.
\end{equation}
 The {\bf canonical two-point amplitude} associated to an elementary normal $\bf j$-group is defined by
\begin{align*}
A_{\mathrm{can}}^\S(x_1,x_2)
:=\mbox{\rm Jac}_{\Phi^{-1}}(e,x_1,x_2)^{1/2}\;\in\;C^\infty(\S^2,\R)\;.
\end{align*}
In the coordinates \eqref{chartS}, it reads
\begin{align}\label{AEXPL}
&A_{\mathrm{can}}^\S(x_1,x_2)=\\
&\big(\cosh a_1\,\cosh a_2\,\cosh(a_1-a_2)\big)^d\big(\cosh 2a_1\,\cosh 2a_2\,
\cosh( 2a_1-2a_2)\big)^{1/2} \;.\nonumber
\end{align}
\end{dfn}

\section{Tempered pair for  elementary normal $\bf j$-groups}
\label{2P}
The aim of this technical section is to prove that the pair $(\S\times\S, S_{\rm can}^\S)$
is tempered, admissible and tame.
We start by splitting the $2d$-dimensional symplectic vector space $(V,\omega^0)$ associated to an elementary normal $\bf j$-group $\S$  
into a direct sum of two Lagrangian subspaces in symplectic duality:
$$
V=\l^\star\oplus\l\;.
$$

The following result establishes temperedness.
\begin{lem}
The pair $(\S\times\S,S_{\mathrm{can}}^\S)$ is tempered. Moreover, the Jacobian of the map
$$
\phi:\S\times\S\to(\s\oplus\s)^\star\,,\quad g\mapsto\Big[X\in\s\oplus\s\mapsto \big(\widetilde X .\,S_{\rm can}^\S\big)(g)\Big]\;,
$$
is proportional to $(A_{\rm can}^{\S})^2$.
\end{lem}
\begin{proof}
Let us fix $\{f_j\}_{j=1}^d$, a basis of $\l^\star$  to which we associate  $\{e_j\}_{j=1}^d$  the symplectic-dual basis 
of $\l$, i.e$.$ it is defined by
 $ \omega^0(f_i,e_j)=\delta_{i,j}$.
We let  $E$ the central element
of the Heisenberg Lie algebra $\h\subset\s$ and $H$ the generator of $\a$ in the one dimensional  split extension
which defines the Lie algebra $\s$:
$$
0\,\to\,\h\,\to\,\s\to\,\a\,\to\,0\;.
$$ 
Accordingly, we consider the following basis of $\s\oplus\s$:
$$
\begin{array}{cccccc}
H_1&:=&H\oplus\{0\}\;,\quad H_2&:=&\{0\}\oplus H\;,\\
f_{j}^1&:=&f_j\oplus\{0\}\;,\quad  f_j^2&:=&\{0\}\oplus f_j\;,\\
e_{j}^1&:=&e_j\oplus\{0\}\;,\quad e_j^2&:=&\{0\}\oplus e_j\;,\\
E_1&:=&E\oplus\{0\}\;,\quad E_2&:=&\{0\}\oplus E\;,
\end{array}
$$
where the index $j$ runs from $1$ to $d=\dim(V)/2$.
From Proposition \ref{proprietes} iii) and with the notation $v=(x,y)\in\l^\star\oplus\l=V$, we see that the left-invariant vector fields on $\S$ are given by:
\begin{equation}\label{INVVECT}
\begin{array}{ccc}
\widetilde{H}&=&\partial_a\,-\,\sum_{j=1}^d(x_j\partial_{x_j}+y_j\partial_{y_j})\,-2t\partial_t\;,\\
\widetilde{f}_j&=&\partial_{x_j}\,-\frac{y_j}{2}\partial_t\;,\\
\widetilde{e}_j&=&\partial_{y_j}+\frac{x_j}{2}\partial_t\;,\\
\widetilde{E}&=&\partial_t\;.
\end{array}
\end{equation}
Thus, we find
\begin{align}\label{comput}
\widetilde{H}_1\,S_{\rm can}^\S&=2t_2\cosh2a_1+2t_1\sinh2a_2
-\omega^0(v_1,v_2)e^{-a_1}\cosh a_2\;,
\\
\widetilde{H}_2\,S_{\rm can}^\S&=-2t_1\cosh2a_2-2t_2\sinh2a_1
-\omega^0(v_1,v_2)e^{-a_2}\cosh a_1\;,\nonumber\\
\widetilde{E}_1\,S_{\rm can}^\S&=-\sinh2a_2\;,\nonumber\\
 \widetilde{E}_2\,S_{\rm can}^\S&=\sinh2a_1\;,
\nonumber\\
\widetilde{f}_j^1\,S_{\rm can}^\S&=y_2^j\cosh a_1\cosh a_2+\tfrac12y_1^j\sinh2a_2\;,\nonumber\\
\widetilde{f}_j^2\,S_{\rm can}^\S&=-y_1^j\cosh a_1\cosh a_2-\tfrac12y_2^j\sinh2a_1\;,\nonumber\\
\widetilde{e}_j^1\,S_{\rm can}^\S&=-x_2^j\cosh a_1\cosh a_2-\tfrac12x_1^j\sinh2a_2\;,\nonumber\\ 
\widetilde{e}_j^2\,S_{\rm can}^\S&=x_1^j\cosh a_1\cosh a_2+\tfrac12x_2^j\sinh2a_1\;\nonumber.
\label{comput}
\end{align}
A computation then shows that the Jacobian of the map 
$\phi:\S\times\S\to(\s\oplus\s)^\star$, underlying Definition \ref{TEMPPAIR}, is given by
\begin{align*}
2^{2d+2}
&\big(\cosh a_1\,\cosh a_2\,\cosh(a_1-a_2)\big)^{2d}
\cosh 2a_1\,\cosh 2a_2\,\cosh 2(a_1-a_2)\\
&\qquad\qquad\qquad= 2^{2d+2}\,{A_{\rm can}^\S}(x_1,x_2)^2\geq2^{2d+2}\;,
\end{align*}
and hence $\phi$ is a global diffeomorphism. It is also clear from Proposition \ref{proprietes} iii),
that the multiplication and inversion maps on $\S\times\S$  are tempered functions in the coordinates
\eqref{comput}. Therefore, the pair
$(\S\times\S,S_{\rm can}^\S)$ is tempered.
\end{proof}

\begin{rmk}
 Note that the formal adjoints of the left invariant vector fields \eqref{INVVECT}, with respect to the inner product of $L^2(\S)$ read:
$$
\widetilde{H}^*=-\widetilde{H}+2d+2\,,\quad \widetilde{f}_j^*=-\widetilde{f}_j\,,\quad
\widetilde{e}_j^*=-\widetilde{e}_j\,,\quad {\widetilde{E}}^*=-\widetilde{E}\;,
$$
so that the assumption \eqref{adjoint} is trivially satisfied.
\end{rmk}

 We will now prove that the tempered pair $(\S\times\S,S_{\rm can}^\S)$  is admissible and tame. For this, we need a decomposition
of the Lie algebra $\s$ and we shall use the following one:
\begin{equation}\label{DECOMP-S}
\s=\bigoplus_{k=0}^3V_k\quad\mbox{\rm where }\quad V_0\;:=\a\,,\quad V_1:=\l^\star\,,\quad V_2:=\l\quad\mbox{\rm and}\quad V_3:=\R E\;.
\end{equation}
Note that both $V_0$ and $V_3$ are of dimension one, while $V_1$ and $V_2$ are $d$-dimensional.
Accordingly, we consider the decompositions of $\s\oplus\s$ given by
$$
\s\oplus\{0\}=\bigoplus_{k=0}^3V_{1,k}\quad\mbox{\rm and}\quad\{0\}\oplus\s=\bigoplus_{k=0}^3V_{2,k}\;,
$$
where the subspaces $V_{ i,k}$, $ i=1,2$, of each factor correspond respectively to the subspaces $V_k$ of $\s$ within the decomposition (\ref{DECOMP-S}).
 We then set:
\begin{equation}
\label{decomposition}
\fV_k:=V_{1,k}\oplus V_{2,k}\quad\mbox{\rm and}\quad\s\oplus\s=\bigoplus_{k=0}^3\fV_{k}\;,
\end{equation}
by which we mean that there are four subspaces involved in the ordered decomposition of 
$\s\oplus\s$. We also let 
$$
\fV^{(k)}:=\bigoplus_{n=0}^k\fV_k\,,\quad k=0,1,2,3\;,
$$
as in \eqref{Vn} and we let $\CU(\fV^{(k)})$ be the unital subalgebra of $\CU(\s\oplus\s)$ generated
by $\fV^{(k)}$ as in \eqref{CUVE}.
Accordingly, we consider the associated tempered coordinates \eqref{TEMPCOORD}:
$$
x_{ i,0}:=\widetilde{H}_ i \,S_{\rm can}^\S\,,\quad x_{ i,1}^j:=\widetilde{f_j^{ i}}\,S_{\rm can}^\S\,,\quad x_{ i,2}^j:=\widetilde{e_j^{ i}}\,S_{\rm can}^\S\,,\quad x_{ i,3}:=\widetilde E_{ i}\,S_{\rm can}^\S\;,
$$
with $i=1,2$, $j=1,\dots,d$ and we use the vector notations:
\begin{align}
\label{TC}
\vec x_0&:=(x_{1,0},x_{2,0})\in\R^2\;,\\
\vec x_1&:=(x_{1,1},x_{2,1}):=\big((x_{1,1}^j)_{j=1}^d,(x_{2,1}^j)_{j=1}^d)\in\R^{2d}\;,\nonumber\\
\vec x_2&:=(x_{1,2},x_{2,2}):=\big((x_{1,2}^j)_{j=1}^d,(x_{2,2}^j)_{j=1}^d)\in\R^{2d}\;,\nonumber\\
\vec x_3&:=(x_{1,3},x_{2,3})\in\R^2\;. \nonumber
\end{align}
According to the notations $(a,v,t)\in\R\times\R^{2d}\times\R\simeq \S$ and $v=(x,y)\in\l^\star\oplus\l=V$, we set
$$
\vec a:=(a_1,a_2)\in\R^2\,,\quad \vec x=(x_1,x_2)\in\R^{2d}\,,\quad \vec y=(y_1,y_2)\in\R^{2d}\,,\quad\vec t:=(t_1,t_2)\in\R^2\;.
$$
We consider the functions
$$
s_{12}:=t_2\sinh2a_1-t_1\sinh2a_2\,,\quad\Omega_{12}:=\omega^0(v_1,v_2)\,,\quad \gamma_{12}:=\cosh a_1\,\cosh a_2\;,
$$
in term of which we have
$$
S_{\rm can}^\S=s_{12}\,+\,\gamma_{12}\,\Omega_{12}\;.
$$
Introducing last
\begin{align}
\label{AB}
A&:=
\begin{pmatrix}
\sinh2a_2&\cosh2a_1\\
-\cosh2a_2&-\sinh2a_1
\end{pmatrix}\;,\\
B&:=
\begin{pmatrix}
-\tfrac{1}{2}\sinh2a_2&-\cosh a_1\,\cosh a_2\\
\cosh a_1\,\cosh a_2&\tfrac{1}{2}\sinh2a_1
\end{pmatrix}\;,\nonumber
\end{align}
\begin{align*}
\vec{\gamma}:=
\begin{pmatrix}e^{-a_1}\cosh a_2,&e^{-a_2}\cosh a_1
\end{pmatrix}\;,\quad
\vec{\delta}:=
\begin{pmatrix}-\sinh2a_2\,,&\sinh2a_1
\end{pmatrix}\;,\nonumber
\end{align*}
the relations given in \eqref{comput} can be summarized as:
\begin{equation}\label{comput2}
\vec x_{3}=\vec\delta\;,\quad \vec{x}_2=B.\vec{x}\;,\quad \vec{x}_1=-B.\vec{y}\;,\quad
\vec{x}_0=2A.\vec{t}-\Omega_{12}\,\vec{\gamma}\;.
\end{equation}

We first treat  the easiest variable $\vec x_{3}$, which lead to multipliers $\alpha_{3}$ that satisfy 
property (ii) of Definition \ref{TEMPADM} with constant $\mu_{3}$:
\begin{lem}
\label{first}
Consider an  element $X\in \CU(\fV_3)$ such that   the associated multiplier $\alpha_X$ is  invertible. 
Then, for every $Y\in\CU(\fV^{(3)})=\CU(\s\oplus\s)$ there exists a positive constant $C_Y$ such that
$$
\big|\widetilde Y \,\alpha_X\big|\leq C_Y\,|\alpha_X|\;.
$$
\end{lem}
\begin{proof}
Note first that $\fV_3$ turns out to be a two-dimensional Abelian Lie algebra. Note also that 
$\alpha_{E_ i}$, $ i=1,2$ is independent of the variable $\vec t$. Thus, given a two-variables 
polynomial $P$, we have
for $X=P(E_1,E_2)\in\CU(\fV_3)$:
$$
\alpha_X=P\big(-\sinh2a_2,\sinh2a_1\big)\;.
$$
It also follows from the explicit expression of the left-invariant vector fields given in \eqref{INVVECT} 
that $\widetilde Y\,\alpha_X=0$ for all $Y\in \CU(\oplus_{k=1}^3\fV_k)$.
 Hence, it suffices to treat the case of 
$Y\in\CU(\fV_0)$. Observe that the restriction of $\widetilde{H}_j$ to functions which depend only on 
$a_j$, equals $\partial_{a_j}$. Thus in this case, we see that $\widetilde Y\,\alpha_X$ is a polynomial of 
the same degree as $P$, but in the variables $e^{\pm a_1}$ and $e^{\pm a_2}$. This is enough to 
conclude when $\alpha_X$ is invertible.
\end{proof}

 Next, we treat the variables $\vec{x}_2$ and  $\vec{x}_1$. We first observe
\begin{lem}\label{RECB}
There exist finitely many matrices $B_{(r)}\in M_2\left(\R[e^{\pm a_1},e^{\pm a_2}]\right)$ such that for 
all integers $N_1$ and $N_2$, the 
elements $\widetilde{H}_1^{N_1}\widetilde{H}_2^{N_2}\,B$ consist of a linear combination of the 
$B_{(r)}$'s, where the matrix $B$ has been defined in \eqref{AB}. The same property holds for the 
matrix $A$.
\end{lem}
\begin{proof}
Set
\begin{align*}D:=
\begin{pmatrix}
-\tfrac{1}{2}\sinh 2a_2&0\\
0&\tfrac{1}{2}\sinh2a_1
\end{pmatrix}
\;,\quad
\Gamma:=\gamma_{12}\begin{pmatrix}0&-1\\
1&0\end{pmatrix}\;,
\end{align*}
and observes that
$$
B=D+\Gamma\;\quad\mbox{and}\quad\;\partial^{2}_{a_ i}\Gamma=\Gamma\,,\quad i=1,2\;.
$$
The derivatives of $B$ therefore all belong to the space generated by  the entries of $D$ and
 $\Gamma$ and by finitely many of their derivatives.
This is enough to conclude since restricted to functions that depend only on the variable $a$, we have $\widetilde H=\partial_a$. The proof for the matrix $A$ is entirely similar.
\end{proof}

We can now deduce what we need for the variables $\vec x_2$ and  $\vec{x}_1$.
\begin{lem}
\label{second}
There exist finitely many tempered functions $\bm_{2,(r)}$ (respectively $\bm_{1,(r)}$) depending on 
the variable $\vec x_{3}$ only, such that for every element $X\in\CU(\fV^{(2)})$ (respectively 
$X\in\CU(\fV^{(1)})$), the element $\widetilde{X}\,\vec x_{2}$ (respectively $\widetilde{X}\,\vec x_{1}$) 
belongs to the space spanned by  $\{\bm_{2,(r)},\,\bm_{2,(r)}\vec x_2\}$ 
(respectively $\{\bm_{1,(r)},\,\bm_{1,(r)}\vec x_1\}$).
\end{lem}
\begin{proof}
This follows from Lemma \ref{RECB} and the expressions (\ref{INVVECT}) for the invariant vector 
fields. Indeed, the latter implies 
that for every $X\in\CU(\bigoplus_{k=1}^2\fV_{k})$ (respectively $X\in\CU(\fV_{1})$) of strictly positive 
homogeneous degree,  $\widetilde{X}\,\vec x_{2}$ (respectively $\widetilde{X}\,\vec x_{1}$) is either 
zero or one of the entries of the matrix $B$.
\end{proof}
\begin{rmk}\label{secondgen}
Note that in view of the expressions (\ref{INVVECT}) and (\ref{comput}) and by symmetry on 
$\vec{x}_1$ and $\vec{x}_2$ the assertion in Lemma \ref{second}
holds for every element $X$ in $\CU(\s\oplus\s)$ for both variables $\vec{x}_1$ and $\vec{x}_2$.
\end{rmk}

Last, we go to the variable $\vec{x}_0$. The next Lemma is proved using 
the same type of arguments as in the proof of Lemma \ref{RECB}.
\begin{lem}\label{RECC}
There exist  finitely many vectors $\gamma_{(r)}\in\R^2[e^{\pm a_1},e^{\pm a_2}]$ such that for all integers $N_1$ and $N_2$, the 
elements $\widetilde{H}_1^{N_1}\widetilde{H}_2^{N_2}\,\gamma$ consist of a linear combination of the $\gamma_{(r)}$'s.
\end{lem}
Observing that $\widetilde H_ i\vec t$ is proportional to $t_ i$ and that $\widetilde H_ i\Omega_{12}=-\Omega_{12}$,
the Lemmas \ref{RECB} and \ref{RECC} then yield the following result.
\begin{lem}
\label{third}
There exist
finitely many matrices $M_{(r)}\in M_2\left(\R[e^{\pm a_1},e^{\pm a_2}]\right)$ and
 finitely many vectors $v_{(s)}\in\R^2[e^{a_1},a^{a_2}]$
such that for all integers $N_1$ and $N_2$, one has
$$
\widetilde{H}_1^{N_1}\widetilde{H}_2^{N_2}\,\vec{x}_0=M_{N_1,N_2}\vec{x}_0\,+\,\Omega_{12}v_{N_1,N_2}\;,
$$
with
$$
M_{N_1,N_2}\,\in\,\Span\{M_{(r)}\}\quad\mbox{and}\quad v_{N_1,N_2}\,\in\,\Span\{v_{(s)}\}\;.
$$
\end{lem}

 The following result is then a  consequence of Lemmas \ref{first}, \ref{second} and \ref{RECC}.
\begin{cor}\label{TEMPAPPROX}
For every $k=0,\dots,3$, there exists a tempered function $0<\bm_k$ with $\partial_{\vec x_{j}}\bm_k=0$ for every $j\leq k$ and such that for every $X\in\CU(\fV^{(k)})$,
there exists $C_X>0$ with
$$
\big|\widetilde{X}\,\vec{x}_k\big|\leq C_X\,\bm_k\,(1+|\vec{x}_k|)\;.
$$
\end{cor}
\begin{rmk}
In fact the function $\bm_0$ above depends on $\vec x_1,\vec x_2, \vec x_3$ and the functions
$\bm_1$ and $\bm_2$ depend on $\vec x_3$ only (and $\bm_3$ is constant as it should be).
\end{rmk}
We are now able to check the admissibility conditions of Definition \ref{TEMPADM}, for the tempered pair $(\S\times\S,S_{\rm can}^\S)$.
\begin{prop}\label{MULTIPLIERS-ELEM}
Define
\begin{align*}
X_0:=1-H^2_1-H_2^2\;,\quad X_1:=1-\sum_{j=1}^d\,\left((f_j^1)^2+(f_j^2)^2\right)\;,\\
 X_2:=1-\sum_{j=1}^d\,\left((e_j^1)^2+(e_j^2)^2\right)\;,\quad
X_3:=1-E_1^2-E_2^2\;.
\end{align*}
Then the corresponding multipliers 
$\alpha_k:=e^{-iS_{\rm can}^\S}\,\widetilde{X}_k\,e^{iS_{\rm can}^\S}$ ,
$k=0,\dots,3$, satisfy conditions (i) and (ii) of Definition \ref{TEMPADM}.
\end{prop}
\begin{proof}
We start by observing the following expression of the multiplier:
$$
\alpha_k=1+|\vec{x}_k|^2-i\beta_k\,,\quad k=0,\dots,3\;,
$$
where
$$
\beta_k:=\widetilde{X}_{1,k}\,x_{1,k}+\widetilde{X}_{2,k}\,x_{2,k}\;,
$$
with obvious notations.
Then we get
$$
\frac{1}{|\alpha_k|^2}=\frac{1}{(1+|\vec{x}_k|^2)^2+\beta_k^2}\leq\frac{1}{(1+|\vec{x}_k|^2)^2}\;,
$$
and the first condition of Definition \ref{TEMPADM} is satisfied for $C_k=1$ and $\rho_k=2$.
Let now $X\in\CU(\fV^{(k)})$ be of strictly positive order. Then, using Sweedler's \eqref{Sweedler}, notations
we get
$$
\widetilde X\,\alpha_k=\sum_{(X)}\big(\widetilde X_{(1)} \vec x_k\big).\big(\widetilde X_{(2)} \vec x_k\big)-
i\widetilde X\,\widetilde{X}_{1,k}\,x_{1,k}-i\widetilde X\,\widetilde{X}_{2,k}\,x_{2,k}\;.
$$
Since $ X_{(1)}, X_{(2)},{X}_{1,k},{X}_{2,k}\in\CU(\fV^{(k)})$,  Corollary \ref{TEMPAPPROX}  yields
$$
|\widetilde{X}\,\alpha_k|\leq C_1\,\bm^2_k(1+|\vec{x}_k|)^2+C_2\,\bm_k(1+|\vec{x}_k|)\;.
$$
As  $1+|\vec x_k|^2\leq |\alpha_k|$, the second condition of Definition \ref{TEMPADM} is satisfied for  $\mu_k=\bm_k(1+\bm_k)$.
\end{proof}

Last, we prove tameness for the pair $(\S\times\S,S_{\mathrm{can}}^\S)$. 
We start by describing
the behavior of the modular weight of the group $\S$:
\begin{lem}
\label{fd-explicit}
Within the chart \eqref{chartS}, we have the following behavior of the
the modular weight $\fd_\S$ of an elementary normal $\bf j$-group $\S$:
$$
\fd_\S\asymp\big[ (a,v,t)\mapsto \cosh a+\cosh2a+|v|(1+e^{2a}+\cosh a)+|t|(1+e^{2a})\big]\;.
$$
\end{lem}
\begin{proof}
Within the decomposition $\R H\oplus V\oplus \R E$ of $\s$, and within the chart \eqref{chartS}
of $\S$,  a quick computation gives
\begin{align}
\label{matrixform}
\Ad_{(a,v,t)}&=
\begin{pmatrix}
\id&0&0\\
-e^a v&e^a&0\\
2e^{2a}t&\tfrac12e^{2a}\omega^0(v,.)&e^{2a}
\end{pmatrix}\;,\\
\Ad_{(a,v,t)^{-1}}&=
\begin{pmatrix}
\id&0&0\\
 v&e^{-a}&0\\
-2t&-\tfrac12e^{-a}\omega^0(v,.)&e^{-2a}
\end{pmatrix}\;.\nonumber
\end{align}
Now, the result follows from the equivalence of the  operator
and Hilbert-Schmidt norms on the finite dimensional vector space $\End(\s)$, together with obvious estimates.
\end{proof}
\begin{cor}
\label{MWP}
The tempered pair $(\S\times\S,S_{\mathrm{can}}^\S)$ is tame.
\end{cor}
\begin{proof}
Combining the last statement of Lemma \ref{fd-semi} with
Lemma \ref{fd-explicit}, we get
$$
\fd_{\S\times \S}(a_1,v_1,t_1;a_2,v_2,t_2)\geq C \big( \cosh a_1+\cosh a_2+|v_1|+|v_2|+|t_1|+|t_2|)
\;,
$$
so that with the relations \eqref{comput2} in mind, we see that there exists $C'>0$ with
$$
\big( |\vec x_0|^2+ |\vec x_1|^2+ |\vec x_2|^2+ |\vec x_3|^2\big)^{1/2}\leq C'\,
\fd_{\S\times \S}^4\;.
$$
According to Definition \ref{tame}, we may set $\mu_\phi=\fd_{\S\times \S}^4$ which is
 tempered since  $\fd_{\S\times \S}$ is by Lemma \ref{fd-semi} and Lemma \ref{fd-tempered}. 
\end{proof}

 We summarize all this by stating the main result of this section:
\begin{thm}
\label{can-temp}
Let $\S$ be an elementary normal $\bf j$-group and let $S_{\rm can}^\S$ be the smooth function on $\S\times\S$ given in Definition \ref{df-Scan}.
Then, the pair $(\S\times\S,S_{\mathrm{can}}^\S)$ is tempered,  admissible and tame.
\end{thm}
\begin{rmk}\label{newdecomp}
From Remark \ref{secondgen} and the above discussion, we observe that setting $\fV_{12}:=\fV_1\oplus\fV_2$ yields a decomposition
into \emph{three} subspaces: $\s\oplus\s=\fV_0\oplus\fV_{12}\oplus\fV_3$ also underlying admissibility but with associated elements
$X_0$, $X_3$ and $X_{12}:=X_1+X_2$. The corresponding multipliers are $\alpha_0$, $\alpha_3$ and $\alpha_{12}=\alpha_1+\alpha_2$.
\end{rmk}

\section{Tempered pairs  for normal $\bf j$-groups} 
\label{2pencore}
Let $\B$ be a normal $\bf j$-group with Lie algebra $\b$. We first observe:
\begin{lem}
\label{mwdi}
Let $\B=\B'\ltimes \S_1$ be a Pyatetskii-Shapiro decomposition with $\S_1$ elementary normal.
Then, there exists $C>0$ such that for every $g_1=(a,v,t)\in\S_1$, 
$g'\in\B'$,  we have:
$$
\fd_{\B}(g_1g')\geq C\big(\fd_{\B'}(g')+ \cosh 2a+ |v|+|t|\big)\;.
$$
\end{lem}
\begin{proof}
By Lemma \ref{fd-semi}, we know that there exists $C_1>0$ such that for all
$g_1\in\S_1$ and $g'\in\B'$, we have,
$$
\fd_{\B}(g_1g')\geq C_1\big(
\fd_{\B'}(g')+\big(1+| \Ad_{g_1}\bR_{g'}|^2+| \bR_{g'^{-1}}\Ad_{g_1^{-1}}|^2
\big)^{1/2}\big)\;,
$$
In the decomposition  $\R H\oplus V\oplus\R E$ of $\s$, the $\B'$-action may be expressed under the following matrix form:
$$
\bR_{g'}=
\begin{pmatrix}
\id&0&0\\
0&A(g')&0\\
0&0&\id
\end{pmatrix}\;,
$$
where $A(g')$ is the matrix in the linear symplectic group
 ${\rm Sp}(V_1,\omega^0_1)$ as defined in Proposition \ref{proprietes}
(iv). From \eqref{matrixform} and setting as usual $g_1=(a,v,t)$, we therefore obtain:
\begin{align*}
\Ad_{g_1}\bR_{g'}&=
\begin{pmatrix}
\id&0&0\\
-e^a v&e^aA(g')&0\\
2e^{2a}t&\tfrac12e^{2a}\omega^0(v,A(g').)&e^{2a}
\end{pmatrix}\;,
\\
\bR_{g'^{-1}}\Ad_{g_1^{-1}}&=
\begin{pmatrix}
\id&0&0\\
 A(g')^{-1}v&e^{-a}A(g')^{-1}&0\\
-2t&-\tfrac12e^{-a}\omega^0(v,.)&e^{-2a}
\end{pmatrix}\;.
\end{align*}
Using once again the equivalence between operator and Hilbert-Schmidt norms on
$\End(\s)$, we deduce for some $C_2>0$:
\begin{align*}
\big(1+| \Ad_{g_1}\bR_{g'}|^2+| \bR_{g'^{-1}}\Ad_{g_1^{-1}}|^2
\big)^{1/2}&\geq C_2\big(\cosh 2a+|v| \cosh a +|t|\big)
\\ &\geq C_2\big(\cosh 2a+ |v|+|t|\big)\;,
\end{align*}
and the proof follows.
 \end{proof}

Now, we let also $\b=\a\oplus\n$ be a decomposition of the Lie algebra of a normal $\bf j$-group
$\B$, with $\n$ the nilradical of $\b$ and $\a$ its orthogonal complement. 
It follows then that $\a$ is an abelian subalgebra, so that $\b=\a\ltimes\n$ and the group $\B$ may be identified to its Lie algebra $\b$ with product
\[ (a,n) \cdot (a',n') = \big(a+a', (e^{-\ad a'}n)\cdot_{\rm CBH} n'\big)\;, \]
where $n\cdot_{\rm CBH}n'$  denotes the Baker-Campbell-Hausdorff  series in the Lie algebra $\n$, which is finite 
since $\n$ is nilpotent. 
\begin{dfn}
\label{adapted}
Let $\{H_j\}_{j=1}^n$ and $\{N_j\}_{j=1}^m$ be bases of $\a$ and $\n$ respectively. The coordinates system
\begin{align*}
\R^{n+m}\to \a\oplus\n\;,\quad
&(a_1,\dots,a_n,n_1,\dots,n_m) \mapsto \\
&\big(\arcsinh(a_1)H_1+\dots+\arcsinh(a_n)H_n,
n_1N_1+\dots+n_mN_m\big)\;,
\end{align*}
are said to be {\bf adapted tempered coordinates} for $\B$.
\end{dfn}

\begin{lem}
\label{lem:mult_inv_temp}
In any adapted tempered coordinates on $\B$, the multiplication and inverse operations are tempered maps $\R^{n+m}\times\R^{n+m}\to\R^{n+m}$ and $\R^{n+m}\to \R^{n+m}$ respectively.
\end{lem}
\begin{proof}

 Let $a_1,\dots,a_n,n_1,\dots,n_m$ be adapted tempered coordinates on $\B$ as in the above definition. Then, since
\[ \sinh(a+a') = \sinh a\cosh a'+\cosh a\sinh a'\;, \] 
the $\{a_i\}$-coordinates of the multiplication of $x,x'\in\R^{n+m}$ read
\begin{align*}
\sinh\big(\arcsinh a_i+\arcsinh a'_i\big)
&= a_i\sqrt{1+{a'_i}^2} + a'_i\sqrt{1+a_i^2}\;,
\end{align*}
so that they clearly are tempered functions in the $a_i,a'_i$ variables.
 For the $\n$ part, recall that there is a decomposition in real root spaces $\n=\bigoplus_\alpha \n_\alpha$ for the adjoint action of $\a$. Now if $n'\in\n_\alpha$,
we have
\begin{align*}
&e^{\ad(\arcsinh (a_1) H_1+\dots+\arcsinh(a_n)H_n)}n' 
= e^{\alpha(H_1)\arcsinh(a_1)+\dots+\alpha(H_2)\arcsinh(a_n)} n' \\
&= \Big(a_1+\sqrt{1+a_1^2}\Big)^{\alpha(H_1)}
\dots
\big(a_n+\sqrt{1+a_n^2}\big)^{\alpha(H_n)} n'\;,
\end{align*}
which is a tempered function in $a_1,\dots,a_n$. As the CBH product in a nilpotent group is polynomial, linearly decomposing $n'_1N_1+\dots n'_mN_m$ along the root space decomposition and using the above computation, we get that the $n_i$ coordinates of the product of $x$ and $x'$ are tempered in all variables.
 For the inverse map, as $(a,n)^{-1}=(-a,-e^{-\ad a}n)$, the above computation also shows the result.
\end{proof}

\begin{lem}
\label{lem:j-algebra-extension-morphism-tempered}
Let $\b=\b'\ltimes \s$ be a Pyatetskii-Shapiro decomposition of a normal $\bf j$-algebra $\b$, with $\s$ an elementary normal $\bf j$-algebra and  with corresponding Lie group decomposition $\B=\B'\ltimes \S$. Denote $\bR:\B'\to \Aut(\S)$ the associated extension homomorphism. Then in any adapted tempered coordinates for $\B'=\a'\oplus\n'$, with $\dim(\a')=n'$, $\dim(\n')=m'$ and
 $\S=\a\oplus\n$ 
 (recall that by Proposition \ref{P-S}, $\n$ is an Heisenberg Lie algebra, thus
 nilpotent), with  $\dim(\a)=1$, $\dim(\n)=m$, $\bR$ is a tempered map $\R^{n'+m'}\times\R^{1+m}\to\R^{1+m}$.
\end{lem}
\begin{proof}
 Let $a_1,\dots,a_{n'},n_1,\dots,n_{m'}$ and $a_{n'+1}, n_{m'+1},\dots,n_{m'+m_1}$ be adapted tempered coordinates for $\B'$ and $\S$ respectively. The group $\B'$ acts trivially on 
$H_{n'+1}$, the generator of  $\a$. Moreover, the coordinates $a_1,\dots,a_{n'+1},n_1,\dots,n_{m'+m_1}$ are adapted tempered coordinates for $\B$. Indeed, one knows \cite[pages 56-57]{PS}  that the infinitesimal action of $H_1,\dots,H_{n'}$ is real semi-simple
with spectrum contained in $\{-\frac{1}{2},0,\frac{1}{2}\}$.
Denote $i':\B'\to \B$ and $i:\S\to \B$ the inclusions seen through the coordinates.
Now by Lemma~\ref{lem:mult_inv_temp}, the map
\begin{align*}
(x',x)\in\B'\times\S\mapsto i'(x') \cdot i(x)\in\B\;,
\end{align*}
is tempered. But the $\n$ part of that product is exactly $\bR_{x'}(x)$ and so, this concludes the proof.
\end{proof}

We are now prepared to state and prove the main result of this chapter.
\begin{thm}
\label{TASP}
Let $\B$ be a normal $\bf j$-group  with  Pyatetskii-Shapiro decomposition 
$\B=\big(\S_N\ltimes\dots\big)\ltimes\S_1$. Parametrizing the elements $g,g'\in\B$ as $g=g_1g_2\dots g_N$ and $g'=g_1'g_2'\dots g_N'$ with $g_i,g_i'\in\S_i$, we define 
\begin{equation*}
 S_{\rm can}^\B:\B\times\B\to\R\,,\quad (g,g')\mapsto\sum_{i=1}^N S_{\rm can}^{\S_i}(g_i,g_i')\;,
\end{equation*}
where $ S_{\rm can}^{\S_i}$ is the canonical phase of $\S_i$ given in Definition \ref{df-Scan}. 
Then the pair $(\B\times\B,S_{\rm can}^\B)$ is tempered  admissible and tame.
\end{thm}
\begin{proof}
 We will use an induction over $N$,  the number of elementary factors in $\B$.
 Accordingly, we set $\B=\B'\ltimes_{\bR}\S$, 
with $\B':=\big(\S_N\ltimes\dots\big)\ltimes\S_2$ and $\S:=\S_1$. 
We then  observe 
 that $\B\times\B=(\B'\times\B')\ltimes_{\bR\times\bR}(\S\times\S)$
and from Lemma \ref{fd-semi} and  Lemma \ref{lem:j-algebra-extension-morphism-tempered},
that the extension homomorphism $\bR\times\bR=:\bR^2$ is tempered within adapted coordinates.
By Theorem \ref{can-temp}, the pair $(\S\times\S,S_{\rm can}^\S)$ is tempered and admissible. By induction hypothesis, the latter also holds for $(\B'\times\B',S_{\rm can}^{\B'})$.
Moreover,  Equations (\ref{comput2})  tell that, in the ``elementary" case of $\S\times\S$, the adapted tempered coordinates and the coordinates associated to the phase function are related to one another through
a tempered diffeomorphism.  By induction hypothesis, the latter also holds for $\B'\times\B'$. 
Obviously,  the extension homomorphism $\bR\times\bR$ is then tempered within the coordinates associated to the phase functions
as well. 
Note that under the parametrization $g=g_1g',h=h_1h'\in\B'$, $g_1,h_1\in\S_1$, $g',h'\in\B'\in\B'$, the multiplication and inverse maps of $\B$ become:
$$
gh=g_1\bR_{g'}(h_1)g'h'\,,\quad g^{-1}=\bR_{{g'}^{-1}}(g_1^{-1}){g'}^{-1}\;,
$$
and similarly for $\B\times\B$.
From this and the  temperedness of the extension homomorphism $\bR\times\bR$, we see that 
temperedness of the multiplication and inversion laws in $\B\times\B$ will immediately follow once 
we will have shown that the map \eqref{TEMPCOORD} is a global diffeomorphism from $\B\times\B$ 
to $(\b\oplus\b)^\star$. We will return to this question while  
examining  admissibility. To this aim, let us set $G_1:=\B'\times\B'$, 
$G_2:=\S\times\S$ and denote respectively 
 by $\g_1$ and $\g_2$ their Lie algebras. Let us also set
 $S_1:=S_{\rm can}^{\B'}$,  $S_2:=S_{\rm can}^\S$,
  and let us assume, by induction hypothesis, that the pair
 $(G_1,S_1)$ is admissible, with associated  decomposition 
 $\g_1=\oplus_{k=0}^{N_1}\mathfrak{W}_k$.
 Let us consider an adapted  basis  of $\g_1$, $\{{}_1w_k\}$, $k=1,\dots,\dim(\g_1)$ with associated 
 coordinates ${}_1(b)_k:=\widetilde{{}_1w_k}.S_1(b)$ on $G_1$.
 Similarly, let us consider the basis $\{{}_2w^r_j\,|\,r=1,2\;;\;j=0,1,2,3\}$ of $\g_2$ adapted to the 
 decomposition (\ref{decomposition}), where, for the values 1 and 2, $j$ consists of a multi-index
 and where $r$ labels the copies of $\S$ in $\S\times\S$. 
 Accordingly, 
 we have the associated coordinate system (\ref{comput2}) on $G_2$ that now reads ${}_2(x)^r_j:=
 \widetilde{{}_2w^r_j}.S_2(x)$.
 
On $G:=G_1\ltimes_{\bR^2}G_2$, with $S(xb):=S_1(b)+S_2(x)$, $x\in G_2$, $b
 \in G_1$, we then compute that:
 $$
 {}_1(xb)_k:=\widetilde{{}_1w_k}.S(xb)={}_1(b)_k\quad\mbox{\rm  and}\quad{}_2(xb)^r_j
 :=\widetilde{{}_2w^r_j}.S(xb)=
 \widetilde{\bR^2_b({}_2w^r_j)}.S_2(x)\;.
 $$
 Hence it suffices to look at the properties of the multipliers within $\g_2$.
 From Pyatetskii-Shapiro's theory, we know that for the values $0$ and $3$ of $j$, the action of $G_1$ is trivial:
 $\bR^2_b({}_2w^r_j)={}_2w^r_j$. Hence:
 $$
 {}_2(xb)^r_j={}_2(x)^r_j\,,\quad\forall j\in\{0,3\}\;.
 $$
 Hence  it suffices to look at the properties of the multipliers within the subspace 
 $V\times V=\fV_1\oplus\fV_2$.
 For $j=1,2$, however, the action is not trivial but stabilizes component-wise   $V\times V$. Accordingly, we set:
 $$
 \bR^2_b({}_2w^r_j)=:\sum_{p=1}^2[\bR^2_b]^p_j\,{}_2w^r_p\,,\quad \forall j\in\{1,2\}\;,
 $$
where, again, $p$ is a multi-index. We therefore have:
\begin{equation}\label{XB2}
 {}_2(xb)^r_j=\sum_{p=1}^2[\bR^2_b]^p_j\,{}_2(x)^r_p\,,\quad\forall j\in\{1,2\}\;.
\end{equation}
 From Pyatetskii-Shapiro's theory, we know that the linear operator $[\bR^2_b]^p_j$ 
 is symplectic, thus it has
  jacobian one.
In particular, this implies that the map
 \eqref{TEMPCOORD} is a global diffeomorphism from $G$ to $\g^\star$. (Hence,  we have completed the
 proof of temperedness of the pair $(\B\times\B,S_{\rm can}^\B)$.)
We now consider the ordered decomposition:
$$
\g=\g_2\oplus\g_1=\Big(\oplus_{j=0}^3\fV_j\Big)\bigoplus\Big(\oplus_{k=0}^{N_1}\mathfrak{W}_k\Big)\;,
$$
where indices occurring on the left ($\g_2$) are considered as lower than the one on  the right ($\g_1$).
 Within this setting, we compute that for every element $X\in\CU(\g_2)$:
$$
{}_{(2)}\alpha_X(xb):=e^{-iS(xb)}\big(\widetilde{X}\,e^{iS}\big)(xb)={}_2\alpha_{\bR^2_b(X)}(x)\;,
$$
where ${}_2\alpha_X:=e^{-iS_2}\big(\widetilde{X}.e^{iS_2}\big)$ denotes the multiplier on $G_2$.
Again, for the extreme values of $j$, we observe that:
$$
{}_{(2)}\alpha_X(xb)={}_2\alpha_{X}(x)\,,\quad\forall X\in\CU(\fV_0\oplus\fV_3)\;,
$$
so in these cases the properties underlying admissibility are trivially satisfied.
For $j=1,2$, we have  with the notation $X_j:=1-\sum_{r=1}^2({}_2w^r_j)^2$ of Proposition \ref{MULTIPLIERS-ELEM}:
$$
\bR^2_b(X_j)=1-\sum_{r=1}^2\left([\bR^2_b]^p_j\,{}_2w^r_p\right)^2\;,
$$
which leads to
$$
{}_{(2)}\alpha_{X_j}(xb)=1\,+\,\sum_{r=1}^2\left({}_2(xb)^r_j\right)^2\,-\,i\,\sum_{r=1}^2\,[\bR_{b_r}]^{p_r}_j\,[\bR_{b_r}]^{p'_r}_j\,
\widetilde{w^r_{p_r}}.{}_2(x)^r_{p'_r}\;,
$$
where $b=(b_1,b_2)\in\B'\times\B'=G_1$. This gives the property $(i)$ of Definition \ref{TEMPADM}.
From the expression (\ref{XB2}) and the structure of the elementary case (Lemma \ref{second}), we 
then observe:
\begin{align}
\label{deg3}
\widetilde{A}.{}_{(2)}\alpha_{X_j}=0\,,\qquad\forall A\in\CU(V\times V)\,,\; \deg(A)\geq 3\;.
\end{align}
Also, setting $-i\beta_{X_j}(xb):=-\,i\,\sum_{r=1}^2\,[\bR_{b_r}]^{p_r}_j\,[\bR_{b_r}]^{p'_r}_j\,
\widetilde{w^r_{p_r}}.{}_2(x)^r_{p'_r}$, we deduce from the expressions (\ref{INVVECT}) and (\ref{comput})  that, for every $A\in V\times V$: $\widetilde{A}.\beta_{X_j}=0$.
From the expression (\ref{XB2}) and setting ${}_{2}(x)_{12}:=({}_{2}(x)_{1},{}_{2}(x)_{2})$, we then deduce that for every $A\in\CU(V\times V)$ of strictly positive degree: 
 \begin{align*}
 \big|\widetilde{A}\big({}_{(2)}\alpha_{X_1}(xb)+{}_{(2)}\alpha_{X_2}(xb)\big) \big|
 =\big|\widetilde{A}|{}_{2}(xb)_{12}|^2\big|&=
\big|\widetilde{\bR^2_b(A)}_x|\bR^2_b({}_{2}(x)_{12})|^2\big|\\
&\leq |\bR^2_b|^{\deg(A)+2}C_A\,\mu_{12}(x)\,|({}_{2}(x)_{12})|^2\;,
\end{align*}
 where the last estimate is obtained from Corollary \ref{TEMPAPPROX}.
Since 
$$
|({}_{2}(x)_{12})|^2=|\bR^2_{b^{-1}}\bR^2_b({}_{2}(x)_{12})|^2\leq
|\bR^2_{b^{-1}}|^2|({}_{2}(xb)_{12})|^2\;,
$$
we then get
\begin{align*}
& \big|\widetilde{A}
 \big({}_{(2)}\alpha_{X_1}(xb)+{}_{(2)}\alpha_{X_2}(xb)\big)\big|
 \leq|\bR^2_b|^{\deg(A)+2}C_A\,\mu_{12}(x)\,|\bR^2_{b^{-1}}|^2|({}_{2}(xb)_{12})|^2 \\
&\qquad\qquad\qquad \leq C_A\,|\bR^2_b|^{\deg(A)+2}\,|\bR^2_{b^{-1}}|^2\mu_{12}(x)\,\left|{}_{(2)}\alpha_{X_1}(xb)+{}_{(2)}\alpha_{X_2}(xb)\right|\;.
\end{align*}
But by \eqref{deg3},  we know that we may assume $\deg(A)\leq2$, hence
\begin{align*}
 \big|\widetilde{A}\big({}_{(2)}\alpha_{X_1}(xb)+{}_{(2)}\alpha_{X_2}(xb)\big)\big|
 \leq C_A\,\fd^6_{G}(b)\,\mu_{12}(x)\,\big|{}_{(2)}\alpha_{X_1}(xb)+{}_{(2)}\alpha_{X_2}(xb)\big|\;.
 \end{align*}
 Defining the element $\mu'_{12}(xb):=\fd^6_{G}(b)\,\mu_{12}(x)$ yields admissibility at the level of  $V\times V$.

Last, tameness also follows by an induction argument,  using Lemma \ref{fd-semi} and Lemma \ref{mwdi}.
\end{proof}

\begin{rmk}
\label{lastS}
Observe that Remarks \ref{LRORNOT}, Lemma \ref{fd-tempered} and  Lemma \ref{mwdi} show that
on the one-variable  Schwartz space  $\CS^{S^\B_{\rm can}}(\B,\CE)$  associated with the two-variable
 tempered, admissible and tame pair $(\B\times\B,S^\B_{\rm can})$ (see Definition \ref{S1P})
one has the same topology associated with
 the equivalent semi-norms
\begin{align*}
f\in\CS^S(\B,\CE)\mapsto\sup_{X\in\,\CU_k(\b)}\sup_{x\in \B}\Big\{\frac
{\fd_\B(x)^n\big\|\widetilde X\,f(x)
\big\|_j}{|X|_k}\Big\}\;,\qquad j,k,n\in\N\;,
\end{align*}
or even with
\begin{align*}
f\in\CS^S(\B,\CE)\mapsto\sup_{X\in\,\CU_k(\b)}\sup_{x\in \B}\Big\{\frac
{\fd_\B(x)^n\big\|\underline X\,f(x)
\big\|_j}{|X|_k}\Big\}\;,\qquad j,k,n\in\N\;.
\end{align*}
 \end{rmk}

\chapter{Non-formal star-products}
\label{NFSP}

 This chapter is devoted to the construction of an  infinite dimensional parameter family
 of non-formal star-products on every negatively curved K\"ahlerian Lie group. In section \ref{SPNG}, we 
 first review
 the construction of one of us of such non-formal star-products, living on a space of distributions. 
 Their covariance properties and their relations with the Moyal product are given there. In 
 section \ref{OIFDF}, we apply the results of chapters \ref{OSIL} and \ref{GeometrySSS}
 to give a proper interpretation of the star-product formula as an oscillatory integral and 
 eventually prove that the Fr\'echet space $\CB(\B)$ becomes a Fr\'echet algebra for all these
 star-products. This is proved in Theorem \ref{achievement}, the main result of this chapter and,
 we believe, and important result in itself.

\section{Star-products on  normal $\bf j$-groups}\label{SPNG}
We consider an  elementary normal $\bf j$-group $\S$ viewed as a symplectic symmetric spaces as in  section \ref{H-group}. We start by recalling the results obtained in
\cite{Bi02, Bi07}.
\begin{dfn} Set $\tilde{\S}:=\{(a,v,\xi)\}=\R\times\R^{2d}\times\R$.
The {\bf twisting map} is the smooth one-parameter family of diffeomorphisms defined as
\begin{equation*}
\phi_\theta:\tilde{\S}\to\tilde{\S}:(a,v,\xi)\mapsto
\Big( a,\sech\left(\tfrac{\theta}{4}\xi\right)v,\tfrac{2}{\theta}\sinh\left(\tfrac{\theta}{2}\xi\right) \Big)
\,,\quad\theta\in\R^*\;.
\end{equation*}
\end{dfn}

Let $\CS(\S)$ be the Euclidean Schwartz space of $\S$,
i.e$.$ the ordinary Schwartz space in the coordinates \eqref{chartS}.  Accordingly, let
$\CS(\S)'$ be the dual space of tempered distributions.
 Let us also denote by
\begin{equation*}
\big(\CF u\big)(a,v,\xi):=\int_{-\infty}^\infty e^{-i\xi t}u(a,v,t)\,{\rm d}t\;,
\end{equation*}
the partial Fourier transform in the $t$-variable. 
For $\gamma>0$, we  let ${\mathcal O}_{C,\gamma}(\R^m)$ be the subset of smooth functions, 
the derivatives of which
 are uniformly polynomially 
bounded:
\begin{align*}
&\mathcal O_{C,\gamma}(\mathbb R^m):=\\
&\big\{f\in C^\infty(\mathbb R^m):\exists r>0:\forall \alpha\in\mathbb N^m:\exists C_\alpha>0,\,|\partial^\alpha f(x)|\leq C_\alpha(1+|x|)^{r-\gamma|\alpha|}\big\}\;.
\end{align*}
 Note that ${\mathcal O}_{C,1}(\R^m)$ is the space of Grossmann-Loupias-Stein symbols, 
traditionally written $\mathcal K(\R^m)$.
\begin{dfn}
\label{gros-gras}
We denote by $ {\bf \Theta}$,  the subspace of $C^\infty(\R,\C)$ constituted by the elements
${\tau}$ such that  $\exp\circ\pm{\tau}$  belong to the space ${\mathcal O}_{C,1}(\R,\C)$,  
and normalized such that $\tau(0)=0$. 
\end{dfn}

Let ${\tau}_0$ be the element of $C^\infty(\tilde\S)$, given by:
$$
{\tau}_0:=\tfrac12\log\, \circ\,{\rm Jac}_{\phi^{-1}_\theta} \;.
$$
Viewed as a function of its last variable only,  ${\bf \tau}_0$ belongs to ${\bf \Theta}$. Indeed, we have:
$$
{\rm Jac}_{\phi^{-1}_\theta}(a,v,\xi)=2^{-d}\frac{\Big(1+\sqrt{1+\tfrac{\theta^2\xi^2}4}\Big)^d}{\sqrt{1+\tfrac{\theta^2\xi^2}4}}\;.
$$

To an element ${\tau}\in{\bf {\bf \Theta}}$, one associates a function on $\tilde\S$
 by 
$[(a,v,\xi)\mapsto\tau(\xi)]$ and, to simplify the notations, we still denote this function
by $\tau$. One then defines a linear injection:
\begin{equation}
\label{eq:T_tau}
T_{\theta,{\tau}}:=\CF^{-1}\circ\exp(\tau_0\,-\tau)\circ(\phi^{-1}_\theta)^\star\circ\CF\;:\CS(\S)\to
\CS(\S)'\;,
\end{equation}
 where, by a slight abuse of notation, we identify a function with the linear operator of 
point-wise multiplication by this function.   We make the following obvious but important observation:
\begin{lem}
\label{bof}
Let $\tau\in\Theta$. Then the inverse of the map $T_{\theta,{\tau}}$, given by
 $$
T_{\theta,{\tau}}^{-1}=\CF^{-1}\circ(\phi_\theta)^\star\circ\exp(-\tau_0+\tau)\circ\CF\;,
$$
defines a continuous linear injection from $\CS(\S)$ to itself. Moreover, 
 $T_{\theta,{\tau}}^{-1}$ extends to a unitary operator on $L^2(\S)$ if and only if $\tau$ is 
 purely imaginary.
\end{lem}
\begin{proof}
As the partial Fourier transform $\CF$ is an homeomorphism from $\CS(\S)$ to $\CS(\tilde\S)$,
and as $\exp(-\tau_0+\tau)$ belongs to ${\mathcal O}_{C,1}(\R,\C)$, a subspace of the Schwartz
multipliers ${\mathcal O}_M(\R,\C)$, it suffices to show that the map $f\mapsto f\circ \phi_\theta$,
is continuous on $\CS(\tilde\S)$. So, let $f\in\CS(\tilde\S)$. Then, as
$$
f\circ \phi_\theta(a,v,\xi)=f\big( a,\sech\left(\tfrac{\theta}{4}\xi\right)v,\tfrac{2}{\theta}\sinh\left(\tfrac{\theta}{2}\xi\right) \big)\;,
$$
we deduce from
$$
\partial_a^{k}(f\circ \phi_\theta)=(\partial_a^{k}f)\circ \phi_\theta\,,\quad
\partial_v^m(f\circ \phi_\theta)=\sech\left(\tfrac{\theta}{4}\xi\right)^{m}(\partial_v^{m}f)\circ \phi_\theta\;,
$$
and from an iterative use of the relation 
$$
\partial_\xi(f\circ \phi_\theta)(a,v,\xi)=-\frac\theta4\frac{\sinh\left(\tfrac{\theta}{4}\xi\right)}
{\cosh\left(\tfrac{\theta}{4}\xi\right)^2} \,v\,(\partial_v f)\circ \phi_\theta(a,v,\xi)+\cosh\left(\tfrac{\theta}{4}\xi
\right)(\partial_\xi f )\circ \phi_\theta(a,v,\xi)\;,
$$
that any (ordinary) derivatives of $f\circ \phi_\theta$ is a finite linear combination  of derivatives 
of $f$ composed with $\phi_\theta$ and with coefficient in the polynomial ring $\R[v,e^{\theta\xi},
e^{-\theta\xi}]$. To conclude  the first claim, we then observe that as a derivative of a Schwartz function is a Schwartz 
function,  any derivatives of $f$ composed with $\phi_\theta$ is  bounded (in absolute value)  by 
$(a^2+|v|^2+\cosh(\theta\xi))^{-n}$, with $n$ arbitrary. For the second claim, observe that
 $$
T_{\theta,{\tau}}^{-1}=\CF^{-1}\circ(\phi_\theta)^\star\circ\,{\rm Jac}_{\phi^{-1}_\theta}^{-1/2}\circ
\exp(\tau)\circ\CF
=\CF^{-1}\circ\,{\rm Jac}_{\phi_\theta}^{1/2}\circ(\phi_\theta)^\star\circ
\exp(\tau)\circ\CF\;,
$$
which entails since $\CF$ and ${\rm Jac}_{\phi_\theta}^{1/2}\circ(\phi_\theta)^\star$ are unitary,
that $T_{\theta,{\tau}}^{-1}$ is unitary if and only if the operator of multiplication by 
$\exp(\tau)$ is unitary, which is equivalent to $\Re(\tau)=0$.
\end{proof}
\begin{rmk}
Denoting by  $T_{\theta,{\tau}}^*$ the formal adjoint of
$T_{\theta,{\tau}}$  with respect to the inner product of $L^2(\S)$, then we have
$T_{\theta,{\tau}}^*=T_{\theta,{-\overline\tau}}^{-1}$ and 
thus, Lemma \eqref{bof} implies that $T_{\theta,{\tau}}^*$ is continuous on $\CS(\S)$ too.
\end{rmk}

Let $\omega^0$ be the standard symplectic structure of $\R^{2d+2}\simeq T^* \R^{d+1}$
 and let  $\star^0_\theta$
be the Moyal  product in its integral form on $\CS(\R^{2d+2})$. Recall that the latter product is
by definition
the composition law of symbols in the Weyl pseudo-differential calculus and that it is  given by
\begin{align*}
 f_1 \star^0_\theta f_2 (x) = \frac{1}{(\pi\theta)^{2(d+1)}} \int_{\R^{2d+2}\times\R^{2d+2}}
e^{\frac{2i}{\theta}S_0(x,y,z)} f_1(y) \,f_2(z) \; {\rm d}y \; {\rm d}z \;,
\end{align*}
where  $S_0(x,y,z):=\omega^0(x,y)+\omega^0(y,z)+\omega^0(z,x)$.
For $\tau\in{\bf \Theta}$, denoting by 
$$
\CE_{\theta,\tau}(\S):=T_{\theta,\tau}\big(\CS(\S)\big)\;,
$$ 
the range subspace of $T_{\theta,\tau}$
in the tempered distribution space $\CS(\S)'$, one has the inclusions
\[
\CS(\S)\,\subset\,\CE_{\theta,\tau}(\S)\,\subset\,C^\infty(\S)\;.
\] 
We consider the
linear isomorphism:
\begin{equation*}
T_{\theta,\tau}^{-1}:\CE_{\theta,\tau}(\S)\to\CS(\S)\;.
\end{equation*}
Identifying $\S\simeq\R^{2d+2}$ by means of the global coordinate system \eqref{chartS}, we transport 
under $T_{\theta,\tau}$ the Moyal product on $\CS(\R^{2d+2})\simeq\CS(\S)$. This yields an 
associative product:
$$
\star_{\theta,\tau}:\CE_{\theta,\tau}(\S)\times \CE_{\theta,\tau}(\S)\to \CE_{\theta,\tau}(\S)\;,
$$
given by
\begin{equation}
\label{SP}
f_1\star_{\theta,\tau}f_2:=T_{\theta,\tau}\big(T_{\theta,\tau}^{-1}(f_1)\star^0_\theta T_{\theta,\tau}^{-1}(f_2)\big)\,,
\quad f_1,f_2\in\CE_{\theta,\tau}(\S)\;.
\end{equation}
The associative algebra $\big(\CE_{\theta,\tau}(\S),\star_{\theta,\tau}\big)$, endowed with the 
Fr\'echet algebra structure transported  under  $T_{\theta,\tau}$ from $\CS(\R^{2d+2})$, satisfies the 
following properties \cite{Bi02,BiMas}:
\begin{thm}
\label{explicit}
Let $\tau\in{\bf {\bf \Theta}}$ and $\theta\neq0$. Then,
\begin{enumerate}
\item[(i)]
For all compactly supported $u,v\in\CE_{\theta,\tau}(\S)$, one has the  integral representation:
\begin{equation}\label{PROD}
u\star_{\theta,\tau}v=\int_{\S\times\S}\,K_{\theta,\tau}(x_1,x_2)\,R^\star_{x_1}(u)\,\,R^\star_{x_2}(v)\,\,
{\rm d}_\S(x_1)\,{\rm d}_\S(x_2)\;,
\end{equation}
where the two-point kernel is given by
\begin{align}
\label{2-point-kernel}
K_{\theta,\tau}(x_1,x_2) &:=
(\pi\theta)^{-2(d+1)}
\,A_{{\theta,\tau}}(x_1,x_2)
\,\exp\big\{\tfrac{2i}{\theta}S_{\mathrm{can}}^\S(x_1,x_2)\big\}\;,
\end{align}
with, in the coordinates\footnote{As usual, we set $x_j=(a_j,v_j,t_j)\in\R^{2d+2}$} \eqref{chartS}:
\begin{align*}
&A_{{\theta,\tau}}(x_1,x_2):=\\
&A_{\rm can}^\S(x_1,x_2)\,\exp\big\{\tau\big(\tfrac2\theta\sinh2a_1\big)+
\tau\big(-\tfrac2\theta\sinh2a_2\big)-\tau\big(\tfrac2\theta\sinh(2a_1-2a_2)\big)\big\} \;,
\end{align*}
and with $S_{\rm can}^\S$ and $A_{\rm can}^\S$ defined in \eqref{SEXPL} and \eqref{AEXPL}.
\item[(ii)]
The product $\star_{\theta,\tau}$ is equivariant under the automorphism group
of the symplectic symmetric space $(\S,s,\omega^\S)$: for all elements $g$ of $\Aut(\S,s,\omega^\S)$ and $u,v\in\CD(\S)$, one has
\[
g^\star( u)\star_{\theta,\tau}g^\star( v)=g^\star(u\star_{\theta,\tau}v)\;.
\]
\end{enumerate}
\end{thm}

\begin{rmk} 
Observe that when $\Im(\tau)\ne 0$, then the amplitude $A_{{\theta,\tau}}$
also contribute to the phase of $K_{\theta,\tau}$. However, as the amplitude is a function
of the variables in $\a\times\a:=\R H\times \R H$ only, in the two dimensional case (i.e$.$
$V=\{0\}$),
 there is no $\tau\in\Theta$ such that
the combined phase coincides with those found in \cite{GGV}. In particular, this means that the upper half-plane is a symplectic 
manifold that supports two non-isomorphic Moyal quantizations.
\end{rmk}

Consider a normal $\bf j$-group decomposed, following Proposition \ref{P-S}, into a semi-direct 
product $\B=\B'\ltimes\S$ where $\S$ is elementary. One knows from Proposition \ref{P-S} and 
\cite{Bi07}
that the extension homomorphism ${\bf R}:\B'\to\Aut(\S)$ underlies a homomorphism from $\B'$ into the isotropy subgroup $\Aut(\S,s,\omega^\S)_e$ at the unit element $e$
of $\S$ viewed as a symmetric space:
$$
{\bf R}:\B'\to\mbox{\rm Sp}(V,\omega^0)\subset\Aut(\S,s,\omega^\S)_e\;,
$$
where $(V,\omega^0)$ is the symplectic vector space attached to $\S$.
In particular, the action of $\B'$ leaves invariant the two-point kernel $K_{\theta,\tau}$ on $\S\times\S$. Iterating the above observation at the level of $\B'$
and  translating the ``extension Lemma" in \cite{BCSV} within the present framework, we obtain:
\begin{prop}
\label{efF}
Let $\B$ be  a normal $\bf j$-group with Pyatetskii-Shapiro decomposition  $\B=(\S_N\ltimes\dots)\ltimes\S_1$ and
fix  $\vec\tau:=(\tau_1,\dots,\tau_N)\in{\bf \Theta}^N$. Parametrizing a group element $g\in\B$ as $g=g_1\dots g_N$, with $g_i\in\S_i$, we consider the two-point kernel on $\B$ given by
\begin{equation}
\label{2P-ext}
K_{\theta,\vec\tau}(g,g'):=K_{\theta,\tau_1}(g_1,g_1')\dots K_{\theta,\tau_N}(g_N,g_N')\;,
\end{equation}
where  $K_{\theta,\tau_i}$  is the  two-points kernel
on $\S_i\times\S_i$, defined in \eqref{2-point-kernel}.
 Then, the bilinear mapping
$$
\star_{\theta,\vec \tau}:=\Big[(u, v)\mapsto\int_{\B\times\B} K_{\theta,\vec\tau}(g,g')\,R^\star_{g}(u)\,
R^\star_{g'}(v)\,{\rm d}_\B(g)\,{\rm d}_\B(g')\Big]\;,
$$
is associative on 
$$
\CE_{\theta,\vec \tau}(\B):=\CE_{\theta,\tau_N}(\S_N)\otimes\dots\otimes\CE_{\theta,\tau_1}(\S_1)\;,
$$
(recall that $\CE_{\theta,\tau_j}(\S_j)$ is nuclear). Moreover, at the level of compactly supported functions, the product  $\star_{\theta,\vec \tau}$ is equivariant under the left-translations in $\B$.
\label{prod:extension}
\end{prop}
\section{An oscillatory integral formula for the star-product }
\label{OIFDF}
In this section, we fix $\B$ a normal $\bf j$-group, with Lie algebra $\b$. We also let
 $\vec\tau\in{\bf\Theta}^N$ be as above ($N$ is the number of elementary components in $\B$) and form the two-point kernel $K_{\theta,\vec\tau}$ 
on $\B\times\B$, defined in \eqref{2P-ext}. Proposition \ref {prod:extension} implies 
 that
 the deformed product
\begin{equation}
\label{PROD-EXT}
u\star_{\theta,\vec \tau}v=\int_{\B\times\B} K_{\theta,\vec\tau}(g,g')\,R^\star_{g}(u)\,
R^\star_{g'}(v)\,{\rm d}_\B(g)\,{\rm d}_\B(g')\;,
\end{equation}
is weakly associative (in the sense of Definition \ref{wass}) and left $\B$-equivariant. The results of chapter \ref{OSIL}  will allow to properly understand the integral in \eqref{PROD-EXT}
as an oscillatory one. As a consequence, we will see that the deformed product extends as a continuous bilinear
and associative map on the function space $\CB(\B,\CA)$, for $\CA$ a Fr\'echet algebra.  
 We start with a simple fact:

\begin{lem}
Let $\B$ be an elementary normal $\bf j$-group and $\vec\tau\in{\bf\Theta}^N$.
 Then the amplitude $
A_{\theta,\vec\tau}$, as given in Proposition \ref{efF},
consists of an element of $\CB^{\mu_\tau}(\B\times\B)$ for a  tempered weight $\mu_\tau$.
\label{A-tau}
\end{lem}
\begin{proof}
Consider first the case where $\B=\S$ is elementary.
Within the notations of section \ref{2P}, we have
$$
|\vec x_3|=|(x_{1,3},x_{2,3})|=\big|\big(-\sinh2a_2,\sinh2a_1\big)\big|=\big(\sinh^22a_2+\sinh^22a_1\big)^{1/2}\;,
$$
so that  the function 
$$
\mu_{\rm can}(x_1,x_2):= \cosh a_1\,\cosh a_2\;,
$$
is a tempered weight. 
As the left invariant vector field $\widetilde H$ on $\S$ restricted to functions of depending on 
the variable $a$ only, coincides
with the partial differentiation operator $\partial_{a}$, we get from  the explicit expression 
$$
A_{\mathrm{can}}^\S(x_1,x_2)=
(\cosh a_1\,\cosh a_2\,\cosh(a_1-a_2))^d
\sqrt{\cosh 2a_1\,\cosh 2a_2\,\cosh 2(a_1-a_2)} \;,
$$
that there exists $\rho>0$ such that
for any $X\in\CU(\s\oplus \s)$, there exists a constant $C_X>0$ with
$$
\big|\widetilde X\,A_{\rm can}^\S\big|\leq C_X\, \mu_{\rm can}^\rho\;.
$$
Hence
$A_{\rm can}^\S\in\CB^{\mu_{\rm can}^{\rho}}(\S\times\S)$.
 Next, since $\tau\in{\bf\Theta}$, we have $\exp\circ\pm\tau \in{\mathcal O}_{C,1}(\R)$.
Thus, there exists $r>0$ such that the $n$-th derivative 
of $\exp\circ\pm\tau (x)$ is bounded by $(1+|x|)^{r-n}$.
Let us denote by  $\deg(\tau)$ such positive number $r$. Since $\exp\circ\pm\tau$  depends 
on the variable $a$ only, among all elements of $\CU(\s\oplus \s)$,  only the powers of $\widetilde H_ i$, $ i=1,2$, give non zero contributions. Therefore, an easy computation shows that
for any $X\in\CU(\s\oplus \s)$, there exists a constant $C_X>0$ with
$$
\big|\widetilde X\,\exp\{\pm\tau\big(\tfrac2\theta\sinh2a\big)\}\big|\leq C_X\, 
(1+|\vec x_3|)^{2\deg(\tau)}\;.
$$
Hence $A_{\theta,\tau}$ belongs to $\CB^{\mu_\tau}(\S\times\S)$ for
 $\mu_\tau=\mu_{\rm can}^{\rho+3\deg(\tau)}$.

 The general case $\B=\B'\ltimes\S$ follows easily by Pyatetskii-Shapiro theory, since only 
 the variables in $V\subset\S$ are affected by the action of $\B'$ and that $A_{\theta,\vec\tau}$
 is independent of these variables.
\end{proof}

We now consider a Fr\'echet algebra $\CA$, with topology underlying a countable family
of sub-multiplicative semi-norms $\{\|.\|_j\}_{j\in\N}$. Combining  Lemma \ref{A-tau}  with   Theorem \ref{TASP} leads us to prove that  
the integral in the expression of the deformed product  \eqref{PROD} can be properly understood as an oscillatory one in the sense of chapter \ref{OSIL}.  In particular, this allows to define the product $\star_{\theta,\vec\tau}$ on $\CB(\B,\CA)$. This is the main result of this chapter.

Before going further, a clarification regarding the notion of temperedness should be made.
Recall that in the framework of a tempered pair $(G,S)$, the notion of temperedness comes from
the global coordinate system associated to the phase function $S$ (see Definition \ref{temp-grp}
and Definition \ref{TEMPPAIR}). Thus, in the case of a normal $\bf j$-group $\B$, this notion is 
a priori only defined on the product group $\B\times\B$. However, we have seen in the proof 
of Theorem \ref{TASP} that the coordinate system \eqref{TEMPCOORD} associated to the phase 
function  $S_{\rm can}^\B$ is related to the adapted tempered coordinates on  $\B\times\B$
(see Definition \ref{adapted}) by a tempered  diffeomorphism. Hence, it  seems natural to
define directly the notion of temperedness on $\B$ by mean of the adapted tempered coordinates.
Having saying that, the important observation is
 that two  functions $\vf_1,\vf_2$  are tempered  on $\B$,
 if and only if $\vf_1\otimes\vf_2$ is tempered on $\B\times\B$.  It implies a great simplification
 of the assumptions  in Theorem
 \ref{OSCKERN}, Proposition \ref{WASS} and Proposition \ref{action-schwartz}, namely the  temperedness at the level of  tensor product of weights can be reduced to
 temperedness at the level of each factor.

\begin{thm}
\label{achievement}
Let $\B$ be a normal $\bf j$-group. Fix
$\vec\tau\in{\bf\Theta}^N$ and let $\underline\mu_1$, $\underline\mu_2$, $\underline\mu_3$ 
be three families of  tempered  (in the sense of adapted tempered coordinates)
 weights on $\B$ of sub-multiplicativity degree $(\underline L_1,
\underline R_2)$,  $(\underline L_{2},\underline R_{2})$, 
 $( \underline L_{3}, \underline R_{3})$. Considering $K_{\theta,\vec\tau}$ the two-point kernel on $\B$ defined in \eqref{2P-ext}, the correspondence
\begin{align*}
\star_{\theta,\vec\tau}: \CB^{\underline\mu_1}(\B,\CA)\times\CB^{\underline\mu_2}(\B,\CA)
&\to\CB^{\underline\nu}(\B,\CA) \quad\mbox{with}\quad \nu_j=\mu_{1,j}^{L_{1,j}}\mu_{2,j}^{L_{2,j}}
\;,\\
(F_1,F_2)&\mapsto \widetilde{\int_{\B\times\B} K_{\theta,\vec\tau}}\,\Big[(x_1,x_2)\mapsto 
R^\star_{x_1}(F_1) R^\star_{x_2}(F_2)\Big]\;,
\end{align*}
is a  continuous bilinear map and is equivariant under the left translations in $\B$ in the 
sense that for all $g\in\B$, we have
$$
L^\star_g(F_1\star_{\theta,\vec\tau}F_2)=(L^\star_gF_1)\star_{\theta,\vec\tau}(L^\star_gF_2)
\quad\mbox{in}\quad  \CB^{\underline\lambda}(\B,\CA)\quad\mbox{ for}\quad
 \underline\lambda=\{\mu_{1,j}^{L_{1,j}R_{1,j}}\mu_{2,j}^{L_{2,j}R_{2,j}}\}\;.
$$
 Moreover, the map $\star_{\theta,\vec\tau}$ is associative in 
the sense that then for every elements $F_j$ in  $\CB^{\underline\mu_j}(\B,\CA)$, $j=1,2,3$,
we have the equality
$$
\big(F_1\star_{\theta,\vec\tau} F_2\big)\star_{\theta,\vec\tau} F_3=F_1\star_{\theta,\vec\tau} 
\big(F_2\star_{\theta,\vec\tau} F_3\big)
\quad\mbox{in}\quad \CB^{\underline\rho}(\B,\CA)\;,
$$
with
$$
\underline\rho=\{\mu_{1,j}^{L_{1,j}^2}{\mu_{2,j}}^{L_{2,j}^2}{\mu_{3,j}}^{L_{3,j}^2}\}_{j\in\N}\;.
$$
In particular, $\big(\CB(\B,\CA),\star_{\theta,\vec\tau}\big)$ is a Fr\'echet algebra.
\end{thm}
\begin{proof}
That the bilinear map $\star_{\theta,\vec \tau}$ (with the domain and image as indicated) is well 
defined and  continuous, follows from Theorem \ref{OSCKERN}
(cf$.$ the above discussion for the condition of temperedness of the weights involved), 
  Theorem \ref{TASP}  
and Lemma \ref{A-tau}.  
Associativity follows from 
associativity in $\CE_{\theta,\vec\tau}(\B)$, which implies weak associativity in the sense
of  Definition \ref{wass} and Proposition \ref{WASS}. So, it remains to prove left $\B$-equivariance.
We first note that by Lemma \ref{SmoothFamily} (ii), the group $\B$ acts on the left continuously 
from  $\CB^{\underline\mu}(\B,\CA)$ to $\CB^{\underline\gamma}(\B,\CA)$, with 
$\underline\gamma=
\{\mu_j^{R_j}\}$ (for any family of weights
$\underline\mu$ of sub-multiplicative degree $(\underline L,\underline R)$).
Also, we have by Lemma \ref{specific-limits} that $F_1\star_{\theta,\vec\tau} F_2=\lim_{n_1,n_2}
F_{1,n_1}\star_{\theta,\vec\tau} F_{2,n_2}$ in $\CB^{\underline\nu}(\B,\CA)$,
with $\underline\nu=\{\mu_{1,j}^{L_{1,j}}\mu_{2,j}^{L_{2,j}}\}$, for any
pair of sequences $\{F_{1,n}\}$ and $\{F_{2,n}\}$ of smooth compactly supported $\CA$-valued 
functions
on $\B$, which converge to $F_1$ and $F_2$, in the topology of 
$ \CB^{\underline{\hat \mu}_1}(\B,\CA)$ and $\CB^{\underline{\hat\mu}_2}(\B,\CA)$ 
for any sequence of weights 
$\underline{\hat\mu}_1$ and $\underline{\hat \mu}_2$ dominating 
$\underline\mu_1$ and $\underline\mu_2$. From continuity of
the left regular action (see Lemma \ref{SmoothFamily} ii) and  left 
$\B$-equivariance at the level of $\CD(\B,\CA)$, we thus  have
$$
L^\star_g\big(F\star_{\theta,\vec\tau} F'\big)=\lim_{n,n'\to\infty}
L^\star_g\big(F_n\star_{\theta,\vec\tau} F'_{n'}\big)=\lim_{n,n'\to\infty}
(L^\star_gF_n)\star_{\theta,\vec\tau} (L^\star_gF'_{n'})\;,
$$ 
where the limits are 
in $\CB^{\underline\lambda}(\B,\CA)$, for 
$ \underline\lambda=\{\mu_{1,j}^{L_{1,j}R_{1,j}}\mu_{2,j}^{L_{2,j}R_{2,j}}\}$.  It remains to find specific approximation 
sequences $\{F_{1,n}\}$ and $\{F_{2,n}\}$, such that $\{L^\star_g F_{1,n}\}$ and 
$\{L^\star_g F_{2,n}\}$
converge to $L^\star_gF_1$ and $L^\star_gF_2$, in the topology of 
$ \CB^{\underline{\hat\gamma}_1}(\B,\CA)$ and $\CB^{\underline{\hat\gamma}_2}(\B,\CA)$
with $\underline{\hat\gamma}_1=\{\hat \mu_{1,j}^{R_{1,j}}\}$ and 
$\underline{\hat\gamma}_2=\{\hat \mu_{2,j}^{R_{2,j}}\}$. 
For this, we observe that the same  construction as in the proof of Lemma \ref{symbols} (viii),  
does the job. Indeed, recall that there, we have constructed the approximation sequence
$\{F_n\}$, by setting for $F\in\CB(\B,\CA)$:
$$
F_n:=e_n\,F\in\CD(\B,\CA)\quad
\mbox{where}\quad e_n:=\int_\B \psi(g)\,R^\star_g(\chi_{C_n})\,{\rm d}_\B (g)
\in\CD(\B)\;,
$$
and $0\leq\psi\in\CD(\B)$, $\int_\B \psi(x)\,{\rm d}_\B(x)=1$, $\{C_n\}$ is an increasing 
sequence of relatively compact open subsets of $\B$ converging to $\B$ and $\chi_{C_n}$ is the characteristic function of $C_n$. Fixing $g\in\B$ and setting
 $C_n^g:=g.C_n$, the sequence  $\{C_n^g\}$ is still an increasing 
sequence of relatively compact open subsets on $\B$ converging to $\B$. Also, as
$$
e_n^g:=L^\star_g(e_n)=\int_\B \psi(g')\,R^\star_{g'}(\chi_{C_n^g})\,{\rm d}_\B (g')
\in\CD(\B)\;,
$$
we deduce  that for all  $j,k\in\N$:
\begin{align*}
\|L^*_g(F_n)-L^*_g(F)\|_{j,k,\underline\gamma}
=\|(1-e_n^g)L^*_g(F)\|_{j,k,\underline\gamma}\quad\mbox{with}\quad \underline{\hat\gamma}
=\{{\hat\mu_j}^{R_j}\}_{j\in\N}\;,
\end{align*}
which, by Lemma \ref{SmoothFamily} (vi), converges to zero as 
$L^*_g(F)\in   \CB^{\underline\gamma}(\B,\CA)$ with $\underline\gamma=\{\mu_j^{R_j}\}$ and 
$\underline{\gamma}\prec\underline{\hat\gamma}$.
\end{proof}

Let $\CS^{S_{\rm can}^\B}(\B,\CA)$ be the one-variable Schwartz space 
associated to the admissible and tame  tempered pair $(\B\times\B,S^\B_{\rm can})$, 
constructed in Definition \ref{S1P}.
The next result  follows immediately from Proposition 
\ref{action-schwartz}.

\begin{prop}
\label{universal-schwartz}
Let $\B$ be a normal $\bf j$-group and 
 $\vec\tau\in{\bf\Theta}^N$. Then, endowed with the multiplication $\star_{\theta,\vec\tau}$,
the space $\CS^{S_{\rm can}^\B}(\B,\CA)$ becomes a Fr\'echet algebra which,
for $\underline\mu$  an arbitrary family of  tempered weights,
acts continuously on   
$\CB^{\underline\mu}(\B,\CA)$,  via
$$
L_{\star_{\theta,\vec\tau}}(F):\vf\mapsto F\star_{\theta,\vec\tau} \vf\,,\qquad F\in
\CB^{\underline\mu}(\B,\CA)\,,\quad
\vf\in\CS^{S_{\rm can}^\B}(\B,\CA)\;.
$$
In particular, $\big(\CS^{S_{\rm can}^\B}(\B,\CA),\star_{\theta,\vec\tau}\big)$ is an ideal  of 
$\big(\CB(\B,\CA),\star_{\theta,\vec\tau}\big)$.
\end{prop}

We now see  that, as expected,  the constant function is an identity for the deformed product.
\begin{prop}
\label{unit}
Let $\B$ be a normal $\bf j$-group. Fix
$\vec\tau\in{\bf\Theta}^N$, $\underline\mu$ a family of tempered weights of sub-multiplicative
degree $(\underline L_,\underline R)$ 
and $F\in\CB^{\underline\mu}(\B,\CA)$. Identifying an element $a\in\CA$ with the function 
$[g\mapsto a]$ in $\CB(\B,\CA)$, we have
$$
a\star_{\theta,\vec\tau}F=aF\,,\,\quad F\star_{\theta,\vec\tau}a=Fa\;,
$$
in $\CB^{\underline\mu}(\B,\CA)$. In particular, if $\CA$ is unital, the element $[g\mapsto 1_\CA]\in
\CB(\B,\CA)$ is the unit
of  $\big(\CB(\B,\CA),\star_{\theta,\vec\tau}\big)$.
\end{prop}
\begin{proof}
Since the constant unit function is a fixed point of the map $T_{\theta,\vec\tau}^{-1}$, for every
$\vf\in\CS^{S_{\rm can}^\B}(\B,\CA)$,
we have:
$$
\vf\star_{\theta,\vec\tau} a=T_{\theta,\vec\tau}\big(T_{\theta,\vec\tau}^{-1}(\vf)\star_\theta^0 a\big)\;,
$$
in $\CS^{S_{\rm can}^\B}(\B,\CA)$.
By Remark \ref{lastS}, we see that the
 transported Schwartz space $\CS^{S_{\rm can}^\B}(\B,\CA)$ is a (dense) subset of the ordinary
Schwartz space $\CS(\b,\CA)$, under the usual identification $\B\simeq\b$. Since 
$T_{\theta,\vec\tau}^{-1}$ preserves the latter space, we see that
$T_{\theta,\vec\tau}^{-1}(\vf)\in\CS(\b,\CA)$. It is well known (see \cite{Ri} for the Fr\'echet algebra 
valued case) that the Moyal product admits the constant function 
as unit element. Thus
$ \vf\star_{\theta,\vec\tau}a=\vf a$ and $ a\star_{\theta,\vec\tau}\vf=a\vf$ for all 
$\vf\in \CS^{S_{\rm can}^\B}(\B,\CA)$ and $a\in\CA$.
Now, consider the injective homomorphism $L_{\star_{\theta,\vec\tau}}$ from 
$\big(\CB^{\underline\mu}(\B,\CA),\star_{\theta,\vec\tau}\big)$ to the algebra of continuous operators acting 
on 
$\CS^{S_{\rm can}^\B}(\B,\CA)$, defined in Proposition \ref{universal-schwartz}. From the previous 
considerations, the associativity of the deformed product and the fact that $\CS^{S_{\rm can}^\B}(\B,
\CA)$ is an ideal of $\CB^{\underline\mu}(\B,\CA)$,
we get
$$
L_{\star_{\theta,\vec\tau}}( F\star_{\theta,\vec\tau}a)=L_{\star_{\theta,\vec\tau}}(Fa)\,,\quad\forall F\in\CB^{\underline\mu}(\B,\CA)\;,
$$
which entails by injectivity that $F\star_{\theta,\vec\tau}a=Fa$ in $\CB^{\underline\nu}(\B,\CA)$,
with $\underline\nu=\{\mu_j^{L_j}\}$.
As $Fa\in \CB^{\underline\mu}(\B,\CA)$, we deduce that the equality $F\star_{\theta,\vec\tau}a=Fa$ 
holds in fact in $\CB^{\underline\mu}(\B,\CA)$.
 The case of $a\star_{\theta,\vec\tau} F$ is entirely similar.
\end{proof}

\chapter{Deformation of Fr\'echet algebras}
\label{DFA}

In this chapter, we  consider  a normal $\bf j$-group $\B$ and a pair $(\CA,\alpha)$,
consisting of a Fr\'echet algebra $\CA$, together with 
a strongly continuous (not necessarily isometric) action $\alpha$ of $\B$ by automorphisms.
For $a\in\CA$, we let $\alpha(a)$ be the $\CA$-valued function on $\B$, defined by 
\begin{equation}
\label{MP}
\alpha(a):=[g\in\B\mapsto \alpha_g(a)\in\CA]
\;.
\end{equation}
Our main goal, achieved in section \ref{DP}, is to show that the formula 
\begin{equation*}
a \star_{\theta,\vec\tau}^\alpha b:=\big(\alpha(a)\star_{\theta,\vec\tau}\alpha(b)\big)(e)\;,
\end{equation*}
equips $\CA^\infty$ with a new noncommutative and associative Fr\'echet algebra structure.
This is proved in Theorem \ref{UDF} as a direct consequence of  Theorem \ref{achievement} and 
of Lemma \ref{embed} which shows that the map \eqref{MP} is continuous from $\CA^\infty$
to $\CB^{\underline\mu}(\B,\CA^\infty)$, for $\underline\mu$ a nontrivial family of   tempered weights
on the group $\B$.

\section{The deformed product}\label{DP}
We do not yet  need  to work with isometric actions.
 We start by  general considerations regarding tempered actions:

\begin{dfn}
\label{temp-action}
A {\bf tempered action}   of a tempered Lie group $G$ on a Fr\'echet algebra $\CA$, 
is given by the data $(\alpha,\underline\mu^\alpha)$ where $\alpha$ is an action of $G$
on $\CA$ and $\underline\mu^\alpha$ is a family of tempered  weights on $G$  such that
for all $j\in\N$, all $a\in\CA$ and all $g\in G$, we have 
$$
\|\alpha_g(a)\|_j\leq \mu^\alpha_j(g)\|a\|_j \;.
$$
\end{dfn}
\begin{rmk}
\label{belcont}
Note that for a tempered action and for $g\in G$ fixed,  $\alpha_g$ acts continuously on
$\CA$. 
\end{rmk}

We denote by $\CA^\infty$ the set of smooth vectors for the action $\alpha$ of $\B$ on $\CA$:
$$
\CA^\infty:=\big\{a\in\CA\;:\; \alpha(a)\in C^\infty(\B,\CA)\big\}\;.
$$
When the action is strongly continuous, $\CA^\infty$ is a dense subspace  of $\CA$.
On this subset, we consider the infinitesimal 
form of the action, given for $X\in\b$ by:
\begin{equation*}
X^\alpha(a):=\frac d{dt}\Big|_{t=0}\,\alpha_{e^{tX}}(a)\,,\qquad a\in \CA^\infty\;,
\end{equation*}
and extended  to the whole universal enveloping algebra $\CU(\b)$, by declaring that the map
$\CU(\b)\to \End(\CA^\infty)$, $X\mapsto X^\alpha$ is an algebra homomorphism.
The subspace $\CA^\infty$ carries a finer   topology  associated with the following set of semi-norms:
$$
\|a\|_{j,X}:=\|X^\alpha( a)\|_j,\qquad a\in\CA^\infty\,,\quad X\in\CU(\b),\,j\in\N\;.
$$
Considering the PBW basis of $\CU(\b)$ associated to an ordered basis of $\b$  as  in \eqref{PBW}, 
one can use only countably many semi-norms to define the 
topology  of $\CA^\infty$. The latter are 
indexed by $(j,k)\in\N^2$,
where $j$ refers to the labeling of the initial family of semi-norms 
$\{\|.\|_j\}_{j\in\N}$ of $\CA$ and $k$ refers
to the labeling of the  filtration $\CU(\b)=\cup_{k\in\N}\CU_k(\b)$ associated to the chosen PBW basis, as 
 defined in \eqref{filter}. In turn, $\CA^\infty$ becomes a Fr\'echet space, for the topology associated with the 
 semi-norms
\begin{equation}
\label{newSN}
\|.\|_{j,k}:\CA^\infty\to[0,\infty)\,,\qquad a\mapsto\sup_{X\in\,\CU_k(\g)}\frac{\|a\|_{j,X}}{|X|_k}=\sup_{X\in\,\CU_k(\g)}\frac{\|X^\alpha(a)\|_j}{|X|_k}\;,
\end{equation}
with $j,k\in\N$ and where $|.|_k$ is the  $\ell^1$-norm of $\CU_k(\b)$ defined in \eqref{norm-CUk}.
As in \eqref{eqeq}, we have
 \begin{align*}
  \|a\|_{j,k}\,\leq \max_{|\beta|\leq k}\, \|a\|_{j, X^\beta}\;,
 \end{align*}
 with $\{X^\beta,\,|\beta|\leq k\}$ the basis $\eqref{PBW}$ of $\CU_k(\b)$.
  Hence the semi-norms  \eqref{newSN} are well defined on $\CA^\infty$.

In the context of a tempered action  on a Fr\'echet algebra $\CA$, we observe
that the restriction of the action to $\CA^\infty$ is also tempered,
 but  never isometric, even if the  action is isometric on $\CA$  unless the group  is Abelian.
  This explains why in our context it is natural to work with tempered actions, 
 rather than with isometric ones.
\begin{lem}
\label{cinqcinq}
Let $(\CA,\alpha,\underline\mu^\alpha)$ be a Fr\'echet algebra endowed with a tempered action
of a tempered Lie group $G$. Then, the restriction of $\alpha$ on $\CA^\infty$ is tempered too,
with:
$$
\|\alpha_g(a)\|_{j,k}\leq C(k)\,
\fd_G(g)^k\,\mu_j^\alpha(g)\,\|a\|_{j,k}\,,\qquad j,k\in\N,\,g\in G,\,a\in\CA^\infty\;.
$$
\end{lem}
\begin{proof}
First remark 
$$
\|\alpha_g(a)\|_{j,k}=\!\!\sup_{X\in\,\CU_k(\g)}\!\!\frac{\|\alpha_g
\big(\big(\Ad_{g^{-1}}(X)\big)^\alpha(a)\big)\|_j}
{|X|_k}\leq \mu_j^\alpha(g) \sup_{X\in\,\CU_k(\g)}\!\!\frac{\|\big(\Ad_{g^{-1}}(X)\big)^\alpha(a)\|_j}
{|X|_k}\;.
$$
As for $X\in\CU_k(\g)$ and $a\in\CA^\infty$, we have
$$
\|X^\alpha(a)\|_j\leq|X|_{k}\sup_{Y\in\,\CU_k(\g)}\frac{\|Y^\alpha(a)\|_j}{|Y|_k}=|X|_{k}\,\|a\|_{j,k}\;,
$$
we get, with $|\Ad_g|_k$ denoting the operator norm of the adjoint action of $G$ on the normed space 
$\big(\CU_k(\g),|.|_k\big)$:
$$
\|\alpha_g(a)\|_{j,k}\leq  \mu_j^\alpha(g) \sup_{X\in\,\CU_k(\g)}\frac{|\Ad_{g^{-1}}(X)|_k}
{|X|_k}\|a\|_{j,k}=  \mu_j^\alpha(g)\,|\Ad_{g^{-1}}|_k\,\|a\|_{j,k}\;,
$$
and one concludes using Lemma \ref{lemtwo}.
\end{proof}

\begin{ex}
{\rm
Applying the former result to $\alpha=R^\star$ and $\CA=\CS^{S_{\rm can}^\B}(\B)$ (which is its own space
of smooth vectors), we see that 
the right-action of $\B$ on   $\CS^{S_{\rm can}^\B}(\B)$ is tempered.}
\end{ex}

The following statement  is the foundation of our construction:
\begin{lem}
\label{embed}
Let $(\alpha,\underline\mu^\alpha)$ be a tempered and strongly continuous action of a Lie group 
$G$ on a 
Fr\'echet algebra $\CA$. Setting then $\underline\nu:=\{\mu_j^\alpha\fd_G^k\}_{j,k\in\N}$, we have 
an equivariant continuous embedding
\begin{align*}
\alpha:\CA^\infty\to\CB^{\underline\nu}(G,\CA^\infty)\;,\qquad
a\mapsto\alpha(a)=[g\in G\mapsto \alpha_g(a)
\in\CA^\infty]\;.
\end{align*}
\end{lem}
\begin{proof}
Note first that for $a\in\CA$ and $g,g_0\in G$, we have
$$
\alpha\big(\alpha_g(a)\big)(g_0)=\alpha_{g_0g}(a)=\big(R^\star_g \alpha(a)\big)(g_0)\;,
$$
and thus $\alpha:a\in\CA\mapsto[g\mapsto\alpha_g(a)]\in C(G,\CA)$ intertwines the actions 
$R^\star$ and $\alpha$. 
Let now $a\in\CA^\infty$ and $X\in\CU(\g)$. By equivariance and strong-differentiability of 
$\alpha$ on $\CA^\infty$, we get 
$$
\widetilde X\alpha(a)=\alpha(X^\alpha a)\;.
$$ 
Since for all $j\in\N$ and  all $a\in\CA$, we have $\|\alpha_g(a)\|_j\leq\mu_j^\alpha(g)\|a\|_j$,   
we deduce that
\begin{align*}
\|\alpha(a)\|_{j,k,\underline\mu^\alpha}=\sup_{X\in\,\CU_k(\g)}\sup_{g\in G}
\frac{\|\widetilde X \alpha_g(a)\|_j}
{\mu^\alpha_j(g)|X|_k}&=\sup_{X\in\,\CU_k(\g)}\sup_{g\in G}\frac{\|\alpha_g(X^\alpha a)\|_j}
{\mu^\alpha_j(g)|X|_k}\\&
\leq\sup_{X\in\,\CU_k(\g)}\frac{\|X^\alpha a\|_j}{|X|_k}=\|a\|_{j,k}\;.
\end{align*}
This analysis shows that the map $\alpha:\CA^\infty\to \CB^{\underline\mu^\alpha}(G,\CA)$ is 
continuous. Now 
we want to take into account the intrinsic topology of $\CA^\infty$ in the target space of the map 
$\alpha$. 
Remark that the topology of $\CB^{\underline\nu}(G,\CA^\infty)$ 
is associated with the countable set of semi-norms
$$
\|F\|_{(j,k),k',\underline\nu}=\sup_{X\in\,\CU_{k'}(\g)}\,\sup_{g\in G}\,\sup_{Y\in\,\CU_{k}(\g)}\,
\frac{\|Y^\alpha\big(\widetilde X F(g)\big)\|_j}{\nu_{j,k}(g)\,|X|_{k'}\,|Y|_k}
\;.
$$
Since $\alpha_{g^{-1}}\circ X^\alpha \circ \alpha_g=(\Ad_{g^{-1}}X)^\alpha$ for all $X\in\CU(\g)$ 
and $g\in G$, we get for $F=\alpha(a)$ and $\underline\nu=\{\mu_j^\alpha\fd_G^k\}_{j,k\in\N}$:
 \begin{align*}
\|\alpha(a)\|_{(j,k),k',\underline\nu}
&=\sup_{X\in\,\CU_{k'}(\g)}\,\sup_{g\in G}\,\sup_{Y\in\,\CU_{k}(\g)}\,
\frac{\|Y^\alpha\big(\widetilde X_g\, \alpha_g(a)\big)\|_j}{\mu^\alpha_j(g)\,\fd_G(g)^k\,|X|_{k'}\,|Y|_k}\\
&=\sup_{X\in\,\CU_{k'}(\g)}\,\sup_{g\in G}\,\sup_{Y\in\,\CU_{k}(\g)}\,
\frac{\|Y^\alpha\big( \alpha_g(X^{\alpha}a)\big)\|_j}{\mu^\alpha_j(g)\,\fd_G(g)^k\,|X|_{k'}\,|Y|_k}\\
&=\sup_{X\in\,\CU_{k'}(\g)}\,\sup_{g\in G}\,\sup_{Y\in\,\CU_{k}(\g)}\,
\frac{\| \alpha_g\big((\Ad_{g^{-1}}Y)^\alpha X^{\alpha}a\big)\|_j}{\mu^\alpha_j(g)\,\fd_G(g)^k\,|X|_{k'}\,|Y|_k}\\
&\leq\sup_{X\in\,\CU_{k'}(\g)}\,\sup_{g\in G}\,\sup_{Y\in\,\CU_{k}(\g)}\,
\frac{\| (\Ad_{g^{-1}}Y)^\alpha X^{\alpha}a\|_j}{\fd_G(g)^k\,|X|_{k'}\,|Y|_k}\\
&\leq\Big(\sup_{g\in G} \frac{|\Ad_{g^{-1}}|_k}{\fd_G(g)^k}\Big)\sup_{X\in\,\CU_{k'}(\g)}\,\sup_{Y\in\,\CU_{k}(\g)}\,
\frac{\| Y^\alpha X^{\alpha}a\|_j}{|X|_{k'}\,|Y|_k}\\
&\leq\Big(\sup_{g\in G} \frac{|\Ad_{g^{-1}}|_k}{\fd_G(g)^k}\Big)\!\Big(\sup_{X\in\,\CU_{k'}(\g)}\,\sup_{Y\in\,\CU_{k}(\g)}
\frac{|YX|_{k+k'}}{|X|_{k'}\,|Y|_k}\Big)\sup_{Z\in\,\CU_{k+k'}(\g)}\!\!
\frac{\| Z^{\alpha}a\|_j}{|Z|_{k+k'}}\\
&\qquad\quad =\Big(\sup_{g\in G} \frac{|\Ad_{g^{-1}}|_k}{\fd_G(g)^k}\Big)\Big(\sup_{X\in\,\CU_{k'}(\g)}\,\sup_{Y\in\,\CU_{k}(\g)}
\frac{|YX|_{k+k'}}{|X|_{k'}\,|Y|_k}\Big)\,\| a\|_{j,{k+k'}}\;,
\end{align*}
and one concludes using Lemma \ref{lemtwo}.
\end{proof}

The next result, although rather obvious, will also play a key role.
\begin{lem}
\label{evaluation}
Let $\CA$ be a Fr\'echet algebra 
and let $\underline\mu$ be a family of tempered weights on a tempered Lie group $G$.
Then, the evaluation map at the unit element, $\CB^{\underline\mu}(G,\CA)\to\CA$, 
$F\mapsto F(e)$, is 
continuous.
\end{lem}
\begin{proof}
Fix $j\in\N$. We have
for any $F\in \CB^{\underline\mu}(G,\CA)$:
$$
\|F(e)\|_j\leq\mu_j(e)\sup_{g\in G}\frac{\|F(g)\|_j}{\mu_j(g)}=
\mu_j(e)\|F\|_{j,0,\underline\mu}\;,
$$
and the result follows immediately.
\end{proof}

 Last, we need to lift the action $\alpha$ from $\CA^\infty$ to 
$\CB^{\underline\nu}(\B,\CA^\infty)$, $\underline\nu=\{\mu_j^\alpha\fd_G^k\}$, and to show that this lift  acts by  automorphisms
of the product $\star_{\theta,\vec\tau}$.
\begin{lem}
\label{hat-alpha}
Let 
$(\alpha,\underline\mu^\alpha)$ be a strongly 
continuous and tempered action of a normal $\bf j$-group $\B$ on a Fr\'echet 
algebra $\CA$ 
and $\underline\mu_1,\underline\mu_2$ be two families of tempered weights 
with sub-multiplicative degree $(\underline L_1,\underline R_1), (\underline L_2,\underline R_2) $.
For $g\in\B$, the map 
$$
\hat\alpha_g:F\mapsto \big[g_0\in\B\mapsto \alpha_g\big(F(g_0)\big)\big]\;,
$$
is continuous on $\CB^{\underline\mu_i}(\B,\CA^\infty)$, $i=1,2$. Moreover, given 
$(\theta,\vec\tau)\in \R^*\times{\bf\Theta}^N$, $\hat\alpha$ defines an action of $\B$ by
automorphisms of the deformed product $\star_{\theta,\vec\tau}$, 
in the sense that for all $F_1\in \CB^{\underline\mu_1}(\B,\CA^\infty)$ 
and $F_1\in \CB^{\underline\mu_2}(\B,\CA^\infty)$, we have for all $g\in\B$:
$$
\hat\alpha_g\big(F\star_{\theta,\vec\tau} F'\big)=\hat\alpha_g(F)\star_{\theta,\vec\tau}
\hat\alpha_g(F')\qquad \mbox{in}\qquad \CB^{\underline\nu}(\B,\CA^\infty)\;,
$$
with
$$
 \underline\nu=\{\mu_{1,j,k}^{L_{1,j,k}}{\mu_{2,j,k}}^{L_{2,j,k}}\}_{j,k\in\N}\;.
$$
\end{lem}
\begin{proof}
For $F\in \CB^{\underline\mu}(\B,\CA^\infty)$, $X,Y\in\CU(\b)$ and $g,g'\in\B$, we have
$$
Y^\alpha\big(\widetilde X\,\hat\alpha_g(F)(g')\big)=\alpha_g\Big((\Ad_{g^{-1}} Y)^\alpha\,\big(
\widetilde X F(g')\big)\Big)\;.
$$
This entails that
\begin{align*}
\|\hat\alpha_g(F)\|_{(j,k),k',\underline\mu}&=\sup_{X\in\,\CU_{k'}(\b)}\,\sup_{g'\in\B}\,\sup_{Y\in\,\CU_{k}(\b)}\,
\frac{\|Y^\alpha\big(\widetilde X\,\hat\alpha_g(F)(g')\big)\|_j}{\mu_{j,k}(g')\,|X|_{k'}\,|Y|_k}\\
&=\sup_{X\in\,\CU_{k'}(\b)}\,\sup_{g'\in\B}\,\sup_{Y\in\,\CU_{k}(\b)}\,
\frac{\|\alpha_g\Big((\Ad_{g^{-1}} Y)^\alpha\,\big(
\widetilde X F(g')\big)\Big)\|_j}{\mu_{j,k}(g')\,|X|_{k'}\,|Y|_k}\\
&\leq C(k)\,\mu_j^\alpha(g)\,\fd_\B(g)^k\sup_{X\in\,\CU_{k'}(\b)}\,\sup_{g'\in\B}\,\sup_{Y\in\,\CU_{k}(\b)}\,
\frac{\|\ Y^\alpha\,\big(\widetilde X F(g')\big)\|_j}{\mu_{j,k}(g')\,|X|_{k'}\,|Y|_k}\\
&\qquad
= C(k)\,\mu_j^\alpha(g)\,\fd_\B(g)^k \|F\|_{(j,k),k',\underline\mu}\;,
\end{align*}
proving the continuity. 

Next, consider  $F_1\in \CB^{\underline\mu_1}(\B,\CA^\infty)$ 
and $F_2\in \CB^{\underline\mu_2}(\B,\CA^\infty)$, together with $\underline{\hat\mu}_1$ and 
$\underline{\hat\mu}_2$, two 
families of tempered weights that dominate respectively $\underline\mu_1$ and $\underline\mu_2$.
Defining $F_{1,n}:=F_1e_n\in\CD(\B,\CA)$ and $F_{2,n}=F_2e_n\in\CD(\B,\CA)$, 
with $e_n\in\CD(\B)$ defined in \eqref{en}, from
$$
\hat\alpha_g(F_{1,n})=\hat\alpha_g(F_1)e_n\,,\quad \hat\alpha_g(F_{2,n})=\hat\alpha_g(F_2)e_n\;,
$$
we deduce from Lemma \ref{SmoothFamily} (viii)  that $\{\hat\alpha_g(F_{1,n})\}$ and 
$\{\hat\alpha_g(F_{2,n})\}$ converges to $\{\hat\alpha_g(F_1)\}$ and 
$\{\hat\alpha_g(F_2)\}$ in the topologies of $\CB^{\underline{\hat\mu}_1}(\B,\CA^\infty)$ 
and $ \CB^{\underline{\hat\mu}_2}(\B,\CA^\infty)$ respectively. 
Thus, we can use  
Lemma \ref{specific-limits} to get the $\hat\alpha$-equivariance at the level of smooth
compactly supported functions from the commutativity of $\hat\alpha$ and $R^\star$:
\begin{align*}
\hat\alpha_g\big(F\star_{\theta,\vec\tau} F'\big)&=\hat\alpha_g\big(\lim_{n,n'\to\infty}F_n\star_{\theta,
\vec\tau}F'_{n'}\big)=\lim_{n,n'\to\infty}\hat\alpha_g\big(F_n\star_{\theta,
\vec\tau}F'_{n'}\big)\\
&=\lim_{n,n'\to\infty}\hat\alpha_g(F_n)\star_{\theta,
\vec\tau}\hat\alpha_g(F'_{n'})=\hat\alpha_g(F)\star_{\theta,
\vec\tau}\hat\alpha_g(F')\;,
\end{align*}
in  $\CB^{\underline\nu}(\B,\CA^\infty)$, with 
$\underline\nu=\{\mu_{1,j,k}^{L_{1,j,k}}{\mu_{2,j,k}}^{L_{2,j,k}}\}$.
\end{proof}

We are now prepared to state   the main result of the first part of this memoir:

\begin{thm}[Universal Deformation Formula of Fr\'echet Algebras]
\label{UDF}
Let  $(\CA,\alpha)$ be a Fr\'echet algebra endowed with 
a tempered and strongly continuous action  of a normal $\bf j$-group $\B$.  Let also
 $\theta\in\R^*$ and $\vec\tau\in{\bf\Theta}^{N}$. 
Then,  $(\CA^\infty,\star_{\theta,\vec\tau}^\alpha)$ is an associative Fr\'echet algebra with 
continuous product.
\end{thm}
\begin{proof}
Let  $\underline\mu^\alpha$ be the 
 family of tempered weights, with sub-multiplicative degree $(\underline L,\underline R)$, associated with the 
 tempered action $\alpha$ as in Definition \ref{temp-action}.
Let $a,b\in\CA^\infty$, then by  Lemma \ref{embed}, 
$\alpha(a),\alpha(b)\in \CB^{\underline\mu}(\B,\CA^\infty)$, where
$\underline\mu=\{\mu_j^\alpha\fd_\B^k\}$. Then, since $\fd_\B$ is sub-multiplicative 
of degree $(1,1)$, Theorem \ref{achievement}
 shows that $\alpha(a)\star_{\theta,\vec\tau}\alpha(b)$ belongs to 
 $\CB^{\underline\nu}(\B,\CA^\infty)$, for $\underline\nu=\{\mu_j^{\alpha 2L_j}\fd_{\B}^{2k}\}$,
and that the map
$$
\CA^\infty\times \CA^\infty\to \CB^{\underline\nu}(\B,\CA^\infty)\,,\quad (a,b)\mapsto
\alpha(a)\star_{\theta,\vec\tau}\alpha(b)\;,
$$
is continuous.
Applying  Lemma \ref{evaluation} for the Fr\'echet algebra $\CA^\infty$ then yields that the
composition of maps
\begin{align*}
&\CA^\infty\times \CA^\infty\to \CB^{\underline\nu}(\B,\CA^\infty)\to \CA^\infty\;,\\
 &(a,b)\mapsto
\alpha(a)\star_{\theta,\vec\tau}\alpha(b)\mapsto \big(\alpha(a)\star_{\theta,\vec\tau}\alpha(b)\big)
(e)=:a\star_{\theta,\vec\tau}^\alpha b\;,
\end{align*}
is continuous. 

It remains to prove associativity. 
With $\hat\alpha$ defined in Lemma \ref{hat-alpha}, we compute for $a,b\in\CA^\infty$ and $g\in\B$:
$$
\alpha\big(a\star_{\theta,\vec\tau}^\alpha b\big)(g)=\alpha_g\big(a\star_{\theta,\vec\tau}^\alpha b\big)=\alpha_g\big(\alpha(a)\star_{\theta,\vec\tau}\alpha (b)(e)\big)
=\hat\alpha_g\big(\alpha(a)\star_{\theta,\vec\tau}\alpha (b)\big)(e)\;.
$$
Using Lemma \ref{hat-alpha}, we deduce the equality in 
$\CB^{\underline\nu}(\B,\CA^\infty)$ (for the value of $\underline\nu$ as indicated above):
$$
\hat\alpha_g\big(\alpha(a)\star_{\theta,\vec\tau}\alpha (b)\big)=\hat\alpha_g(\alpha(a))\star_{\theta,\vec\tau}\hat\alpha_g(\alpha (b))\;.
$$
As a short computation  shows,  for $a\in\CA$ and $g\in\B$, we have 
$\hat\alpha_g\big(\alpha(a)\big)=L^\star_{g^{-1}}\big(\alpha(a)\big)$. Thus, using the equivariance of the product $\star_{\theta,\vec\tau}$ under the left regular action, as stated in Theorem \ref{achievement}, we get the equalities 
$$
\hat\alpha_g(\alpha(a))\star_{\theta,\vec\tau}\hat\alpha_g(\alpha (b))=L^\star_{g^{-1}}(\alpha(a))\star_{\theta,\vec\tau}L^\star_{g^{-1}}(\alpha (b))
=L^*_{g^{-1}}\big(\alpha(a)\star_{\theta,\vec\tau}\alpha (b)\big)\;,
$$
 in  $\CB^{\underline\lambda}(\B,\CA^\infty)$, for 
 $\underline\lambda=\{\mu_j^{\alpha2L_jR_j}\fd_\B^{2k}\}$. Evaluating this equality at the unit
 element, yields, by Lemma \ref{evaluation}, the equality in $\CA^\infty$ (remember that $g\in\B$
 is fixed):
$$
\alpha\big(a\star_{\theta,\vec\tau}^\alpha b\big)(g)=L^*_{g^{-1}}\big(\alpha(a)\star_{\theta,\vec\tau}
\alpha (b)\big)(e)=\big(\alpha(a)\star_{\theta,\vec\tau}\alpha (b)\big)(g)\;.
$$
Hence,  the functions $\alpha\big(a\star_{\theta,\vec\tau}^\alpha b\big)$
and $\alpha(a)\star_{\theta,\vec\tau}\alpha (b)$ coincide.
This implies for $a,b,c\in\CA^\infty$:
$$
a\star_{\theta,\vec\tau}^\alpha\big(b\star_{\theta,\vec\tau}^\alpha c\big)=
\big(\alpha(a)\star_{\theta,\vec\tau}\alpha\big(b\star_{\theta,\vec\tau}^\alpha c\big)\big)(e)
=
\big(\alpha(a)\star_{\theta,\vec\tau}\big(\alpha(b)\star_{\theta,\vec\tau}\alpha (c)\big)\big)(e)\;,
$$
and the associativity of $\star_{\theta,\vec\tau}^\alpha$ on $\CA^\infty$ follows
from associativity of $\star_{\theta,\vec\tau}$ on the triple Cartesian product of the 
space $\CB^{\underline\nu}(\B,\CA)$, as stated in Theorem \ref{achievement}.
\end{proof}
\begin{rmk}
 Contrarily to the $\R^{2d}$-action case treated in \cite{Ri}, in the non-Abelian situation the 
original action is no longer an automorphism of the deformed product $\star_{\theta,\vec\tau}$ on 
$\CA^\infty$. This can be understood as  the chief reason 
 to introduce the whole oscillatory integrals machinery in chapter \ref{OSIL} and also 
  to consider the spaces $\CB^{\underline\mu}(\B,\CA)$.
\end{rmk}

To conclude this section, we establish   a formula 
 for the deformed product 
$\star_{\theta,\vec\tau}^\alpha$ on $\CA^\infty$, which in some sense, is more natural.
It will also clarify an important point, namely that the universal deformation of the
algebra\footnote{Recall that $C_{ru}(\B)$ denotes the $C^*$-algebra of right uniformly continuous 
and  bounded functions on $\B$, endowed with the sup-norm.}
$\CA=C_{ru}(\B)$, for the action $\alpha=R^\star$ coincides with 
$\big(\CB(\B),\star_{\theta,\vec\tau}\big)$.

\begin{prop}
\label{OIFSP-prop}
Let $(\alpha,\underline\mu^\alpha)$ be a strongly continuous and tempered action of a normal
$\bf j$-group $\B$ on a Fr\'echet algebra $\CA$. Then, for
 $a,b\in\CA^\infty$ and $\theta\in\R^*$, $\vec\tau\in{\bf\Theta}^N$, we have
\begin{equation}
\label{OIFSP}
a\star_{\theta,\vec\tau}^\alpha b=\widetilde{\int_{\B\times\B} K_{\theta,\vec\tau}}\, \big(\alpha(a)
\otimes\alpha(b)\big)\;,
\end{equation}
where we denote
$$
\alpha(a)
\otimes\alpha(b):\B\times\B\to\CA^\infty:\quad (x,y)\mapsto \alpha_x(a)\alpha_y(b)\;.
$$
\end{prop}
\begin{proof}
Since for $a\in\CA^\infty$,   the element
$\alpha(a)$ belongs to $ \CB^{\underline\mu}(\B,\CA^\infty)$,
$\underline\mu=\{\mu_j^\alpha\fd_\B^k\}$, by Lemma \ref{embed} and the Leibniz rule,
 we get that 
$$
\alpha(a)\otimes\alpha(b)\in \CB^{\underline\mu\otimes\underline\mu}
(\B\times\B,\CA^\infty)\;,
$$
which shows that the right hand side of \eqref{OIFSP} is indeed well defined. Next, by 
construction we have
$$
\alpha(a)\star_{\theta,\vec\tau}\alpha( b)=\widetilde{\int_{\B\times\B} K_{\theta,\vec\tau}}\Big(
{\mathcal R}\otimes {\mathcal R} \big(\alpha(a),\alpha(b)\big)\Big)\in \CB^{\underline\lambda}(\B,\CA^\infty)
\;,
$$
with
$$
\underline\lambda=\{\mu_j^{\alpha2L_jR_j}\fd_\B^{2k}\}\;,
$$
where the map ${\mathcal R}\otimes {\mathcal R} $ has been defined in Lemma \ref{genout}. Now,
using Lemma \ref{specific-limits}, we get with the element $e_n\in\CD(\B)$ defined in 
\eqref{en}, $n\in\N$, the equality in $\CB^{\underline\lambda}(\B,\CA^\infty)$:
$$
\alpha(a)\star_{\theta,\vec\tau}\alpha( b)=\lim_{n,m\to\infty}\int_{\B\times\B} K_{\theta,\vec\tau}(x,y)\, 
R^\star_x\big(e_n\alpha(a)\big)\,R^\star_y\big(e_m\alpha(b)\big)\,{\rm d}_\B(x)\,{\rm d}_\B(y)\;.
$$
By Lemma \ref{evaluation}, we know that the evaluation at the neutral element 
is continuous
from $\CB^{\underline\lambda}(\B,\CA^\infty)$ to $\CA^\infty$. Thus  we get
\begin{align*}
a \star_{\theta,\vec\tau}^\alpha b&=\big(\alpha(a)\star_{\theta,\vec\tau}\alpha( b)\big)(e)\\
&=
\lim_{n,m\to\infty}\int_{\B\times\B} K_{\theta,\vec\tau}(x,y)\, 
e_n(x)\alpha_x(a)\,e_m(y)\alpha_y(b)\,{\rm d}_\B(x)\,{\rm d}_\B(y)\;,
\end{align*}
one then concludes using Proposition \ref{uniqueness}.
\end{proof}

\begin{cor}
Let $\theta\in\R^*$ and $\vec \tau\in{\bf\Theta}^N$. For
$\CA=C_{ru}(\B)$ and $\alpha=R^\star$,
we have 
 $$
 \big(\CA^\infty,\star_{\theta,\vec\tau}^\alpha\big)=\big(\CB(\B),\star_{\theta,\vec\tau}\big)\;.
 $$
\end{cor}
\begin{proof}
By Lemma \ref{symbols} (ii), the set of smooth vectors in $C_{ru}(\B)$ for the right-regular action is
$\CB(\B)$. By the Proposition above, their algebraic structures coincide too.
\end{proof}
\section{Relation with the fixed point algebra}\label{RFPA}
Under slightly more restrictive conditions on the tempered action $(\alpha,\underline\mu^\alpha)$,
we exhibit a relationship between the deformed Fr\'echet algebra 
$(\CA^\infty,\star^\alpha_{\theta,\vec\tau})$ and a fixed point subalgebra of 
$\big(\CB^{\underline\mu^\alpha}(\B,\CA),\star_{\theta,\vec\tau}\big)$. This construction is very similar
to the construction of isospectral deformations of Connes and Dubois-Violette \cite{CDV}.
So, throughout this paragraph, we still assume that $\CA$ is a Fr\'echet algebra carrying a strongly continuous and tempered action $\alpha$ of a normal $\bf j$-group $\B$. But now, we further 
assume that the action is {\em almost-isometric}. By this, we mean that there exists
a family of tempered weights $\underline\mu^\alpha$ such that for $a\in\CA$ and all $g\in\B$, we have
$$
\|\alpha_g(a)\|_j=\mu_j^\alpha(g)\,\|a\|_j\;.
$$
\begin{ex}
For  any Lie group $G$, take $\CA=L^p(G)$,  $p\in[1,\infty)$.
Then, on this Banach space (not algebra), the right regular action 
 is almost isometric with  associated weight given by $\Delta_G^{1/p}$.
\end{ex}

Also, to simplify the discussion below,   we assume that each weight
$\mu_j^\alpha$ is sub-multiplicative.
We start with the simple observation that for any family of tempered weights $\underline\mu$, 
the extended action (defined in Lemma \ref{hat-alpha})
$\hat\alpha$ on  $\CB^{\underline\mu}(\B,\CA)$, commutes with the left regular action $L^\star$. 
This leads us to defined the commuting composite action 
$\beta:=\hat\alpha\circ L^\star=L^\star\circ\hat\alpha$, explicitly given by:
\begin{equation*}
 \big(\beta_g F\big)(g_0):=\alpha_g\big(F(g^{-1}g_0)\big)\,,\quad g,g_0\in\B,\,\quad 
 F\in \CB^{\underline\mu}(\B,\CA)\;.
\end{equation*}
Note also that by Lemma \ref{SmoothFamily} (ii) and Lemma \ref{hat-alpha}, for fixed $g\in\B$,
$\beta_g$ sends continuously $\CB^{\underline\mu}(\B,\CA)$ to $\CB^{\underline\nu}(\B,\CA)$,
with $\underline\nu=\{\mu_j^{R_j}\}$. Thus, in present context 
of sub-multiplicative weights, $\beta_g$ is continuous on $\CB^{\underline\mu}(\B,\CA)$.
Now observe that for an almost isometric action of a Lie group $G$ on a Fr\'echet algebra
$\CA$, the map  $\alpha:a\mapsto[g\mapsto\alpha_g(a)]$, is an isometric
embedding of $\CA^\infty$ into
$\CB^{\underline\mu^\alpha}(G,\CA)$. Indeed for all $j,k\in\N$, we have
\begin{align}
\label{isometric}
\|\alpha(a)\|_{j,k,\underline\mu^\alpha}&=\sup_{X\in\,\CU_k(\g)}\,\sup_{g\in G}\frac{\|\big(\widetilde X\, \alpha(a)\big)(g)\|_j}
{\mu_j^\alpha(g)\,|X|_k}\nonumber\\
&=\sup_{X\in\,\CU_k(\g)}\,\sup_{g\in G}\frac{\| \alpha_g\big(X^\alpha a\big)\|_j}
{\mu_j^\alpha(g)\,|X|_k}
=\sup_{X\in\,\CU_k(\g)}\,\frac{\|X^\alpha a\|_j}{|X_k|}=\|a\|_{j,k}\;.
\end{align}
By $\big(\CB^{\underline\mu^\alpha}(\B,\CA)\big)^\beta$, we denote the closed subspace of 
$\CB^{\underline\mu^\alpha}(\B,\CA)$ of fixed points for the action $\beta$.
It is then immediate to see that the image of $\CA^\infty$ under $\alpha$ lies
inside $\big(\CB^{\underline\mu^\alpha}(\B,\CA)\big)^\beta$. 
Reciprocally,  an element $F$ of $\big(\CB^{\underline\mu^\alpha}(\B,\CA)\big)^\beta$,  satisfies 
$F(g)=\alpha_g\big(F(e)\big)$ for all $g\in\B$,
i.e$.$  $F=\alpha(a)$ with $a:=F(e)\in\CA$. But by our 
 assumption of almost-isometry and \eqref{isometric}, we have 
 $\|a\|_{j,k}=\|F\|_{j,k,\underline\mu^\alpha}$, for all $j,k\in\N$,
and thus $a=F(e)$ has to be smooth. This proves that 
$\alpha:\CA^\infty\to\big(\CB^{\underline\mu^\alpha}(\B,\CA)\big)^\beta$ is an isomorphism
of Fr\'echet spaces, which is isometric  for each seminorm.  
Moreover, the map $\alpha$ is an algebra homomorphism.
Indeed,    by the arguments given in the proof of
Theorem \ref{UDF}, applied  to the case of an almost-isometric action with sub-multiplicative 
weights $\underline\mu^\alpha$,  for all $a,b\in\CA^\infty$ we have the equality
$$
\alpha\big(a\star_{\theta,\vec\tau}^\alpha b\big)=
\alpha(a)\star_{\theta,\vec\tau}\alpha (b)\quad\mbox{in}\quad
 \big(\CB^{\underline\mu^\alpha}(\B,\CA)\big)^\beta\;,
$$
which also shows that $\big(\CB^{\underline\mu^\alpha}(\B,\CA)\big)^\beta$
is an algebra for $\star_{\theta,\vec\tau}$.
In summary, we have proved the following:

\begin{prop}
\label{UDF2}
Let $\theta\in\R^*$, $\vec\tau\in{\bf\Theta}^N$ and let $(\CA,\alpha,\underline\mu^\alpha,\B)$ 
be a Fr\'echet algebra endowed 
with a strongly continuous tempered
and almost-isometric action
of a normal $\bf j$-group $\B$.
Then, we have  an isometric isomorphism of Fr\'echet algebras:
$$
\big(\CA^\infty,\star^\alpha_{\theta,\vec\tau}\big)\simeq
\Big(\big(\CB^{\underline\mu^\alpha}(\B,\CA)\big)^\beta,\star_{\theta,\vec\tau}\Big)\;.
$$
\end{prop}
\begin{rmk}
We stress that the assumption of sub-multiplicativity for the family of weights 
$\underline\mu^\alpha$,
associated with the tempered action $\alpha$,
is in fact irrelevant in the previous result. However it is unclear to us whether a
similar statement holds without the  assumption of almost-isometry. 
\end{rmk}

\section{Functorial properties of the deformed product}\label{FPDP}

To conclude with the deformation theory at the level of Fr\'echet algebras, 
we establish  some functorial properties. We come back to 
the general setting of a strongly continuous and tempered action
$(\alpha,\underline\mu^\alpha)$ of a normal $\bf j$-group $\B$ on a Fr\'echet
algebra $\CA$ (i.e$.$ we no longer assume that the action is almost isometric). 
We start with the question of 
algebra homomorphisms.

\begin{prop}
\label{HOM}
Let $(\CA,\{\|.\|_j\},\alpha)$, $(\CF,\{\|.\|'_j\},\beta)$ be two Fr\'echet algebras endowed with 
strongly continuous and tempered actions of a normal
$\bf j$-group $\B$ by automorphisms. Let also  
 $T:\CA\to\CF$ be a continuous  homomorphism which 
  intertwines the actions $\alpha$ and $\beta$. Then for any $\theta\in\R^*$ and $\vec\tau\in{\bf\Theta}^N$,  the map 
$T$ restricts to
a homomorphism from $(\CA^\infty,\star^\alpha_{\theta,\vec\tau})$ to 
$(\CF^\infty,\star^\beta_{\theta,\vec\tau})$. 
\end{prop}
\begin{proof}
Since by assumption $T\circ\alpha=\beta\circ T$, we get for any $P\in\CU(\b)$ that
$T\circ P^\alpha=P^\beta\circ T$, which entails that $T$ restricts to
a continuous map from $\CA^\infty$ to $\CF^\infty$.  The remaining part of the statement
follows then by Lemma \ref{continuous-maps}.
\end{proof}

 Next, we prove that if a Fr\'echet algebra is endowed with
a continuous  involution, then the latter will also define a continuous  involution for the deformed product,
under the mild condition\footnote{In Lemma \ref{continvol} we will see how to
suppress this extra condition.} that  $\bar\tau(-a)=\tau(a)$. 
Indeed, the latter implies that
$$
\overline{K_{\theta,\tau}}(x_1,x_2)=K_{\theta,\tau}(x_2,x_1)\;,
$$
so by  Lemma \ref{specific-limits}, we get:
\begin{prop}
Let $(\CA,\alpha)$ be a Fr\'echet algebras endowed with a 
strongly continuous tempered action  of a normal
$\bf j$-group $\B$. Assuming that  for $\theta\in\R^*$ and $\vec\tau\in{\bf\Theta}^N$, we have $\bar\tau_j(-a)=\tau_j(a)$, $j=1,\dots,N$, then any continuous involution of $\CA^\infty$
is a continuous involution of $(\CA^\infty,\star^\alpha_{\theta,\vec\tau})$ too.
\end{prop}

In a similar way, we deduce from Lemma \ref{specific-limits} that the deformation 
is ideal preserving:

\begin{prop}
\label{IDEAL}
Let $(\CA,\alpha)$ be a Fr\'echet algebras endowed with  a
strongly continuous tempered action  of a normal
$\bf j$-group $\B$ and $\theta\in\R^*$, $\vec\tau\in{\bf\Theta}^N$.
If $\CI$ is a closed $\alpha$-invariant  ideal of $\CA$, then $\CI^\infty$ is a
closed  ideal of $(\CA^\infty,\star^\alpha_{\theta,\vec\tau})$.
\end{prop}

We now examine the consequence of the fact that the constant function is the unit of $\big(\CB(\B),\star_{\theta,\vec\tau}\big)$. 

\begin{prop}
Let $(\CA,\alpha)$ be a Fr\'echet algebras endowed with a
strongly continuous and tempered action of a normal
$\bf j$-group $\B$ and $\theta\in\R^*$, $\vec\tau\in{\bf\Theta}^N$.
If $a\in\CA^\infty$ is fixed by the action $\alpha$, then for
$b\in\CA^\infty$, we have
$$
a\star^\alpha_{\theta,\vec\tau} b=ab,\qquad b\star^\alpha_{\theta,\vec\tau} a=ba\;.
$$
\end{prop}
\begin{proof}
This is a consequence of Proposition \ref{unit} together with the defining relation of the deformed 
product:
$$
a\star^\alpha_{\theta,\vec\tau} b=\big(\alpha(a)\star_{\theta,\vec\tau} \alpha(b)\big)(e)=\big(a\star_{\theta,\vec\tau} \alpha(b)\big)(e)=\big(a \alpha(b)\big)(e)=ab\;.
$$
The second equality is entirely similar.
\end{proof}

Next, we study the question of the existence of a bounded approximate unit for the Fr\'echet algebra
$(\CA^\infty,\star_{\theta,\vec\tau}^\alpha)$. We recall that a Fr\'echet algebra
$(\CA,\{\|.\|_j\})$ admits a bounded approximate unit if there exists a net 
$\{e_\lambda\}_{\lambda\in\Lambda}$ of elements of $\CA$ such that for any $a\in\CA$, the nets
$\{ae_\lambda\}_{\lambda\in\Lambda}$  and $\{e_\lambda a\}_{\lambda\in\Lambda}$  converges 
to $a$ and such that for each $j\in\N$, there exists $C_j>0$ such that for every 
$\lambda\in\Lambda$, we have $\|e_\lambda\|_j\leq C_j$.

\begin{prop}
\label{BAU}
Let $(\CA,\alpha)$ be a Fr\'echet algebra endowed with a strongly continuous and tempered 
action of a normal $\bf j$-group $\B$ and such that $\CA$ admits a bounded approximate unit.
Then for any  $\theta\in\R^*$, $\vec\tau\in{\bf\Theta}^N$, the Fr\'echet algebra 
$(\CA^\infty,\star^\alpha_{\theta,\vec\tau})$ 
admits a bounded approximate unit too.
\end{prop}
\begin{proof}
Let $\{f_\lambda\}$ be a net of bounded approximate units   for  $\CA$, let also
$0\leq \psi\in\CD(\B)$ be of $L^1$-norm one and define 
$$
e_\lambda:=\int_\B \psi(g)\,\alpha_g(f_\lambda)\,{\rm d}_\B( g)\;.
$$
Observe that even if $\{f_\lambda\}$ is not smooth, $\{e_\lambda\}$ is.
Indeed, for all  $X\in\CU(\b)$, we have 
$$
X^\alpha e_\lambda=\int_\B \underline X\psi(g)\,\alpha_g(f_\lambda)\,{\rm d}_\B( g)\;,
$$
and we get for the semi-norms 
defining the topology of $\CA^\infty$, with $\underline\mu^\alpha$ the family of tempered weights
associated to the temperedness of the action $\alpha$:
\begin{align*}
\|e_\lambda\|_{j,k}=\sup_{X\in\,\CU_k(\b)}\frac{\|X^\alpha e_\lambda\|_j}{|X|_k}&\leq
\sup_{X\in\,\CU_k(\b)}\int_\B |\underline X\psi|(g)\,\|\alpha_g(f_\lambda)
\|_j\,{\rm d}_\B( g)
\\&\leq\sup_{X\in\,\CU_k(\b)}\int_\B \frac{|\underline X\psi|(g)}{|X|_k}\,\mu^\alpha_j(g)\,{\rm d}_\B( g)
\times\|f_\lambda\|_j\;.
\end{align*}
Hence, the net $\{e_\lambda\}$ belongs to $\CA^\infty$ and  is  semi-norm-wise bounded in 
$\lambda\in\Lambda$
as $\|f_\lambda\|_j$ is. Next, we show that it is indeed an
 approximate unit for $\CA^\infty$: Since $\int \psi=1$, we first note that for any $a\in \CA$
$$
e_\lambda a-a=\int_\B \psi(g)\,\big(\alpha_g(f_\lambda)a-a\big)\,{\rm d}_\B( g)=
\int_\B \psi(g)\,\alpha_{g}\big(f_\lambda\,\alpha_{g^{-1}}(a)-\alpha_{g^{-1}}(a)\big)\,{\rm d}_\B(g)\;,
$$
which gives
$$
\|e_\lambda a-a\|_j\leq \int_\B \psi(g)\,\mu_j^\alpha(g)\,\|f_\lambda\,\alpha_{g^{-1}}(a)-\alpha_{g^{-1}}(a)\|_j\,{\rm d}_\B(g)\;,
$$
which converges to zero because $\|f_\lambda\,\alpha_{g^{-1}}(a)-\alpha_{g^{-1}}(a)\|_j$ does by assumptions and because $\psi$ is compactly supported. The general case is treated recursively exactly as in the proof of Lemma \ref{symbols} (viii).
Hence,  $\CA^\infty$ (with its original algebraic structure) admits a bounded approximate 
unit too. Now, we will  prove that a bounded 
approximate unit for $\CA^\infty$  is also a bounded approximate unit for 
$(\CA^\infty,\star^\alpha_{\theta,\vec\tau})$. 
 So, let  $\{e_\lambda\}$ be any bounded approximate unit for $\CA^\infty$. First observe that
 if we view the product $\star_{\theta,\vec\tau}$ as a bilinear map
$$
\star_{\theta,\vec\tau}:\CB(\B)\times\CB^{\underline\mu}(\B,\CA^\infty)\to
\CB^{\underline\nu}(\B,\CA^\infty)\;,\quad\underline\mu=\{\mu_j^\alpha\fd_\B^k\}_{j,k\in\N}\,,\;
\underline\nu=\{\mu_j^{\alpha L_j}\fd_\B^k\}_{j,k\in\N}\;,
$$
a slight adaptation of the arguments of Proposition \ref{unit} shows that for all $a\in\CA^\infty$:
$$
1\star_{\theta,\vec\tau}\alpha(a)=\alpha(a)\;,
$$
where $1$ denotes the unit element of $\CB(\B)$.
 Combining this with Proposition \ref{OIFSP-prop} gives the equality in $\CA^\infty$:
$$
e_\lambda \star_{\theta,\vec\tau}^\alpha a-a=\widetilde{\int_{\B\times\B} K_{\theta,\vec\tau}}\, \big(\alpha(e_\lambda)
\otimes\alpha(a)-1\otimes\alpha(a)\big)\;,
$$
where
\begin{align*}
&\alpha(e_\lambda)
\otimes\alpha(a)-1\otimes\alpha(a)\\&\qquad
:=\big[(x,y)\in\B\times\B\mapsto \alpha_x(e_\lambda)
\alpha_y(a)-\alpha_y(a)\big]\in \CB^{\underline\mu\otimes\underline\mu}(\B\times\B,\CA^\infty)\;.
\end{align*}
But by Proposition \ref{weak-conv}, we know that given $(j,k)\in\N^2$, there exist
positive integers $\vec r\in\N^{4N}$, such that (where  the differential operator $\bD_{\vec r}$ is
given in \eqref{bD}), we have 
$$
e_\lambda \star_{\theta,\vec\tau}^\alpha a-a=\int_{\B\times\B} K_{\theta,\vec\tau}(x,y)\, 
\bD_{\vec r}\big(\alpha_x(e_\lambda)
\alpha_y(a)-\alpha_y(a)\big)\,{\rm d}_\B(x)\,{\rm d}_\B(y)\;,
$$
with the integral   being absolutely convergent for the semi-norm $\|.\|_{j,k}$ of $\CA^\infty$. 
Now, take an increasing sequence $\{C_n\}_{n\in\N}$ of relatively compact open subsets in $G$, 
such that $\lim_n C_n=G$ and fix $\varepsilon>0$. By the absolute convergence 
in  the semi-norm $\|.\|_{j,k}$ of the integral above
and since the net $\{e_\lambda\}_{\lambda\in\Lambda}$ is bounded in the semi-norm $\|.\|_{j,k}$,
there exists $n\in\N$ such that
\begin{align*}
\Big\|\int_{\B\times\B\setminus C_n} K_{\theta,\vec\tau}(x,y)\, 
\bD_{\vec r}\big(\alpha_x(e_\lambda)
\alpha_y(a)-\alpha_y(a)\big),{\rm d}_\B(x)\,{\rm d}_\B(y)\Big\|_{j,k}\leq\varepsilon\;.
\end{align*}
Moreover, since $\{e_\lambda\}$ is an approximate unit for $\CA^\infty$, from a compactness
argument, we deduce that for any $n\in\N$, we have
\begin{align*}
\lim_\lambda\Big\|\int_{ C_n} K_{\theta,\vec\tau}(x,y)\, 
\bD_{\vec r}\big(\alpha_x(e_\lambda)
\alpha_y(a)-\alpha_y(a)\big)\,{\rm d}_\B(x)\,{\rm d}_\B(y)\Big\|_{j,k}=0\;.
\end{align*}
This concludes the proof as the arguments for $a\star_{\theta,\vec\tau}^\alpha e_\lambda$ are
similar.
\end{proof}

Lastly, we show that the deformation associated with a normal $\bf j$-group
coincides with the iterated deformations of each of its elementary normal $\bf j$-subgroups.

\begin{prop}
\label{defofdef}
Let $\B$ be a normal $\bf j$-group with Pyatetskii-Shapiro decomposition $\B=\B'\ltimes\S$, where
$\B'$ is a normal $\bf j$-group and  $\S$  is an elementary normal $\bf j$-group. 
Let $\CA$ be a Fr\'echet 
algebra endowed with a strongly continuous and tempered action $(\alpha,\underline\mu^\alpha)$ of 
$\B$. Denote by 
$\alpha^{\B'}$
(respectively by $\alpha^\S$) the restriction of $\alpha$ to $\B'$ (respectively to $\S$).
For $\CC$ a subspace of $\CA$, denote by $\CC^\infty_\B$ (respectively by $\CC^\infty_{\B'}$, 
$\CC^\infty_\S$) the set of smooth vectors in $\CC$ for the action of $\B$ (respectively of $\B'$, $\S$).
 Then, for $\theta\in\R^*$ and
 $\vec\tau=(\vec\tau',\tau_1)\in{\bf\Theta}^{N+1}$ 
 ($N$ is the number of elementary factors in $\B'$), we have
$$
\big((\CA^\infty_\S,\star^{\alpha^\S}_{\theta,\tau_1})_{\B'}^\infty,\star^{\alpha^{\B'}}_{\theta,\vec\tau'}\big)=(\CA^\infty_\B,\star_{\theta,\vec\tau}^\alpha)\;.
$$
\end{prop} 
\begin{proof}
Observe that being the restrictions of a strongly continuous and tempered action, 
the action $\alpha^\S$ of $\S$ on $\CA$ is
 also strongly continuous and tempered.  But the action $\alpha^{\B'}$ of $\B'$ on $\CA_\S^\infty$ is
 also strongly continuous (which is rather obvious) and tempered.
 To see that, note  that  for $g'\in\B'$ and $a\in \CA_\S^\infty$,
we have
\begin{align*}
\|\alpha^{\B'}_{g'}(a)\|_{j,k}=\sup_{X\in\,\CU_k(\s)}\frac{\|X^{\alpha^\S}\, \alpha^{\B'}_{g'}(a)\|_j}{|X|_k}
&=
\sup_{X\in\,\CU_k(\s)}\frac{\| \alpha^{\B'}_{g'}\big((\Ad_{g^{\prime-1}}X)^{\alpha^\S}a\big)\|_j}
{|X|_k}\\
&\leq \mu_j^\alpha(g')\,\sup_{X\in\,\CU_k(\s)}\frac{\| (\Ad_{g^{\prime-1}}X)^{\alpha^\S}a\|_j}
{|X|_k}\;.
\end{align*}
As $\B'$ acts on $\S$ by conjugation, it acts on $\CU_k(\s)$ and by Lemma \ref{lemtwo}, we deduce that
\begin{equation}
\label{eq-c}
\|\alpha^{\B'}_{g'}(a)\|_{j,k}\leq C(k)\, \mu_j^\alpha(g')\,\fd_\B(g')^k\,\|a\|_{j,k}\;.
\end{equation}
By Lemma \ref{fd-semi}, we have
$$
\fd_\B\big|_{\B'}\asymp \fd_{\B'}+\big[g'\mapsto (1+|\bR_{g'}|^2+|\bR_{g'^{-1}}|^2)\big]\;,
$$
and by Lemma \ref{lem:j-algebra-extension-morphism-tempered}, the extension homomorphism
$\bR$ is tempered. Hence,
 the action  $\alpha^{\B'}$ of $\B'$ on $\CA_\S^\infty$ is tempered with associated family
of tempered weights given by $\{\mu_j^\alpha\fd_\B^k\big|_{\B'}\}_{(j,k)\in\N^2}$.
Note also that the subspace of smooth vectors for $\B$ coincides with the subspace
of smooth vectors for $\B'$ within the subspace of smooth vectors for $\S$,
i.e$.$
$$
\CA_\B^\infty=(\CA_\S^\infty)_{\B'}^\infty\;.
$$
Indeed, the inclusion $\CA_\B^\infty\subset(\CA_\S^\infty)_{\B'}^\infty$ is clear since 
$a\in (\CA_\S^\infty)_{\B'}^\infty$ if and only if for all $X'\in\CU(\b')$, all $X\in\CU(\s)$
and all $j\in\N$, we have
$$
\|{X'}^{\alpha^{\B'}} X^{\alpha^\S}\,a\|_j<\infty\;,
$$
and $X'X\in\CU(\b)$.
But this also gives the reversed inclusion since $[\b',\s]\subset\s$, any element
of $\CU(\b)$ can be written as a finite sum of elements of the form $X'X$, with $X'\in\CU(\b')$
and $X\in\CU(\s)$. 

Next, we show that the action $\alpha^{\B'}$ of $\B'$ is by automorphisms  on the
 deformed Fr\'echet algebra
$(\CA^\infty_\S,\star^{\alpha^\S}_{\theta,\tau_1})$. First, by Proposition \ref{OIFSP-prop}
and Lemma \ref{specific-limits}, we get with the elements $e_n\in\CD(\S)$ defined in 
\eqref{en}, $n\in\N$, and for
 $a,b\in\CA^\infty_\S$:
\begin{align*}
a\star_{\theta,\tau_1}^{\alpha^\S} b&=\widetilde{\int_{\S\times\S} K_{\theta,\tau_1}}\, \big(\alpha(a)
\otimes\alpha(b)\big)\\
&=\lim_{n,m\to\infty}\int_{\S\times\S} K_{\theta,\tau_1}(x,y)\, 
e_n(x)\alpha_x^\S(a)\,e_m(y)\alpha_y^\S(b)\,{\rm d}_\S(x)\,{\rm d}_\S(y)\;.
\end{align*}
Observe also that \eqref{eq-c} shows that for $g'\in\B'$ fixed, the operator $\alpha^{\B'}_{g'}$ is continuous on 
$\CA^\infty_\S$. From this and the absolute 
convergence of the integrals in the product $\star_{\theta,\tau_1}^{\alpha^\S}$ at
the level of compactly supported functions, 
we deduce that
for $a,b\in \CA^\infty_\S$ and $g'\in \B'$:
\begin{align*}
&\alpha^{\B'}_{g'}\big(a\star_{\theta,\tau_1}^{\alpha^\S}\, b\big)
=\lim_{n,m\to\infty}\alpha^{\B'}_{g'}\Big(\int_{\S\times\S} K_{\theta,\tau_1}(x,y)\, 
e_n(x)\alpha_x^\S(a)\,e_m(y)\alpha_y^\S(b)\,{\rm d}_\S(x)\,{\rm d}_\S(y)\Big)\\
&=\lim_{n,m\to\infty}\int_{\S\times\S} K_{\theta,\tau_1}(x,y)\, 
e_n(x)\alpha_{g'x}(a)\,e_m(y)\alpha_{g'y}(b)\,{\rm d}_\S(x)\,{\rm d}_\S(y)\\
&=\lim_{n,m\to\infty}\int_{\S\times\S} K_{\theta,\tau_1}(x,y)\, 
e_n(x)\alpha^\S_{\bR_{g'}(x)}\big(\alpha^{\B'}_{g'}(a)\big)\,
e_m(y)\alpha^S_{\bR_{g'}(y)}\big(\alpha^{\B'}_{g'}(b)\big)
\,{\rm d}_\S(x)\,{\rm d}_\S(y)\;.
\end{align*}
Remember that
$$
\bR\in{\rm Hom}\big(\B',\Aut(\S,s,\omega^\S)\cap {\rm Sp}(V,\omega^0)\big)\;,
$$
where $(V,\omega^0)$ is the symplectic vector space attached to $\S$.
Using the invariance of the two-point kernel and  of the left Haar measure under
the action of ${\rm Sp}(V,\omega^0)=\Aut(\S)\cap\Aut(\S,s,\omega^\S)$, we get:
\begin{align*}
&\alpha^{\B'}_{g'}\big(a\star_{\theta,\tau_1}^{\alpha^\S}\,b\big)\\
&=
\lim_{n,m\to\infty}\int_{\S\times\S} K_{\theta,\tau_1}(x,y)\, 
e_n(x)\alpha^\S_{x}\big(\alpha^{\B'}_{g'}(a)\big)\,e_m(y)\alpha^S_{y}\big(\alpha^{\B'}_{g'}(b)\big)
\,{\rm d}_\S(x)\,{\rm d}_\S(y)\\
&=
\alpha^{\B'}_{g'}(a)\star_{\theta,\tau_1}^{\alpha^\S}\,\alpha^{\B'}_{g'}( b)\;.
\end{align*}
Thus, both Fr\'echet algebras 
$\big((\CA^\infty_\S,\star^{\alpha^\S}_{\theta,\tau_1})_{\B'}^\infty,\star^{\alpha^{\B'}}_{\theta,\vec
\tau'}\big)$ and $(\CA^\infty_\B,\star_{\theta,\vec\tau}^\alpha)$ are well defined and their 
underlying sets coincide. It remains to show that their  algebraic structures coincide too.
But this follows from  Proposition \ref{Fubini} as the extension homomorphism
$\bR$ of $\B=\B'\ltimes_\bR\S$ is tempered.
\end{proof}

\chapter{Quantization of polarized symplectic symmetric spaces}
\label{QPSSS}
 In the previous chapters, we defined a deformation of every Fr\'echet algebra that 
admits a strongly continuous (and tempered)
 action of a normal {\bf j}-group.
In particular, the method applies to every $C^*$-algebra which the group acts on. However, in that
 $C^*$-case, our procedure does not yet yield,  at this stage, pre-$C^*$-structures on the deformed 
 algebras. To cure this problem (in the case of an isometric action) 
 we will represent 
our deformed algebras by bounded operators on an Hilbert space. The present chapter consists in 
defining these Hilbert space representations.

 The construction relies on defining a unitary representation of the normal {\bf j}-group 
$\B$ at hand. This unitary representation
is obtained as the tensor product of irreducible unitary representations of the symplectic symmetric 
spaces (cf$.$ (\ref{SSs})) underlying the elementary factors $\S_j$
in the Pyatetskii-Shapiro decomposition $\B=\big(\S_N\ltimes\dots\big)\ltimes\S_1$ of the normal 
{\bf j}-group (cf$.$ (\ref{AES})).

 At the level of an elementary factor $\S$, the unitary representation Hilbert space will be 
defined through a variant of Kirillov's 
orbit method when viewing the symplectic symmetric space $\S$ as a polarized co-adjoint orbit of a 
central extension of its transvection group  (i.e$.$ the subgroup of the automorphism group 
of the symplectic symmetric space,
generated by products of even numbers of geodesic symmetries--see Proposition 
\ref{TRANSVECTIONGRP} below).
However, the construction applies to a much more general situation than the one of elementary 
normal {\bf j}-groups: the situation
of what we call ``elementary local symplectic symmetric spaces". In particular, for these 
spaces, we will obtain an explicit formula for the composition of symbols (see Proposition \ref{3PK}) 
which will be our main tool to investigate  the problem of $C^*$-deformations in the next sections. 
We therefore opted, within the present chapter \ref{QPSSS},  to start with presenting this more 
general situation in the sections
\ref{PSSS} to \ref{The three-point kernel}, and then to later pass to the particular case of elementary 
normal {\bf j}-groups in chapter \ref{QKLG}.

Once the elementary case is treated, one then needs to pass to the case of a normal 
{\bf j}-group $\B$. In that case, as already mentioned above, the Hilbert space $\CH$ will
consist in the tensor product of the Hilbert spaces $\CH_j$ representing its elementary factors 
$\S_j$. In order to define the $\B$-action on the tensor product Hilbert space $\CH$, we will need to 
represent every elementary group factor $\S_j$ on each $\CH_k$ for every 
$k$ less than or equal to $j$. The issue here is that the transvection group  of $\S_k$ does 
not  generally  contain $\S_j$ as a subgroup.
However, we will show that $\S_j$ injects into the subgroup of the automorphism group 
$\Aut(\S_k)$ of $\S_k$ 
that preserves Kirillov's polarization (see Proposition \ref{oulala}) on $\S_k$. The latter property will 
be shown, already in the more general case
of local symplectic symmetric spaces, to be sufficient to extend the action of $\S_k$ on $\CH_k$ to 
an action
of $\S_j\ltimes\S_k$ on $\CH_j\otimes\CH_k$ (see section \ref{Extensions  of polarization 
quadruples} below). Iterating this result will then lead to our definition of the unitary representation 
of  $\B$ on $\CH$ and its associated symbol calculus.

\quad

At the level of an elementary normal $\bf j$-group $\S$, the above mentioned quantization map 
(allowing to represent our algebras by bounded operators), realizes the program of the construction
of `covariant Moyal quantizers' in the sense of 
\cite{CGV} (see also  \cite[Section 3.5]{GBFV}). Such an object is of the 
following nature: Let $(M,\omega)$ be a symplectic manifold and let $G$ be a Lie subgroup
of the symplectomorphisms group ${\rm Symp}(M,\omega)$, together with $(\CH,U)$, a Hilbert 
space carrying a unitary representation of $G$. Then, a covariant Moyal quantizer, is a family
$\{\Omega(x)\}_{x\in M}$ of densely defined self-adjoint operators on $\CH$ such that
\begin{align}
U(g)\,\Omega(x)\,U(g)^\star&=\Omega(g.x)\,,\qquad \forall (g,x)\in G\times M\;,
\label{c1}\\
{\rm Tr}\big(\Omega(x)\big)&=1\,,\qquad \forall x\in M\;,
\label{c2}\\
{\rm Tr}\big(\Omega(x)\,\Omega(y)\big)&=\delta_x(y)\,,\qquad \forall (x,y)\in M\times M
\label{c3}\;,
\end{align}
where the traces have to be understood in the distributional sense and where $\delta_x(y)$ is a shorthand 
for the reproducing kernel of the Liouville measure on $M$.
Associated to a Moyal quantizer, there are quantization and dequantization  maps,
respectively given, on suitable domains and with $d\mu$ the Liouville measure on $M$, by
$$
\Omega:f\mapsto \int_M f(x)\,\Omega(x)\,d\mu(x)\;,
$$
and by
$$
\sigma:A\mapsto {\rm Tr}\big(\Omega(x) \,A\big)\;.
$$
Then, condition \eqref{c1} ensures that both quantization and dequantization maps are 
$G$-equivariant. Condition \eqref{c2} is a normalization condition and
condition \eqref{c3} says that $\Omega$ and $\sigma$ are formal inverses of each other. 
At the non-formal level, condition \eqref{c3} together with self-adjointness of the quantizers,
implies that $\Omega$ is a unitary operator from $L^2(M)$ to ${\mathcal L}^2(\CH)$,
 the Hilbert space of Hilbert-Schmidt operators on $\CH$. Transporting the algebraic structure
 from ${\mathcal L}^2(\CH)$ to $L^2(M)$, one therefore obtains a non-formal 
 and $G$-equivariant star-product at the level
 of square integrable functions:
 $$
 f_1\star f_2:=\sigma\big(\Omega(f_1)\,\Omega(f_2)\big)\;.
 $$
 This non-formal star-product turns to be a tri-kernel product, with distributional kernel
computable as it is
 given by the trace of three quantizers:
  $$
 f_1\star f_2(x)=\int_{M\times M} f_1(x_1)\,f_2(x_2)\,{\rm Tr}\big(\Omega(x)\,\Omega(x_1)\,
 \Omega(x_2)\big)\,d\mu(x_1)\,d\mu(x_2)\;.
 $$
 
 To construct such a quantizer for an elementary local symplectic symmetric spaces, we 
 essentially follow the pioneer ideas of Unterberger \cite{Un84} and recast them in a much 
 more general geometric context. Note that, however, in the work of Unterberger the condition
 \eqref{c3} does not hold but, as we will see in the present chapter,
  is restorable by a minor modification of his construction. This condition is 
 also called the `traciality property' as it eventually allows to prove that:
 $$
 \int_M   f_1\star f_2(x)\,d\mu(x)= \int_M   f_1(x)\, f_2(x)\,d\mu(x)\,,\qquad\forall f_1,f_2\in L^2(M)\;.
 $$
 Traciality is a property not shared by every quantization map (for example its fails for the 
 coherent-state quantization) and  proved to be fundamental in the work of Gracia-Bond\'ia 
 and V\'arilly. For instance, in \cite{CGV}, they were able to construct (and to prove uniqueness of)  an
 equivariant quantization map on the (regular) orbits of the Poincar\'e group in $3+1$
 dimensions solely from the axioms \eqref{c1},  \eqref{c2},  \eqref{c3}.

\section{Polarized symplectic symmetric spaces}\label{PSSS}

In this section, we introduce a particular class of symplectic symmetric spaces (called hereafter 
``polarized") for which 
we will be able to define a symmetric space variant of Kirillov's orbit method (see the next 
section). Within the present section,
after defining polarized symplectic symmetric spaces, we associate to every such space an 
algebraic object (its ``polarization quadruple") on which we will later base our unitary 
representation.

\quad

Let $(M,s,\omega)$ be a symplectic symmetric space and let $\Aut(M,s,\omega)$ be its
automorphism group and $\mathfrak{aut}(M,s,\omega)$ its derivation algebra (see Definition \ref{def:SSS}
and Proposition \ref{autSSS}).
Choose  a base point $o$ in $M$. Then the conjugation by the  
symmetry at $o$ yields an
involutive automorphism of the automorphism group:
\begin{align}
\label{sigma}
\sigma:\Aut(M,s,\omega)\to\Aut(M,s,\omega)\,,\quad g\mapsto s_o\circ g\circ s_o=:\sigma(g)\;.
\end{align}

The following result is a simple adaptation to the symplectic situation of a standard fact for general symmetric spaces \cite{Lo}:
\begin{prop}\label{TRANSVECTIONGRP}
The smallest subgroup $G(M)$ of $\Aut(M,s,\omega)$ 
that is stable under $\sigma$ and that acts transitively on $M$ is a Lie subgroup of $\Aut(M,s,\omega)$.
It coincides with group generated by products of an even number of symmetries:
$$
G(M)=\mbox{\rm gr}\{s_x\circ s_y\;|\;x,y\in M\}\;.
$$
The group $G(M)$ is called the {\bf transvection group} of $M$.
\end{prop}

We now come to the notion of polarized symplectic symmetric spaces:
\begin{dfn}
A symplectic symmetric space $(M,\omega,s)$ is said to be {\bf polarizable} if it admits a 
$G(M)$-invariant Lagrangian tangent distribution. A choice of such a transvection-invariant 
distribution  $\tilde{W}\subset TM$
determines a {\bf polarization} of $M$, in which case one speaks about a {\bf polarized} 
symplectic symmetric space.
\end{dfn}
The infinitesimal version of the notion of symplectic symmetric space is given in the two following 
definitions\footnote{Involutions $\sigma$ are denoted the same way either at the Lie 
group or Lie algebra level.}  (see \cite{Bi95, BCG}):
\begin{dfn}\label{SILA}
A {\bf symplectic involutive Lie algebra} (shortly a  ``siLa") is a triple $(\g,\sigma,\varpi)$  where $(\g,\sigma)$ is an {\bf involutive Lie algebra} (shortly an  ``iLa") i.e$.$
 $\g$ is a finite dimensional
real Lie algebra and $\sigma$ is an involutive automorphism of $\g$, and, where 
$\varpi\in\bigwedge^2\g^\star$ is a 
Chevalley two-cocycle on $\g$ (valued in the trivial representation on $\R$) such that,
denoting by 
$$
\g=\k\oplus\p\;,
$$
the $\pm1$-eigenspace decomposition of $\g$ associated to the involution 
$\sigma=:\id_\k\oplus(-\id_\p)$,
the cocycle $\varpi$ contains $\k$ in its radical and restricts to $\p\times\p$ as a non-degenerate two-form.
 A {\bf morphism} between two such siLa's is a Lie algebra homomorphism which
 intertwines both involutions and  two-cocycles.
\end{dfn}
\begin{dfn}\label{TST}
A {\bf transvection symplectic triple} is a siLa $(\g,\sigma,\varpi)$ where
\begin{enumerate}
\item[(i)] the action of $\k$ on $\p$ is faithful, and
\item[(ii)] $[\p,\p]=\k\;.$
\end{enumerate}
\end{dfn}

\begin{lem}
Every symplectic symmetric space determines a transvection symplectic triple.
\end{lem}
\begin{proof}
Let $(M,s,\omega)$ be a symplectic symmetric space. Define a 
transvection symplectic triple $(\g,\sigma,\varpi)$ as follows:
 $\g$ is the Lie algebra of  the transvection group $G(M)$, $\sigma$ is the
 restriction to $\g$  of the
differential of the involution 
 \eqref{sigma}   and $\varpi=\pi_{\star e}^\star\omega_o$ is the 
pullback of $\omega$ at the base point $o$ by the differential $\pi_{\star e}:\g\to T_oM$ 
of the projection
$\pi:G(M)\to M$, $g\mapsto g.o$. See \cite{Bi95} for the details.
\end{proof}
Defining the notion of isomorphism in the obvious way,  one 
knows from \cite{Bi95, BCG} the following 
result. It is a symplectic adaptation
of the classical analogue for general affine symmetric spaces \cite{Lo}.

\begin{prop}  
\label{correspondence}
The correspondence  $(M,s,\omega)\mapsto(\g,\sigma,\varpi)$ described above  induces a bijection
between the isomorphism classes of simply connected symplectic symmetric spaces and the 
isomorphism classes of transvection
symplectic triples.
\end{prop}

 In the above statement, the ``reverse direction" (i.e$.$ $(\g,\sigma,\varpi)\mapsto(M,s,\omega)$) is obtained as follows. One starts by considering the (abstract) connected simply connected Lie group $G$ whose Lie algebra is $\g$. Then the automorphism $\sigma$
of $\g$ uniquely determines an automorphism of $G$, again denoted by $\sigma$. The connected component $K$ of the subgroup
of $G$ constituted by the elements that are fixed by $\sigma$ is then automatically a closed subgroup of $G$. The coset space $M:=G/K$
is then naturally a connected smooth manifold which is also simply connected consequently to the  connectedness of $K$ and simple connectedness of $G$. One then check that the formulae
\begin{equation}\label{SSF}
s_{gK}(g'K):=g\sigma(g^{-1}g')K\;,
\end{equation}
define a structure of symmetric space $(M,s)$ in the sense of Loos. The Lie algebra of $K$  then coincides with $\k$, implying a natural isomorphism $\p\to T_K(M)$. The latter yields a Lie algebra homomorphism from the original iLa $(\g,\sigma)$ to the  transvection iLa of $(M,s)$. This homomorphism turns out to be an isomorphism due to the above conditions (i) and (ii) imposed on $(\g,\sigma)$.
 Note that from the above exposition, one extracts
\begin{lem}
Every iLa determines a simply connected symmetric space.
\end{lem}

Note also that two non-isomorphic iLa's (when non-transvection) could determine the same symmetric space. More precisely:
\begin{dfn}\label{SAMESS}
Let $(\g_j,\sigma_j)$ ($j=1,2$) be two  iLa's with associated simply connected symmetric spaces  denoted by $(M_j,s_j)$ respectively. One says that they {\bf determine the same simply connected symmetric space} if $M_1$ and $M_2$ are isomorphic as symmetric spaces.
\end{dfn}
Of course, one has an analogous definition in the symplectic case.

\vspace{3mm}

 Now given a siLa $(\g,\sigma,\varpi)$, the two-cocycle $\varpi$ is generally not exact, or equivalently, the symplectic action of the 
transvection group on the symplectic symmetric space
is not Hamiltonian  (see \cite{Bi95,BCG}). However, it is always possible to centrally extend the transvection group in such a 
way that the extended group acts 
on $M$ in a Hamiltonian way. The associated moment mapping then is a symplectic 
equivariant covering onto a co-adjoint orbit, in accordance with the classical general result for Hamiltonian homogeneous symplectic spaces \cite{S}.
The situation we consider in the present article concerns such \emph{non-exact} transvection triples 
underlying polarized symplectic symmetric spaces. 
\begin{lem}
\label{pol}
Let $(M,s,\omega)$ be a symplectic symmetric space, polarized by a transvection-invariant 
Lagrangian distribution $\tilde{W}\subset TM$.  These data correspond (via the correspondence 
of Proposition \ref{correspondence}) to a $\k$-invariant Lagrangian subspace $W$ in $\p$.
\end{lem}
\begin{proof} Under  the linear isomorphism 
$\pi_{\star e}|_{\p}:\p\to T_oM$ the subspace $\tilde{W}_o$
of $T_oM$ corresponds to a Lagrangian subspace $W$ of the symplectic vector space
$(\p,\varpi)$.
\end{proof}

According to the previous Lemma, we use the following terminology:
\begin{dfn}
A siLa  $(\g,\sigma,\varpi)$ is called {\bf polarized} if it
is endowed with $W$, a $\k$-invariant  Lagrangian subspace of $(\p,\varpi)$.
\end{dfn}

Let $(\g,\sigma,\varpi)$ be a non-exact 
transvection symplectic triple  (i.e$.$ the Chevalley two-cocycle $\varpi$ is not exact)  polarized by a Lagrangian subspace
$W\subset \p$. Let us consider $\fD$, the  algebra of $W$-preserving symplectic endomorphisms of
 $\p$. Note that the faithfulness condition (i) of Definition \ref{TST} implies the inclusion 
$\k\subset\fD$.
The vector space $\fD\oplus\p$ then naturally carries a structure of Lie algebra (containing $\g$)
 that underlies
a siLa. We centrally extend the latter in order to define a new siLa:
$$
\fL:=\fD\oplus\p\oplus\R Z\;,
$$
 with table given by
$$
[X,Y]_\fL:=[X,Y]+\varpi(X,Y)\,Z\,,\quad [X,Z]_\fL:=0\,,\quad \forall X,Y\in \fD\oplus\p\;,
$$
where $[.,.]$ denotes the Lie bracket in $\fD\oplus\p$ and where we have extended the 2-form 
$\varpi$ on $\p$
to a 2-form on the entire Lie algebra $\fD\oplus\p$ by zero on $\fD$. 
\begin{lem}\label{EXACTSILA}
Let $(\g,\sigma,\varpi)$ be a non-exact polarized  transvection symplectic triple.
Within the notations given above,
consider the element $\xi\in\fL^\star$ defined by
$$
\langle\xi,Z\rangle=1\,,\quad\xi\big|_{\fD\oplus\p}=0\;.
$$
Define moreover:
$$
\tilde{\fD}:=\fD\oplus\R Z\,,\quad\sigma_\fL:=\id_{\tilde{\fD}}\oplus(-\id_\p)\;,
$$
and
$$
 \tilde\g=\g\oplus\R Z\,,\quad \tilde\k=\k\oplus\R Z\,,\quad \tilde\sigma:=\id_{\tilde{\k}}
\oplus(-\id_\p)\;.
$$
Then, the triples $(\fL,\sigma_\fL,\delta\xi)$ and $(\tilde\g,\tilde\sigma,\delta\xi\big|_{\tilde \g})$ are 
exact siLa's.
\end{lem}
\begin{proof}
We give the proof for the first triple  $(\fL,\sigma_\fL,\delta\xi)$ only, the second case being handled in 
a similar way.
Since for $X,Y\in\p$, we have $[X,Y]\in\k\subset \fD$, we get
$$
\delta\xi(X,Y)=\langle \xi,[X,Y]_\fL\rangle=\langle \xi,[X,Y]+\varpi(X,Y)Z\rangle=\varpi(X,Y)\;,
$$
which at once proves closedness, non-degeneracy and $\tilde\fD$-invariance of the 2-form 
$\delta\xi$ on $\p$.
\end{proof}
\begin{rmk}
The exact siLa's $(\fL,\sigma_\fL,\delta\xi)$ and $(\tilde\g,\tilde\sigma,\delta\xi\big|_{\tilde \g})$
and the non-exact siLa $(\g,\sigma,\varpi)$  all three determine the same simply connected
symplectic symmetric space.
\end{rmk}
\begin{dfn}\label{POLARIZATION}
Given an exact siLa $(\g,\sigma,\varpi)$ (i.e$.$ $\varpi=\delta\xi$ for an element $\xi\in\g^\star$), 
by a {\bf polarization} affiliated to  
$\xi$, we mean   a $\sigma$-stable Lie subalgebra $\b$ of $\g$ containing $\k$ and maximal  for
 the property of being isotropic with respect to the two-form $\varpi$. 
\end{dfn}

 The following statement is classical (see e.g$.$ \cite{Lo}).
 \begin{lem}\label{ILGSS}
Let $(G,\sigma)$ be a connected  involutive Lie group i.e$.$ a connected Lie group $G$ equipped with an involutive automorphism $\sigma$. Let us denote by $K$ the connected component  of the subgroup of $G$ constituted by the $\sigma$-fixed elements.
Then $K$ must be closed and Formula (\ref{SSF}) defines a structure of symmetric space on the quotient manifold $M=G/K$.
\end{lem}
Note that in this slightly more general situation, $M$ need not  necessarily be simply connected. 

Within the above setting let us denote by  $(\g,\sigma)$ the involutive Lie algebra  associated to $(G,\sigma)$. Let us furthermore assume that it underlies an exact siLa with polarization $\b$ as in Definition \ref{POLARIZATION}.
In that context, $M$ automatically becomes a polarized symplectic symmetric space. Denote  by $B:=\exp\{\b\}$ the analytic (i.e$.$ connected) Lie subgroup of $G$ with 
Lie algebra $\b$. One has $K\subset B$ and we will always assume $B$ to be closed in $G$. Since
 $B$ is stable under $\sigma$, the coset space $G/B$ admits the following natural family 
 of involutions $\underline\sigma$:
\begin{equation}\label{SSACTION}
M\times G/B\to G/B\,,\quad
(gK,g_0B)\mapsto \underline{\sigma}_{gK}(g_0B):=g\sigma(g^{-1}g_0)B\;.
\end{equation}

\begin{dfn}
\label{def:PQ}
With the same notations as above, 
the quadruple $(G,\sigma,\xi,B)$ is called a {\bf polarization quadruple}. Its infinitesimal version
$(\g,\sigma,\xi,\b)$ is called the associated {\bf infinitesimal polarization quadruple}.
A {\bf morphism} between two polarization quadruples $(G_j,\sigma_j,\xi_j,B_j)$, $j=1,2$, is 
defined as a  Lie group homomorphism
$$
\phi: G_1\to G_2\;,
$$
that intertwines the involutions, such that $\phi( B_1)\subset B_2$ and such that, denoting again 
by $\phi$ its differential at the unit element, one has
$\phi^\star\xi_2=\xi_1$.
\end{dfn}
The map (\ref{SSACTION}) corresponds to  an `action' of the symmetric
space $M= G/K$ on the manifold $G/B$. The following result is a consequence of immediate computations.
\begin{lem}\label{MGB}
Let $(G,\sigma,\xi,B)$ be a polarization quadruple. Then,
the following properties hold for all $x$ and $y$ in $M$:
$$
\underline{\sigma}_x^2=\Id_{G/B}\,,\quad \underline{\sigma}_x\circ\underline{\sigma}_y\circ\underline{\sigma}_x= 
\underline{\sigma}_{s(x,y)}\;.
$$
Moreover, we have the $G$-equivariance property:
$$
g\circ\underline{\sigma}_x\circ g^{-1}=\underline{\sigma}_{g.x}\,,
\qquad\forall g\in G\,,\,\forall x\in M\;.
$$
\end{lem}
We also observe that under mild conditions on the modular functions
of $G$ and $B$, a $G$-invariant and $M$-invariant measure
always exists  on the manifold $G/B$:
\begin{lem}
\label{G-M-invariant}
Let $(G,\sigma)$ be an involutive Lie group and $B$ a $\sigma$-stable closed subgroup of $G$
such that the modular function of $B$ coincides with the restriction to $B$ of the modular function
of $G$.
Then, there exists a (unique up to normalization) Borelian measure ${\rm d}_{G/B}$ on 
the manifold $G/B$ which is both invariant under $G$ and under the action of $M=G/K$ given in  
\eqref{SSACTION}.
\end{lem}
\begin{proof}
Under the closedness  condition of $B$ and under the coincidence assumption for the modular 
functions, it is well known  that   there exists a (unique up to normalization) 
$G$-invariant Borelian measure ${\rm d}_{G/B}^o$ on  $G/B$. Now, define 
${\rm d}_{G/B}:={\rm d}_{G/B}^o+\underline\sigma_K^\star\, {\rm d}_{G/B}^o$.
As $\underline\sigma_K$ is an involution, the latter measure is $\underline\sigma_K$-invariant.
Moreover, from 
$L^\star_g\circ \underline\sigma_K^\star=\underline\sigma_K^\star\circ L^\star_{\sigma(g)}$ on
$G/B$ for all $g\in G$, we deduce that ${\rm d}_{G/B}$ is also $G$-invariant. By uniqueness of
${\rm d}_{G/B}^o$, the latter is a multiple of the former. To conclude with the $M$-invariance, 
it suffices
to observe  that for all $g\in G$, we have $\underline\sigma_{gK}^\star=L_g^\star\circ
\underline\sigma_K^\star\circ L_{g^{-1}}^\star$.
\end{proof}

We end this section by constructing two canonical exact  polarization quadruples out of
a non-exact transvection triple. We omit the proof which is immediate.
\begin{prop}
\label{pol-quad}
Let $(\g,\sigma,\varpi)$ be a non-exact  transvection symplectic triple polarized by $W\subset \p$. 
Within the context of Lemma \ref{EXACTSILA}, we set $\fB:=\tilde{\fD}\oplus W$ and 
$\b:=\tilde\k\oplus W$.
Then the quadruples $(\fL,\sigma_\fL,\xi,\fB)$ and $(\tilde \g,\tilde\sigma,\xi\big|_{\tilde\g},\b)$ 
are polarization quadruples.
\end{prop}

The latter observation leads us to introduce the following terminology:
\begin{dfn}
\label{TRANSFULL}
The polarization quadruple $(\fL,\sigma_\fL,\xi,\fB)$ associated to a
non-exact polarized transvection triple $(\g,\sigma,\varpi)$, is called the (infinitesimal) 
{\bf full polarization quadruple}.
The sub-quadruple $(\tilde{\g},\tilde{\sigma},\xi\big|_{\tilde\g},\b)$  is called the associated  (infinitesimal) 
{\bf transvection quadruple}.

 In the sequel, we will denote by $\L$ the connected, simply connected Lie group with 
Lie algebra $\fL$ and we will consider the connected Lie subgroup $\tilde{G}$ of $\L$ tangent to 
$\tilde{\g}$. We will denote by $\bB$ (respectively $B$) 
the connected Lie subgroup of $\L$ (respectively of $\tilde G$) associated to $\fB$ (respectively to $\b$).
\end{dfn} 
 We summarize the present section by the following 
\begin{prop}
\begin{enumerate}
\item[(i)] Every symplectic symmetric space uniquely determines a transvection symplectic triple (cf$.$ Definition \ref{TST}).
\item[(ii)] Every siLa (cf$.$ Definition \ref{SILA}) uniquely determines a simply connected symplectic symmetric space.
\item[(iii)] In the simply connected  \underline{and} transvection case, the correspondences mentioned in items (i) and (ii) are inverse to one another.
\item[(iv)] Every non-exact transvection symplectic triple $(\g,\sigma,\varpi)$ (cf$.$ Definition \ref{TST} with non-exact two-cocycle) uniquely determines a pair of exact polarization quadruples: the associated infinitesimal full polarization quadruple $(\fL,\sigma_\fL,\xi,\fB)$ and its transvection sub-quadruple $(\tilde{\g},\tilde{\sigma},\xi\big|_{\tilde\g},\b)$ (cf$.$ Proposition \ref{pol-quad} and Definition \ref{TRANSFULL}).
\item[(v)] The three siLa's involved in item (iv): $(\g,\sigma,\varpi)$, $\big(\tilde{\g},\tilde{\sigma},\delta
\xi\big|_{\tilde\g}\big)$ and 
$\left(\fL,\sigma_\fL,\delta\xi\right)$ all  determine the same simply connected polarized symplectic symmetric space (cf$.$ Definition \ref{SAMESS}).
\end{enumerate}
\end{prop}

\section{Unitary representations of symmetric spaces}
\label{KIRILLOV}

 In this section, we fix $(\g,\sigma,\varpi)$ a non-exact polarized transvection triple,
to which  we   associate a  (connected and simply connected) transvection quadruple $(\tilde{G},\tilde\sigma,\xi,B)$, according to the
construction underlying Definition \ref{TRANSFULL}.
We start with the following
 {\em pre-quantization} condition: in the sequel we will  always assume that
  the character $\xi|_\b:\b\to\R$ exponentiates to $B$
 as a unitary character
 $$
 \chi:B\to U(1)\,,\qquad b\mapsto\chi(b)\;.
 $$
By this we mean that we assume the existence of a Lie group homomorphism $\chi$ whose differential at the identity coincides with $\xi|_\b$.
Note that then, the character is automatically fixed by the restriction to $B$ of the involution:
$$
\sigma^\star\chi=\chi\;.
$$
Of course, the pre-quantization condition is satisfied when the group $B$ is exponential,
as it will be the case for Pyatetskii-Shapiro's elementary normal $\bf j$-groups:

\begin{equation}
\label{unitary-char}
\chi(b):=e^{i\langle\xi,\log(b)\rangle}\,,\quad b\in B\;.
\end{equation}

 \begin{lem}
 Let $(\tilde \g,\tilde\sigma,\xi,\b)$ be the transvection quadruple of 
 a non-exact  trans\-vection triple $(\g,\sigma,\varpi)$ such that $B$ is exponential. Then,
the pre-quantization condition is satisfied.
\end{lem}
\begin{proof}
Since $B$ is exponential, by the BCH formula,   the statement will follow 
from $\xi\big([\b,\b]_{\tilde g}\big)=0$.
By construction of $\xi$ (see Definition \ref{EXACTSILA}), this will follow if the $Z$-component
of    $[\b,\b]_{\tilde g}$ vanishes. But the latter reads $\varpi(\b,\b)Z$ which reduces to zero by 
Definition \ref{POLARIZATION} of a polarization quadruple and by Proposition \ref{pol-quad} which shows that 
$(\tilde g,\tilde\sigma,\xi,\b)$ is indeed  
a polarization quadruple. 
\end{proof}

We then form the line bundle:
$$
E_\chi:=\tilde{G}\times_\chi\C\to \tilde G/B\;,
$$
and consider the associated induced representation of $\tilde{G}$ on the smooth sections
 $\Gamma^\infty(E_\chi)$. We will denote the latter representation by   $U_\chi$.  Identifying 
 as usual $\Gamma^\infty(E_\chi)$
with the space of $B$-equivariant functions:
\begin{align}
\label{sections}
\Gamma^\infty(E_\chi)&\simeq C^\infty(\tilde{G})^B\nonumber
\\&:=\big\{\hat{\varphi}\in C^\infty(\tilde{G})\;|\;
\hat{\varphi}(gb)
=\overline{\chi}(b)\hat{\varphi}(g),\;\forall b\in B,\,\forall g\in \tilde{G}\big\}\;,
\end{align}
the representation $U_\chi$ is  given by the restriction to $C^\infty(\tilde{G})^B$ of the 
left-regular representation:
\begin{equation}
\label{U}
[U_\chi(g)\tilde\varphi]^\wedge(g'):=\hat{\varphi}(g^{-1}g')\,,\quad \forall \tilde\varphi\in
\Gamma^\infty(E_\chi) \;.
\end{equation}
 We endow the line bundle $E_\chi$ with the Hermitian structure, defined in terms of the 
identification \eqref{sections} by:
$$
h_{gB}\big(\tilde\varphi_1,\tilde\varphi_2\big):=\overline{\hat\varphi_1(g)}\,\hat\varphi_2(g)\,,\qquad 
\forall\tilde\varphi_1,\tilde\varphi_2\in\Gamma^\infty(E_\chi)\,,\quad gB\in \tilde G/B\;.
$$ 
We make the assumption that the modular function of $B$ coincides with the restriction
to $B$ of the  modular function of $\tilde G$. By Lemma \ref{G-M-invariant}, this condition
implies  the existence of a $\tilde G$-invariant and $\underline \sigma_{\tilde K}$-invariant
Borelian  measure 
$ {\rm d}_{ \tilde G/B}$ on $\tilde G/B$.  Here,  $\underline \sigma_{\tilde K}:\tilde G/B\to \tilde G/B$, $gB\mapsto\tilde \sigma(g)B$ is the involutive 
diffeomorphism given in \eqref{SSACTION} for the involutive pair 
$(\tilde G,\tilde \sigma)$ and subgroup $B$, underlying the  transvection quadruple 
$(\tilde{G},\tilde\sigma,\xi,B)$. 
We then let $\CH_\chi$ be the Hilbert space completion of $\Gamma_c^\infty(E_\chi)$ 
 for the inner product:
$$
\langle\tilde\varphi_1,\tilde\varphi_2\rangle
:=\int_{ \tilde G/B}\,h_{gB}(\tilde\varphi_1,\tilde\varphi_2)
\,{\rm d}_{ \tilde G/B}(gB)\;.
$$
Of course, the induced representation $U_\chi$ of $\tilde{G}$ then naturally acts on 
$\CH_\chi$ by unitary operators.
Now observe  that the  $\tilde{\sigma}$-invariance of character $\chi$, 
implies that the pull back under
$\tilde{\sigma}$ of an equivariant function is again equivariant. Therefore, we get
a linear involution:
\begin{align*}
\Sigma: \CH_\chi\to \CH_\chi\,,\quad 
[\Sigma\tilde\varphi]^\wedge:=\tilde{\sigma}^\star\hat{\varphi}\;.
\end{align*}
Also,  the
 $\underline{\sigma}_{\tilde K}$-invariance of  the measure $ {\rm d}_{ \tilde G/B}$  implies:
\begin{align*}
\langle\Sigma\tilde\varphi_1,\Sigma\tilde\varphi_2\rangle&=\int_{ \tilde G/B}\,
h_{gB}(\Sigma\tilde\varphi_1,
\Sigma\tilde\varphi_2)
\,{\rm d}_{ \tilde G/B}(gB)\\
&=
\int_{ \tilde G/B}\,h_{\underline\sigma_{\tilde K}(gB)}
(\tilde\varphi_1,\tilde\varphi_2)\,{\rm d}_{ \tilde G/B}(gB)
=\langle\tilde\varphi_1,\tilde\varphi_2\rangle\;,
\end{align*}
for all $\tilde\vf_1,\tilde\vf_2\in\CH_\chi$, 
showing that $\Sigma$ is not only involutive but also self-adjoint. Thus the element $\Sigma$ 
belongs to $\CU_{sa}(\CH_\chi)$, the collection of   unitary and self-adjoint operators on $\CH_\chi$.
 When composed with the representation $U_\chi$ of $\tilde G$, the operator $\Sigma$
satisfies the following  properties, whose proofs  consist in direct computations:

\begin{prop}
\label{propO}
 Let $(M,s,\omega)$ be the  polarized 
symplectic symmetric space associated to a  transvection quadruple $(\tilde{G},\tilde\sigma,\xi,B)$ (cf$.$ Lemmas \ref{ILGSS} and \ref{MGB}).
Assume that  the modular function of $B$ coincides with the restriction to $B$ of the modular function of
$\tilde G$.
Then the map
$$
\tilde G\to\CU_{sa}(\CH_\chi)\,,\quad g\mapsto U_\chi(g)\,\Sigma\, U_\chi(g)^*\;,
$$
is constant on the left cosets of $\tilde K$ in $\tilde G$. The corresponding 
mapping:
$$
\Omega:M=\tilde G/\tilde K\to\CU_{sa}(\CH_\chi)\,,\quad g\tilde K\mapsto\Omega(g\tilde K):=
U_\chi(g)\,\Sigma\, U_\chi(g)^*\;,
$$
defines a {\bf unitary representation} of the symmetric space $M=\tilde G/\tilde K$ in the sense that, for all $x$, $y$
in $M$ and $g$ in $\tilde G$, the following representative properties hold:
 \begin{align}
 \label{RP}
 \Omega(x)^2&=\Id_{\CH_\chi}\;,\nonumber\\
\Omega(x)\,\Omega(y)\,\Omega(x)&=\Omega(s_xy)\;,\\
U_\chi(g)\,\Omega(x)\, U_\chi(g)^*&=\Omega(g.x)\;.\nonumber
\end{align}
\end{prop}

\begin{dfn}
The pair $(\CH_\chi,\Omega)$ is called the {\bf unitary representation of $(M,s)$ induced
by the character} $\chi$ of $B$.
\end{dfn}

We are now ready to define our prototype of quantization map on a polarized 
symplectic symmetric space:
\begin{dfn}
\label{QM}
Let $(M,s,\omega)$ be a the  polarized 
symplectic symmetric space.
 Denote by $L^1(M)$ the space of integrable functions on $M$ with respect to the $\tilde G$-invariant 
 (Liouville) measure ${\rm d}_M$.
Denote by $\CB(\CH_\chi)$ the space of bounded linear operators on 
the Hilbert space $\CH_\chi$.
Consider the $\tilde G$-equivariant continuous linear map:
$$
\Omega: L^1(M)\to\CB(\CH_\chi)
\,,\quad f\mapsto\Omega(f):=\int_M\,f(x)\,\Omega(x)\,{\rm d}_M(x)\;.
$$
The latter is called the {\bf quantization map} of $M$ induced by the transvection quadruple 
$(\tilde{G},\tilde\sigma,\xi,B)$.
\end{dfn}
\begin{rmk}
 From $\|\Omega(x)\|=\|\Sigma\|=1$ (the norm here is the uniform norm on $\CB(\CH_\chi)$), 
we get the obvious estimate 
$\|\Omega(f)\|\leq\|f\|_1$, from which the continuity of the quantization map follows.
Also, from Proposition \ref{propO} and from the $\tilde G$-invariance of ${\rm d}_M$, the covariance 
property at the level of the quantization map reads:
\begin{equation*}
U(g)\,\Omega(f)\,U(g)^*=\Omega (^g\!f)\,,\qquad \forall f\in L^1(M)\,,\;\forall g\in \tilde G\;,
\end{equation*}
where $^g\!f:=[g_0\tilde K\mapsto f(g^{-1}g_0\tilde K)]$.
The latter equivariance property, 
 under the full group of automorphisms of the symplectic symmetric space, is an important difference 
 between the present 
 ``Weyl-type"
construction and the classical
coherent-state-quantization approach. For instance, as it appears already in the flat case of 
$\R^{2n}$, 
holomorphic coherent-state-quantization (i.e$.$ Berezin-Toeplitz in that case) yields a equivariance group that is isomorphic to $U(n)\ltimes\C^n$,
while Weyl  quantization is equivariant under the full automorphism group of $\R^{2n}$ i.e$.$  
$\mbox{\rm Sp}(n,\R)\ltimes\R^{2n}$. 
Another essential difference is unitarity (see Proposition \ref{Berezin-prop} below).
Lastly, the quantization defined above is in 
general
not positive, i.e$.$ for $0\leq f \in L^1(M)$, $\Omega(f)$ is not necessarily a positive operator. 
But since $\Omega(f)^*=\Omega(\overline f)$,  it maps real-valued functions to self-adjoint 
operators.  
\end{rmk}

\begin{rmk}
The quantization map of Definition \ref{QM} is a generalization of the Weyl quantization,
from the point of view of symmetric spaces. Moreover, we will see that for the symmetric
space underlying a two-dimensional elementary normal $\bf j$-group (i.e$.$ for the affine
group of the real line)  this construction coincides with  Unterberger's
Fuchs calculus  \cite{Un84}.
\end{rmk}

Our next step is to 
introduce a functional parameter in the construction of the quantization map.
There are several reasons for doing this, among which there is one of a purely analytical nature: 
obtaining a quantization map which is a unitary operator from the Hilbert space of square 
integrable symbols, $L^2(M)$, to the Hilbert space of 
Hilbert-Schmidt operators on $\CH_\chi$, denoted by $\mathcal L^2(\CH_\chi)$.
This unitarity property will enable us to define a non-formal $\star$-product on $L^2(M)$
in a straightforward way. 
\begin{dfn}
\label{Sigma-m}
Identifying  a  Borelian function  ${\bf m}$ on $\tilde G/B$ with the  operator on 
$\CH_\chi$ of point-wise 
multiplication by this function, we let
$$
\Sigma_{\bf m}:={\bf m}\circ\Sigma\;.
$$
\end{dfn}

When $\bm$ is locally essentially bounded,
 the family of operators $U_\chi(g)\Sigma_{\bf m}U_\chi(g)^*$, $g\in \tilde G$,  can be defined
 on the common 
domain $\Gamma^\infty_c(E_\chi)$. Note however that the latter family of operators is
 not necessarily constant 
on the left cosets of $\tilde K$ in 
$\tilde G$ and unless $\underline{\sigma}_{\tilde K}^\star\,{\bf m}=
{\bf m}^{-1}$, one loses the involutive property for $\Sigma_{\bf m}$ (but  one always keeps the
$\tilde G$-equivariance). 
 Also, these operators are bounded on $\CH_\chi$, if and only if the function ${\bf m}$ is essentially 
 bounded. But we will see in Theorem \ref{unitarity} that  in order to obtain a unitary quantization map
 we are forced to consider such unbounded $\Sigma_\bm$'s.
 We mention 
 a simple self-adjointness criterion, interesting on its own. 
 
 \begin{lem}
 \label{SA}
 Let  ${\bf m}$ be a locally essentially bounded Borelian function  on $\tilde G/B$ such that 
 $\underline \sigma_{\tilde K}^\star\, \overline{\bf  m}={\bf m}$. 
  Define 
$$
\tilde\Omega_{\bf m}(g):=U_\chi(g)\,\Sigma_{\bf m}\,U_\chi(g)^*\,,\quad g \in \tilde G\;,
$$ 
on the domain
$$
B_g:=\big\{\varphi\in  \CH_\chi:|{\bf m}_g|\,\varphi\in  \CH_\chi\big\}\quad
 \mbox{where}\quad {\bf m}_g:=[g_0B\mapsto {\bf  m}(g^{-1}g_0B)]\,,\quad g\in \tilde G\;.
 $$
Then, $\tilde\Omega_{\bf m}(g)$ is self-adjoint on $ \CH_\chi$. Moreover,  $\Gamma^\infty_c(E_\chi)$ 
is a common core for all $\tilde\Omega_{\bf m}(g)$'s, $g\in \tilde G$.
  \end{lem}
\begin{proof}
Note first  that the formal adjoint of  $\Sigma_{\bf m}$ is 
$\Sigma_{ \underline \sigma_{\tilde K}^\star\, \overline{\bf  m}}$. Therefore,  when 
$\underline \sigma_{\tilde K}^\star\, \overline{\bf  m}={\bf m}$, 
the operator $\Omega_{\bf m}(x)$ is symmetric on 
$\Gamma^\infty_c(E_\chi)$. Next, we remark that as both $\Sigma$ and $U_\chi(g)$, $g\in \tilde G$, preserve  
$\Gamma^\infty_c(E_\chi)$, we get  for $\vf\in \Gamma^\infty_c(E_\chi)$:
$$
\tilde\Omega_{\bf m}(g)^2\vf=|{\bf m}_g|^2\vf\;,
$$
that is, $\tilde\Omega_{\bf m}(g)$ squares on $\Gamma^\infty_c(E_\chi)$ to a multiplication operator. 
Since $\tilde\Omega_{\bf m}(g)$ is symmetric on  the space of smooth compactly supported sections
 of $E_\chi$, the latter entails that
$$
\|\tilde\Omega_{\bf m}(g)\vf\|=\| |{\bf m}_g|\vf\|\,,\quad\forall \vf\in \Gamma^\infty_c(E_\chi)\;,
$$
and thus $\tilde\Omega_{\bf m}(g)$ is well defined on  $B_g$. Then, the same computation as above, 
shows that  $\tilde\Omega_{\bf m}(g)$ is also symmetric on its domain $B_g$. 
Observe  that $B_g$ is
complete in the graph norm, given by
$ \|\psi\|^2 + \||{\bf m}_g|\psi\|^2$, and that
$ \Gamma^\infty_c(E_\chi)$ is dense in $B_g$ for this norm. Thus $\tilde\Omega_{\bf m}(g)$, with
domain $B_g$, is a closed operator. Clearly
$B_g \subset {\rm dom}\big( \tilde\Omega_{\bf m}(g)^*\big)$, since $\tilde\Omega_{\bf m}(g)$ is 
symmetric on $B_g$.

Choose an increasing sequence 
of relatively compact open sets $\{C_n\}_{n\in\mathbb N}$ in $\tilde G/B$,
 converging to 
$\tilde G/B$.   For $n\in\mathbb N$, let $\chi_n$ be the 
indicator function of $C_n$. Then of course
$\chi_n  \vf \in B_g$ for all $\vf\in\CH_\chi$.  Note also that for $\vf\in B_g$, we have by definition of 
${\bf m}_g
$ and from the relation $\underline \sigma_{\tilde K}^\star\, {\bf  m}=\overline{\bf m}$:
$$
\tilde\Omega_{\bf m}(g)\vf={\bf m}_g \Omega(x)\vf=\Omega(x)\overline{\bf m}_g\vf\,,\quad 
x=g\tilde K\in M\;.
$$
Thus  for $g\in \tilde G$, $x=g\tilde K\in M$,
 $\psi \in {\rm dom}\big( \tilde\Omega_{\bf m}(g)^*\big)$, $\vf\in\CH_\chi$ and using the fact that 
 $\chi_n{\bf m}_g$ is essentially  bounded (i.e$.$  the associated multiplication operator is
 bounded), we get 
 \begin{align*}
 \langle \vf,\chi_n \tilde\Omega_{\bf m}(g)^*\psi\rangle=\langle\tilde \Omega_{\bf m}(g)\chi_n\vf,\psi
 \rangle=
 \langle\Omega(x)\overline{\bf m}_g\chi_n\vf,\psi\rangle=\langle\vf, \chi_n{\bf m}_g\Omega(x)\psi
 \rangle\;.
 \end{align*}
 Using the
monotone convergence theorem, we  obtain
\begin{align*}
\| \tilde\Omega_{\bf m}(g)^*\psi\|&=\lim_{n\to\infty}\| \chi_n\tilde\Omega_{\bf m}(g)^*\psi\|=
\lim_{n\to\infty}\sup_{\|\vf\|=1}\big|\langle \vf,\chi_n \tilde\Omega_{\bf m}(g)^*\psi\rangle\big|
\\
&=\lim_{n\to\infty}\sup_{\|\vf\|=1}\big|\langle\vf, \chi_n{\bf m}_g\Omega(x)\psi\rangle\big|
=
\lim_{n\to\infty}\|\chi_n{\bf m}_g\Omega(x)\psi\|\\&=
\|{\bf m}_g\Omega(x)\psi\|
=\|\Omega(x)\overline{\bf m}_g\psi\|=\|\overline{\bf m}_g\psi\|=\||{\bf m}_g|\psi\|\;,
\end{align*}
so that necessarily $\psi \in B_g$. Thus ${\rm dom}\big( \tilde\Omega_{\bf m}(g)^*\big)= B_g $, 
as required.
Note lastly that $\Gamma_c^\infty(E_\chi)$ being dense in each
$B_g$ for the graph norm, it is a common core for all the $ \tilde\Omega_{\bf m}(g)$, which
are therefore essentially selfadjoint on that domain.
\end{proof}
\begin{rmk}
\label{TRANS=FULL}
At this early stage of the construction, it is important to observe 
 that our representation of $M$ (Proposition \ref{propO})
 and the associated quantization map (Definition \ref{QM}) could have been equally defined
 starting with the full polarization quadruple $(\L,\sigma_\L,\xi,\bB)$  (see Definition 
 \ref{TRANSFULL}) of  a non-exact polarized transvection triple $(\g,\sigma,\varpi)$, 
 instead of the transvection  quadruple $(\tilde{G},\tilde{\sigma},\xi\big|_{\tilde G},B )$.
In particular, all the results of sections \ref{Locality and the one-point phase}, 
\ref{Unitarity and midpoints}, \ref{The star-product} and
\ref{The three-point kernel} can be thought as arising from the full quadruple.
We will make great use of this observation in
section  \ref{Extensions  of polarization quadruples}.
\end{rmk}

\section{Locality and the one-point phase}
\label{Locality and the one-point phase}

We next pass to the notion of locality in the context of transvection quadruples, out of which we will
 be able to give an explicit expression of the operators  $\Omega_\bm(x)$ on $\CH_\chi$.
\begin{dfn}\label{GLOBAL}
Within the notations of Definition \ref{TRANSFULL}, we say that the polarized symplectic 
symmetric space $(M,s,\omega)$,  associated to a non-exact polarized transvection triple 
$(\g,\sigma,\varpi)$,  is {\bf local} whenever there exists a subgroup $Q$ of $\tilde{G}$ 
such that:
\begin{enumerate}
\item[(i)] The map
\begin{equation*}
Q\times B\to\tilde{G}\,,\quad(q,b)\mapsto qb\;.
\end{equation*}
is a  global diffeomorphism. In particular, $Q$ is  closed as a subgroup of $\tilde{G}$.
\item[(ii)] For all $q\in Q$ and $b\in B$, one has
$$
\bC_q(\tilde{\sigma}(b)b^{-1})\;\in\;B\;,
$$
where $\bC_g(g'):=gg'g^{-1}$ denotes the conjugate  action of $\tilde{G}$ on itself. 
\item[(iii)] For every $q\in Q$, setting 
$\tilde{\sigma}q=:\left(\tilde{\sigma}q\right)^Q\,\left(\tilde{\sigma}q\right)^B\,$
relatively to the global decomposition  $\tilde{G}=Q.B$, one has:
$$
\chi\big((\tilde{\sigma}q)^B\big)=1\;.
$$
\end{enumerate}
\end{dfn}

For a local  symplectic symmetric space, the identification $Q\simeq \tilde G/B$
allows to transfer the symmetric space structure of the former to the latter:
\begin{lem}
\label{WhenLocal}
Let $(M,s,\omega)$ be a local  symplectic symmetric space. Then:
\begin{enumerate}
\item[(i)] The mapping:
$$
\underline{s}:Q\times Q\to Q\,,\quad(q,q')\mapsto \underline{s}_q(q')
:=q\left(\tilde{\sigma}(q^{-1}q')\right)^Q\;,
$$
defines a left-invariant structure of symmetric space on the Lie group $Q$. 
\item[(ii)] Moreover, the global diffeomorphism
$$
Q\to\tilde{G}/B\,,\quad q\mapsto qB\;,
$$ 
intertwines  the symmetry $\underline{s}$ with the involution $\underline{\sigma}$,
defined in \eqref{SSACTION} for the  transvection quadruple $(\tilde G,\tilde\sigma,\xi,B)$:
$$
\underline{s}_qq'\;\mapsto\;\underline{\sigma}_{q\tilde K}(q'B)\,,\qquad\forall q,q'\in Q\;.
$$
\item[(iii)] Under the identification $Q\simeq\tilde G /B$ given above,
the $(\tilde G,\underline{\sigma}_{\tilde K})$-invariant measure ${\rm d}_{\tilde G/B}$
on ${\tilde G/B}$ 
constructed in Lemma
\ref{G-M-invariant}, becomes a  $(\tilde G,\underline{s}_e)$-invariant measure ${\rm d}_Q$ 
on the Lie group $Q$,
which  is also  a left-invariant Haar measure on $Q$.
\item[(iv)] Last, we have an  isomorphism of Hilbert spaces
$\CH_\chi\simeq L^2(Q)$ induced  by the $\tilde G$-equivariant isomorphism:
\begin{equation*}
C^\infty(\tilde G)^B\to C^\infty(Q)\,,\quad\hat{\varphi}\mapsto{\varphi}:=\hat{\varphi}\big|_Q\;,
\end{equation*}
under which, we have 
$\Sigma=\underline s_e^\star$.
\end{enumerate}
\end{lem}
\begin{proof} Item (ii) follows from a direct check implying in turn the left-$Q$-equivariance of 
$\underline{s}$ from the left-$\tilde G$-equivariance 
of $\underline{\sigma}$. 
The fact that $\underline{s}_q$ fixes $q$ isolatedly  is a consequence of the following 
observation.
Considering the linear epimorphism $p_{\star \tilde K}:\p\to\q=T_eQ\simeq T_B(\tilde G/B)$ 
tangent to the projection 
$p:M\simeq\tilde G/\tilde K\to \tilde G/B\simeq Q$, $g\tilde K\mapsto gB$, one observes that for every $X\in\p$:
$$
\underline{\sigma}_{\tilde K\star B}(p_{\star \tilde K}X)=(p\circ s_{\tilde K})_{\star \tilde K}(X)
=-p_{\star\tilde K}X\;.
$$
Hence $\underline{\sigma}_{\tilde K\star B}=-\id_{\q}$ and (i) follows.
Last, (iii) and (iv) are immediate consequences of the $Q$-equivariant identification 
$Q\simeq \tilde G/B$.
 \end{proof}

From now on, we will always make the identification $\CH_\chi\simeq L^2(Q)$,
under which we can derive the action
of the individual operators $\Omega(x)$, $x\in M$. For this, we need a preliminary result:
\begin{lem}
\label{RTY}
Let $(M,s,\omega)$ be a local symplectic symmetric space and
 let ${\hat{\varphi}}\in C^\infty(\tilde{G},\C)^B$ be a $B$-equivariant function.
Then, for all $q,q_0\in Q$ and $b\in B$,  one has:
$$
L_{qb}^\star\circ\tilde{\sigma}^\star\circ L^\star_{(qb)^{-1}}\,{\hat{\varphi}}(q_0)=\bE(q_0^{-1}q
b)\,
{\hat{\varphi}}\big(\underline{s}_qq_0\big)\;,
$$
with
\begin{equation}\label{E}
\bE(qb):=\overline{\chi}\left(\bC_q\big(\tilde{\sigma}(b)b^{-1}\big)\right)\;,
\end{equation}
\end{lem}
\begin{proof}
A direct computation yields:
$$
L_{qb}^\star\circ\tilde{\sigma}^\star\circ L_{(qb)^{-1}}^\star\,{\hat{\varphi}}(q_0)=
{\hat{\varphi}}(qb\tilde{\sigma}(b^{-1})\tilde{\sigma}(q^{-1}q_0))=
{\hat{\varphi}}(q\tilde{\sigma}(q^{-1}q_0)\bC_{\tilde{\sigma}(q_0^{-1}q)}(b\tilde{\sigma}(b^{-1})))\;.
$$
Under the assumption of locality (Definition \ref{GLOBAL}), we have $\bC_{\tilde{\sigma}(q_0^{-1}q)}(b\tilde{\sigma}(b^{-1}))\;\in\;B$.
The $\tilde{\sigma}$-invariance of $\chi$ and item (iii) of Definition \ref{GLOBAL} then yield the 
formula.
\end{proof}
\begin{rmk}
We call the function $\bE$ in \eqref{E} the {\em one-point phase}. Observe that the latter is
well defined thanks to the second condition in the
assumption of locality   (Definition \ref{GLOBAL}).
\end{rmk}
\begin{cor}
\label{expr}
Let $(M,s,\omega)$ be a local  symplectic symmetric space.
For  $\vf\in L^2(Q)$ and $x=qb\tilde K\in M$, $q\in Q$, $b\in B$, we have
$$
\Omega(x)\,\vf(q_0)=\bE(q_0^{-1}q
b)\,
{{\varphi}}\big(\underline{s}_q(q_0)\big)\;,
$$
where $\bE$ is the phase defined in \eqref{E} and $\underline s$ is the symmetry of the Lie group
$Q$ constructed in Lemma \ref{WhenLocal}.
\end{cor}

\section{Unitarity and midpoints for elementary spaces}
\label{Unitarity and midpoints}

In addition to locality (Definition \ref{GLOBAL}), we will  assume further 
conditions on the structure
of our polarized symplectic symmetric space $(M,s,\omega)$,
 which will enable us to give an explicit 
expression of the three-point
kernel associated to a WKB-quantization of $M$ as well as to prove the triviality of the 
associated Berezin transform  (see Definition 
\ref{Berezin} below). Recall that the notion of  midpoint map on a symmetric
space is given in Definition \ref{def-midpoint}.

\begin{dfn}\label{ELEMENTARY}
A local symplectic symmetric space $(M,s,\omega)$ is called {\bf elementary} when, 
within the context of the section \ref{KIRILLOV}, the following additional conditions
are satisfied:
\begin{enumerate}
\item[(i)] The symmetric space $(Q,\underline s)$ is solvable and  admits a 
(necessary unique) midpoint
map. (By a result  proved in  \cite{YVthese}, this implies that $Q$ is exponential.)
\item[(ii)] There exists an exponential Lie subgroup $\Y$ of $B$ normalized by 
$Q$ and such that the semi-direct product 
$$
\S:=Q\ltimes\Y\;\subset\;\tilde{G}\;,
$$
acts simply transitively on $M$.
\item[(iii)] Denoting by $\fY$ the Lie algebra of $\Y$, there exists a global
diffeomorphism $\Psi:Q\to\q$ such that
$$
\langle\xi\,,\,(\Ad_{q^{-1}}-\Ad_{(\underline{s}_eq)^{-1}})y\rangle
=\langle\xi\,,\,[\Psi(q)\,,\,y]\rangle\,,\quad\forall\,y\in\fY\,,\forall q\in Q\;.
$$
\item[(iv)] The Lie algebras $\fY$ and $\q$ are 
Lagrangian subspaces of the Lie algebra $\s$ of $\S$ that are in symplectic duality with respect to 
the
evaluation at the unit element $e$ of $\S$ of the symplectic structure transported from $M$ to $\S$ via the diffeomorphism $\S\to M=\tilde{G}/\tilde{K}$, $x\mapsto x\tilde{K}$. 
\end{enumerate}
\end{dfn}

\begin{rmk}
\label{Rc1}
Let $(M,s,\omega)$ be an
elementary symplectic symmetric space.
\begin{enumerate}
\item[(i)]  From now on, we always make the $\S$-equivariant
identification:
\begin{equation*}
\S=Q\ltimes\Y\to M\,,\qquad qb\mapsto qb\tilde K\;.
\end{equation*}
 Observe that under this identification, 
the $\tilde G$-invariant Liouville measure ${\rm d}_M$ on $M$ is a left Haar measure on $\S$,
which under the parametrization $g=qb$, $q\in Q$, $b\in\Y$, is 
proportional to any product of left invariant Haar
measures on $Q$ and on $\Y$.
We simply denote the latter by ${\rm d}_\S$.
\item[(ii)] For ${\bf m}$ a locally essentially bounded Borelian function  on $Q$ and
$x\in\S$ , the operator 
$\tilde \Omega_\bm(x)$ given in Lemma \ref{SA} will  be simply
 denoted by $ \Omega_\bm(x)$. This is coherent with the identification above and with 
 the notation of Proposition \ref{propO} when $\bm=1$. Moreover,
 the family $\{ \Omega_\bm(x)\}_{x\in\S}$ satisfies the first axiom, \eqref{c1}, of a covariant Moyal
 quantizer, as given at the very beginning of this chapter.
 \item[(iii)] Since $Q$ normalizes $\Y$, the restriction 
to $\S=Q\ltimes \Y$ of the representation $U_{\chi}$  of $\tilde G$ on 
$\CH_\chi\simeq L^2(Q)$ given in \eqref{U},  reads:
$$
U_\chi(qb)\vf(q_0)=\chi\big(\bC_{q_0^{-1}q}(b)\big)\,\vf(q^{-1}q_0)\;.
$$
\end{enumerate}
\end{rmk}

In the elementary case,
we  observe the following relation between the one-point phase $\bE$
and the diffeomorphism $\Psi:Q\to\q$:
\begin{lem}
\label{E-Psi}
Let  $(M,s,\omega)$ be an elementary symplectic symmetric space. Then, for $q\in Q$ and $y\in\fY$, we have:
$$
\bE(q^{-1}e^y)=\exp\{i\langle\xi,[\Psi(q),y]\rangle\}\;.
$$
\end{lem}
\begin{proof}
By definition, we have for $q\in Q$ and $b\in \Y$:
$$
\bE(q^{-1}b)=\overline\chi\big(\tilde\sigma(\bC_{\tilde\sigma(q^{-1})}(b))\,\bC_{q^{-1}}(b^{-1})\big)\;.
$$
Next, we write
$$
\tilde\sigma(q^{-1})=\tilde\sigma(q)^{-1}=\big(\tilde\sigma(q)^Q\;\tilde\sigma(q)^B\big)^{-1}
=\big(\tilde\sigma(q)^B\big)^{-1}\,\big(\underline s_eq\big)^{-1}\;,
$$
to get
$$
\bC_{q^{-1}}(b^{-1})\,\tilde\sigma\big(\bC_{\tilde\sigma(q^{-1})}(b)\big)
=\big(\tilde\sigma\big(\tilde\sigma(q)^{B}\big)^{-1}\big)\;
\tilde\sigma\big(\bC_{(\underline s_eq)^{-1}}(b)\big)\;\tilde\sigma\big(\tilde\sigma(q)^B\big)
\;\big(\bC_{q^{-1}}(b^{-1})\big)\;.
$$
Since $Q$ normalizes $\Y$ and $B$ is $\tilde\sigma$-stable, we  observe 
that each of the four factors in the right hand side above, belong to $B$. Thus, we can split
$\bE(q^{-1}b)$ according to this decomposition, to get
$$
\bE(q^{-1}b)=\chi\big(\bC_{q^{-1}}(b)\big)\;
\overline\chi\big(\bC_{(\underline s_eq)^{-1}}(b)\big)\;.
$$ 
The result follows from the definition of the diffeomorphism $\Psi$ and the character $\chi$.
\end{proof}

Now, using Corollary \ref{expr}, Definition \ref{Sigma-m} and Lemma \ref{SA}, under the identification $\S\simeq M$, we note
\begin{prop}
\label{expression}
Let $(M,s,\omega)$ be an elementary symplectic symmetric space. 
Let $q\in Q$, $b\in\Y$ and let
$\bm$ be an essentially locally bounded Borelian  function on
$ Q$. Then the densely defined (on 
$B_{qb}=B_q$--see Definition \ref{Sigma-m} and Lemma \ref{SA}) operator 
$\Omega_\bm(qb)$, acts  as:
$$
{\Omega}_\bm(qb){\varphi}(q_0)=\bm(q^{-1}q_0)\, \bE(q_0^{-1}qb)
\,{\varphi}(\underline{s}_{q}q_0)\;,
\quad \forall \vf\in B_q\subset \CH_\chi\;,\quad\forall q_0\in Q\;.
$$
\end{prop}
\begin{cor}
\label{kernel-op}
Let $(M,s,\omega)$ be an elementary symplectic symmetric space,
 $\bm$ be an essentially locally bounded  Borelian  function on
$Q$ and  $f\in\CD(\S)$. Then the operator $\Omega_\bm(f)$ defined by
\begin{align*}
\Omega_\bm(f): \CD(Q)&\to \CD'(Q)\;,\\
\vf\mapsto \Omega_\bm(f)\vf:=\Big[\psi\in\CD(Q)&\mapsto
\int_{Q\times \S}\psi(q_0)\,f(qb)\,\big({\Omega}_\bm(qb){\varphi}\big)(q_0)\,{\rm d}_\S(qb)\,
{\rm d}_Q(q_0)\Big]\;,
\end{align*}
has a distributional kernel given by
\begin{align*}
&\Omega_{\bf m}(f)[q_0,q]=\\&
{\bf m}\big(\mid(e, q_0^{-1} q)^{-1}\big)\,
\big|{\rm Jac}_{(\underline s^e)^{-1}}\big|(q_0^{-1}q)\,
\int_\Y\, f(\mid(q_0, q)b)
\, \bE\big(\mid(e, q_0^{-1} q)b\big)\,{\rm d}_\Y(b)\;.
\end{align*}
\end{cor}
\begin{proof}
Observe first that under  the decomposition $\S=Q\ltimes \Y$,  the left Haar measure  
on $\S$ coincides with the product of left Haar measures  on $Q$ and $\Y$:
$$
{\rm d}_\S(qb)={\rm d_Q}(q)\,{\rm d_\Y}(b)\,,\quad\forall q\in Q\,,\;\forall b\in \Y\;.
$$
For $f\in\CD(\S)$  and any Borelian ${\bf m}$, it is clear that $\Omega_{\bf m}(f)$ defines a 
continuous operator from $\CD(Q)$ to $\CD'(Q)$ and acts as:
\begin{align*}
\Omega_{\bf m}(f)\vf(q_0)
&=\int_{Q\ltimes\Y} \, f(qb)\,{\bf m}(q^{-1}q_0)\, \bE(q_0^{-1}qb)\,\vf\big(\underline s_q(q_0)\big)\, 
{\rm d}_Q(q)\,{\rm d}_\Y(b)\,,\quad \vf\in \CD(Q)\;.
\end{align*}
For any $q_0\in Q$, we set
$q'(q):=\underline{s}_{q}(q_0)$ and
we get from the defining property of the midpoint map that
$q=\mid(q_0, q')$.
 Now observe that left-translations (in the group $Q$)  are automorphisms of  
the symmetric space $(Q,\underline s)$. Indeed, for all $q_0,q,q'$ in $Q$, we have
$$
L_{q_0}(\underline{s}_q(L_{q_0^{-1}}q'))=q_0q\left(\tilde{\sigma}(q^{-1}q_0^{-1}q')\right)^Q=\underline{s}_{q_0q}(q')\;.
$$
Hence, by Remark \ref{MDPT}, we get
$$
\mid(q_0, q')=q_0\,\mid(e,q_0^{-1}q')=L_{q_0}\circ (\underline s^e)^{-1}\circ L_{q_0^{-1}}(q')\;,
$$
the invariance of the Haar measure ${\rm d}_Q$ under left translation gives:
$$
\big|{\rm Jac}_{\mid(q_0,.)}\big|(q')
=\big|{\rm Jac}_{(\underline s^e)^{-1}}\big|(q_0^{-1}q')\;.
$$
Therefore, a direct computation shows that:
\begin{align*}
\Omega_{\bf m}(f)\vf(q_0)=
\int_{Q\ltimes\Y} f(\mid(q_0, q)b)\,
&{\bf m}(\mid(e, q_0^{-1} q)^{-1})\,\big|{\rm Jac}_{(\underline s^e)^{-1}}\big|(q_0^{-1}q)
\\
&\times\bE\big(\mid(e, q_0^{-1} q)b\big)\,\vf(q)\, {\rm d}_Q(q)\,{\rm d}_\Y(b)\;,
\end{align*}
and the result follows by identification.

\end{proof}
\begin{rmk}
\label{verifc2}
 As a consequence of the preceding corollary, we deduce  that for 
an elementary symplectic symmetric space, the second axiom, \eqref{c2},
of a covariant Moyal quantizer is satisfied when 
$\bm(e)|{\rm Jac}_{(\underline s^e)^{-1}}|(e)=1$. Indeed,
 for $f\in\CD(\S)$, we deduce  that
 \begin{align*}
 {\rm Tr}\big[\Omega_\bm(f)\big]&=
 \int_Q \Omega_\bm(f)[q,q]\,{\rm d}_Q(q)\\&=
 \bm(e)\big|{\rm Jac}_{(\underline s^e)^{-1}}\big|(e)\int_{Q\times\Y} f(qb)\,\bE(b)\,{\rm d}_\Y(b)\,
 {\rm d}_Q(q)=\int_{\S} f(x)\,{\rm d}_\S(x),
 \end{align*}
 since $\bE(b)=1$ for all $\b\in\Y$ by Definition \ref{ELEMENTARY} (iii) and Lemma \ref{E-Psi}.
\end{rmk}

Our next aim is  to understand the geometrical conditions  on the functional parameter
$\bm$ necessary for the quantization map to extend to a unitary operator from $L^2(M)$
to the Hilbert space of Hilbert-Schmidt operators on $\CH_\chi$. 
For this, we introduce the following
specific function
on $Q$:
\begin{dfn}
\label{trucs}
For $(M,s,\omega)$  an elementary symplectic symmetric space,
define the function $\bm_0$ on $Q$ as:
\begin{equation}
\label{m0}
{\bf m_0}(q):=\left|\mbox{\rm Jac}_{\underline s^e}(q^{-1})
\,\mbox{\rm Jac}_{\Psi}(q)\right|^{1/2}\;.
\end{equation}
\end{dfn}

\begin{rmk}
\label{unitary-not-SA}
 Observe that both  $|\mbox{\rm Jac}_{\Psi}|$ and  $|\mbox{\rm Jac}_{\underline s^e}|$
are  $\underline s_e^\star$-invariant. 
Indeed,  we have $\Psi\circ\underline s_e =-\Psi$  and 
$\underline s^e\circ \underline s_e=\underline s_e\circ \underline s^e$, hence the claim follows from 
$|\mbox{\rm Jac}_{\underline s_e}|=1$  (cf$.$ Lemma \ref{WhenLocal} (iii)).
However, $\bm_0$ need not be $\underline s_e^\star$-invariant, as 
$|\mbox{\rm Jac}_{\underline s^e}|$ need not be invariant under the inversion map on $Q$. 
Thus, from Lemma
\ref{SA}, the  quantization map  $\Omega_{\bm_0}$, need not send real functions to self-adjoint
operators.
\end{rmk}

We can now state one of the main results of this chapter:

\begin{thm}
\label{unitarity}
Let $(M,s,\omega)$ be an elementary symplectic symmetric space 
and $\bm$ be a Borelian  function on
$Q$ which is (almost everywhere) dominated by ${\bf m_0}$ and assume that $\Y$
is Abelian. 
 Then the quantization map: 
 $$
 \Omega_{\bf m}:f\mapsto \Omega_{\bf m}(f):=\int_\S f(qb)\,\Omega_\bm(qb)\,{\rm d}_\S(qb)\;,
 $$
  is a bounded operator from $L^2(\S)$ to ${\mathcal L}^2(\CH_\chi)$ with
  $$
  \|\Omega_\bm\|\leq \| \bm/{\bm_0}\|_\infty\;.
  $$
  Moreover, $\Omega_\bm$  is a unitary operator  if and only if $|{\bf m}|={\bf m_0}$. 
  \end{thm}
\begin{proof}
Recall that a linear operator 
$T:\CD(Q)\to \CD'(Q)$ extends to a Hilbert-Schmidt operator on $L^2(Q)$,
 if and only  its distributional kernel belongs to
 $L^2(Q\times Q)$. 
 In this case,  its Hilbert-Schmidt norm coincides with the $L^2$-norm of its kernel.
Thus by Corollary \ref{kernel-op}, we deduce that if $f\in\CD(\S)$, then
 the square of the Hilbert-Schmidt norm of $\Omega_{\bf m}(f)$ reads
\begin{align*}
&\|\Omega_{\bf m}(f)\|^2_2
=\\
&\int_{Q^2\times\Y^2} |{\bf m}|^2(\mid(e, q_0^{-1} q)^{-1})\, \big|{\rm Jac}_{(\underline s^e)^{-1}}\big|^2(q_0^{-1}q)
\, 
\overline f(\mid(q_0, q)b)\, f(\mid(q_0, q)b_0)\,
\\
&\qquad\qquad\qquad\quad
\times\overline \bE\big(\mid(e, q_0^{-1} q)b\big)\,\bE\big(\mid(e, q_0^{-1} q)b_0\big)\,
{\rm d}_Q(q)\,{\rm d}_\Y(b) \,{\rm d}_Q(q_0)\,
{\rm d}_\Y(b_0)\;.
\end{align*}
Performing the change of variables $\mid(q_0, q)\mapsto q$ (the inverse of the one we performed
in the proof of Corollary \ref{kernel-op}) and using the relation between the function $\bE$
and the diffeomorphism $\Psi:Q\to q$ given in Lemma \ref{E-Psi}, we get the following expression
for $\|\Omega_{\bf m}(f)\|^2_2$:
\begin{align*}
&
\int_{Q^2\times\Y^2} |{\bf m}|^2(q^{-1}q_0)\,\big|{\rm Jac}_{(\underline s^e)^{-1}}\big|
\big(\underline s^e(q_0^{-1}q)\big)\,
 \overline f(qb)\, f(qb_0)\,
  e^{i\langle\xi,[\Psi(q^{-1} q_0),\log(b)-\log(b_0)]\rangle}\\
  &\qquad\qquad\qquad\qquad\qquad\qquad\qquad\qquad\quad\qquad\qquad\times
 {\rm d}_Q(q)\,{\rm d}_\Y(b) \,{\rm d}_Q(q_0)\,
{\rm d}_\Y(b_0)\;,
\end{align*}
and the latter can be rewritten as
$$
\int_{Q^2\times\Y^2} \frac{|{\bf m}|^2(q_0)}{\big|{\rm Jac}_{\underline s^e}\big|
\big(q_0^{-1}\big)}
 \overline f(qb) f(qb_0)
  e^{i\langle\xi,[\Psi( q_0),\log(b)-\log(b_0)]\rangle}
 {\rm d}_Q(q){\rm d}_\Y(b) {\rm d}_Q(q_0)\,
{\rm d}_\Y(b_0)\;,
$$
where in the last line, we used left-invariance of the Haar measure on $Q$. Setting  
$w_0=\Psi(q_0)\in\q$ and ${\rm d}w_0$ the Lebesgue 
measure on $\q$, we get
\begin{align*}
\|\Omega_{\bf m}(f)\|^2_2&=
\int_{Q\times\q\times\Y^2} \frac{|{\bf m}|^2(\Psi^{-1}(w_0)) \,\overline f(qb)\, f(qb_0)}{\big|{\rm Jac}_{\underline s^e}\big|
\big(\Psi^{-1}(w_0)^{-1}\big)\big|{\rm Jac}_{\Psi}\big|
\big(\Psi^{-1}(w_0)\big)}\,
  e^{i\langle\xi,[w_0,\log(b)-\log(b_0)]\rangle}\\
   &\qquad\qquad\qquad\qquad\qquad\qquad\qquad\qquad\qquad\quad\times
  {\rm d}_Q(q)\,{\rm d}_\Y(b) \,{\rm d}w_0\,
{\rm d}_\Y(b_0)\\
&=
\int_{Q\times\q\times\Y^2} \frac{|{\bf m}|^2(\Psi^{-1}(w_0))}{{\bf m_0}^2(\Psi^{-1}(w_0))}
\, \overline f(qb)\, f(qb_0)\,
e^{-i\langle\xi,[w_0,\log(b)-\log(b_0)]\rangle}\\
 &\qquad\qquad\qquad\qquad\qquad\qquad\qquad\qquad\qquad\quad\times {\rm d}_Q(q)\,{\rm d}_\Y(b) \,{\rm d}w_0\,
{\rm d}_\Y(b_0)\\
&=
\int_{Q\times\q} \frac{|{\bf m}|^2(\Psi^{-1}(w_0))}{{\bf m_0}^2(\Psi^{-1}(w_0))}
\, \left|\int_\Y f(qb)\, \,
e^{-i\langle\xi,[w_0,\log(b)]\rangle}\,{\rm d}_\Y(b)\right|^2\,{\rm d}_Q(q) \,
{\rm d}w_0\;.
\end{align*}
We will next use the relation (which follows from the construction of $\xi\in\tilde\g^\star$):
$$
\langle\xi,[X,Y]\rangle=\varpi(X,Y)\;,\quad\forall X,Y\in\tilde \g\;,
$$
and the fact that  $\q$ and $\fY$ are Lagrangian subspaces  in symplectic duality 
 (see Definition \ref{ELEMENTARY}).   Now, since
$\Y$ is Abelian and exponential, the exponential map $\fY\to\Y:y\mapsto\exp(y)$ is a measure preserving diffeomorphism 
(the vector space $\fY$ being endowed with a normalized Lebesgue measure ${\rm d}y$) inducing the isometric linear identification:
$$
\exp^\star: L^2(\Y)\to L^2(\fY)\;.
$$
Within this set-up, the (isometrical) Fourier transform reads:
$$
\CF:L^2(\Y)\to L^2(\q)\,,\quad\varphi\mapsto\Big[w_0\mapsto\int_\Y\varphi(b)e^{-i\langle\xi,[w_0,\log(b)]
\rangle}\,{\rm d}_\Y(b)\Big]\;,
$$  
where, again, the vector space $\q$ is endowed with a normalized Lebesgue measure ${\rm d}w_0$.
We therefore observe that
$$
\|\Omega_{\bf m}(f)\|^2_2=\int_{Q\times\q} \frac{|{\bf m}|^2(\Psi^{-1}(w_0))}{{\bf m_0}^2(\Psi^{-1}(w_0))}
\, \left|\CF(L^\star_qf)(w_0)\right|^2\,{\rm d}_Q(q) \,
{\rm d}w_0\;,
$$
where we set $L^\star_qf:b\mapsto f(qb)$.
Hence, we see that if  
$|{\bf m}|\leq C\, {\bf m_0}$, we have:
\begin{align*}
\|\Omega_{\bf m}(f)\|^2_2&\leq C^2\,\int_{Q\times\q}
\, \left|\CF(L^\star_qf)(w_0)\right|^2\,{\rm d}w_0\,{\rm d}_Q(q)\\&\quad=C^2\,\int_{Q\times\Y}
\, \left|\,L^\star_qf(b)\right|^2\,{\rm d}_\Y(b)\,{\rm d}_Q(q)\\
&\quad=C^2\|f\|_2^2\,,\quad\forall f\in\CD(\S)\;.
\end{align*}
i.e$.$
$$
\|\Omega_{\bf m}(f)\|_2\;\leq C\;\|f\|_2\,,\quad\forall f\in\CD(\S)\;.
$$
By density,  we deduce that $\Omega_{\bf m}(f)$
is Hilbert-Schmidt for all $f\in L^2(\S)$ and with equality of norms 
if and only if $|{\bf m}|= {\bf m_0}$. In this case, a  similar computation shows that 
for any $f_1,f_2\in L^2(\S)$, we have
\begin{equation*}
\Tr\big[\Omega_{\bf m}(f_1)^*\,\Omega_{\bf m}(f_2)\big]=\int_\S \,\overline f_1(qb)\,f_2(qb)\,
{\rm d}_\S(qb)\;,
\end{equation*}
which terminates the proof.
\end{proof}

\begin{rmk}
\label{wht}
 Let $(M,s,\omega)$ be an elementary symplectic symmetric space 
and $\bm$ be an essentially bounded function on $Q$. By Lemma
\ref{SA}, we know that when $\underline s_e^\star\, \overline{\bf  m}={\bf m}$, then the family
$\{\Omega_{\bm_0}(x)\}_{x\in \S}$ consists of self-adjoint operators.
By  Remarks \ref{Rc1} (ii) and \ref{verifc2} we also know (under a mild   normalization condition)
that the first two axioms, \eqref{c1} and \eqref{c2},  of a covariant Moyal quantizer (as defined at the
beginning of this chapter) are satisfied. 
Now, observe that the content of the previous Lemma can be summarized as follow:
$$
\Tr\big[\Omega_{\underline s_e^\star\overline {\bf m}} (x)\,\Omega_{\bf m} (y)\big]=\delta_x(y)\quad
\Longleftrightarrow\quad |{\bf m}|={\bf m_0}\;.
$$
Hence,  when $\bm_0$ is $\underline s_e$-invariant, then the family 
$\{\Omega_{\bm_0}(x)\}_{x\in \S}$ satisfies also the third axiom
 \eqref{c3}.
\end{rmk}

\section{The $\star$-product as the composition law of symbols}
\label{The star-product}

\begin{dfn}
Let $(M,s,\omega)$ be an elementary symplectic symmetric space such that
 $\bm/\bm_0\in L^\infty(Q)$ and such that $\Y$ is Abelian.  Then let 
\begin{align*}
\sigma_\bm:{\mathcal L}^2(\CH_\chi)\to L^2(\S)\;,
\end{align*}
be  the adjoint of the quantization 
map $\Omega_\bm$. We call  the latter the {\bf symbol map}.
\label{Symbol-Map}
\end{dfn}

 Recall that the defining property of the symbol map is
$$
\langle f, \sigma_\bm[A]\rangle=\Tr\big[\Omega_\bm(f)^*\,A\big]\,,\qquad
 \forall A\in{\mathcal L}^2(\CH_\chi)\,,\quad\forall f\in L^2(\S)\;.
$$
Hence, the symbol map   is formally given  by 
\begin{align}
\label{symbol-map}
\sigma_{\bf m}[A](x)=\Tr\big[A\,\,\Omega_{\underline s_e^\star\,\bf \overline  m}(x)\big]\,,
\quad x\in\S\;.
\end{align}
Here again, the trace on the right hand side is understood in the distributional sense on
 $\CD(\S)$. Note however that when $\bm$ is essentially bounded, this expression for the
 symbol map genuinely holds on ${\mathcal L}^1(\CH_\chi)$, the ideal of trace-class 
 operators on $\CH_\chi$. When $\bm$ is only locally essentially bounded, it also holds
 rigorously on the dense subspace of  ${\mathcal L}^2(\CH_\chi)$, consisting of finite linear
 combinations of rank one operators $|\vf\rangle\langle\psi|$ with $\psi\in\CD(Q)$
 and $\vf\in \CH_\chi$ arbitrary.

\begin{dfn}
\label{Berezin}
Let $(M,s,\omega)$ be an elementary symplectic symmetric space.
Assuming  that $\bm/\bm_0\in L^\infty(Q)$, we then set
$$
B_\bm:L^2(\S)\to L^2(\S)\;,\quad f\mapsto \sigma_\bm\circ\Omega_\bm(f)\;,
$$
and we call this linear operator the {\bf Berezin transform} of the quantization map 
$\Omega_{\bm}$.
\end{dfn}
\begin{rmk}
The Berezin transform measures the obstruction for the symbol map to be inverse
of the quantization map. Said differently, the unitarity of the quantization map on $L^2(\S)$
is equivalent to the triviality of the associated Berezin transform. 
\end{rmk}
\begin{prop}
\label{Berezin-prop}
Let $(M,s,\omega)$ be an elementary symplectic symmetric space  with
$\Y$ Abelian
and let $\bm$ a Borelian function on $Q$ such that
$\bm/\bm_0\in L^\infty(Q)$ and such that $\Y$ is Abelian. Then 
\begin{enumerate}
\item[(i)] The Berezin transform is a positive and bounded operator on $L^2(\S)$,
with
$$
\|B_\bm\|\leq  \|\bm/\bm_0\|_\infty^2\;.
$$
\item[(ii)] The Berezin transform is a kernel operator with distributional kernel given
by 
$$
B_\bm[x_1,x_2]=\Tr\big[\Omega_{\underline s_e^\star\,\overline \bm}(x_1)\,\Omega_\bm(x_2)\big]\;,
$$
and the latter can be identified with
$$
B_\bm[x_1,x_2]=
\delta_{q_1}(q_2)\times
\int_{\q}  \frac{{|\bf m|^2}}{{\bf m_0^2}}\big(\Psi^{-1}(w)\big)
e^{i\langle\xi,[w,\log b_1-\log b_2]\rangle} \,{\rm d}w\;,
$$
where $x_j=q_jb_j\in\S$,  $j=1,2$.
\end{enumerate}
\end{prop}
\begin{proof}
 Note that by construction $B_\bm=\Omega_\bm^*\circ\Omega_\bm$, yielding positivity and boundedness
from boundedness of $\Omega_\bm$. The operator norm  estimate comes  from those 
 of Theorem \ref{unitarity}. The second claim
comes from the computation done in the proof of this Theorem.
\end{proof}

For $\bm=\bm_0$, the symbol map $\sigma_{\bm}$ is the inverse of 
the quantization map $\Omega_\bm$. In particular,
the associated Berezin transform is trivial.
 A $\tilde G$-equivariant associative product  $\star_{\bm_0}$ on $L^2(\S)$ 
 is then defined:
 \begin{align}
 \label{prod-op-ptv}
f_1\star_{\bm_0}f_2:=\sigma_{\bm_0}\big[\Omega_{\bm_0}(f_1)\,\Omega_{\bm_0}(f_2)\big]
\,,\quad\forall f_1,f_2\in L^2(\S)\;.
\end{align}
We deduce from \eqref{symbol-map} that the
product \eqref{prod-op-ptv} is a three-point kernel product, 
with distributional kernel given by the 
operator trace of a product of three $\Omega$'s:
\begin{align*}
f_1\star_{\bm_0}f_2(x)=\int_{\S\times\S}f_1(y)\,f_2(z)\,\Tr\big[\Omega_{\underline s_e^\star\bm_0}(x)\,
\Omega_{\bm_0}(y)
\,\Omega_{\bm_0}(z)\big]\,{\rm d}_\S(y)\,\,{\rm d}_\S(z)\;.
\end{align*}
We will return to the explicit form of the three-point kernel and its geometric interpretation in the next 
section. 

 We now come to an important point. Putting together  Remark \ref{unitary-not-SA}
and Theorem \ref{unitarity}, we see that in general the quantization map $\Omega_{\bm}$
need not be unitary and involution preserving (the complex conjugation on $L^2(\S)$
and the adjoint on ${\mathcal L}^2(\CH_\chi)$) at the same time. However, in most
cases (e.g$.$ for elementary normal $\bf j$-groups), the function $\bm_0$ is 
$\underline s_e$-invariant, which implies that the complex conjugation is an involution of 
the Hilbert algebra $\big(L^2(\S),\star_{\bm_0}\big)$:
\begin{prop}
Let $(M,s,\omega)$ be an elementary symplectic symmetric space with $\Y$ Abelian.
Assuming further that $\underline s_e^\star\bm_0=\bm_0$, then for all $f_1,f_2\in L^2(\S)$,
we have:
$$
\overline{f_1\star_{\bm_0} f_2}=\overline{f_2}\star_{\bm_0}\overline{f_1}\;.
$$
\end{prop}

Next, we pass to a possible approach\footnote{This is however not the approach we will follow for
the symmetric spaces underlying elementary normal $\bf j$-groups--see Proposition
 \ref{Relation-Weyl}.}
 to define a $\star$-product for the quantization
map $\Omega_\bm$ in the more general context of an arbitrary function $\bm$ (and without
the assumption that $\Y$ is Abelian).

\begin{dfn}
Let $(M,s,\omega)$ be a polarized, local and elementary symplectic symmetric space
and
fix $\bm$ a  locally essentially bounded Borelian function  on $Q$. We then let
$L^2_\bm(\S)$, be the Hilbert-space of classes of measurable functions on $\S$
for which the norm underlying the following scalar product is finite\footnote{We do not exclude the
possibility that  $L^2_\bm(\S)$ be trivial.}:
\begin{align*}
\langle f_1,f_2\rangle_{\bf m}:=\int_{Q\times\q}  \frac{|{\bf m}|^2}{{\bf m_0}^2}\big(\Psi^{-1}(w)\big)
\,&\Big(\int_\Y  \overline f_1(qb_1)e^{-i\langle\xi,[w,\log b_1]\rangle}\, {\rm d}_\Y(b_1) \Big)\\
&\times\Big(\int_\Y f_2(qb_2)\,e^{i\langle\xi,[w,\log b_2]\rangle}{\rm d}_\Y(b_2) \Big)
\,{\rm d}_Q(q) \,{\rm d}w\;.
\end{align*}
\end{dfn}
\begin{rmk}
Formally, we have 
$$
\langle f_1,f_2\rangle_{\bf m}=\langle f_1,B_\bm f_2\rangle=\Tr\big[
\Omega_\bm(f_1)^*\Omega_\bm(f_2)\big]\;,
$$
where  $\langle.,.\rangle$ denotes the inner product of $L^2(\S)$ and $\Tr$ is the operator trace 
on $\CH_\chi$.
\end{rmk}

Repeating the computations done in the proof of Theorem \ref{unitarity}, we deduce following
extension of the latter:
\begin{prop}
\label{generic-prod}
Let $(M,s,\omega)$ be an elementary symplectic symmetric space
and
let $\bm$ be a locally essentially bounded Borelian function  on $Q$. 
\begin{enumerate}
\item[(i)] The quantization map $\Omega_\bm$ is a unitary operator from $L^2_\bm(\S,{\rm d}_\S)$
to ${\mathcal L}^2(\CH_\chi)$,
\item[(ii)] Associated to the quantization map $\Omega_\bm$, there is a deformed  product 
$\star_\bm$ on $L^2_\bm(\S,{\rm d}_\S)$, which is formally given by:
$$
f_1\star_\bm f_2= B_{\bm}^{-1}\circ\sigma_\bm\big[\Omega_\bm(f_1)\Omega_\bm(f_2)\big]\;,
$$
\item[(iv)] 
 $\big(L^2_\bm(\S),\star_\bm\big)$ is  a Hilbert algebra
 and the complex conjugation is an involution when $\underline s_e^\star \overline \bm=\bm$.
\end{enumerate}
\end{prop}

\section{The three-point kernel} 
\label{The three-point kernel} 
The aim of this section is to compute the distributional three-point kernel 
$\Tr\big[\Omega_{\bm_0}(x)\,\Omega_{\bm_0}(y)\,\Omega_{\bm_0}(z)\big] $
of the product $\star_{\bm_0}$ given in \eqref{prod-op-ptv}. We start with two preliminary results
extracted from  \cite{YVthese}:

\begin{thm}
\label{3points}
Let $(M,s,\omega)$ be an elementary symplectic symmetric space.
Given three points $q_0,q_1,q_2$ in $Q$, the equation 
$$
\underline{s}_{q_2}\underline{s}_{q_1}\underline{s}_{q_0}(q)=q\;,
$$
admits a unique solution $q\equiv q(q_0,q_1,q_2)\in Q$.  In particular, this yields a well-defined
map
$$
Q^3\to Q^3:(q_0,q_1,q_2)\mapsto(q,\underline{s}_{q_0}(q),\underline{s}_{q_1}\underline{s}_{q_0}(q))\;.
$$
The latter is a global diffeomorphism called 
the {\bf medial triangle map} whose inverse is given by:
$$
\Phi_Q:Q^3\to Q^3\,,\quad (q_0,q_1,q_2)\mapsto \big(\mid(q_0,q_1),\mid(q_1,q_2),
\mid(q_2,q_0)\big)\;.
$$
\end{thm}

\begin{figure}[htb]
\centering
\begin{tikzpicture}[scale=2.4]
\coordinate (O) at (0,0) ; 
\coordinate (A) at (90:1cm) ;
\coordinate (B) at (-150:1cm) ;
\coordinate (C) at (-20:1cm) ; 
\coordinate (K) at (-85:0.225cm) ; 
\coordinate (L) at (35:0.263cm) ;
\coordinate (M) at (150:0.268cm) ; 
\coordinate (P) at (-85:2.4cm) ; 
\coordinate (Q) at (35:1.5cm) ;
\coordinate (R) at (150:2cm) ; 
\draw (A) node[above=1pt] {$x_1$} arc(0:-60:1.732cm) ; 
\draw (B) node[below left] {$x_2$} arc(119.6:70.4:2.175cm) ; 
\draw (C) node[below right] {$x_3$} arc(-103.5:-186.5:1.237cm) ; 
\draw (K) node[below=2pt] {$\mid(x_2,x_3)$} ;
\draw (L) node[above right] {$\mid(x_3,x_1)$} ;
\draw (M) node[above left] {$\mid(x_1,x_2)$} ; 
\foreach \pt in {A,B,C,K,L,M}  \draw (\pt) node {$\mathord{\scriptstyle\bullet}$} ;
\end{tikzpicture}
\caption{The medial triangle of three points}
\label{fg:medial-triangle}
\end{figure}
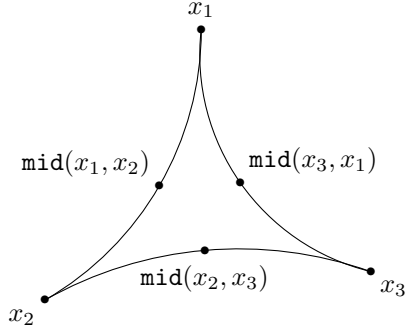

\begin{lem}\label{3Pinter} Let us consider the map
$$
\nu:Q^3\to Q^3:(q_0,q_1,q_2)\mapsto(\mid(e, q_0^{-1} q_1)\,,\,\mid(e, q_1^{-1} q_2)\,,\,\mid(e, q_2^{-1} q_0))\;.
$$
Then
$$
\nu\,\circ\,\Phi^{-1}_Q (q_0,q_1,q_2)=(\,q^{-1}q_0\,,\,(\underline s_{q_0} q)^{-1}q_1
\,,\,(\underline s_{q_1}\underline s_{q_0} q)^{-1}q_2\,)\;,
$$
where $q$ is the unique fixed point of $\underline{s}_{q_2}\underline{s}_{q_1}\underline{s}_{q_0}$.
\end{lem}
\begin{proof} We have $\mid(e,q_0^{-1}q_1)=q_0^{-1}\mid(q_0,q_1)$ and similarly for the other points.
Therefore setting 
$$
(q_0',q_1',q_2'):= \big(\mid(q_0,q_1),\mid(q_1,q_2),
\mid(q_2,q_0)\big)
=\Phi_Q(q_0,q_1,q_2)\;,
$$
we get
$$
\nu(q_0,q_1,q_2)=(q_0^{-1}q_0',q_1^{-1}q_1',q_2^{-1}q_2')\;.
$$
A quick look at the above figure leads to $\underline{s}_{q'_2}\underline{s}_{q'_1}\underline{s}_{q'_0}(q_0)=q_0$ and the assertion immediately follows.
\end{proof}

\begin{thm}
\label{3PK}
Let $(M,s,\omega)$ be an elementary symplectic symmetric space with
$\Y$ Abelian.
Assume that ${\Jac}_{\underline s^e}$ is invariant under the inversion map on $Q$. 
Let $q\equiv q(q_0,q_1,q_2)$ denote the unique solution of the equation 
$\underline{s}_{q_2}\underline{s}_{q_1}\underline{s}_{q_0}(q)=q$ in $Q$. 
Set 
$$
\CJ(q_0,q_1,q_2):=|{\Jac}_{(\underline s^e)^{-1}}|(q_0^{-1}q_1)\;
|{\Jac}_{(\underline s^e)^{-1}}|(q_1^{-1}q_2)\;|{\Jac}_{(\underline s^e)^{-1}}|(q_2^{-1}q_0)\;.
$$
Then, one has:
$$
\underline{s}_e^\star\overline{\bm_0}=\bm_0\;,
$$
and   the kernel
 of the product $\star_{\bm_0}$ \eqref{prod-op-ptv} is given by:

\begin{align*}
K_{\bm_0}(x_0,x_1,x_2)&=
\CJ^{1/2}(\Phi_Q^{-1}(q_0,q_1,q_2))\;\big|{\rm Jac}_{\Phi_Q^{-1}}\big|(q_0,q_1,q_2)\\
&\qquad\times\,|\mbox{\rm Jac}_{\Psi}|^{1/2}\big(q_0^{-1}q\big)\,
|\mbox{\rm Jac}_{\Psi}|^{1/2}\big(q_1^{-1}\underline s_{q_0} q\big)\,
|\mbox{\rm Jac}_{\Psi}|^{1/2}\big(q_2^{-1}\underline s_{q_1}\underline s_{q_0} q\big)\\
&\qquad\qquad\times
\bE\big(q^{-1}x_0\big)\, 
\bE\big((\underline s_{q_0} q)^{-1}x_1\big)\, 
\bE\big((\underline s_{q_1}\underline s_{q_0} q)^{-1}x_2\big)\;.
\end{align*}
\end{thm}
\begin{proof}
Under the assumption that ${\Jac}_{\underline s^e}$ is invariant under the inversion map on $Q$,  
$\bm_0$ is $\underline s_e$-invariant and reads:
$$
\bm_0(q)=\left|\mbox{\rm Jac}_{\underline s^e}(q)
\,\mbox{\rm Jac}_{\Psi}(q)\right|^{1/2}\;.
$$
Accordingly, for $f\in L^2(\S)$, we get from Corollary \ref{kernel-op}
 the following expression for the operator kernel
of $\Omega_{\bm_0}(f)$:
\begin{align*}
\Omega_{\bf m_0}(f)[q_0,q]=&
|\mbox{\rm Jac}_{\Psi}|^{1/2}\big(\mid(e, q_0^{-1} q)^{-1}\big)\,
\big|{\rm Jac}_{(\underline s^e)^{-1}}\big|^{1/2}(q_0^{-1}q)\\
&\qquad\times\int_\Y\, f(\mid(q_0, q)b)
\, \bE\big(\mid(e, q_0^{-1} q)b\big)\,{\rm d}_\Y(b)\;.
\end{align*}
Since for $f\in L^2(\S)$, $\Omega_{\bm_0}(f)$ is Hilbert-Schmidt, 
the product of three $\Omega_{\bm_0}(f)$'s is a fortiori trace-class and thus
we can employ the formula:
\begin{align*}
&\Tr\Big[\Omega_{\bf m_0}(f_0)\Omega_{\bf m_0}(f_1)\Omega_{\bf m_0}(f_2)\Big]
=\\
&\int_{Q^3} \Omega_{\bf m_0}(f_0)[q_0,q_1]\,\Omega_{\bf m_0}(f_1)[q_1,q_2]\,
\Omega_{\bf m_0}(f_2)[q_2,q_0]\,{\rm d}_Q(q_0)\,{\rm d}_Q(q_1)\,{\rm d}_Q(q_2)\;.
\end{align*}
Using Theorem \ref{3points} and the formula above for the kernel of $\Omega_{\bm_0}(f_j)$,
we see that the above trace  equals: 
\begin{align*}
&\int_{Q^3\times\Y^3}
|\mbox{\rm Jac}_{\Psi}|^{1/2}\big(\mid(e, q_0^{-1} q_1)^{-1}\big)\,
|\mbox{\rm Jac}_{\Psi}|^{1/2}\big(\mid(e, q_1^{-1} q_2)^{-1}\big)\\
&\times|\mbox{\rm Jac}_{\Psi}|^{1/2}\big(\mid(e, q_2^{-1} q_0)^{-1}\big)
\,\CJ^{1/2}(q_0,q_1,q_2)\, f_0(\mid(q_0, q_1)b_0)\, f_1(\mid(q_1, q_2)b_1)\\
&\times f_2(\mid(q_2, q_0)b_2)
\bE\big(\mid(e, q_0^{-1} q_1)b_0\big)\, 
\bE\big(\mid(e, q_1^{-1} q_2)b_1\big)\, 
\bE\big(\mid(e, q_2^{-1} q_0)b_2\big)\\
&\times{\rm d}_Q(q_0)\,{\rm d}_Q(q_1)\,{\rm d}_Q(q_2)
\,{\rm d}_\Y(b_0)\,{\rm d}_\Y(b_1)\,{\rm d}_\Y(b_2)\;.
\end{align*}
Performing the change of variable $(q_0,q_1,q_2)\mapsto \Phi_Q^{-1}(q_0,q_1,q_2)$
and setting $x_j:=q_jb_j\in\S$, $j=0,1,2$,
we get
\begin{align*}
&\int_{\S^3}
|\mbox{\rm Jac}_{\Psi}|^{1/2}\big(q_0^{-1}q\big)\,
|\mbox{\rm Jac}_{\Psi}|^{1/2}\big(q_1^{-1}\underline s_{q_0} q\big)\,
|\mbox{\rm Jac}_{\Psi}|^{1/2}\big(q_2^{-1}\underline s_{q_1}\underline s_{q_0} q\big)
\CJ^{1/2}(\Phi_Q^{-1}(q_0,q_1,q_2))\\
&\times\big|{\rm Jac}_{\Phi_Q^{-1}}\big|(q_0,q_1,q_2)
\bE\big(q^{-1}x_0\big)\, 
\bE\big((\underline s_{q_0} q)^{-1}x_1\big)\, 
\bE\big((\underline s_{q_1}\underline s_{q_0} q)^{-1}x_2\big)\,
 f_0(x_0)\, f_1(x_1)\\
&\times f_2(x_2)\,{\rm d}_\S(x_0)\,{\rm d}_\S(x_1)\,{\rm d}_\S(x_2)
\;,
\end{align*}
where $q\equiv q(q_0,q_1,q_2)$ is the unique solution of the equation 
$\underline{s}_{q_2}\underline{s}_{q_1}\underline{s}_{q_0}(q)=q$ (see  Lemma \ref{3Pinter}). 
The result then follows by identification.
\end{proof}

Last, using Lemma \ref{E-Psi}, we deduce the following expression for the phase in the 
 kernel of the product $\star_{\bm_0}$:
\begin{cor}
\label{phase}
Write $K_{\bm_0}=A_{\bm_0}\,e^{-iS}$ for the three-point kernel given in Proposition \ref{3PK}. Then we have
for $x_j=q_jb_j\in\S$, $q_j\in Q$, $b_j\in\Y$, $j=0,1,2$:
\begin{align*}
&S(x_0,x_1,x_2)=\\
&\langle\xi,[\Psi(q_0^{-1}q),\log b_0]\rangle+
\langle\xi,[\Psi(q_1^{-1}\underline s_{q_0} q),\log b_1]\rangle+
\langle\xi,[\Psi(q_2^{-1}\underline s_{q_1}\underline s_{q_0}q),\log b_2]\rangle\;,
\end{align*}
where $q\equiv q(q_0,q_1,q_2)$ is the unique solution of the equation 
$\underline{s}_{q_2}\underline{s}_{q_1}\underline{s}_{q_0}(q)=q$ (see Theorem \ref{3points}).
\end{cor}

\section{Extensions  of polarization quadruples}
\label{Extensions  of polarization quadruples}

We first  observe  that, given two polarization quadruples 
$({G}_j,\sigma_j,\xi_j,B_j)$, $j=1,2$ (in the general sense of Definition \ref{def:PQ}),
 a {morphism} $\phi$ between them yields a $G_1$-equivariant intertwiner:
\begin{equation}
\label{INTERTWINER}
\phi^\star:C^\infty(G_2)^{B_2}\to C^\infty(G_1)^{B_1}\;,
\end{equation}
such that
$$
 U_{\chi_1}(g_1)
\phi^\star\hat{\varphi_2}=\phi^\star\big(U_{\chi_2}(\phi(g_1))\hat{\varphi_2}\big)\;.
$$
Note that the condition $\phi^\star\xi_2=\xi_1$ (cf$.$ Definition \ref{def:PQ})
implies $\chi_1(b_1)=\chi_2(\phi(b_1))$ in view of (\ref{unitary-char}). 

 In the context of the  transvection and full   polarization quadruples, we observe:
\begin{lem}\label{UNISO}
Let  $(M,s,\omega)$ be an elementary symplectic 
symmetric space,  associated to a non-exact polarized transvection triple 
$(\g,\sigma,\varpi)$. Consider $(\L,\sigma_\L,\xi,\bB)$ and
$(\tilde{G},\tilde{\sigma},\xi\big|_{\tilde G},B )$ and  the 
full and transvection polarization quadruples as in Definition \ref{TRANSFULL}. 
Then, the intertwiner (\ref{INTERTWINER}) corresponding to the 
injection $\tilde{G}\to\L$ is a linear isomorphism.
\end{lem}
\begin{proof} The injection $j:\tilde{G}\to\L$ induces a global diffeomorphism 
$\tilde{G}/B\to\L/\bB$, $gB\mapsto g\bB$.  Indeed, the map $\tilde{G}/\tilde{K}\to\L/\tilde{\bD}:g\tilde{K}\mapsto g\tilde{\bD}$
is an identification. Considering the natural projections $\tilde{G}/\tilde{K}\to\tilde{G}/B$ and 
$\L/\tilde{\bD}\to\L/\bB$, one observes that the diagram
$$
\begin{array}{ccc}
\tilde{G}/\tilde{K}&\longrightarrow&\L/\tilde{\bD}\\
\downarrow&&\downarrow\\
\tilde{G}/B&\stackrel{\phi}{\longrightarrow}&\L/\bB
\end{array}\;,
$$
where $\phi(gB):=g\bB$, is commutative. In particular, $\phi$ is surjective. Examining its differential
proves that is also a submersion. The space $Q=\tilde{G}/B$, being exponential, has trivial fundamental
group. The map $\phi$ is therefore a diffeomorphism.
Also the restrictions $C^\infty(\tilde{G})^B\to C^\infty(Q)$
and $C^\infty(\L)^\bB\to C^\infty(Q)$ are linear isomorphisms and one observes that 
$j^\star\hat{\varphi}|_Q=\hat{\varphi}|_Q$.
 \end{proof}

Note that when the modular function of $\bB$ coincides with the restriction to $\bB$ of the modular
function of $\L$, then there exists a $\L$-invariant measure on $\L/\bB$. From the isomorphism
$\L/\bB\simeq \tilde G/ B\simeq Q$, we see that the later is a left-Haar measure on $Q$. 
Hence, under the 
assumption above, we deduce that the left-Haar measure ${\rm d}_Q$ is also $\L$-invariant. This,
together 
with  the above Lemma, yields:
\begin{lem}
\label{UNISO2}
In the setting of  Lemma \ref{UNISO} and
when  the modular function of $\bB$ coincides with the restriction to $\bB$ of the modular
function of $\L$, the injection
$\tilde{\fD}\to\fL$ induces a unitary representation ${\mathcal R}:\tilde{\bf D}\to\CU(\CH_\chi)$ of the 
corresponding analytic subgroup $\tilde{\bf D}\subset\L$ on the representation space $\CH_\chi$
associated to  the transvection quadruple $(\tilde{G},\tilde{\sigma},\xi,B)$.
\end{lem}

Consider now two Lie groups $G_j$, $j=1,2$, with unitary representations $(U_j,\CH_j)$, 
together with a Lie group homomorphism
$$
\rho:G_1\to G_2\;,
$$
and form the associated semi-direct product
$G_1\ltimes_{\bR}G_2$,
where
\begin{equation}
\label{hom}
\bR_{g_1}(g_2):=\bC_{\rho(g_1)}(g_2)=\rho(g_1)g_2\rho(g_1)^{-1}\,,\quad\forall g_j\in G_j\;.
\end{equation}
We deduce the representation homomorphism:
\begin{align}
\label{rephom}
{\mathcal R}:G_1\to\CU(\CH_2)\,,\quad g_1\mapsto U_2(\rho(g_1))\;.
\end{align}
Within this setting, we first observe:
\begin{lem}
\label{OBS}
Parametrizing an element $g\in G_1\ltimes_{\bR} G_2$ as $g=g_1.g_2$,
the map
$$
U:G_1\ltimes_{\bR}G_2\to\CU(\CH_1{\otimes}\CH_2)\,,\quad
g\mapsto U_1(g_1)\otimes
{\mathcal R}(g_1)U_2(g_2)\;,
$$
defines a unitary representation of $G_1\ltimes_{\bR}G_2$ on the tensor product Hilbert space 
$\CH:=\CH_1{\otimes}\CH_2$.
\end{lem}
\begin{proof} Let $g_j,g_j'\in G_j$, $j=1,2$. Then, on the first hand:
\begin{align*}
U(g_1g_2.g_1'g_2')=U\big(g_1g_1'\bR_{{g_1'}^{-1}}(g_2)g_2'\big)&= U_1(g_1g_1')\otimes
{\mathcal R}(g_1g_1')\,U_2\big(\bR_{{g_1'}^{-1}}(g_2)g_2'\big)\\& =U_1(g_1g_1')\otimes
{\mathcal R}(g_1)\,U_2(g_2)\,{\mathcal R}(g_1')\,U_2(g_2')\;,
\end{align*}
while on the second hand:
\begin{align*}
U(g_1g_2) \,U(g_1'g_2')&=\big(U_1(g_1)\otimes{\mathcal R}(g_1)\,U_2(g_2)\big)\,\big(U_1(g_1')
\otimes{\mathcal R}(g_1')\,U_2(g_2')\big)\\
&=U_1(g_1g_1')\otimes{\mathcal R}(g_1)\,U_2(g_2)\,{\mathcal R}(g_1')\,U_2(g_2')\;,
\end{align*}
and the proof is complete.
\end{proof}

Consider Lastly two full polarization quadruples $(\L_j,\sigma_j,\xi_j,\bB_j)$, $j=1,2$,  
associated to two  local symplectic symmetric spaces $(M_j,s_j,\omega_j)$.
Let also  $(U_{\chi_j},\CH_{\chi_j},\Omega_j)$ be the unitary representation of $\L_j$ and of the 
representation
of the symplectic symmetric space  $M_j=\tilde G_j/\tilde K_j=\L_j/\tilde\bD_j$ (see Remark 
\ref{TRANS=FULL}).
Finally, let $\bf K$ be a  Lie subgroup  of $\L_1$ that acts  transitively on $M_1$, and 
consider the associated
 ${\bf K}$-equivariant  diffeomorphism:
 $$
 \vf: {\bf K}/({\bf K}\cap \tilde \bD_1)\to M_1\;.
 $$ 
Now, given  a Lie group homomorphism $\rho:{\bf K}\to \tilde \bD_2\subset \L_2$, we can form the 
semi-direct product  ${\bf K}\ltimes_{\bR} \L_2$, according to \eqref{hom}. 
Under these conditions, we have the global 
identification:
\begin{align}
\label{trick}
\big( {\bf K}\ltimes_{\bR} \L_2\big)/\big(({\bf K}\cap \tilde \bD_1)\ltimes_\bR  \tilde \bD_2\big)&\to  
M_1\times M_2\;,\\
(g_1.g_2)({\bf K}\cap \tilde \bD_1)\ltimes_\bR \tilde \bD_2&\mapsto\big(\vf(g_1{\bf K}\cap
 \tilde\bD_1),g_2
\tilde \bD_2\big)\;.\nonumber
\end{align}

\begin{prop}\label{EXTENSION2}
 Let  $U$ be the unitary representation of ${\bf K}\ltimes_{\bR} \L_2$ on 
 $\CH:=\CH_{\chi_1}\otimes\CH_{\chi_2}$ constructed in Lemma \ref{OBS}.
Then under the conditions displayed above,
\begin{enumerate}
\item[(i)] the map
$$
{\bf \Omega}: {\bf K}\ltimes_{\bR} \L_2\to\CU_{sa}(\CH)\,,\quad
g\mapsto U(g)\circ\left(\Sigma_1\otimes\Sigma_2\right)\circ U(g)^*\;,
$$
is constant on the left cosets of $({\bf K}\cap \tilde \bD_1)\ltimes_\bR \tilde \bD_2$ in 
${\bf K}\ltimes_{\bR} \L_2$. 
\item[(ii)] For every $g_1\in{\bf K}$ and $g_2\in \L_2$, one has\footnote{Warning: the reverse order in the group elements.}
$$
{\bf \Omega}(g_2.g_1)=\Omega_1(g_1)\otimes\Omega_2(g_2)\;.
$$
\item[(iii)] Under the identification \eqref{trick}, the quotient map
$$
{\bf \Omega}: M_1\times M_2\to\CU_{sa}(\CH)\;,
$$
is $ {\bf K}\rtimes_{\bR} \L_2$-equivariant.
\end{enumerate}
\end{prop}
\begin{proof}
We start by checking item (ii). Observe first that
\begin{align*}
U(g_2g_1)=U(g_1\bR_{g_1^{-1}}g_2)&=U_{\chi_1}(g_1)\otimes{\mathcal R}(g_1)U_{\chi_2}
(\bR_{g_1^{-1}}g_2)\\&=U_{\chi_1}(g_1)\otimes U_{\chi_2}(g_2){\mathcal R}(g_1)\;.
\end{align*}
Now, since $\rho$ is $\tilde \bD_2$-valued, we have for all $g_1\in{\bf K}$:
$$
{\mathcal R}(g_1)\Sigma_2{\mathcal R}(g_1)^{*}=U_{\chi_2}(\rho(g_1))\Sigma_2U_{\chi_2}(\rho(g_1))^*=\Sigma_2\;,
$$
 (cf$.$ Lemma \ref{EXACTSILA}).
Moreover, one has
\begin{align*}
{\bf \Omega}(g_2g_1)&=\left(U_{\chi_1}(g_1)\otimes U_{\chi_2}(g_2){\mathcal R}(g_1)\right)\circ
\left(\Sigma_1\otimes\Sigma_2\right)\circ\left(U_{\chi_1}(g_1)^*\otimes {\mathcal R}(g_1)^*U_{\chi_2}(g_2)^*
\right)\\
&=U_{\chi_1}(g_1)\Sigma_1 U_{\chi_1}(g_1)^*\otimes U_{\chi_2}(g_2){\mathcal R}(g_1)\Sigma_2  
{\mathcal R}(g_1)^*U_{\chi_2}(g_2)^*\\&=\Omega_1(g_1)\otimes\Omega_2(g_2)\;.
\end{align*}
This implies (ii) and (i) consequently.
 Regarding item (iii), one observes at the level of ${\bf K}\ltimes_{\bR} \L_2$ that
$$
{\bf \Omega}(gg')=U(g)\,{\bf \Omega}(g')\,U(g)^{*}\;,
$$
which is enough to conclude.
\end{proof}
\begin{rmk}
\label{FP}
In the same manner as in Definition \ref{Sigma-m}, given a Borelian function $\bm$
on the product manifold $(\L_1/\tilde\bD_1)\times(  \L_2/\tilde\bD_2)$, we may define for
$g\in \bK\ltimes_\bR \L_2$:
$$
{\bf \Omega}_\bm(g):= U(g)\circ\bm\circ\left(\Sigma_1\otimes\Sigma_2\right)\circ U(g)^*\;.
$$
\end{rmk}

 Of course, the above procedure can be iterated, namely one observes:
\begin{prop}
\label{EXTPROP}
Let  $(\L_j,\sigma_j,\xi_j,\bB_j)$,  $j=1,\dots,N$, be $N$  full  polarization quadruples, 
associated to $N$  elementary symplectic symmetric spaces $(M_j,s_j,\omega_j)$
satisfying the extra conditions of 
coincidence of the modular function on $\bB_j$ with the restriction to $\bB_j$ of the
modular function of $\L_j$, according to Lemmas \ref{UNISO} and \ref{UNISO2}.
 For every $j=1,\dots,N-1$, consider a subgroup ${\bK}_j$ that acts transitively
on $M_j$ together with a Lie group homomorphism $\rho_j:{\bK}_j\to \tilde \bD_{j+1}$.
 Set  ${\bK}_N:=\L_N$,  
assume that for every such $j$, the subgroup $\rho_j({\bK}_j)$ normalizes 
${\bK}_{j+1}$ in $\L_{j+1}$ and
denote by $\bR_{j}$ the corresponding homomorphism from ${\bK}_j$ to 
$\Aut({\bK}_{j+1})$.
Then, iterating the procedure described in Proposition \ref{EXTENSION2} yields a map
$$
{\bf \Omega}: M_1\times M_2\times\dots\times M_N\to\CU_{sa}(\CH)\;,
$$
into the self-adjoint unitaries on the product Hilbert space $\CH:=\CH_{\chi_1}\otimes\dots
\otimes\CH_{\chi_N}$ 
that is equivariant under the natural action of the Lie group
$$
\Big(\dots\Big(\big({\bK}_1\ltimes_{\bR_1}{\bK}_2\big)\ltimes_{\bR_2}{\bK}_3\Big)
\ltimes_{\bR_3}\dots\Big)\ltimes_{\bR_{N-1}}\L_N\;.
$$
\end{prop}

This  `elementary' tensor product construction for the quantization map on direct products
of polarized symplectic symmetric spaces (but with covariance under
semi-direct products of subgroups of the covariance group of each piece) allows
to transfer most of the results of the previous  sections. For notational convenience, 
we formulate all that follows in the context of two elementary pieces, i.e$.$ in the context of Proposition
\ref{EXTENSION2} rather than in the context of Proposition \ref{EXTPROP}.

So in all that follows, we assume we are given two  elementary  symplectic symmetric spaces
$(M_j,s_j,\omega_j)$, $j=1,2$ (see Definition \ref{ELEMENTARY}). 
We also let $\S_j=Q_j\ltimes \Y_j$, $j=1,2$, be the subgroups of
$\tilde G_j$ (and thus of $\L_j$) that acts simply transitively on $M_j$. 
We also assume that we are given a homomorphism $\rho:\S_1\to \tilde \bD_2$ (i.e$.$ the role
of $\bK$ in Proposition \ref{EXTENSION2} is played by $\S_1$). In this particular context, 
the identification \eqref{trick} becomes:
$$
\S_1\ltimes_{\bR}\S_2\to M_1\times M_2\,,\quad g_1.g_2\mapsto (g_1\tilde\bD_1,g_2\tilde\bD_2)\;.
$$ 
We also let $\bm_0:=\bm^1_0\otimes\bm_0^2$ be the smooth function on
$Q_1\times Q_2$, where $\bm_0^j$ is the function on $Q_j$ given in \eqref{m0}.
Combining Proposition \ref{EXTENSION2} with the results of sections 
\ref{Locality and the one-point phase}, 
\ref{Unitarity and midpoints}, \ref{The star-product} and \ref{The three-point kernel}, 
we eventually obtain:

\begin{thm}
\label{uiop}
Let $(M_j,s_j,\omega_j)$, $j=1,2$, be two elementary  symplectic symmetric spaces and
consider an homomorphism $\rho:\S_1\to \tilde \bD_2$. Within
the notations given above, we have:
\begin{enumerate}
\item[(i)] Identifying $\CH_{\chi_1}\otimes \CH_{\chi_2}$ with $L^2(Q_1\times Q_2)$
 and parametrizing an element $g\in\S_1\ltimes_{\bR}\S_2$ as 
$g=q_2b_2q_1b_1$, $q_j\in Q_j$, $b_j\in\Y_j$, we have for $\vf\in\CD(Q_1\times Q_2)$:
\begin{align*}
&{\bf \Omega}_{\bm_0}(g){\varphi}(\bar q_1,\bar q_2)=\\&\bm_0^1(q_1^{-1}\bar q_1)\,
\bm_0^2(q_2^{-1}\bar q_2)\,
  \bE^{\S_1\ltimes\S_2}\big(\bar q_1^{-1}q_1b_1,\bar q_2^{-1}q_2b_2)
\,{\varphi}(\underline{s_1}_{q_1}\bar q_1,\underline{s_2}_{q_2}\bar q_2)\;,
\end{align*}
where 
$$
  \bE^{\S_1\ltimes\S_2}\big(q_1b_1,q_2b_2):=
    \bE^{\S_1}\big(q_1b_1)\,\bE^{\S_2}(q_2b_2)\,,\quad \forall q_j\in Q_j,\, \forall b_j\in\Y_j\;,
    $$
    and $\bE^{\S_j}$, $j=1,2$, is the one-point phase attached to each elementary symplectic
    symmetric space $M_j$ as given in Lemma \ref{RTY}.
\item[(ii)] Moreover, when $\Y_1$ and $\Y_2$ are Abelian, the map
\begin{align*}
{\bf \Omega}_{\bm_0}: L^2(\S_1\ltimes_{\bR}\S_2)&\to 
{\mathcal L}^2\big(\CH_{\chi_1}\otimes \CH_{\chi_2}\big)\,,\\ f&\mapsto \int_{\S_1\ltimes_{\bR}\S_2}
f(g)\,{\bf \Omega}_{\bm_0}(g)\,{\rm d}_{\S_1\ltimes_{\bR}\S_2}(g)\;,
\end{align*}
is unitary and $\S_1\ltimes_{\bR}\S_2$ equivariant.
\item[(iii)] Denoting by ${\bf \sigma}_{\bm_0}$ the adjoint of ${\bf \Omega}_{\bm_0}$, the associated
deformed product:
$$
f_1\star_{\bm_0} f_2:={\bf \sigma}_{\bm_0}\big[{\bf \Omega}_{\bm_0}(f_1)\,
{\bf \Omega}_{\bm_0}(f_2)\big]\;,
$$
takes on $\CD(\S_1\ltimes_{\bR}\S_2)$ the expression
$$
\int_{(\S_1\ltimes_{\bR}\S_2)^2} K_{\bm_0}^{\S_1\ltimes_{\bR}\S_2}
(g,g',g'')\,f_1(g')\,f_2(g'')\,
{\rm d}_{\S_1\ltimes_{\bR}\S_2}(g')\,{\rm d}_{\S_1\ltimes_{\bR}\S_2}(g'')\;,
$$
where the three-points kernel $K_{\bm_0}^{\S_1\ltimes_{\bR}\S_2}$ is given, with $g=g_2g_1$,
 $g'=g'_2g'_1$ and  $g''=g''_2g''_1$  by
$$
K_{\bm_0}^{\S_1\ltimes_{\bR}\S_2}(g,g',g''):=K_{\bm_0^1}^{\S_1}(g_1,g'_1,g''_1)\,
K_{\bm_0^2}^{\S_2}(g_2,g'_2,g''_2)\;,
$$
with $K_{\bm_0}^{\S_j}$, $j=1,2$, as given in Proposition \ref{3PK}.
\end{enumerate}
\end{thm}

\chapter{Quantization of K\"ahlerian Lie groups}
\label{QKLG}

The aim of this chapter is two-fold. First, 
 we establish that the symplectic symmetric space
$(\S,s,\omega^\S)$ associated with an elementary normal $\bf j$-group $\S$ (see section
\ref{H-group}) underlies an elementary local and polarized symplectic symmetric space,
in the sense of Definitions \ref{TRANSFULL}, \ref{GLOBAL} and \ref{ELEMENTARY}. 
Using the whole construction of chapter \ref{QPSSS}, we will
then be able to construct a $\B$-equivariant quantization map for any normal $\bf j$-group
$\B$. Second, we will show that the composition law of symbols associated to this quantization map,
coincides exactly with the left-$\B$-equivariant star-products on $\B$ that we have considered in
chapter \ref{NFSP}.

\section{The transvection quadruple  of an elementary normal $\bf j$-group}
\label{H-group2}

  We start by  describing the non-exact  transvection siLa $(\g,\sigma,\varpi)$ underlying the
symplectic symmetric space structure $(\S,s,\omega^\S)$ of an elementary normal $\bf j$-group, as 
described in section \ref{H-group}.

The transvection group $G$ of $\S$ is the connected and simply connected Lie group
 whose  Lie  algebra $\g$ is a one-dimensional split extension of two copies of the 
Heisenberg algebra:
\begin{equation}
\label{frak-g}
\g:=\a\ltimes_\rho(\h\oplus\h)\;,
\end{equation}
where, again, $\a=\R H$ and  the extension homomorphism is given by $\rho:=\rho_\h\oplus(-\rho_\h)
\in{\rm Der}(\h\oplus\h)$, 
with $\rho_\h$  defined in \eqref{split}. 
The involution $\sigma$ of $\g$ is given by
\begin{equation}
\label{Invol}
\sigma\big(aH+ (X\oplus Y)\big):=(-aH)+ (Y\oplus X)\,,\qquad \forall a\in\R\,,\quad \forall X,Y\in\h\;.
\end{equation}
One has the associated  ($\pm1$)-eigenspaces decomposition:
\begin{equation*}
\g=\k\oplus\p\,,\qquad\k:=\h_+\quad\mbox{ and }\quad \p:=\a\oplus\h_-\;,
\end{equation*}
where for every subspace $F\subset\h$, we set 
\begin{equation}
\label{pm}
F_\pm:=\{X\oplus (\pm X),\;X\in F\}\subset\h\oplus\h\;,
\end{equation}
and for every element $X\in\h$ we let 
$X_\pm:=\tfrac{1}{2}(X\oplus(\pm X))\in\h\oplus\h$.
Last, we define $\varpi\in\Lambda^2\g^\star$ by
\begin{equation}
\label{2-form}
\varpi(H,E_-)=2\quad\mbox{and}\quad \varpi(v_-,v'_-)=\omega^0(v,v')\,,\quad
\forall v,v'\in V\;,
\end{equation}
and by zero everywhere else on $\g\times\g$. Note that  $\varpi$ is $\k$-invariant and
its restriction to $\p$ is non-degenerate. This implies that $\varpi$ is a 
Chevalley two-cocycle (see \cite{Bi95}). 
Also, from $[H,\h_-]=\h_+=\k$, we deduce that $[\p,\p]=\k$
and clearly the action of $\k$ on $\p$ is   faithful. 
Thus, in terms of  Definition \ref{TST}, 
we have proved the following:
\begin{prop}
The siLa $(\g,\sigma,\varpi)$ defined by \eqref{frak-g}, \eqref{Invol} and \eqref{2-form},
is a transvection symplectic triple.
\end{prop}

Consider now $(\tilde\g,\tilde\sigma,\delta\xi)$ the exact siLa constructed out of the
non-exact siLa $(\g,\sigma,\varpi)$ as in Lemma \ref{EXACTSILA}.
Recall that $\tilde\g$ is the one-dimensional central extension of $\g$ with generator
$Z$ and table
$$
[X,Y]_{\tilde\g}=[X,Y]_\g+\varpi(X,Y)\,Z\,,\quad\forall X,Y\in\g\;.
$$
The involution $\tilde\sigma$ equals ${\rm id}_{\tilde \k}\oplus(-{\rm id}_\p)$, where
$\tilde\k=\k\oplus\R Z$ and $\xi\in\tilde\g^\star$ is defined by $\langle\xi,Z\rangle=1$
and $\xi\big|_{\g}=0$.
Accordingly, we set $ \tilde G=\exp\{\tilde\g\}$ and $\tilde K=\exp\{\tilde\k\}$.
We identify $\tilde \g$ with $\tilde G$ via the global chart:
 \begin{equation}
 \label{tildeG}
aH+v_1\oplus v_2+t_1E\oplus t_2 E+\ell Z\mapsto \exp\{aH\}\exp\{v_1\oplus v_2+t_1E\oplus t_2 E
+\ell Z\}\;,
 \end{equation}
 where $a,t_1,t_2,\ell\in\R$ and $v_1,v_2\in V$. The group law of $\tilde G$ in these
 coordinates then reads:
\begin{align}
\label{prod}
& (a,v_1,v_2,t_1,t_2,\ell)(a',v_1',v_2',t_1',t_2',\ell')=
\Big(a+a',e^{-a'}v_1+v_1',e^{a'}v_2+v_2',\nonumber\\
&\hspace{1,7cm}e^{-2a'}t_1+t_1'+\tfrac12e^{-a'}\omega^0(v_1,v_1'),
 e^{2a'} t_2+t_2'+\tfrac12e^{a'}\omega^0(v_2,v_2'),\nonumber\\
&\vspace{2cm}\ell+\ell'
+(e^{-2a'}-1)t_1+(e^{2a'}-1)t_2+\tfrac12\omega^0(e^{-a'}v_1-e^{a'}v_2,v_1'-v_2')\Big)\;,
\end{align}
and the inversion map is given by:
\begin{align*}
 &(a,v_1,v_2,t_1,t_2,\ell)^{-1}=\\
&  \big(-a,-e^av_1,-e^{-a}v_2,-e^{2a}t_1,- e^{-2a}t_2, -\ell
-(e^{2a}-1)t_1-(e^{-2a}-1)t_2\big)\;.
 \end{align*}
  Moreover the involution $\tilde{\sigma}$ admits the following expression:
 $$
 \tilde{\sigma} (a,v_1,v_2,t_1,t_2,\ell)=(-a,v_2,v_1,t_2,t_1,\ell)\;.
 $$
Under the parametrization of $\tilde G$ given above,
we consider the following global coordinates system  on $\tilde G/\tilde K$:
 \begin{equation}
 \label{ChartS}
\tilde G/\tilde K\to\R^{2d+2},\quad (a,v_1,v_2,t_1, t_2,\ell ) \widetilde K\mapsto
 \big(a,v_1-v_2,t_1-t_2-\tfrac12\omega^0(v_1,v_2)\big)\;.
 \end{equation}
From the formula $s_{g\tilde K}(g'\tilde K)=g\tilde\sigma(g^{-1}g')\tilde K$ for the symmetry on 
$\tilde G/\tilde K$, we deduce the following  isomorphism of symplectic symmetric spaces: 
 \begin{prop}
 Under the identifications $\S\simeq\R^{2d+2}\simeq \tilde G/\tilde K$ associated with the charts 
 \eqref{chartS} and \eqref{ChartS}, the symplectic symmetric space $\tilde G/\tilde K$ underlying
 the exact siLa  $(\tilde\g,\tilde\sigma,\delta\xi)$ defined above, is isomorphic to the symplectic
 symmetric space $(\S,s,\omega^\S)$ underlying an elementary normal $\bf j$-group $\S$,
 as given in section \ref{PStheory}.
 \end{prop}

Next, we need to  endow $(\S,s,\omega^\S)$ with a structure of 
polarized symplectic symmetric space. From Lemma \ref{pol}, it suffices to specify $W$, a
$\k$-invariant Lagrangian subspace  of $\p$. 
To this aim, we again consider the splitting of the $2d$-dimensional symplectic 
vector space $(V,\omega^0)$ 
into a direct sum of two Lagrangian subspaces in symplectic duality:
$$
V=\l^\star\oplus\l\;.
$$
Relatively to this decomposition and within the notation \eqref{pm}, we define:
\begin{align*}
W:=\l_-\oplus\R E_-\subset\g\;.
\end{align*}
Following then Proposition \ref{pol-quad},
we let
$$
\b:=\tilde\k\oplus W=\k\oplus\R Z\oplus W\:,
$$
be the  polarization Lie algebra. 
 Accordingly with  the terminology introduced in Definition \ref{TRANSFULL},
we call $(\tilde\g,\tilde\sigma,\xi,\b)$ the  transvection quadruple of the
symplectic symmetric space $(\S,s,\omega^\S)$.

 Regarding the question of existence of $\tilde G$-invariant measures on the 
 homogeneous spaces  $\tilde G/\tilde K$ and $\tilde G/B$, we first observe the following
 fact:
\begin{lem}
Let $(\g=\k\oplus\p,\sigma)$ be the involutive Lie algebra associated to a solvable, simply connected, oriented symmetric space $M=G/K$ such that $[\p,\p]=\k$. Then both 
$G$ and $K$ are unimodular Lie groups.
\end{lem}
\begin{proof} Under the orientation hypothesis, let $\nu_\p$ denote a $K$-invariant volume element on $\p\simeq T_KM$. 
Since $\k=[\p,\p]\subset[\g,\g]$, under the solvability assumption,  the Lie algebra $\k$ is nilpotent. Therefore
for every $Z\in\k$, one has $\Tr\left(\ad_Z|_\k\right)=0$ and there exists an $\ad$-invariant volume element
$\nu_\k$ on $\k$. The volume element $\nu_\k\wedge\nu_\p$ on $\g$ is therefore $\k$-invariant. It is also $\ad_\p$-invariant
since, due to the iLa condition, for every $X\in\p$, the element $\ad_X$ is trace-free. The simple-connectedness of $M$ implies the connectedness
of $K$. The latter is therefore unimodular, as well as $G$.
\end{proof}

 \begin{rmk}
 \label{DoExist}
The above lemma implies  that $ \tilde G$, $\tilde K$ and $B:=\exp\{\b\}$
 are all unimodular  ($\b$ is nilpotent). In particular,  there exist 
  $\tilde G$-invariant measures on the homogeneous
 spaces  $\tilde G/\tilde K$ and $\tilde G/B$.
 \end{rmk}

Last, we need to specify the local and elementary structures underlying
the polarized symplectic symmetric space $(\S,s,\omega^\S)$, as introduced in Definition
\ref{GLOBAL} and Definition \ref{ELEMENTARY} respectively.  We first note:
\begin{lem}
\label{OK}
Let
$$
\q:=\a\oplus(\l^\star\oplus0)\quad\mbox{and}\quad\fY:=(\l\oplus0)\oplus\R\big((E\oplus0)+Z\big)\;.
$$
Then $\q$
is a Lie subalgebra of $\tilde\g$ supplementary to $\b$ and 
$\fY$ is an Abelian Lie subalgebra of $\b$ which is normalized by  $\q$. Moreover, 
 the associated semi-direct product $\q\ltimes\fY$
is naturally isomorphic to the Lie algebra $\s$ and
induces the vector space decomposition
$\tilde{\g}=\s\oplus\tilde{\k}$. 
\end{lem}
\begin{proof}
First observe that for all $X\in\h$, one has $X\oplus0=X_-+X_+\in\h\oplus\h$ and
therefore
\begin{align*}
\varpi(H,E\oplus0)&=\varpi(H,E_-+E_+)=\varpi(H,E_-)=2\;,\\
\varpi(v\oplus0,v'\oplus0)&=\varpi(v_-+v_+,v'_-+v'_+)
=\varpi(v_-,v'_-)=\omega^0(v,v')\,,\quad \forall v,v'\in V\;.
\end{align*}
The fact that $\q$ is a Lie subalgebra follows from 
$$
[H,\l^\star\oplus 0]_{\tilde \g}=[H,\l^\star\oplus 0]+\varpi(H,\l^\star\oplus0) Z=\l^\star\oplus 0\;,
$$
and
$$
[\l^\star\oplus 0,\l^\star\oplus 0]_{\tilde{\g}}=\omega^0(\l^\star,\l^\star)\big(E\oplus0+Z\big)=0\;.
$$
 Next, observe that $\fY$ is Abelian:
 $$
 [\fY,\fY]_{\tilde \g}=[\l\oplus 0,\l\oplus0]_{\tilde \g}+[\l\oplus 0,\R\big((E\oplus0)+Z\big)]_{\tilde \g}=
 \omega^0(\l,\l)\big((E\oplus0)+Z\big)=0\;.
 $$
To see that $\q$ normalizes $\fY$, let $x\in\l$ and $t\in\R$. Then one has
\begin{align*}
&[H,x\oplus0+t(E\oplus0+Z)]_{\tilde \g}\\&=[H,x\oplus0]+\varpi(H,x_-)Z+t[H,E\oplus0]+t\varpi(H,E_-)Z
\\&=x\oplus0+2t(E\oplus0+Z)\;.
\end{align*}
Similarly, for all $y\in\l^\star$, one has:
$$
[y\oplus0,x\oplus0+t(E\oplus0+Z)]_{\tilde\g}=[y\oplus0,x\oplus0]
+\omega^0(y,x)Z=\omega^0(y,x)(E\oplus0+Z)\;.
$$
The rest of the statement is immediate. 
\end{proof}

\begin{rmk}
\label{normalaswell}
Neither $\q$ nor $\fY$ are $\tilde\sigma$-stable. However, since $\tilde\sigma(\q)=\a\oplus(0\oplus \l^\star)$ and $[0\oplus \l^\star,\fY]=0$, one sees that $\tilde\sigma(\q)$ normalizes $\fY$ as well.
\end{rmk}

\begin{lem}
Equipped  with the subgroup $Q=\exp\{\q\}$ of $\tilde G$, the  symplectic symmetric 
space $(\S,s,\omega^\S)$ is local in the sense of Definition \ref{GLOBAL}.
\end{lem} 
\begin{proof}
Note first that
$$
\b=\k\oplus W\oplus\R Z=\h_+\oplus \l_-\oplus\R E_-\oplus\R Z=
\l^\star_+\oplus(\l\oplus\l)\oplus(\R E\oplus\R E)\oplus\R Z\;.
$$
Thus,  under the parametrization \eqref{tildeG} of $\tilde G$, we have
\begin{equation}
\label{B}
B=\big\{(0,n\oplus m_1,n\oplus m_2,t_1,t_2,\ell)\,:\,m_1,m_2\in\l\,,
\;n\in\l^\star\,,\;t_1,t_2,\ell\in\R\big\}\;,
\end{equation}
and
\begin{equation*}
Q=\big\{(a, n,0,0,0,0)\,:\,
\;n\in\l^\star\,,\;a\in\R\big\}\;.
\end{equation*}
Thus for $q=(a, n,0,0,0,0)$ and $b=(0,n'\oplus m_1',n'\oplus m_2',t_1',t_2',\ell')\in B$, we have
using \eqref{prod}:
$$
q.b=\big(a,(n+n')\oplus m_1',n'\oplus m_2',t_1'+\tfrac12\omega^0(n,m_1'),t_2',\ell'
+\tfrac12\omega^0(n,m_1'-m_2')\big)\;,
$$
from which we deduce that the map 
$$
Q\times B\to\tilde G\;,\quad (q,b)\mapsto q.b\;,
$$
is a global diffeomorphism (i.e$.$  the first condition of Definition \ref{GLOBAL} is satisfied).
Note that identifying $B$ with $\b$, one has
$$
b\tilde{\sigma}(b^{-1})=2\,b_\p\;,
$$
where we set $b=:b_{\tilde{\k}}+ b_\p$ according to the vector space decomposition $\b=\tilde{\k}\oplus(\b\cap\p)$.

For the second condition, observe that as $\b\cap\p=\l_-\oplus\R E_-$, we get
$$
[\a,\b\cap\p]=[\a,\l_-]\oplus[\a,\R E_-]=\l_+\oplus\R E_+\subset\h_+=\k\subset\b\;.
$$
To check the last condition, consider  $q= (a, n,0,0,0,0)\in Q$, with $a\in\R$, $n\in\l^\star$. We
then have
\begin{align}
\label{QB}
\tilde\sigma q= (-a,0, n,0,0,0)=(-a,-n,0,0,0,0)(0,n,n,0,0,0)=(\tilde\sigma q)^Q\,
(\tilde\sigma q)^B\;,
\end{align}
since $(0,n,n,0,0,0)\in\l^\star_+\subset\h_+\subset\b$. Thus,
$\chi\big((\tilde\sigma q)^B\big)=1$
since for $b=(0,n\oplus m_1,n\oplus m_2,t_1,t_2,\ell)\in B$, we have
$\chi(g)=e^{i\ell}$.
\end{proof}

Next, we come to the symmetric space structure of the group $Q$:
\begin{lem}
\label{Qs}
In the global chart:	
\begin{equation}
\label{chartQ}
\q\simeq\a\oplus\l^\star \to Q\,,\quad (a,n)\mapsto \exp\{aH\}\exp\{n\oplus0\}\;,
\end{equation}
the left invariant symmetric space structure $\underline s$ on $Q$ described in  Lemma
\ref{WhenLocal} reads:
\begin{equation}
\label{SSQ}
\underline s_{(a,n)}(a',n')=\big(2a-a',2\cosh(a-a')\,n-n'\big)\;.
\end{equation}
Moreover, the symmetric space $(Q,\underline s)$ admits a midpoint map, which in the coordinates
above, is given by:
\begin{equation*}
\mid\big((a,n),(a_0,n_0)\big)=\Bigl(\frac{a+a_0}2,\frac{n+n_0}{2}\sech(\tfrac{a-a_0}2)\Bigr)\;.
\end{equation*}
\end{lem}
\begin{proof}
By definition we have $\underline s_qq'=q\tilde \sigma\big(q^{-1}q'\big)^Q$ and the formula
for the symmetry follows easily from \eqref{prod} and \eqref{QB}. 
The  formula for the midpoint map comes from a direct computation of the inverse diffeomorphism
of the partial map $\underline s^q:=[Q\ni q'\mapsto \underline s_{q'}q\in Q]$.
\end{proof}

\begin{rmk}
\label{N}
Setting $A:=\exp\{\a\}$ and $N:={\exp}\{\l^\star\oplus 0\}$ we have the global decomposition
$Q=AN$ and for $q=an\in Q$, the symmetry at the neutral element reads 
$\underline s_eq=a^{-1}n^{-1}$.
Also,  the global chart
\begin{equation}
 \label{ChartQ}
 \widetilde G/B\to\R^{d+1},\quad (a,n_1\oplus m_1,n_2\oplus m_2,t_1,t_2,\ell ) B
 \mapsto (a,n_1-n_2)\;,
 \end{equation}
 $a,t_1,t_2,\ell\in\R$, $n_1,n_2\in\l^\star$, $m_1,m_2\in\l$, 
 identifies  $\tilde G/B$ with $Q$ via the coordinate system \eqref{chartQ}.
\end{rmk}  

\begin{lem}
The Abelian subgroup  $\Y:=\exp\{\fY\}$ of $\tilde G$,
endows the local symplectic symmetric space $(\S,s,\omega^\S)$ with an elementary
structure in the sense of Definition \ref{ELEMENTARY}.  
\end{lem}
\begin{proof}
We already know by Lemma \ref{Qs}, that the left invariant symmetric space $(Q,\underline s)$ 
admits a midpoint map. We also know by Lemma \ref{OK} that $\Y$ is normalized by
$Q$ and that $\S$ is isomorphic to $Q\ltimes\Y$. But we need to know
that $Q\ltimes\Y$ acts simply transitively on  the symmetric space
$\tilde G/\tilde K$. For this,  
let $g=(a,v,0,t,0,t)\in Q\ltimes\Y$ and $g'=(a',v_1',v_2',t_1',t_2',\ell')\in\tilde G$.
Then we get 
\begin{align*}
gg'=&
\big(a+a',e^{-a'}v+v_1',v_2',e^{-2a'}t+t_1'+\tfrac12e^{-a'}\omega^0(v,v_1'),t_2',\\
&\qquad\qquad\qquad\qquad\qquad\qquad\qquad \ell'+e^{-2a'}t+\tfrac12\omega^0(e^{-a'}v,v_1'-v_2')
\big)\;,
\end{align*}
and thus in the chart \eqref{ChartS} of $\tilde G/\tilde K$, we get:
$$
gg'\tilde K\mapsto \big(a+a',e^{-a'}v+v_1'-v_2',e^{-2a'}t+t_1'-t_2'-\tfrac12\omega^0(v_1',v_2')+\tfrac12\omega^0(e^{-a'}v,v_1'-v_2')\big)\;.
$$
This  means that under the identification $\S\simeq\tilde G/\tilde K$, $Q\ltimes\Y\simeq\S$
acts by left translations and the second condition of Definition \eqref{ELEMENTARY}
is verified. For the third condition, note that under the parametrization \eqref{tildeG}
of $\tilde G$, we have
\begin{align*}
\Y=\big\{(0,m,0,t,0,t)\,:\, m\in\l\,,\;t\in\R\big\}\;.
\end{align*}
Take $q=(a,n,0,0,0,0)\in Q$, $a\in\R$, $n\in\l^\star$ and $b=e^y=(0,m,0,t,0,t)\in\Y$, $m\in\l$, $t\in\R$.
A computation then shows that
$$
\langle\xi,\big(\Ad_{q^{-1}}-\Ad_{(\underline s_e q)^{-1}}\big)y\rangle
=2t\sinh2a+2\cosh a\,\omega^0(n,m)\;,
$$
which entails that 
\begin{align}
\label{Psi-explicit}
\Psi(a,n)=\big(2\sinh2a,2n\cosh a\big)\;.
\end{align}
The last condition follows from \eqref{2-form}.
\end{proof}
\begin{rmk}
\label{tildeG-actions}
Parametrizing $\tilde G$ as in \eqref{tildeG}, $\S\simeq\tilde G/\tilde K$ as in \eqref{ChartS}
and $Q \simeq\tilde G/B$ as in  \eqref{ChartQ},
we have the following expression for the action of $\tilde G$ on $\S$:
\begin{align*}
&(a,v_1,v_2,t_1,t_2,\ell).(a',v',t')\\
&=\big(a+a',e^{-a'}v_1-e^{a'}v_2+v',
e^{-2a'}t_1-e^{2a'}t_2+t'-\tfrac12\omega^0(v_1,v_2)\big)\;,
\end{align*}
and on $Q$:
$$
(a,n_1\oplus m_1,n_2\oplus m_2,t_1,t_2,\ell).(a',n')=\big(a+a',e^{-a'}n_1-e^{a'}n_2+n'\big)\;.
$$
\end{rmk}
\begin{rmk}
\label{fd-explicit2} 
From similar methods than those leading to Lemma \ref{fd-explicit}, we deduce that  we have:
$$
\fd_Q\asymp[(a,n)\in Q\mapsto \cosh a+|n|(1+e^a)]\;.
$$
\end{rmk}
From the Remark above 	and in analogy with Remark \ref{lastS}, we define the Fr\'echet valued Schwartz space of
 $Q$, denoted $\CS(Q,\CE)$, as the set
of smooth functions such that all left (or right) derivatives decrease faster than any power 
of $\fd_Q$. The latter space is Fr\'echet for the semi-norms:
\begin{align*}
f\in\CS(Q,\CE)\mapsto\sup_{X\in\,\CU_k(\q)}\sup_{x\in Q}\Big\{\frac
{\fd_Q(x)^n\big\|\widetilde X\,f(x)
\big\|_j}{|X|_k}\Big\}\;,\qquad j,k,n\in\N\;,
\end{align*}
or even for
\begin{align*}
f\in\CS(Q,\CE)\mapsto\sup_{X\in\,\CU_k(\q)}\sup_{x\in Q}\Big\{\frac
{\fd_Q(x)^n\big\|\underline X\,f(x)
\big\|_j}{|X|_k}\Big\}\;,\qquad j,k,n\in\N\;.
\end{align*}

\section{Quantization of  elementary normal $\bf j$-groups}\label{Quantization of  elementary normal groups}

In this section, we specialize the different ingredients of our quantization map
  in the case of the  elementary  symplectic
symmetric space $(\S,s,\omega^\S)$ determined by an elementary normal $\bf j$-group.
We also  (re)introduce a real parameter $\theta$
in the definition of the character \eqref{unitary-char}:
\begin{align*}
\chi_\theta(b):=\exp\{\tfrac i\theta\langle\xi,\log(b)\rangle\}\,,\quad b\in B\,,\quad \theta\in\R^*\;,
\end{align*}
which is globally defined as $B$ is exponential.
By Lemma \ref{G-M-invariant} and Remark \ref{DoExist}, the Haar measure ${\rm d}_\S$ 
on $\S$ (respectively ${\rm d}_Q$ on $Q$) is invariant under
 both $ s_e^\star$ (respectively $\underline s_e^\star$) and $\tilde G$.
Observe that under the  parametrization \eqref{B} of the group  $B$, we have 
$\chi_\theta(b)=\exp\{\tfrac i\theta\ell\}$.
Note  that within the chart \eqref{chartQ}, any left-invariant Haar measure
${\rm d}_Q$ on $Q$ is a  multiple of the Lebesgue measure on $\q$
(these facts are transparent in Equations \eqref{SSs}, \eqref{SSQ} and in Remark
\ref{tildeG-actions}). 
Also, within the chart \eqref{chartS}, any left invariant Haar measure ${\rm d}_\S$ on $\S$
is a multiple of the Lebesgue measure on $\q\ltimes\fY$.  By Remark \ref{Rc1} (iii), the restriction 
to $\S=Q\ltimes \Y$ of the induced representation $U_{\chi_\theta}$ (that we denote by 
$U_\theta$ from now on) of $\tilde G$ on $L^2(Q,{\rm d}_Q)$ reads within the charts \eqref{chartS}
 on $\S$ and \eqref{chartQ} on $Q$:
\begin{align*}
&U_\theta(a,v,t)\psi(a_0,n_0)\\
&=\exp\Big\{\tfrac{i}\theta\big( e^{2(a-a_0)}t+
\omega^0(\tfrac12e^{a-a_0}n-n_0,e^{a-a_0}m)\big)\Big\}\,
\psi\big(a_0-a,n_0-e^{a-a_0}n\big)\;,
\end{align*}
where $(a,v,t)\in\S$ with $a,t\in\R$ and $v=n\oplus m\in\l^\star\oplus\l=V$
and $(a_0,n_0)\in Q$ with $a_0\in\R$ and $n_0\in\l^\star$.
\begin{rmk}
In accordance  with the  notations of earlier chapters, from now on,  we make explicit
the dependence on  the parameter $\theta\in\R^*$ in all the objects we are considering.
For instance, we now set $\Omega_{\theta,\bm}$ instead $\Omega_{\bm}$ for the quantization map,
 $\star_{\theta,\bm}$ instead of $\star_{\bm}$ for the associated composition product, 
 $K_{\theta,\bm}$ instead of $K_{\bm}$ for its three-points kernel, $\bE_\theta$ instead of
 $\bE$ for the one-point phase etc.
\end{rmk}
 
\begin{lem}
\label{S1POINT}
Within the coordinates \eqref{chartS} on $\S$ and under the decomposition
$v=n\oplus m\in\l^\star\oplus\l=V$, the one-point phase $\bE_\theta$ of Lemma \ref{RTY} reads
\begin{align*}
\bE_\theta(a,v,t)=\exp\big\{-\tfrac{2i}\theta\big(t\sinh 2a+\omega^0(n,m)\cosh^2a\big)\big\}\;.
\end{align*}
\end{lem}
\begin{proof} Recall that from Lemma \ref{E-Psi}, we have:
$$
\bE_\theta(q^{-1}e^y)=\exp\big\{\tfrac{i}{\theta}\langle\xi,[\Psi(q),y]\rangle\big\}\;.
$$
Which from (\ref{Psi-explicit}) becomes:
$$
\bE_\theta(q^{-1}e^y)=\exp\big\{\tfrac{2i}{\theta}\left(t\sinh2a+\omega^0(n,m)\cosh a\right)\big\}\;,
$$
with $q=(a,n,0,0,0,0)\in Q$, $a\in\R$, $n\in\l^\star$ and $b=e^y=(0,m,0,t,0,t)\in\Y$, $m\in\l$, $t\in\R$.
Since $q^{-1}=(-a,-e^{a}n,0,0,0,0)$, one gets
$$
\bE_\theta(qb)=\exp\big\{-\tfrac{2i}{\theta}\left(t\sinh2a+\omega^0(n,m)e^a\cosh a\right)
\big\}\;.
$$
 One concludes by observing that the coordinates $qb$ are related to the coordinates \eqref{chartS} through
$$
(a,v,t)=(a,n,0,0,0,0).(0,m,0,t-\tfrac{1}{2}\omega^0(n,m),0,t-\tfrac{1}{2}\omega^0(n,m))\;.
$$
\end{proof}

From \eqref{SSQ} and \eqref{Psi-explicit}, we observe:
 $$
\left|\mbox{\rm Jac}_{\underline s^e}(a,n)\right|={2^{d+1}\cosh^{d}a}\,,\qquad 
\left|\mbox{\rm Jac}_{\Psi}(a,n)\right|={2^{d+2}\cosh2a\cosh^{d}a}\;,
$$
so that the element $\bm_0$ given in \eqref{m0}  reads:
$$
\bm_0(a,n)=2^{d+2}\,\cosh^{1/2}2a\,\cosh^{d}a\;.
$$
From this, we deduce:
\begin{prop} 
\label{Omega-explicit}
Parametrizing $\S$ as in \eqref{chartS} and $Q$ as in \eqref{chartQ},
we have the following expression for the action on $L^2(Q,{\rm d}_Q)$ of the unitary quantizer
$ \Omega_{\theta,\bm_0}(x)$, $x\in\S$, associated with the
polarized symplectic symmetric space underlying an elementary normal $\bf j$-group $\S$:
\begin{align*}
\label{Omega-explicit}
& \Omega_{\theta,\bm_0}(a,v,t)\psi(a_0,n_0)
 =2^{d+1}\cosh(2a-2a_0)^{1/2}\cosh(a-a_0)^d\\
 &\times\exp\Big\{\frac{2i}\theta\big( \sinh(2a-2a_0)t
 +\omega^0(\cosh(a-a_0)n-n_0,\cosh(a-a_0)m)\big)\Big\}\\
&\times\psi\big(2a-a_0,2\cosh(a-a_0)n-n_0\big)\;.
\end{align*}
\end{prop}
\begin{rmk}
Observe that $\left|\mbox{\rm Jac}_{\underline s^e}\right|$ is $\underline s_e$-invariant, as it is
an even function of the variable $a$ only. Thus, by Lemma \ref{SA} and
Remark \ref{unitary-not-SA}, we deduce that
the unitary quantization map $ \Omega_{\theta,\bm_0}$ is also
 compatible with the natural involutions of its source and range spaces 
(the complex conjugation on $L^2(\S,{\rm d}_\S)$ and the adjoint on ${\mathcal L}^2(L^2(Q,
{\rm d}_Q)$). In 
particular, it sends real-valued functions to self-adjoint operators.
\end{rmk}

 Our next result is one of the key steps of this chapter: it renders
transparent the link  between chapters \ref{NFSP} and \ref{DFA} and chapters 
\ref{QPSSS} and \ref{QKLG}.
For this, we need the explicit  expression of the tri-kernel $K_{\theta,\bf m_0}$
of the product \eqref{prod-op-ptv} for an elementary normal $\bf j$-group $\S$.
First, observe that the unique solution of the equation $\underline s_{q_2}
\circ \underline s_{q_1}\circ \underline s_{q_0}(q)=q$, $q,q_0,q_1,q_2\in Q$, as given 
in Lemma \ref{3points}, reads
$$
q=\big(a_0-a_1+a_2,\cosh(a_0-a_1)n_2-\cosh(a_2-a_0)n_1+
\cosh(a_1-a_2)n_0\big)\;.
$$
From \cite{Qi,BCSV}, we extract
\begin{lem} Within the notations of Proposition \ref{3PK}, we have
$$
\CJ=\left|\mbox{\rm Jac}_{\Phi_Q}\right|\;.
$$
\end{lem}
Then, Proposition \ref{3PK} and a straightforward computation gives
$K_{\theta,\bf m_0}=A_{\bf m_0} \,e^{\tfrac{2i}\theta S}$ with:
\begin{align*}
A_{\bf m_0}(x_1,x_2,x_3)&={\bf m_0}(a_1-a_2)\,{\bf m_0}(a_2-a_3)\,{\bf m_0}(a_3-a_1)
\\&=2^{3d+3}\,\cosh(2a_1-2a_2)^{1/2}\,\cosh(a_1-a_2)^{d}\,\cosh(2a_2-2a_3)^{1/2}\,\\
&\qquad\qquad\times\cosh(a_2-a_3)^{d}\,\cosh(2a_3-2a_1)^{1/2}\,\cosh(a_3-a_1)^{d}
\;,
\end{align*}
and
\begin{align*}
S(x_1,x_2,x_3)&=
\sinh(2a_2-2a_1)t_3+
\sinh(2a_1-2a_3)t_2+\sinh(2a_3-2a_2)t_1\\&\quad
+\cosh(a_1-a_2)\,\cosh(a_2-a_3)\,\omega^0(v_3,v_1)\\&\quad+
\cosh(a_2-a_3)\,\cosh(a_3-a_1)\,\omega^0(v_1,v_2)\\
&\quad+
\cosh(a_3-a_1)\,\cosh(a_1-a_2)\,\omega^0(v_2,v_3)\;.
\end{align*}
By identification, we thereby obtain:

\begin{prop}
For  $g_1,g_2,g_3\in\S$ and $\theta\in\R^*$, we have
$$
K_{\theta, \bm_0}(g_1,g_2,g_3)=K_{\theta,0}(g_1^{-1}g_2,g_1^{-1}g_3)\;,
$$
 where the three-point kernel  on the left hand side of the above equality is given in 
 Proposition \ref{3PK} and the two-point kernel on the right hand side is given in Theorem 
 \ref{explicit} for $\tau=0$. In particular, the products $\star_{\theta,\bm_0}$ and $\star_{\theta,0}$
 coincide on $L^2(\S,{\rm d}_\S)$.
\end{prop}

 Before giving the link between  the generic  kernels $K_{\theta, \bm}$
and $K_{\theta,\tau}$, hence a fortiori between the generic
 products $\star_{\theta,\bm}$ and $\star_{\theta,\tau}$,
we will give the relation between our quantization map and Weyl's. Denote by $\Omega^0$, the Weyl quantization map of $\S$ in the  Darboux chart
\eqref{chartS}. For  a function $f$ on $\S$, $\Omega^0(f)$ is an operator on 
$L^2(\R^{d+1})\simeq L^2(Q,{\rm d}_Q)$ given   (up to a normalization constant) by
\begin{align*}
&\Omega^0(f)\psi(a_0,n_0):=\\
&C(\theta,d)\int_{\R^{2d+2}} e^{-\frac{i}\theta(2(a_0-a)t+\omega^0(n_0-n,m))}\,
 f\big(\tfrac{a+a_0}2,\tfrac{n+n_0}2,m,t\big)\,\psi(a,n)\,{\rm d}a\,{\rm d}n\,{\rm d}m\,{\rm d}t\;.
\end{align*}
Then, recall that for $\tau\in{\bf\Theta}$ (see Definition \ref{gros-gras}), the inverse 
$T_{\theta,\tau}^{-1}$ of the map \eqref{eq:T_tau} is continuous on  the `flat' Schwartz 
space\footnote{By this we mean the ordinary Schwartz space in the global chart \eqref{chartS}.}
$\CS(\S)$. As the Weyl quantization maps continuously Schwartz functions to trace-class operators,
we deduce  that $\Omega^0\circ T_{\theta,{\tau}}^{-1}$ is well defined and continuous from 
$\CS(\S)$ to ${\mathcal L}^1\big(L^2(Q,{\rm d}_Q)\big)$. From this, we get
\begin{prop}
\label{Relation-Weyl} 
Let $\tau\in\Theta$ (see Definition \ref{gros-gras}). Then, as
continuous operators from\footnote{Observe that 
$\CS^{S_{\rm can}}(\S)\subset \CS(\S)$.} $\CS(\S)$
 to the trace ideal ${\mathcal L}^1\big(L^2(Q,{\rm d}_Q)\big)$, we have 
\begin{equation}
\label{transport-Om}
\Omega^0\circ T_{\theta,\tau}^{-1}=\Omega_{\theta,\bm}\,,
\quad\mbox{where}\quad  {\bf m}(a,n)=\bm_0(a,n)\,\exp\big\{{\tau\big(\tfrac2\theta\sinh2a\big)}\big\}\;.
\end{equation}
\end{prop}
\begin{proof}
By density, it suffices to show that for $f\in \CS(\S)$, the operators  $\Omega_{\theta,\bm}(f)$ and 
$\Omega^0\big(T_{\theta,\tau}^*(f)\big)$
coincide on $\CD(Q)$.
Note then that for $f\in\CS(\S)$, we have
\begin{align*}
&T_{\theta,\tau}^{-1}(f)(a,v,t)
=\\&
2\pi\int_{\R^2} \frac{\cosh^{1/2}\big(\tfrac{\theta\xi}2\big)}{\cosh^d\big(\tfrac{\theta\xi}4\big)}\,
\,e^{\tau\big(\tfrac{2}{\theta}\sinh\big(\tfrac{\theta\xi}{2}\big) \big)}
f\Big(a,\sech\big(\tfrac{\theta\xi}{4}\big)v,t'\Big)\,
e^{i\xi t-i\tfrac2\theta\sinh\big(\tfrac{\theta\xi}2\big)t'}\,{\rm d}\xi\,{\rm dt'}\;.
\end{align*}
Hence, for $\psi\in\CD(Q)$, we get
\begin{align*}
&\Omega^0\big(T_{\theta,\tau}^{-1}(f)\big)\psi(a_0,n_0)
=\\
&C_1(d)\int_{\R^{2d+4}} 
 \frac{\cosh^{1/2}\big(\tfrac{\theta\xi}2\big)}{\cosh^d\big(\tfrac{\theta\xi}4\big)}
 \,e^{\tau\big(\tfrac{2}{\theta}\sinh\big(\tfrac{\theta\xi}{2}\big) \big)}
  \,e^{-\frac{i}\theta(2(a_0-a)t+\omega^0(n_0-n,m))}\,
  e^{i\xi t-i\tfrac2\theta\sinh\big(\tfrac{\theta\xi}2\big)t'}\\
&\times f\Big(\tfrac{a+a_0}2,\sech\big(\tfrac{\theta\xi}{4}\big)(n+n_0)/2,
\sech\big(\tfrac{\theta\xi}{4}\big)m,t'\Big)
\,\psi(a,n)\,{\rm d}\xi\,{\rm dt'}\,{\rm d}a\,{\rm d}n\,{\rm d}m\,{\rm d}t\;.
\end{align*}
Performing the change of variable $a\mapsto2a-a_0$, 
$n\mapsto 2\cosh(\tfrac{\theta\xi}4)n-n_0$ and
$m\mapsto \cosh(\tfrac{\theta\xi}4)m$, we get
\begin{align*}
&\Omega^0\big(T_{\theta,\tau}^{-1}(f)\big)\psi(a_0,n_0)
=\\&
C_2(d)\int_{\R^{2d+4}} 
 \cosh^{1/2}\big(\tfrac{\theta\xi}2\big)\cosh^d\big(\tfrac{\theta\xi}4\big)
 \,e^{\tau\big(\tfrac{2}{\theta}\sinh\big(\tfrac{\theta\xi}{2}\big) \big)}
  \,f(a,n,m,t')\,e^{i\xi t-i\tfrac2\theta\sinh\big(\tfrac{\theta\xi}2\big)t'}\\
 &
\times e^{-\frac{2i}\theta\big(2(a_0-a)t+\omega^0\big(n_0-\cosh\big(\tfrac{\theta\xi}
4\big)n,\cosh\big(\tfrac{\theta\xi}4\big)m\big)\big)}\
\,\psi\big(2a-a_0,2\cosh(\tfrac{\theta\xi}4)n-n_0\big)\\
&\times{\rm d}\xi\,{\rm dt'}\,{\rm d}a\,{\rm d}n\,{\rm d}m\,
{\rm d}t\;.
\end{align*}
Integrating out the $t$-variable yields a factor $\delta\big(\xi-\tfrac4\theta(a_0-a)\big)$ and thus we get
\begin{align*}
&\Omega^0\big(T_{\theta,\tau}^{-1}(f)\big)\psi(a_0,n_0)
=\\
&C_3(d)\int_{\R^{2d+2}}  \,
 \cosh^{1/2}2(a-a_0)\cosh^d(a-a_0)\,
 \,e^{\tau\big(\tfrac{2}{\theta}\sinh\big(2(a-a_0)\big) \big)}
 \,f(a,n,m,t')\\
 &
\times e^{\tfrac{2i}\theta\big(\sinh\big(2(a-a_0)\big)t'+\omega^0\big(\cosh(a-a_0)n-n_0,\cosh(a-a_0)m\big)\big)}\\
&\times\psi\big(2a-a_0,2\cosh(a-a_0)n-n_0\big)\,{\rm d}a\,{\rm d}n\,{\rm d}m\,{\rm d}t'\;,
\end{align*}
which by Proposition \ref{Omega-explicit} coincides with  $\Omega_{\theta,\bm}(f)\psi$.
\end{proof}

From the above result  and the defining relation \eqref{SP} for the product $\star_{\theta,\tau}$
for a generic element $\tau\in{\bf\Theta}$, we then deduce:
\begin{prop}
\label{ouf}
To every $\tau\in{\bf\Theta}$, associate a right-$N$-invariant\footnote{
See Remark \ref{N} for the definition of the abelian subgroup $N$.} function $\bm$ on $Q$ as in 
\eqref{transport-Om}. Then, the three point kernel K$_{\theta, \bm}$ of the product 
$\star_{\theta,\bm}$ 
defined in Proposition \ref{generic-prod} (ii), is related to the two-point kernel $K_{\theta,\tau}$
given in  Theorem  \ref{explicit} via:
$$
K_{\theta, \bm}(g_1,g_2,g_3)=K_{\theta,\tau}(g_1^{-1}g_2,g_1^{-1}g_2)\,,\quad \forall 
g_1,g_2,g_3\in\S\;.
$$
 In particular, the product $\star_{\theta,\bm}$ is well defined on $\CB(\S)$  and 
 coincide with $\star_{\theta,\tau}$.
\end{prop}

\begin{rmk}
\label{indicate}
From now on, to indicate that a right-$N$-invariant
borelian  function $\bm$ on $Q$  is associated to an element $\tau\in{\bf\Theta}$, 
as in \eqref{transport-Om}, we just write $\bm\in{\bf\Theta}(\S)$.
\end{rmk}

\begin{rmk}
Observe that, considering the `medial triangle'  three-point function $\Phi_\S$  given in Proposition
\ref{LB} (ii), and writing $K_{\theta,\bf m}=A_{\bf m} \,e^{\tfrac{2i}\theta S}$, for $\bm$ a
right-$N$-invariant function on $Q$, we have:
\begin{align*}
 A_{\bf m}(x_1,x_2,x_3)&={\bf m}(a_1-a_2)\,\frac{{\bf m}_0^2(a_2-a_3)}
{\overline{\bf m}(a_2-a_3)}\,{\bf m}(a_3-a_1)\\
&=\big|\mbox{\rm Jac}_{\Phi_\S^{-1}}(x_1,x_2,x_3)\big|^{1/2}
\frac{{\bf m}}{\bm_0}(a_1-a_2)\,\frac{{\bf m}_0}
{\overline{\bf m}}(a_2-a_3)\,\frac{{\bf m}}{\bm_0}(a_3-a_1)\;.
\end{align*}
The expression above for the amplitude (in the case where $\bm$ is right-$N$-invariant), could also
be derived from the explicit expression  for 
 the Berezin transform (see Definition \ref{Berezin}). Indeed, in the 
present  situation, its distributional kernel in coordinates \eqref{chartS}, reads:
 \begin{align*}
&B_{\theta,\bm}[x_1,x_2]\\
&=\delta(a_1-a_2)\,\delta(n_1-n_2)\,\delta(m_1-m_2)\int_\R
\frac{|\bm|^2}{\bm_0^2}\big(\tfrac12\arcsinh(\tfrac12 a_0)\big)\,e^{\frac{i\theta}2 a_0(t_1-t_2)}
\,{\rm d}a_0\;.
\end{align*}
By standard Fourier-analysis arguments, we deduce:
$$
B_{\theta, \bf m}=\frac{|{\bf m}|^2}{{\bf m_0}^2}
\big(\tfrac12\arcsinh(\tfrac i\theta\partial_t)\big)\;.
$$
\end{rmk}

Finally, let us discuss  the question of the involution for the generic
product $\star_{\theta,\bf m}$. Since in general, the formal adjoint of $\Omega_{\theta,\bm}(x)$
on $L^2(Q,{\rm d}_Q)$ is $\Omega_{\theta,\underline s_e^\star \overline\bm}(x)$,
we deduce
$$
\overline{f_1\star_{\theta,\bf m} f_2}=\overline f_2\star_{\theta,\underline s_e^\star \overline\bm} 
\overline{f_1}\;.
$$
Hence, we obtain that the natural involution for the product $\star_{\theta,\bf m}$ is
\begin{equation}
\label{invol}
\ast_{\theta,\bm}:f\mapsto \frac{{\bf m}}{\underline s_e^\star \overline\bm }
\big(\tfrac12\arcsinh(\tfrac i\theta\partial_t)\big) \overline f\;.
\end{equation}
\begin{rmk}
From Theorem \ref{unitarity} and the previous expression for the involution, we observe
that the element $\bm_0$ defined in \eqref{m0} is uniquely determined by the requirement
that the associated quantization map is both unitary and involution preserving.
\end{rmk}
\section{Quantization of  normal $\bf j$-groups}
\label{QNG}
 Consider now a normal $\bf j$-group $\B=(\S_N\ltimes_{\bR^{N-1}}\dots)\ltimes_{\bR^1}\S_{1}$ 
 with associated extension morphisms 
 \begin{align}
 \label{pia}
 \bR^j\in{\rm Hom}\big((\S_N\ltimes\dots)\ltimes\S_{j+1},{\rm Sp}(V_j,\omega^0_j)\big)\,,
 \qquad j=1,\dots,N-1\;,
 \end{align}
  as in \eqref{AES}. 
 We wish to apply Proposition \ref{EXTPROP} to this situation. For this, 
 recall that for $\S$ an elementary normal $\bf j$-group viewed as an elementary symplectic 
 symmetric space, $\fD$ denotes the Lie algebra of  $W$-preserving symplectic endomorphisms of 
 $\p$ where the Lagrangian subspace $W$ has been chosen to be $\l_-\oplus\R E_-
 \simeq \l\oplus\R E$.

 \begin{prop}
 \label{ouaip}
Denote by $\fD_0$ the stabilizer Lie subalgebra in $\fsp(V,\omega^0)$ of the Lagrangian subspace $\l$.
\begin{enumerate}
\item[(i)] Let $\mbox{\rm Sym}(\l)$ be the space of endomorphisms of $\l$ that are symmetric 
with respect to a given Euclidean scalar product on $\l$. Let also
 $$
 \eta:\End(\l)\times{\rm Sym}(\l)\to{\rm Sym}(\l)\,,\quad (T,S)\mapsto T\circ S+S\circ T^t\;.
 $$
  Then, endowing $\mbox{\rm Sym}(\l)$
with the structure of an Abelian Lie algebra, one has the isomorphism:
$$
\fD_0\simeq\End(\l)\ltimes_\eta\mbox{\rm Sym}(\l)\;.
$$
\item[(ii)] The Lie algebra $\fD_0$ contains an Iwasawa component of $\fsp(V,\omega^0)$.
\item[(iii)] Letting $\fD_0$ trivially act on the central element $E$ of $\h$
induces an isomorphism: 
$$
\fD\simeq \fD_0\ltimes\h\;.
$$
\end{enumerate}
\end{prop}
\begin{proof}
Item (i) is immediate from an investigation at the matrix form level. 

Item (iii) follows from the fact that
the derivation algebra of the non-exact polarized
transvection symplectic triple  $(\g,\sigma,\varpi)$  underlying 
an elementary normal ${\bf j}$-group $\S$ admits the symplectic Lie algebra $\fsp(V,\omega^0)$ as Levi-factor \cite{Bi07}.

Item (ii) follows from a
dimensional argument combined with Borel's conjugacy Theorem of maximal solvable subgroups in complex simple Lie groups. 
Indeed, on the first hand, the dimension of the Iwasawa factor of $\fsp(V,\omega^0)$ equals 
$\dim\fsp(V,\omega^0)-\dim\u(d)$
that is $2d+\frac{2d(2d-1)}{2}-d^2=d(d+1)$ with $2d:=\dim V=2\dim\l$. On the other hand, the dimension of the Borel factor in $\End(\l)$
equals $d+\frac{d(d-1)}{2}$ which equals $\dim\mbox{\rm Sym}(\l)$. Hence $\fD_0$ contains a maximal solvable Lie subalgebra of dimension $2(d+\frac{d(d-1)}{2})=d(d+1)$.
Borel's Theorem then yields the assertion since $\fsp(V,\omega^0)$ is totally split. \end{proof}

 From \cite{Bi07}, we observe that the full polarization quadruple of $\S$ underlies the Lie group 
$\L=\mbox{\rm Sp}(V,\omega^0)\ltimes\tilde{G}$. Hence:

\begin{cor}
Let $\S$ be an elementary normal $\bf j$-group viewed as an elementary symplectic symmetric
space (see Definition \ref{ELEMENTARY}) and let $(\L,\sigma_\L,\xi,\bB)$ be the associated  full 
polarization quadruple  (see Definition \ref{TRANSFULL}). Then we have the global decomposition
$\L=Q\bB$ and moreover $\Delta_\L\big|_{\bB}=\Delta_\bB$.
\end{cor}

We can now prove the conditions needed to apply
Proposition \ref{EXTPROP}:
\begin{prop}
\label{oulala}
Let $\B=(\S_N\ltimes_{\bR^{N-1}}\dots)\ltimes_{\bR^1}\S_{1}$ be a normal $\bf j$-group,
 to which one associates the full 
polarization quadruples $(\L_j,\sigma_{\L_j},\xi_j,\bB_j)$, $j=1,\dots,N$,   of the $\S_j$'s. Then, there exists an homomorphism $\rho_j:(\S_N\ltimes\dots)\ltimes\S_{j}\to\tilde\bD_{j-1}$
whose image  normalizes $\S_{j-1}$ in $\L_{j-1}$ and such that  the extension
homomorphism $\bR_j$ constructed in \eqref{hom}, coincides
with the extension
homomorphism $\bR^j$ underlying the Pyatetskii-Shapiro's decomposition     \eqref{pia}.
\end{prop}
\begin{proof}
Firstly, by Pyatetskii-Shapiro's theory \cite{PS}, one knows that the action of  $(\S_N\ltimes\dots)\ltimes\S_{j}$ on
 $\S_{j-1}$ factors through a 
solvable subgroup of ${\rm Sp}(V_{j-1},\omega_{j-1}^0)$. Setting $(AN)_{j-1}$ the Iwasawa factor of
${\rm Sp}(V_{j-1},\omega_{j-1}^0)$, we thus get an homomorphism:
$$
\tilde\rho_j:(\S_N\ltimes\dots)\ltimes\S_{j}\to (AN)_{j-1}\;.
$$
But Proposition \ref{ouaip} (ii) asserts that $(AN)_{j-1}$ is a subgroup of $\exp\{\fD_{0,j-1}\}$,
where $\fD_{0,j-1}$ is the stabilizer Lie subalgebra in $\fsp(V_{j-1},\omega_{j-1}^0)$ of the 
Lagrangian subspace $\l_{j-1}$.   
Combining this   with the isomorphism of Proposition \ref{ouaip} (iii),  yields another homomorphism:
$$
\hat\rho_j:(AN)_{j-1}\to \bD_{j-1}\subset \tilde\bD_{j-1}\;.
$$
Hence $\rho_j:=\hat\rho_j\circ \tilde\rho_j$ is the desired homomorphism.
Now, observe that by  \cite[ Proposition 2.2 item (i)]{Bi07},  the group $\S_{j-1}$ viewed as a 
subgroup of $\L_{j-1}$,  is normalized by $\mbox{\rm Sp}(V_{j-1},\omega_{j-1}^0)$ for the action 
given in Proposition  \ref{proprietes} (iv) and that this action is precisely the one  associated with the extension homomorphism $\bR^j$ in the decomposition \eqref{pia}.
Thus, all that remains to do is to prove that the extension homomorphisms $\bR^j$ and $\bR_j$
coincide. Here,  $\bR_j(g):=\bC_{\rho_j(g)}\in\Aut(\S_{j-1})$, $g\in\S_j$, is the extension 
homomorphism constructed in \eqref{hom}. But that $\bR^j=\bR_j$ follows from a very general
fact about homogeneous spaces.  Namely, observe that if $M=G/K$, then action of the isotropy 
$K\times M\to M$, $(k,gK)\to kg K$, lifts to $G$ as the restriction to $K$ of the conjugacy action
$K\times G\to G$, $(k,g)\mapsto kgk^{-1}$ (indeed: $kgK=kgk^{-1}K$). 
\end{proof}

From this, we deduce that Proposition \ref{EXTPROP} and Theorem \ref{uiop}
are valid in the case of a normal $\bf j$-group. Moreover, we also deduce that the associated
product $\star_{\theta,\bm}$ coincides with $\star_{\theta,\vec \tau}$ of Proposition \ref{efF}:
\begin{prop}
\label{ouf2}
Let $\B$ be a normal $\bf j$-group. 
To every $\vec\tau\in{\bf\Theta}^N$, we associate a function $\bm$ on $Q_N\times\dots Q_1$  by $\bm=\bm_N\otimes
\dots\otimes\bm_1$ where $\bm_j$ is related to $\tau_j$ as in 
\eqref{transport-Om}. Then, the three-point kernel K$_{\theta, \bm}$ of the product $\star_{\theta,\bm}$ 
defined in Theorem  \ref{uiop} (iii), is related with the two-point kernel $K_{\theta,\tau}$
given in  Proposition  \ref{efF} \eqref{2P-ext} via:
$$
K_{\theta, \bm}(g_1,g_2,g_3)=K_{\theta,\vec \tau}(g_1^{-1}g_2,g_1^{-1}g_2)\,,\quad \forall 
g_1,g_2,g_3\in\B\;.
$$
 In particular, the product $\star_{\theta,\bm}$ is well defined on $\CB(\B)$  and 
 coincide with $\star_{\theta,\vec\tau}$.
\end{prop}
\begin{rmk}
\label{indicateB}
To indicate that a  function $\bm=\bm_N\otimes
\dots\otimes\bm_1$ on $Q_N\times\dots Q_1$ is related to elements $\tau_j\in\bf\Theta$ as in 
\eqref{transport-Om}, we just write $\bm\in{\bf\Theta}(\B)$.
\end{rmk}

We also quote the following extension of Proposition \ref{Relation-Weyl}:
\begin{prop}
\label{Relation-Weyl2} 
Let $\B$ be a normal $\bf j$-group. For any  $\bm\in{\bf\Theta}(\B)$,
 the quantization map $\Omega_{\theta,\bm}$ is a
continuous operator from $\CS(\B):=\CS(\S_N)\otimes\dots\otimes\CS(\S_1)$
 to ${\mathcal L}^1\big(L^2(Q_N,{\rm d}_{Q_N})\otimes\dots
\otimes L^2(Q_1,{\rm d}_{Q_1})\big)$.
\end{prop}

Our last result concerns the representation homomorphism \eqref{rephom}.
\begin{lem}
\label{OM}
Let $\B=\B'\ltimes_\bR\S_1$, $\S_1=Q_1\ltimes\Y_1$.
Then, the restriction of the  homomorphism $
\CR:\B'\to \CU(L^2(Q_1))$
(see \eqref{rephom})  underlying the one
$\rho=:\B'\to \tilde\bD_1$ of Proposition \ref{oulala}, defines a tempered action
(see Definition \ref{temp-action}) of $\B'$ on
$\CS(Q_1)$.
\end{lem}
\begin{proof}
Parametrizing $\S_1\simeq\s_1=\{(a,n,m,t)\}$, $Q_1\simeq q_1=\{(a,n)\}$ as usual, 
the matrix of the Pyatetskii-Shapiro
extension homomorphism $\bR_{g'}$ is expressed under the form
\begin{equation*}
\bR_{g'}=\begin{pmatrix} 0&0&0&0\\
0&\bR_+(g')&0 &0\\0& \bR_-(g')& (\bR_+(g')^T)^{-1}&0\\
0&0&0&0 \end{pmatrix}\;,
\end{equation*}
where (by Borel's conjugacy theorem) the element $\bR_+(g')$ is (conjugated to) an upper triangular matrix form
in $\End(\l^\star)$. The determinant of this element then consists (as a small induction argument shows) in
a product of powers of the modular functions of the elementary factors of $\S_2,...,\S_N$ of $\B'$:
$$
\det(\bR_+)=\Pi_{j=2}^N\Delta_{\S_j}^{n_j}\;.
$$
One deduces  from a careful examination
of the construction of $\rho$ in Proposition \ref{oulala}, that $\CR$
 is explicitly given by:
\begin{equation*}
\CR(g')\varphi(a,n)=\det(\bR_+(g'))^{-\frac12}
\varphi(a,\bR_+(g')^{-1}n)\;.
\end{equation*}
 Temperedeness results in 
a direct computation using that the semi-norms of $\CS(Q_1)$ 
are given by (see the discussion right after Remark \ref{fd-explicit2}):
\begin{align*}
\|f\|_{k,n}:=\sup_{X\in\,\CU_k(\q_1)}\sup_{x\in Q_1}\Big\{\frac
{\fd_{Q_1}(x)^n\big|\widetilde X\,f(x)
\big|_j}{|X|_k}\Big\}\;,
\end{align*}
 that  left invariant vector fields associated with $H$ the 
 generator of $\a$ and $\{f_j\}_{j=1}^d$ be a basis of $\l^\star$,
  read in the coordinates \eqref{chartQ}:
 $$
 \widetilde H =\partial_a-\sum_{j=1}^dn_j\partial_{n_j}\,,\qquad \widetilde f_j=\partial_{n_j}\,,\quad j=1,\dots,d\;,
 $$  
together with the behavior of the modular weight $\fd_{Q_1}$ given in Remark \ref{fd-explicit2}).
\end{proof}

\chapter{Deformation of $C^*$-algebras}\label{SCSTAR}

Throughout this chapter, we consider a $C^*$-algebra $A$, endowed with an 
{\em isometric and strongly
continuous} action $\alpha$ of a normal $\bf j$-group $\B$. 
In chapter \ref{DFA}, we have seen how to deform the
 {\em Fr\'echet algebra} $A^\infty$, consisting of smooth vectors for the action $\alpha$. 
 Our goal here is  to construct a  $C^*$-norm on
  $(A^\infty,\star_{\theta,\bm}^\alpha)$, in order to get, after completion, a deformation theory at the
  $C^*$-level. We stress that from now on, the isometricity
assumption of the action is fundamental. The way we will define this $C^*$-norm is based on
 the pseudo-differential calculus introduced in the previous two chapters. 

\quad

The basic ideas of the construction can be summarized
as follow. Consider $\CH$, a separable Hilbert space carrying a faithful representation of $A$. We will thereof 
identify $A$ with its image in  $\CB(\CH)$. Let $\CS^{S_{\rm can}^\B}(\B,A)$ 
be the $A$-valued one-point Schwartz space associated to the tempered pair 
$(\B\times\B,S_{\rm can}^\B)$ as given in Definition  \ref{SCHH}. Since 
this space is a subset of the
 flat $A$-valued Schwartz space of $\B$, Proposition \ref{Relation-Weyl2}
shows that for every $\bm\in{\bf\Theta}(\B)$, the map
\begin{equation}
\label{TheMap}
f\in \CS^{S_{\rm can}^\B}(\B,A)\mapsto\Omega_{\theta,\bm}^\B(f):=\int_\B  \Omega_{\theta,\bm}^\B(x)
\otimes f(x)\,
{\rm d}_\B(x)\;,
\end{equation}
is well defined and takes values\footnote{Given a Hilbert space $\CH$, $\CK(\CH)$ denotes the  $C^*$-algebra of compact operators.} in 
$\CK\big(\CH_\chi\big)\otimes A$, where $\CH_\chi:=L^2(Q_1)\otimes\dots\otimes 
L^2(Q_N)$. 
 Then, the main step is to extend the map \eqref{TheMap}, from   $\CS^{S_{\rm can}^\B}(\B,A)$ to 
 $\CB(\B,A)$. As for 
 $a\in A^\infty$, the $A$-valued function $\alpha(a):=[g\in\B\mapsto\alpha_g(a)\in A]$ belongs to 
 $\CB(\B,A)$, we will
 define a new norm on $A^\infty$ by setting
 \begin{equation*}
 \|a\|_{\theta,\bm}:=\big\|\Omega_{\theta,\bm}^\B\big(\alpha(a)\big)\big\|\;,
  \end{equation*}
  where the norm on the right hand side above, denotes the  $C^*$-norm of
  $\CB\big(\CH_\chi\otimes\CH\big)$.
  This will eventually be achieved by proving a non-Abelian $C^*$-valued version of the Calder\'on-Vaillancourt 
  Theorem
 in the context of the present pseudo-differential calculus. 
 This theorem will be proved using wavelet analysis
 and oscillatory integrals  methods.

\section{Wavelet analysis}
\label{WA}
Let $\B$ be a normal $\bf j$-group, with Pyatetskii-Shapiro decomposition  
$$
\B=(\S_N\ltimes\dots)\ltimes\S_{1}\;,
$$
where the $\S_j$'s, $j=1,\dots,N$, are elementary normal $\bf j$-groups. 
Recall that our choice of  parametrization  is:
$$
\S_N\times\dots\times\S_1\to\B\,,\qquad (g_N,\dots,g_1)\mapsto g_1\dots g_N\;.
$$

\begin{rmk}
\label{invariant-parametrization}
Observe that the extension 
homomorphism at each step, $\bR^j$, being valued in ${\rm Sp}(V_j,\omega^0_j)$,  it preserves
any left invariant Haar measure ${\rm d}_{\S_j}$ on $\S_j$:
\begin{equation}
\label{FG}
(\bR^j_{g'})^\star {\rm d}_{\S_j}={\rm d}_{\S_j}\,,\quad\forall g'\in (\S_N\ltimes\dots)\ltimes\S_{j-1}\;.
\end{equation}
This implies that the product of left invariant
Haar measures ${\rm d}_{\S_1}\otimes\dots\otimes{\rm d}_{\S_N}$ defines a left invariant Haar
measure on $\B$ under both parametrizations $g=g_1\dots g_N$ or $g=g_N\dots g_1$ of 
$g\in\B$. 
\end{rmk}

The aim of this section is to construct a weak resolution of the identity on the tensor
product  Hilbert space
$L^2(Q_N)\otimes\dots
\otimes L^2(Q_1):=\CH_\chi$
from  a suitable family of coherent states for  $\B$, that we now introduce.

\begin{dfn}
\label{CoheState}
Let $\B$ be a normal $\bf j$-group.
Given a mother wavelet $\eta\in\CD(Q_N\times\dots\times Q_1)$, let
 $\{\eta_x\}_{x\in \B}$ be the family of coherent states defined by
$$
\eta_x:=U_\theta(x)\eta\,,
\quad x\in\B\;,
$$
where $U_\theta$ is the unitary representation of $\B$ on 
$\CH_\chi$ constructed  in Lemma \ref{OBS} for the morphism
underlying Proposition \ref{oulala}.
\end{dfn}
Observe that in the elementary case, we have:
\begin{align}
\label{bE02}
\eta_x(q_0)=\bE^0_\theta(q^{-1}q_0b)\,\eta(q^{-1}q_0)\,,\quad x=qb\in\S\,,\quad q_0\in Q\;,
\end{align}
where the  phase $\bE^0_\theta$ is defined by
\begin{align}
\label{bE0}
\bE^0_\theta(x):=\chi_\theta\big(\bC_{q^{-1}}(b)\big)\,,\quad x=qb\in\S\;.
\end{align}
In the generic case, setting $\B=\B'\ltimes_\bR\S_1$ with $\B'$ a normal $\bf j$-group
and $\S_1$ an elementary normal $\bf j$-group,  for $\eta=\eta'\otimes\eta^1$,
$\eta'\in\CD(Q_N\times\dots\times Q_2)$, $\eta^1\in\CD(Q_1)$ 
and parametrizing $g\in\B$ as $g=g'.g_1$, $g'\in\B'$,
$g_1\in\S_1$, we have (see Lemma \ref{OBS}):
\begin{align}
\label{bE01}
\eta_g=\eta'_{g'}\otimes\CR(g')\eta^1_{g_1}\;,
\end{align}
where  $\CR:\B'\to U(L^2(Q_1))$ is the homomorphism underlying \eqref{rephom}.
Parametrizing now  $g\in\B$ as $g=g_1.g'$, $g'\in\B'$,
$g_1\in\S_1$, we find
\begin{align}
\label{bE03}
\eta_g=\eta'_{g'}\otimes\big( \CR(g')\eta^1\big)_{g_1}\;,
\end{align}
 
\begin{prop}
\label{first-type}
Let $\B$ be a normal $\bf j$-group, $\CE$ a complex Fr\'echet space  and  
 $\eta\in\CD(Q_N\times\dots\times Q_1)$. Then,
 the map
 \begin{align*}
& \CF^\eta: \,L^\infty(Q_N\times\dots\times Q_1,\CE)\to L^\infty(\B,\CE)\;,\\
 & f\mapsto 
 \Big[x\in\B\mapsto\int_{Q_1\times\dots\times Q_N} f(q_N,\dots,q_1)\,
 \eta_x(q_N,\dots,q_1){\rm d}_{Q_N}(q_N)\dots{\rm d}_{Q_1}(q_1) \in\CE\Big]\;,
 \end{align*}
restricts as	a continuous map $\CF^\eta:
\CS(Q_N\times\dots\times Q_1,\CE)\to \CS^{S_{\rm can}}(\B,\CE)$. 
\end{prop}
\begin{proof}
Assume first  that $\B=\S$ is an elementary normal $\bf j$-group.   
We denote by $\bE^0_\theta$ the element of $C^\infty(\S)$
given in  \eqref{bE0}, so that with $x=qb\in\S$, $q\in Q$, $b\in\Y$, we have
for every  $f\in\CS(Q,\CE)$:
$$
 \big(\CF^\eta f\big)(x)=\int_Q f(q_0)\bE^0_\theta(q^{-1}q_0b)\eta(q^{-1}q_0){\rm d}_Q(q_0)
 =\int_Q f(qq_0)\bE^0_\theta(q_0b)\eta(q_0)\,{\rm d}_Q(q_0)\;.
$$
  Decomposing as usual $\q=\a\oplus\n$, we let $H$ be the 
 generator of $\a$ and $\{f_j\}_{j=1}^d$ be a basis of $\n$. Then, from the expressions
 given in \eqref{INVVECT}, we see that the associated left invariant vector fields 
  read in the coordinates \eqref{chartQ}:

 $$
 \widetilde H =\partial_a-\sum_{j=1}^dn_j\partial_{n_j}\,,\qquad \widetilde f_j=\partial_{n_j}\,,\quad j=1,\dots,d\;.
 $$  
 Moreover, in the chart \eqref{chartS} of $\S$, with $x=(a,n\oplus m,t)$, the function $\bE^0_\theta$ 
 takes the following form:
 $$
 \bE^0_\theta(x)=\exp\big\{\tfrac i\theta\big(e^{-2a}t-e^{-a}\omega^0(n,m)\big)\big\}\;. 
  $$

Hence, defining
 $$
i\widetilde H  \,\bE^0_\theta=: \alpha \, \bE^0_\theta\quad\mbox{and}\quad -\sum_{j=1}^d\widetilde f_j^2\, \bE^0_\theta=:\beta \, \bE^0_\theta\;,
$$
a simple computation gives
$$
\alpha(x)=\tfrac {2}\theta(e^{-2a}t-e^{-a}\omega^0(n,m))\,,\qquad \beta(x)=\theta^{-2}e^{-2a}
|m|^2\;,
$$
where $|m|^2=\sum_{j=1}^d\omega^0(f_j,m)^2$.  Moreover, it is easy to see that both
 $\alpha$ and $\beta$ are eigenvectors of $\widetilde H$ with eigenvalue $-2$ 
 and that $\widetilde f_j \beta=0$.
 Hence setting $\widetilde P:= 1-\sum_{j=1}^d\widetilde f_j^2$, we get by  integration by parts on the 
 $q_0$-variable and with $k, k'\in\N$ arbitrary:
 \begin{align*}
& \big(\CF^\eta f\big)(x)=\\&\int_Q\bE^0_\theta(q_0b)(1-\widetilde  H_{q_0}^{2k'})
 \left[\frac{\widetilde 
 P^k_{q_0}\big[f(qq_0)\,\eta(q_0)\big] }{\big(1+ \alpha(q_0b)^{2k'}+\alpha^{(k')}(q_0b)\big)\big(1+\beta(q_0b)\big)^k}\right]\,
 {\rm d}_Q(q_0)
\end{align*}
\noindent where
$$
\alpha^{(k')}\;:=\;\sum_{r=1}^{2k'-1}c^{k'}_r\,\alpha^r\quad(c^{k'}_r\in\C)\;.
$$
 This easily entails that
\begin{align*}
 \big\|\big(\CF^\eta f\big)(x)\big\|_j\leq C(k, k')\, \int_Q\frac{\|(1-\widetilde H_{q_0}^{2k'})
 \widetilde P^k_{q_0}\left[f(qq_0)\,
\eta(q_0)\right]\|_j}{\big(1+ \alpha(q_0b)^{2k'}\big)\big(1+\beta(q_0b)\big)^k}
\,{\rm d}_Q(q_0)\;.
\end{align*}
By left invariance of $\widetilde H$ and $\widetilde P$, we get up to a redefinition
of $f\in\CS(Q,\CE)$ and of $\eta\in\CD(Q)$:
$$
\| \big(\CF^\eta f\big)(x)\|_j\leq C(k, k')\int_Q \frac{\|f(qq_0)\|_j\,
|\eta(q_0)|}{\big(1+ \alpha(q_0b)^{2k'}\big)\big(1+\beta(q_0b)\big)^k}\,{\rm d}_Q(q_0)\;.
$$
 Now,  given any tempered weight $\mu$ on $\S$, one therefore has
\begin{align*}
&\int_\S\mu(x)\,\|\big(\CF^\eta f\big)(x)\|_j{\rm d}_\S(x)\leq\\& \qquad C(k, k')
\int_{Q^2\times\Y} \frac{\mu(qb)\,\|f(qq_0)\|_j\,
|\eta(q_0)|}{\big(1+\alpha(q_0b)^{2k'}\big)\big(1+\beta(q_0b)\big)^k}\,{\rm d}_Q(q_0)
{\rm d}_Q(q)\,{\rm d}_\Y(b)\;.
\end{align*}
Observing that for every $X\in\CU(\q)$, the element $\widetilde{X}_q[qb\mapsto\bE^0_\theta(qb)]$ only 
depends on the variable $\bC_{q^{-1}}(b)=:\exp\Ad_{q^{-1}}y$, where $b=e^y\in\Y$,
changing the variable following $b_0:=:\exp(y_0):=\bC_{q_0^{-1}}(b)$ ($y_0\in\fY$) yields
\begin{align*}
&\int_{Q^2\times\Y} \frac{\mu(qb)\,\|f(qq_0)\|_j\,
|\eta(q_0)|}{\big(1+\alpha(q_0b)^{2k'}\big)\big(1+\beta(q_0b)\big)^k}\,{\rm d}_Q(q_0)
{\rm d}_Q(q)\,{\rm d}_\Y(b)=\\
&\int_{Q^2\times\Y} \frac{\mu(qq_0b_0q_0^{-1})\,\left|\det\Ad_{q_0}|_\fY\right|\,\|f(qq_0)\|_j\,
|\eta(q_0)|}{\big(1+\alpha(b_0)^{2k'}\big)\big(1+\beta(b_0)\big)^k}\,{\rm d}_Q(q_0)
{\rm d}_Q(q)\,{\rm d}_Q(b_0)\;.
\end{align*}
Writing $y_0=(m_0,t_0)\in\l\times\R$ as before, we observe that $\alpha(b_0)=\frac{2t_0}{\theta}$, $\beta(b_0)=\theta^{-2}|m_0|^2$ and ${\rm d}_\Y(b_0)={\rm d}y_0$.
Therefore this last integral exists as soon as $k$ and $k'$ are large enough in front of the polynomial growth of $\mu$. Since for every 
$Y\in\CU(\s)$, one has $\widetilde{Y}_x[x\mapsto\eta_x]=\left(dU_\theta(Y)\eta\right)_x$, the left-invariant derivatives of $\CF^\eta f$ follow the exact same treatment.

\noindent Passing to non-elementary case, 
we first observe that by Lemma \ref{OM},  the restricted action $\CR$ of the tempered
Lie group $\B'$ on $\CS(Q_1)$ is  tempered in the sense of Definition \ref{temp-action}. Therefore, in the above discussion, replacing $\eta$ by $\CR(g')\eta$ yields constants $C(k,k')$ that now depend polynomially on the tempered weights $\{\mu^{\CR}_{k,j}\in C^\infty(\B')\}$ associated to the tempered action $\CR$ (within the  notations of  Definition  \ref{temp-action}).

\noindent Now, within the conventions of \eqref{bE01} and considering $f:=f'\otimes f_1\in\CS(Q_N\times...\times Q_2,\CE){\otimes}\CS(Q_1)$ and $\eta:=\eta'\otimes\eta^1\in\CD(Q_N\times...\times Q_2){\otimes}\CD(Q_1)$, one has:
$$
\left(\CF^\eta f\right)(g_1g')\;=\;\big(\CF^{\eta'}f'\big)(g')\,\big(\CF^{\CR(g')\eta^1}f_1\big)(g_1)\;.
$$
Let $d\CR$ the infinitesimal form of $\CR$, that is for $\eta^1\in\CD(Q_1)$, it is
defined  by $d\CR(X)\eta^1:=\tfrac d{dt}|_{t=0}\CR(e^{tX})\eta^1$ for 
$X\in\b'$ and extended as an algebra morphism to $\CU(\b')$. Note that for any $X\in\CU(\b')$ and $g'\in\B'$, we have
$$
\widetilde X_{g'}\,\CR(g')\eta^1=\CR(g')\big(d\CR(X)\eta^1\big).
$$
Now, for every tempered weight $\mu=\mu'\otimes\mu_1$  on $\B$ and $X=X'\otimes X_1\in\CU(\b)$ (obvious notations),  denoting
$$
\tilde\eta^1:= d\CR(X'_{(2)})\eta^1\in\CD(Q_N\times\dots\times Q_2)\;,
$$
we then have
\begin{align}
\label{ber}
&\int_\B\mu(g)\,\|\widetilde{X}.\left(\CF^\eta f\right)(g)\|_j{\rm d}_\B(g)\\
\leq\sum_{(X')}\int_{\B'}\mu'(g')&\int_{\S_1}\mu_1(g_1)\,
\|\widetilde{X'_{(1)}}_{g'}\big(\CF^{\eta'}f'(g')\big)\widetilde{X'_{(2)}}_{g'}\widetilde{X_1}_{g_1}\big(\CF^{\CR(g')\eta^1}f_1(g_1)\big)\|_j{\rm d}g_1\,{\rm d}g'\nonumber\\
&\hspace{-2cm}=\sum_{(X')}\int_{\B'}\mu'(g')
\|\widetilde{X'_{(1)}}_{g'}\big(\CF^{\eta'}f'(g')\big)\|_j\int_{\S_1}\mu_1(g_1)\big|
\widetilde{X_1}_{g_1}\big(\CF^{\CR(g')\tilde\eta^1)}f_1(g_1)\big)\big|{\rm d}g_1{\rm d}g'\nonumber
\end{align}
From what we have proven in the elementary case and by induction hypothesis on $N$, the
map
$$
\B'\to\R_+,\qquad g'\mapsto\int_{\S_1}\mu_1(g_1)\big|
\widetilde{X_1}_{g_1}\big(\CF^{\CR(g')\tilde\eta^1)}f_1(g_1)\big)\big|{\rm d}g_1\;,
$$
is temperate. Hence, the integral \eqref{ber}
exists. We conclude by nuclearity of the Schwartz spaces $\CS(Q_N\times\dots\times Q_1,\CE)$.
\end{proof}

We then deduce the following consequence:
\begin{cor}
\label{ONE}
Let $\B$ be a normal $\bf j$-group and  
$\eta\in\CD(Q_N\times\dots\times Q_1)$. Then,
 the maps $[\B\ni x\mapsto \langle\eta,\eta_x\rangle]$ and $[\B\ni x\mapsto \langle\eta,\eta_{x^{-1}}\rangle]$ belong to $L^1(\B)$. 
\end{cor}
\begin{proof}
This follows from Proposition \ref{first-type} with $\CE=\C$ since  $\langle\eta,\eta_{x}\rangle
=\CF^{\overline\eta}\big(\eta\big)(x)$ and $\langle\eta,\eta_{x^{-1}}\rangle
=\overline{\langle\eta,\eta_{x}\rangle}$.
\end{proof}

The next result is probably well known\footnote{It can be viewed as a Banach space 
valued version of 
Schur's test Lemma.}  but since we are unable to locate it
in the literature and since we use it several  times, we deliver a proof.
\begin{lem}
\label{CS}
Let $(X,\mu)$ be a $\sigma$-finite measure
space and $\CH$ a separable Hilbert space.
Consider an element 
$$
K\in L^\infty\big(X\times X,\mu\otimes\mu;\CB(\CH)\big)\;,
$$
such that
$$
c_1^2:=\sup_{x\in X}\int_X\big\|K(x,y)\big\|_{\CB(\CH)}d\mu(y)<\infty\,,\quad
 c_2^2:=\sup_{y\in X}\int_X\big\|K(x,y)\big\|_{\CB(\CH)}d\mu(x)<\infty\;.
$$
Then, the associated kernel operator is bounded on $L^2(X,\mu;\CH)\simeq L^2(X,\mu)\otimes\CH$ with operator norm not exceeding $ c_1\,c_2$.
\end{lem}
\begin{proof}
Let $T_K$ be the operator associated with the kernel  $K$.  
 For vectors  $\Phi,\Psi\in L^2(X,\mu;\CH)\cap L^1(X,\mu;\CH)$, 
 we have $\big|\langle\Phi,T_K\Psi\rangle\big|<\infty$.
Moreover,  the Cauchy-Schwarz inequality gives
\begin{align*}
\big|\langle\Phi,T_K\Psi\rangle\big|
&= \Big|\int_{X\times X}\big\langle\Phi (x) ,K(x,y)\,\Psi(y)\big\rangle_\CH\,d\mu(y)\,d\mu(x)
\Big|\\&\leq
\int_{X\times X}\|\Phi (x)\|_\CH \,\|K(x,y)\|_{\CB(\CH)}\,\|\Psi(y)\|_\CH\,d\mu(y)\,d\mu(x)\\
&\leq\Big(\int_{X\times X}\|\Phi (x)\|_\CK ^2\,\|K(x,y)\|_{\CB(\CH)}\,d\mu(y)\,
d\mu(x)\Big)^{1/2}\\
&\qquad\qquad\times\Big(\int_{X\times X}
\|K(x,y)\|_{\CB(\CH)}\,\|\Psi(y)\|_\CH^2\,d\mu(y)\,d\mu(x)\Big)^{1/2}\\
&\leq c_1\,c_2\,\|\Phi\|_{L^2(X,\mu;\CH)} \,\|\Psi\|_{L^2(X,\mu;\CH)} \;,
\end{align*}
and the claim follows immediately. 
\end{proof}

We are now able to prove that the family of coherent states $\{\eta_x\}_{x\in\B}$ provides a  weak resolution of the identity.
\begin{prop}
\label{WRI}
Let $\B$ be a normal $\bf j$-group, $\CH$ a separable Hilbert space and  
$\eta\in\CD(Q_N\times\dots\times  Q_1)\setminus\{0\}$.
Then, for all   $\Phi,\Psi\in\CH_\chi\otimes\CH$ 
the following relation holds:
\begin{align}
\label{cons}
\langle \Psi,\Phi\rangle_{\CH_\chi\otimes\CH}= C^\B(\eta)^{-1}\int_\B \big\langle\langle\eta_x,\Psi
\rangle_{\CH_\chi},\langle\eta_x,\Phi\rangle_{\CH_\chi}\big\rangle_\CH\;{\rm d}_\B( x)
\;,
\end{align}
 where $\langle \eta_x,\Psi\rangle_{\CH_\chi}$ is the vector in $\CH$ is defined by:
$$
\big\langle\vf,\langle \eta_x,\Psi\rangle_{\CH_\chi}\big\rangle_\CH:=
\big\langle\eta_x\otimes\vf,\Psi\big\rangle_{\CH_\chi\otimes\CH}\,,\quad\forall \vf\in\CH\;,
$$
and 
$$
C^\B(\eta):=(2\pi\theta)^{\dim(\B)/2}\|\Delta_{Q_N\times\dots\times Q_1}\,\eta\|_2^2\;,
$$
where $\Delta_{Q_N\times\dots\times Q_1}=
\Delta_{Q_N}\otimes\dots\otimes\Delta_{ Q_1}$, with $\Delta_{Q_j}$ the modular function of $Q_j$.
\end{prop}
\begin{proof}
We first demonstrate that for $\Phi\in\CH_\chi\otimes\CH$, the map 
$[x\in\B\mapsto\langle\eta_x,\Phi\rangle_{\CH_\chi}\in\CH]$ belongs to $L^2(\B,\CH)$. 
To see this, let $\{\B_j\}_{j\in\N}$ be an increasing sequence of relatively compact subsets of $\B$, 
which converges  to $\B$. For each $j\in\N$, we define the operator 
$$
T^\eta_j: L^2(\B_j,\CH)\to\CH_\chi\otimes\CH\,,\quad F\mapsto\int_{\B_j}\eta_x\otimes F(x)\,{\rm d}_\B(x)\;.
$$
Clearly, each $T^\eta_j$ is bounded:
\begin{align*}
\|T^\eta_j F\|_{\CH_\chi\otimes\CH}&\leq \int_{\B_j}\|\eta_x\|_{\CH_\chi}\|F(x)\|_\CH\,{\rm d}_\B(x)\\
&\leq\|\eta\|_{\CH_\chi}
 {\rm meas}(\B_j)^{1/2}\Big(\int_{\B_j}\|F(x)\|_\CH^2\,{\rm d}_\B(x)\Big)^{1/2}\\&\quad=
\|\eta\|_{\CH_\chi}
 {\rm meas}(\B_j)^{1/2}\|F\|_ {L^2(\B,\CH)}\;.
 \end{align*}
To see that  the family $\{T^\eta_j\}_{j\in\N}$ is in fact uniformly bounded, note that the adjoint of 
$T^\eta_j$ reads:
$$
{T^\eta_j}^*:\CH_\chi\otimes\CH\to L^2(\B_j,\CH)\,,\quad \Phi\mapsto
 \big[x\in\B_j\mapsto\langle\eta_x,\Phi\rangle_{\CH_\chi}\in\CH\big]\;.
$$ 
Hence for $F\in L^2(\B_j,\CH)$ we get
$$
\big|T^\eta_j\big|^2 F(x)=\int_{\B_j}\langle\eta_x,\eta_{y}\rangle\,F(y)\,{\rm d }_\B(y)\;,
$$
that is $|T^\eta_j|^2=S^\eta_j\otimes{\rm Id}_\CH$, where $S^\eta_j\in\CB(L^2(\B_j))$ is a kernel operator with  
kernel $K^\eta_j(x,y)= \langle\eta_x,\eta_{y}\rangle$. Applying Lemma
\ref{CS}, a simple change of variable gives 
$\|S^\eta_j\|\leq\|[x\mapsto  \langle\eta,\eta_{x}\rangle]\|_1:=C$ which is finite by Lemma \ref{ONE} 
and of course, is uniform in $j\in\N$. Finally, since 
$$
\int_{\B_j}\|\langle\eta_x,\Phi\rangle_{\CH_\chi}\|_\CH^2\,{\rm d}_\B(x)
=\|{T^\eta_j}^* \Phi\|_{L^2(\B_j,\CH)}^2\leq C\|\Phi\|_{\CH_\chi\otimes\CH}^2\;,
$$
taking the limit $j\to\infty$ gives
$$
\int_{\B}\|\langle\eta_x,\Phi\rangle_{\CH_\chi}\|_\CH^2\,{\rm d}_\B(x)
\leq C\|\Phi\|_{\CH_\chi\otimes\CH}^2\;,
$$
as needed.
The rest of the proof is computational.
Assume first that $\B=\S$ is elementary. In this case, 
for $\Phi,\Psi\in \CH_\chi\otimes\CH$ and $\eta\in\CD(Q)$, we have in 
chart \eqref{chartS}  and from \eqref{bE02}:
\begin{align*}
& \int_\S\big\langle \langle\eta_x,\Psi\rangle_{\CH_\chi},\langle\eta_x,\Phi\rangle_{\CH_\chi}\big
\rangle_\CH\;{\rm d}_\S(x)
 =\\
 &\int_{\R^{4d+4}}\langle\Psi(a_0,n_0),\Phi(a_1,n_1)\rangle_\CH\,
\eta\big(a_0-a,n_0-e^{a-a_0}n\big)\,\overline \eta\big(a_1-a,n_1-e^{a-a_1}n\big)\\
&\times \exp\Big\{\frac{i}\theta\big( (e^{-2a_0}-e^{-2a_1})e^{2a}\big(t+\omega^0(n,m)\big)-e^{a}
\omega^0(n_0e^{-a_0}-n_1e^{-a_1},m)\big)\Big\}\\&\times
\,{\rm d}a\,{\rm d}n\,{\rm d}m\,{\rm d}t\,{\rm d}a_0\,{\rm d}n_0
\,{\rm d}a_1\,{\rm d}n_1\;.
\end{align*}
Integrating out the $t$-variable yields a factor $2\pi\theta\,e^{-2a+2a_0}\delta(a_0-a_1)$, 
the former expression then
becomes
\begin{align*}
& 2\pi\theta\int_{\R^{4d+2}}\langle\Psi(a_0,n_0),\Phi(a_0,n_0)\rangle_\CH\,
\eta\big(a_0-a,n_0-e^{a-a_0}n\big)\,\overline \eta\big(a_0-a,n_1-e^{a-a_0}n\big)\\
&\qquad\qquad\qquad\times
 e^{-2a+2a_0} \exp\Big\{-\frac{i}\theta\big( e^{a-a_0}\omega^0(n_0-n_1,m)\big)\Big\}
 \,{\rm d}a\,{\rm d}n\,{\rm d}m\,{\rm d}a_0\,{\rm d}n_0\,{\rm d}n_1\;.
\end{align*}
Integrating  the $m$-variables, yields a factor $(2\pi\theta)^de^{-da+da_0}\delta(n_0-n_1)$ and we get,
up to a constant:
\begin{align*}
&\int_{\R^{2d+2}}\langle\Psi(a_0,n_0),\Phi(a_0,n_1)\rangle_\CH\,
\eta\big(a_0-a,n_0-e^{a-a_0}n\big)\,\overline \eta\big(a_0-a,n_0-e^{a-a_0}n\big)\\
&\qquad\qquad\qquad\qquad\qquad\qquad
\qquad\qquad\qquad\times e^{-(d+2)(a-a_0)}\,{\rm d}a\,{\rm d}n\,{\rm d}a_0\,{\rm d}n_0\;,
\end{align*}
which after an affine change of variable and restoring the constants, gives
\begin{align*}
(2\pi\theta)^{d+1}\langle\Psi,\Phi\rangle_{\CH_\chi\otimes\CH}\;\int_{\R^{d+1}}
|\eta(a,n)|^2\, e^{(2d+2)a}\,{\rm d}a\,{\rm d}n\;,
\end{align*}
which is all we needed  since $\Delta_Q(a,n)=e^{(d+1)a}$.

\noindent The case of a generic normal $\bf j$-group $\B$ is  treated with the same argument lines
as in Proposition \ref{first-type}: set $\B'\ltimes\S_1$, with $\S_1$ elementary normal and assume
that the relation \eqref{cons} holds for $\B'$.
With the notations of \eqref{bE01}, we have 
 for all $\varphi=\varphi'\otimes\varphi^1,\psi=\psi'\otimes\psi^1\in\CH_\chi$ and omitting
 the constant $C^\B(\eta)^{-1}=C^{\B'}(\eta')^{-1}C^{\S_1}(\eta^1)^{-1}$):
\begin{align*}
&\int_\B\langle \varphi,\eta_x\rangle\langle\eta_x,\psi\rangle{\rm d}_\B  (x)\\&
=\int_{\B'\times\S_1}
\langle\varphi',\eta'_{g'}\rangle\langle\varphi^1,\CR(g')\eta^1_{g_1}\rangle
\langle\CR(g')\eta^1_{g_1},\psi^1\rangle
\langle\eta'_{g'},\psi'\rangle{\rm d}_{\B'}(g')\,{\rm d}_{\S_1}(g_1)\\
&=\int_{\B'}
\langle\varphi',\eta'_{g'}\rangle\Big(\int_{\S_1} 
\langle\CR({g'}^{-1})\varphi^1,\eta^1_{g_1}\rangle
\langle\eta^1_{g_1},\CR({g'}^{-1})\psi^1\rangle {\rm d}_{\S_1}(g_1)\Big)
\langle\eta'_{g'},\psi'\rangle{\rm d}_{\B'}(g')\\
&=\int_{\B'}
\langle\varphi',\eta'_{g'}\rangle
\langle\CR({g'}^{-1})\varphi^1,\CR({g'}^{-1})\psi^1\rangle
\langle\eta'_{g'},\psi'\rangle{\rm d}_{\B'}(g')\\
&=\langle\varphi^1,\psi^1\rangle{\hspace{-1mm}}\int_{\B'}{\hspace{-2mm}}
\langle\varphi',\eta'_{g'}\rangle
\langle\eta'_{g'},\psi'\rangle{\rm d}g'=\langle\varphi^1,\psi^1\rangle
\langle\varphi',\psi'\rangle=\langle\varphi,\psi\rangle.
\end{align*}
\end{proof}
\begin{rmk}
In the following, we will absorb the constant
$C(\eta)^{1/2}$ of Proposition \ref{WRI}
in a redefinition of the mother wavelet $\eta\in\CD(Q_N\times\dots\times Q_1)$.
\end{rmk}

\begin{rmk}
Other types of weak resolution of the identity can be constructed in this setting.
For instance, setting $\tilde\eta_x:=\Omega^\B_{\bm_0}(x)\eta$, where 
$\eta$ is arbitrary in $ \CH_\chi$,
we have from the unitarity of the quantization map $\Omega^\B_{\bm_0}$:
$$
\langle \psi,\phi\rangle= \|\eta\|^{-2}\int_\B \langle\psi,\tilde\eta_x\rangle
\langle\tilde\eta_x,\phi\rangle\;{\rm d}_\B( x)
\;,
$$
for all   $\phi,\psi\in\CH_\chi$. Similarly, let $W^\eta_{x,y}$ be
the Wigner function on $\B$, associated to a pair of wavelets $\eta_x,\eta_y$:
$$
W^\eta_{x,y}(z):=\sigma_{\bm_0}\big[|\eta_x\rangle\langle \eta_y|\big](z)=
\langle\eta_y|\Omega_{\bm_0}(z)|\eta_x\rangle\,,\quad x,y,z\in\B\;.
$$
 Then, these  Wigner functions
may be used to construct a weak resolution of the identity on $L^2(\B)$:
For all $f_1,f_2\in L^2(\B)$, we have
$$
\langle f_1,f_2\rangle=\|\eta\|^{-4}\int_{\B\times\B} \langle f_1,W^\eta_{x,y}\rangle\,
\langle W^\eta_{x,y},f_2\rangle\,{\rm d}_\B(x)\,
{\rm d}_\B(y)\;.
$$
\end{rmk}

\begin{prop}
\label{CV-pre-born}
Let $\B$ be a normal $\bf j$-group, $\CH$ be a separable Hilbert space and $A$ be a densely defined operator on 
$ \CH_\chi\otimes\CH$, whose domain contains the 
(algebraic) tensor product $\CD(Q_N\times\dots\times Q_1)\otimes\CH$.
  For $x,y\in\B$ and $\eta\in \CD(Q_N\times\dots\times Q_1)$ define  
the element $\langle \eta_x,A\eta_y\rangle_{\CH_\chi}$  of $\CB(\CH)$ by means of the quadratic form
$$
\CH\times\CH\to\C\,,\quad(\phi,\psi)\mapsto \langle\eta_x\otimes\phi, A \,
\eta_y\otimes\psi\rangle_{\CH_\chi\otimes\CH}\;.
$$
Assuming further that
$$
\sup_{y\in \B}\;\int_\B \|\langle \eta_x,A\eta_y\rangle_{\CH_\chi}\|_{\CB(\CH)}\,{\rm d_\B}(x)<\infty
\,,\quad \sup_{x\in \B}\;\int_\B \|\langle \eta_x,A\eta_y\rangle_{\CH_\chi}\|_{\CB(\CH)}\,{\rm d_\B}(y)
<\infty\;,
$$
then $A$ extends to a bounded operator on $\CH_\chi\otimes\CH$.
\end{prop}
\begin{proof}
Since $\eta_x$ is smooth and compactly supported, our assumption about the domain of $A$ ensures 
that $\langle \eta_x,A\eta_y\rangle_{\CH_\chi}$  is well defined as an element of $\CB(\CH)$. Thus, Lemma \ref{CS} applied to $(X,\mu)=(\B,{\rm d}_\B)$, yields that the operator
$\tilde A$ on $L^2(\B,\CH)$  given by
$$
\tilde A F(x):=\int_\B\langle \eta_x,A\eta_y\rangle_{\CH_\chi}\,F(y)\,{d\mu}(y)\;,
$$
and is bounded, with
\begin{align*}
\|\tilde A\|&\leq 
\Big(\sup_{y\in \B}\;\int_\B \|\langle \eta_x,A\eta_y\rangle_{\CH_\chi}\|_{\CB(\CH)}\,{\rm d_\B}(x)
\Big)^{1/2}\\&\qquad\qquad\qquad\times
\Big(\sup_{x\in \B}\;\int_\B \|\langle \eta_x,A\eta_y\rangle_{\CH_\chi}\|_{\CB(\CH)}\,{\rm d_\B}(y)\Big)^{1/2}<\infty\;.
\end{align*}
For $\Phi\in\CH_\chi\otimes\CH$ define the $\CH$-valued function on $\B$:
 $\tilde\Phi:=[x\in\B\mapsto \langle\eta_x,\Phi\rangle_{\CH_\chi}\in\CH]$.  By Proposition \ref{WRI} we know that  $\tilde\Phi$  belongs to $L^2(\B,\CH)$ with $\|\Phi\|_{\CH_\chi\otimes\CH}
 =\|\tilde\Phi\|_{L^2(\B,\CH)}$.
Take now $\Phi,\Psi\in \,{\rm dom}\, A$. In this case, we can use twice the resolution of the identity to get
\begin{align*}
\langle\Phi,A\Psi\rangle_{\CH_\chi\otimes\CH}&=\int_{\B\times\B} \big\langle\langle\eta_x,\Phi\rangle_{\CH_\chi},\langle\eta_x,A\eta_y\rangle_{\CH_\chi}\langle\eta_y,\Psi\rangle_{\CH_\chi}\big\rangle_{\CH}\,{\rm d_\B}(x)\,{\rm d_\B}(y)\\&=
\langle\tilde\Phi,\tilde A\tilde\Psi\rangle_{L^2(\B,\CH)}\;.
\end{align*}
Therefore, we conclude that
\begin{align*}
\big|\langle\Phi,A\Psi\rangle_{\CH_\chi\otimes\CH}\big|&=\big|\langle\tilde\Phi,\tilde A\tilde\Psi\rangle_{L^2(\B,\CH)}\big|\\&\leq \|\tilde\Phi\|_{L^2(\B,\CH)}\,\|\tilde\Psi\|_{L^2(\B,\CH)}\,\|\tilde A\|=\|\Phi\|_{\CH_\chi\otimes\CH}\,\|\Psi\|_{\CH_\chi\otimes\CH}\,\|\tilde A\|<\infty\;,
\end{align*}
and the result follows immediately.
\end{proof}

\section{A tempered pair from the one-point phase}
  Let $\S$ be an elementary normal $\bf j$-group. Consider the one-point phase
 $\bE_\theta$, defined in (\ref{E}), and given by
\begin{equation}
\label{bS}
\bE_\theta(qb)= \overline{\chi}_\theta\left(\bC_q(b^{-1}\tilde{\sigma}b)\right)=:e^{\tfrac{i}\theta {\bf S}(qb)}\;.
\end{equation}
Recall that from Lemma \ref{S1POINT}, we have in the coordinates \eqref{chartS} and
up to a global constant factor:
$$
\bS(a,n\oplus m,t)=t\sinh 2a+\omega^0(n,m)\cosh^2a\;.
$$

The aim of  this section it to prove that the pair $(\S,\bf S)$,  is tempered, 
admissible and tame.
For this, we consider the following  decomposition
of the Lie algebra $\s$ (i.e$.$ the one we used in Equation \eqref{DECOMP-S}):
\begin{equation*}
\s=\bigoplus_{k=0}^3V_k\quad\mbox{\rm where }\quad V_0\;:=\a\,,\quad V_1:=\l^\star\,,\quad V_2:=\l
\quad\mbox{\rm and}\quad V_3:=\R E\;.
\end{equation*}
As usual, we us fix $\{f_j\}_{j=1}^d$, a basis of $\l^\star$  to which we associate  $\{e_j\}_{j=1}^d$  the 
symplectic-dual basis of $\l$,  defined by $\omega^0(f_i,e_j)=\delta_{i,j}$. Associated to the 
decomposition $v=n\oplus m\in\l^\star\oplus\l= V$, we get coordinates
$$
n_j:=\omega^0(n,e_j)\,,\quad m_j:=\omega^0(f_j,m)\,,\qquad j=1,\dots,d\;.
$$
From the expressions \eqref{INVVECT} of the 
left invariant vector fields  of $\S$, in the chart \eqref{chartS}, 
 we  get the following coordinates system on $\S$:
\begin{align}
\label{again}
x_0:=\widetilde H\, \bS=2e^{-2a}t-(1+e^{-2a})\omega^0(n,m)\,,&\quad
x_1^j:=\widetilde f_j \,\bS=(1+e^{-2a})m_j\;,\\
 x_2^j:=\widetilde e_j\,\bS=(1+e^{2a})n_j\,,&\quad x_3:=\widetilde E\,\bS=\sinh2a\;.
 \nonumber
\end{align}
We then deduce:
\begin{lem}
\label{phi}
The pair $(\S,\bf S)$ is tempered in the sense of Definition \ref{TEMPPAIR}. 
Moreover, the Jacobian of the map
$$
\phi:\S\to\s^\star\,,\;g\mapsto\big[\s\to\R\,,\;X\in\s\mapsto 
\big(\widetilde X \,{\bf S}\big)(g)\big]\;,
$$
is proportional to $ \bm_0^2\times\Delta_\S^{-1/(d+1)}$.
\end{lem}

The following Lemma is actually all that we need to prove admissibility 
(in the sense of Definition \ref{TEMPADM}) of the tempered pair $(\S,\bf S)$:
\begin{lem}
\label{lem:temp}
For every $k\in\{0,1,2,3\}$, there exists a tempered function $\bm_k>0$ with $\partial_{x_{j}}\bm_k=0$ for every $j\leq k$ and such that for every $X\in\CU(V^{(k)})$,
there exists $C_X>0$ with
$$
\big|\widetilde{X}\,{x}_k\big|\leq C_X\,\bm_k\,(1+|{x}_k|)\;.
$$
\end{lem}
\begin{proof}
From the computations, for $k\in\N^*$ and $i,j=1,\dots,d$:
\begin{align*}
&\widetilde H^k\, x_0=(-1)^k\big(2^{2k+1} e^{-2a}t-2^{k}(1+2^ke^{-2a})\omega^0(n,m)\big)\;,&\\
&\widetilde H^k\, x_1^j=(-1)^k\big(1+3^ke^{-2a}\big)m_j\,,\quad \widetilde f _i \,x_1^j=0\;,&\\
&\widetilde H^k\, x_2^j=\big((-1)^k+e^{2a}\big)n_j\,,\quad \widetilde f _i \,x_2^j=(1+e^{2a})\delta_i^j\,,\quad \widetilde e_i \,x_2^j=0\;,&\\
&\widetilde H^k\, x_3=2^{k+1}\begin{cases} \cosh2a\,,\quad &\mbox{$k$  even}\\
\sinh2a\,,\quad &\mbox{$k$ odd}\end{cases}\,,\quad  \widetilde f _i \,x_3=0\,,\quad\widetilde e_i \,x_3=0\,,\quad\widetilde E\,x_3=0\;,&
\end{align*}
and  elementary estimates, we obtain:    
\begin{align*}
|\widetilde X x_0|&\leq C^0_X (1+|x_0|)(1+|x_1|\,|x_2|)\,,\quad\forall X\in\CU(V_0)\;,\\
|\widetilde X x_1^j|&\leq C^1_X (1+|x_1^j|)\,,
\quad\quad\quad\quad\quad\quad\forall X\in\CU(V^{(1)})\;,\\
|\widetilde X x_2^j|&\leq C^2_X (1+|x_2^j|)(1+|x_3|)\,,\quad\forall X\in\CU(V^{(2)})\;,\\
|\widetilde X x_3|&\leq C^3_X (1+|x_3|)\,,\quad\quad\quad\quad\forall X\in\CU(\s)\;,
\end{align*}
 and the claim follows with $\bm_0(x)=(1+|x_1|\,|x_2|)$, $\bm_1(x)=1$, $\bm_2(x)=(1+|x_3|)$, $\bm_3(x)=1$.
\end{proof}

Repeating the   arguments of  the proof of  Proposition \ref{MULTIPLIERS-ELEM}, we 
deduce admissibility for the tempered pair $(\S,\bS)$. 
\begin{lem}
\label{bel}
Define
$$
X_0:=1-H^2\in\CU(V_0)\;,\quad X_1:=1-\sum_{j=1}^d\,f_j^2\in\CU(V_1)\;,
$$
$$
 X_2:=1-\sum_{j=1}^d\,e_j^2\in\CU(V_2)\;,\quad
X_3:=1-E^2\in\CU(V_3)\;.
$$
Then the corresponding multipliers $\alpha_k:=\bE^{-1}\,\widetilde{X}_k\,\bE$ satisfy conditions (i) 
and (ii) of Definition \ref{TEMPADM}, with $\rho_k=2$ and the $\mu_k$'s are given by the $\bm_k$'s
of Lemma \ref{lem:temp}.
\end{lem}

Lastly, we observe that tameness (see Definition \ref{tame}) follows from Lemma \ref{fd-explicit}
and arguments very similar to those of Corollary \ref{MWP}.
We then summarize all this by stating the main result of this section:

\begin{thm}
\label{TP1P}
Let $\S$ be an elementary normal $\bf j$-group and let $\bS\in C^\infty(\S)$ be as given in \eqref{bS}.
Then the pair $(\S,\bf S)$ is  tempered, admissible and tame. 
\end{thm}

\begin{rmk}
For $\B$
a generic   normal $\bf j$-group $\B$, we could also define a one-point tempered pair, by setting
\begin{equation}
\label{bEB}
\bE_\theta^\B:=\exp\{\tfrac{2i}\theta\bS^\B\}:\B\to\mathbb U(1)\,,\qquad {\bf S}^\B:\B\to\R\,,\quad g\mapsto\sum_{j=1}^N {\bf S}^{\S_j}(g_j)\;,
\end{equation}
where ${\bf S}^{\S_j}$ is the one-point phase \eqref{bS} of each elementary factor of $\B$ in the parametrization
 $g=g_1\dots g_N\in\B$, relative to a Pyatetskii-Shapiro decomposition.
 Then temperedness and admissibility will follow from arguments very similar than those of 
 Theorem \ref{TASP}.
\end{rmk}

\begin{rmk}
The one-point Schwartz space $\CS^{S_{\rm can}}(\S)$ associated with the two-point pair
$(\S\times\S,S_{\rm can})$ as given in Definition \ref{S1P}, coincides with the one-point Schwartz space $\CS^{\bS}(\S)$ associated 
with the one-point par $(\S,\bS)$.
\end{rmk}

\section{Extension of the oscillatory integral}

Associated to a tempered, admissible and tame pair $(G,S)$,
 we have constructed in section \ref{OITP} a continuous linear map  for any Fr\'echet space 
 $\CE$ and  any element $\bm\in\CB^\mu(G)$ (with $\mu$ a tempered weight on $G$):
 $$
 \widetilde{\int_G \bE\,\bm}:\CB^{\underline\mu}(G,\CE)\to\CE\;,
 $$
 which extends the ordinary integral on $\CD(G,\CE)$ and that we called
 the {\em oscillatory integral}.

 The aim of the present section is to explain how for the tempered pair 
 $(\S,\bS)$ of Theorem \ref{TP1P}, one can enlarge the domain of definition
 of the oscillatory integral. For this let $\CE$ be a Fr\'echet space, $\underline\mu=\{\mu_j\}_{j\in\N}$
 a family of tempered weights on $\S$   and  $\nu$ be a fixed 
tempered and $\Y$-right-invariant  weight\footnote{We may view $\nu$ as a function on $Q$.} 
on $\S$. 
 Let us then consider the following subspace  of $C^\infty(\S,\CE)$:
 \begin{align*}
&\CB^{\underline\mu,\nu}(\S,\CE):=
\Big\{F\in C^\infty(\S,\CE)\;:\;\forall\,(j,X,Y)\in\N\times\CU(\q)
\times \CU(\fY)\,,\,\exists \,C\;:\;\nonumber\\
&\qquad\qquad\qquad\qquad\qquad\qquad\qquad\qquad
\|\widetilde{X}\,\widetilde Y\,F(qb)\|_j\leq C\,\nu(q)^{\deg(X)}\,\mu_j(qb)
\Big\}\;.
\end{align*}
This space may be understood as a variant of the symbol space $\CB^{\underline\mu}(\S,\CE)$, 
where a
specific dependence of the family of weights $\underline\mu$ in the degree of the derivative
is allowed. We endow the latter space with the following set of semi-norms:

 \begin{equation}
\label{norms3}
 \|F\|_{j,k_1,k_2,\underline\mu,\nu}:=\sup_{X\in\,\CU_{k_1}(\q)}\,\sup_{Y\in\,\CU_{k_2}(\fY)}\,
 \sup_{qb\in \S}\Big\{\frac{\|\widetilde{X}\widetilde{Y}F(qb)\|_j}
 {\mu_j(qb)\,\nu(q)^{k_1}\,|X|_{k_1}\,|Y|_{k_2}}\Big\}\;,
 \end{equation}
where $j,k_1,k_2\in\N$, and  $\bigcup_{k\in\N}\CU_{k}(\q)$, $\bigcup_{k\in\N}\CU_{k}(\fY)$ are the filtrations of
$\CU(\q)$ and $\CU(\fY)$ associated to the choice of PBW basis as explained in \eqref{PBW}. 
 As expected, the space $\CB^{\underline\mu,\nu}(\S,\CE)$ is  Fr\'echet  for the topology 
induced by the
semi-norms \eqref{norms3} and most the properties of Lemma  \ref{SmoothFamily} remain true.

\begin{lem} 
\label{SmoothFamilybis}
Let $(\S,\CE)$ be as above, let $\underline\mu$, $\underline\rho$ and $\underline{\hat\mu}$ be three
 families of weights  on
 $\S$ and let $\nu$, $\lambda$ and $\hat\nu$ be three right-$\Y$-invariant weights on $\S$.
\begin{enumerate}
\item[(i)] The space $\CB^{\underline\mu,\nu}(\S,\CE)$ is Fr\'echet.
\item[(ii)] The bilinear map:
\begin{equation*}
\CB^{\underline\mu,\nu}(\S)\times\CB^{\underline\rho,\lambda}(\S,\CE)\to
\CB^{\underline\mu.\underline\rho, \nu.\lambda}(\S,\CE)\,,\quad (u,F)\mapsto[g\in \S\mapsto u(g)\,F(g)\in\CE]\;,
\end{equation*}
is  continuous.
\item[(iii)]  If there exists $C>0$ such that $\underline\mu\leq C\underline{\hat\mu}$ and
$\nu\leq C\hat\nu$, then 
$\CB^{\underline\mu,\nu}(\S,\CE)\subset \CB^{\underline{\hat\mu},\hat\nu}(\S,\CE)$ continuously.
\item[(iv)] Assume that $\underline\mu\prec\underline\mu$ and
$\nu\prec\hat\nu$. 
Then, the closure of  $\CD(\S,\CE)$
in $\CB^{\underline{\hat\mu},\hat\nu}(\S,\CE)$ contains $\CB^{\underline{\mu},\nu}(\S,\CE)$.
 In particular,
 $\CD(\S,\CE)$ is a dense subset of $\CB^{\underline{\mu},\nu}(\S,\CE)$ for the induced topology of 
 $\CB^{\underline{\hat\mu},\hat\nu}(\S,\CE)$.
\end{enumerate}
\end{lem}
\begin{proof}
The first assertion follows from the fact that a countable projective limit of Fr\'echet spaces is Fr\'echet
and that $\CB^{\underline\mu,\nu}(\S,\CE)$ can be realized as the countable projective
 limit of the family
of Banach spaces underlying the norms $\sum_{i=0}^j \sum_{l_1=0}^{k_1}\sum_{l_2=0}^{k_2} 
\|.\|_{i,l_1,l_2,\underline\mu,\nu}$.
The proof of all the other statements are identical to their counter-parts in Lemma  
\ref{SmoothFamily}.
\end{proof}

We are now able to prove our extension result for the oscillatory integral associated to
the admissible,  tempered and tame pair $(\S,\bS)$:
\begin{thm}
\label{OIEXT}
Let $\underline\mu$ be family of tempered weights on $\S$, $\nu$ a $\Y$-right-invariant 
tempered weight on $\S$
and $\bm$ an element of $\CB^\lambda(\S)$ for another tempered weight $\lambda$ on $\S$. 
Let also $\bD_{\vec r}$, $\vec r\in\N^4$, be the differential operator constructed in \eqref{bD}.
Then for all $j\in\N$, there exist
 $\vec r_j\in\N^{4}$, $C_j>0$ and $k_j,l_j\in\N$, such that  for every element 
 $F\in \CB^{\underline\mu,\nu}(\S,\CE)$, we have
\begin{equation*}
\int_\S\|\bD_{\vec r} \,\bm(g) F(g)\|_j\,{\rm d}_\S(g)\leq C_j\,\|F\|_{j,k_j,l_j,\underline\mu,\nu}\;.
\end{equation*}
Consequently, the oscillatory integral constructed in Definition \ref{OI} for the tame
and admissible tempered pair $(\S,\bS)$, originally defined in $\CB^{\underline\mu}(\S,\CE)$, 
extends as 
a continuous map:
$$
\widetilde{\int_\S \bm\,\bE}:\CB^{\underline\mu,\nu}(\S,\CE)\to \CE\;.
$$
\end{thm}
\begin{proof}
The proof is very similar to those of Proposition \ref{PROPIP}, so we focus on the 
differences due to the particular behavior at infinity of an element of $\CB^{\underline\mu,\nu}(\S,\CE)$.

By Lemma \ref{phi}, the  Radon-Nikodym 
derivative of the left Haar measure on $\S$ with respect to the Lebesgue measure on $\s^\star$, is 
bounded by a polynomial of order $2d+4$ in the coordinate $x_3$. 
For each $j\in\N$,  the weight $\mu_j$ is also bounded by a polynomial in $x_0,x_1,x_2,x_3$. 
Now, observe that by construction of the operator $\bD_{\vec r}$ in \eqref{bD}, 
we have for any $\vec r=(r_0,r_1,r_2, r_3)\in\N^4$, with 
$K_1= 2r_0+2r_1$, $K_2=2r_2+2r_3$ and with the notations given 
in \eqref{bDPhi}:
\begin{align}\label{LDOM}
|\bD_{\vec r} \,F|
&\leq|\Psi_0|\,|\Psi_{1,0}|\,|\Psi_{2,1,0}|\,|\Psi_{3,2,1,0}|\,\big|\widetilde{X}'_{3,2,1,0}\,F\big|
\nonumber\\
&\leq \,C\,|\Psi_0|\,|\Psi_{1,0}|\,|\Psi_{2,1,0}|\,|\Psi_{3,2,1,0}|\,\mu_j\,\nu^{2r_0+2r_1}
\|F\|_{j,K_1,K_2,\underline\mu,\nu}\;.
\end{align}
This will gives the estimate we need, if we prove that the function in front of 
$\|F\|_{j,K_1,K_2,\underline\mu,\nu}$ in \eqref{LDOM}
 is integrable for a suitable choice of $\vec r\in\N^{4}$.
We  prove  a stronger result, namely that given $\vec R\in\N^{4}$, there exists 
$\vec r\in\N^4$ such that
\begin{align*}
&|\Psi_0|\,|\Psi_{1,0}|\,|\Psi_{2,1,0}|\,|\Psi_{3,2,1,0}|\,\nu^{2r_0+2r_1}\\
&\leq\frac{C}{
(1+|x_0|)^{R_0}(1+|x_1|)^{R_1}(1+|x_2|)^{R_2}(1+|x_3|)^{R_3}}\;.
\end{align*}
From Corollary \ref{ESTIMATION} and Lemma \ref{bel}, we obtain the following estimation:
\begin{align*}
&|\Psi_0|\,|\Psi_{1,0}|\,|\Psi_{2,1,0}|\,|\Psi_{3,2,1,0}|\\
&\leq
C\,\frac{(1+|x_1||x_2|)^{2r_0^2}}{(1+|x_0|)^{2r_0}}\frac1{(1+|x_1|)^{2r_1}}
\frac{(1+|x_3|)^{2r_2(r_0+r_1+r_2)}}{(1+|x_2|)^{2r_2}}\frac1{(1+|x_3|)^{2r_3}}\;.
\end{align*}
Lastly  (this is the main difference with the proof of Proposition \ref{PROPIP}),
note that $\nu$, the tempered function on $Q$, can be bounded by $|x_2|^{p_2}|x_3|^{p_3}$
for some integers $p_2,p_3$. Hence $|\Psi_0|\,|\Psi_{1,0}|\,|\Psi_{2,1,0}|\,|\Psi_{3,2,1,0}|\,\nu^{2r_0+2r_1}$ is smaller than
\begin{align*}
&C\,
(1+|x_0|)^{-2r_0}(1+|x_1|)^{-2r_1+2r_0^2}(1+|x_2|)^{-2r_2+2r_0^2+2p_2(r_0+r_1)}\\
&\qquad\qquad\qquad\qquad\qquad\qquad\quad\;\times
(1+|x_3|)^{-2r_3+2r_2(r_0+r_1+r_2)+2p_3(r_0+r_1)}\;,
\end{align*}
and the claim follows.
\end{proof}

\section{A  Calder\'on-Vaillancourt type  estimate}

For $j=1,\dots,N$, fix $\bm_j$ a $\Y_j$-right-invariant tempered weight on $\S_j$
(that we identify in a natural manner as a function on $Q_j$), in the sense of Definition
\ref{temp-grp} for the tempered pair $(\S_j,\bS^{\S_j})$ underlying Theorem \ref{TP1P}.
 Let also $A$ be a $C^*$-subalgebra of $\CB(\CH)$, with $\CH$ a separable 
Hilbert space.
Our aim here is to prove that for $F\in\CB(\B,A)$,
  the operator $\Omega_{\theta,\bm}(F)$, defined via
a suitable quadratic form  
on $\CH_\chi\otimes\CH$, is bounded\footnote{Observe that this property holds
for $F\in \CS^{S_{\rm can}^\B}(\B,A)$, by Proposition \ref{Relation-Weyl2}  as 
$\CS^{S_{\rm can}^\B}(\B,A)\subset \CS(\B;A)$.}. 
We start by proceeding formally, in order 
to explain our global strategy. Also, to simplify the notations, we assume
first that $\B=\S$ is elementary.
So let
$\Phi,\Psi\in  \CH_\chi\otimes\CH$. Using twice the resolution of the identity of
Proposition \ref{WRI}, we write
\begin{align*}
&\langle\Phi,\Omega_{\theta,\bm}(F)\Psi\rangle_{\CH_\chi\otimes\CH}=\\
&\qquad\qquad
\int_{\S\times\S}\big\langle\langle\eta_x,\Phi\rangle_{\CH_\chi},\langle\eta_x,\Omega_{\theta,\bm}(F)\eta_y\rangle_{\CH_\chi}\langle\eta_y,\Psi\rangle_{\CH_\chi}\big\rangle_\CH\,{\rm d}_\S(x)\,{\rm d}_\S(y)\;.
\end{align*}
Next, we use the $\S$-covariance of the pseudo-differential calculus to get
$$
\langle\eta_x,\Omega_{\theta,\bm}(F)\eta_y\rangle_{\CH_\chi}=\langle\eta_{y^{-1}x},\Omega_{\theta,\bm}(L^\star_{y^{-1}}F)\eta\rangle_{\CH_\chi}\;.
$$
Then, we exchange the integrals over $\S$ and $Q$ and expand the scalar product of $\CH_\chi$
to obtain:
\begin{align*}
\langle\eta_{x},\Omega_{\theta,\bm}(F)\eta\rangle_{\CH_\chi}= \int_{\S\times Q} F(q_0qb)\, \overline{\eta_x}(q_0)\,\bm(q^{-1})\,\bE(qb)\,\eta\big(q_0\underline s^eq\big)\,{\rm d}_Q(q_0)\,{\rm d}_\S(qb)\;.
\end{align*}
Given $F\in\CB(\S,A)$, this suggests to define the function
\begin{equation}
\label{Salakis}
F^\eta:\S\times Q\to A\,,\qquad (qb,q_0)\mapsto F(q_0qb)\, \eta\big(q_0\underline s^e q\big)\;,
\end{equation}
so that  with $\hat\bm(qb):=\bm(q^{-1})$ and with the notations of Proposition \ref{first-type}, we will have
$$
\langle\eta_{x},\Omega_{\theta,\bm}(F)\eta\rangle_{\CH_\chi}=\CF^{\overline \eta}\Big(\int_\S\bE(y)\, 
\hat\bm(y)\,F^\eta(.,y)\,{\rm d}_\S(y)\Big)(x)\;.
$$
Consequently, we obtain
\begin{align*}
&\langle\Phi,\Omega_{\theta,\bm}(F)\Psi\rangle_{\CH_\chi\otimes\CH}=\\
&\int_{\S\times\S}\big\langle\langle\eta_x,\Psi\rangle_{\CH_\chi},
\CF^{\overline \eta}\Big(\int_\S\bE(z)\, \hat\bm(z)\,\big(L^\star_{y^{-1}}F\big)^\eta(.,z)\,
{\rm d}_\S(z)\Big)(y^{-1}x)\big\langle \eta_y,\Psi\rangle_{\CH_\chi}\big\rangle_\CH\\
&\times{\rm d}_\S(x)\,{\rm d}_\S(y)\;.
\end{align*}
Surprisingly, this is the right hand side of the (formal) equality above which gives rise
to a well defined and bounded quadratic form on $\CH_\chi\otimes\CH$, once the 
integral sign in the middle is replaced by an oscillatory one in the sense of
Theorem \ref{OIEXT} for the tempered pair $(\S,\bS)$.

Coming back to the case of a generic normal $\bf j$-group $\B$, 
the most important step is to understand 
the properties of the corresponding map $F\mapsto F^\eta$ given in \eqref{Salakis}. 

\begin{lem}
\label{F-eta}
Let $A$ be a $C^*$-algebra, $\B$ be a normal $\bf j$-group with Pyatetskii-Shapiro decomposition  
$\B=(\S_N\ltimes\dots)\ltimes\S_{1}$ and 
 $\eta\in\CD(Q_N\times\dots\times Q_1)$. Then the map
\begin{align*}
  F\mapsto F^\eta&:=
\Big[q_Nb_N\in\S_N\mapsto\Big[q_{N-1}b_{N-1}\in\S_{N-1}\mapsto
 \dots  \Big[ q_1b_1\in \S_1\mapsto \\
 &\quad\Big[(q'_N,\dots,q'_1)\in 
  Q_N\times\dots\times Q_1  \mapsto
  F(q'_1q_1b_1\dots q'_{N-1}q_{N-1}b_{N-1}q'_Nq_Nb_N)\\
  &\quad\quad\quad\quad\quad\quad\quad\quad\quad\quad\quad
  \quad\quad\quad\quad\times
   \eta\big(q'_N\underline s^eq_N,\dots,q'_1\underline s^eq_1 \big)\in A\Big]\Big]
   \dots\Big]\Big]\;,
\end{align*}
is continuous  from $\CB(\B,A)$ to  
$$
\CB^{\underline\mu_N,\nu_N}\big(\S_N,
\CB^{\underline\mu_{N-1},\nu_{N-1}}
\big(\S_{N-1},\dots\CB^{\underline\mu_1,\nu_1}\big(\S_1,
\CS(Q_N\times\dots\times Q_1,A)\big)\dots\big)\big)\;,
$$ 
where for $j=N,\dots1$, we have settled:
$$
\nu_j:=\fd_{Q_j}^2\,,\quad \underline\mu_j:=\big\{\fd_{\S_j}^{n_j(k_{j-1},l_{j-1};\dots;k_1,l_1;k,l)}\big\}
_{(k_{j-1},l_{j-1};\dots;k_1,l_1;k,l)\in\N^{2j}}\;,
$$
where $(k_{j-1},l_{j-1};\dots;k_1,l_1;k,l)\in\N^{2j}$ labels the semi-norms of the space
$$
\CB^{\underline\mu_{j-1},\nu_{j-1}}
\big(\S_{j-1},\dots\CB^{\underline\mu_1,\nu_1}\big(\S_1,
\CS(Q_N\times\dots\times Q_1,A)\big)\dots\big)\;,
$$
and the exponent $n_j(k_{j-1},l_{j-1};\dots;k_1,l_1;k,l)\in\N$ is linear in its arguments.
\end{lem}
\begin{proof}
For notational convenience, we assume that $\B$ contains only two elementary factors, i.e$.$ 
$\B=\S_2\ltimes\S_1$ with $\S_1,\S_2$ elementary normal $\bf j$-groups. This is enough to
understand the global mechanism and the proof for a 
generic normal
$\bf j$-group with an arbitrary number of elementary factors will then follow by induction, 
without essential supplementary  difficulties. In this simplified  case, we have to prove that the map
\begin{align*}
  F\mapsto F^\eta&:=\big[q_2b_2\in\S_2\mapsto\big[q_{1}b_{1}\in\S_{1}
  \mapsto\big[(q'_2,q'_1)\in Q_2\times Q_1  \mapsto\\
  &\qquad\qquad\qquad\qquad\qquad
  F(q'_1q_1b_1 q'_{2}q_{2}b_{2})\,
   \eta\big(q'_2\underline s^eq_2,q'_1\underline s^eq_1 \big)\in A\big]\big]\big]\;,
\end{align*}
is continuous  from $\CB(\B,A)$ to  
\begin{align}
\label{the-space}
\CB^{\underline\mu_2,\nu_2}\big(\S_2,\CB^{\underline\mu_1,\nu_1}
\big(\S_1,\CS(Q_2\times Q_1,A\big)\big)\;.
\end{align}
By the discussion following Remark \ref{fd-explicit2}, 
it is clear that one may regard $\CS(Q_2\times Q_1,A)$ as a Fr\'echet space for the topology 
induced by the following countable set of semi-norms:
$$
\|f\|_{k,j}:=\sup_{Z^2\in\,\CU_k(\q_2)}\sup_{Z^1\in\,\CU_k(\q_1)}
\sup_{(q_2,q_1)\in Q_2\times Q_1}\Big\{\frac{\fd_{Q_2}(q_2)^j\fd_{Q_1}(q_1)^j
\|\underline Z^2_{q_2}\,\underline Z^1_{q_1}\,
 f(q_1,q_2)\|}{|Z^2|_k|Z^1|_k}\Big\}\;,
$$
i.e. we may use right-invariant vector fields instead of left-invariant one since they are related 
by tempered functions with tempered inverses. 
Note then that  the natural Fr\'echet topology of 
the space \eqref{the-space} 
is associated with the following countable family of semi-norms (indexed by $(k_2,l_2,k_1,l_1,k,j)
\in\N^6$):
\begin{align*}
&\Phi\!\mapsto\!\!\!
\sup_{X\in\,\CU_{k_2}(\q_2)}\!\sup_{X'\in\,\CU_{l_2}(\fY_2)}\!\sup_{q_2b_2\in\S_2}\!
\sup_{Y\in\,\CU_{k_1}(\q_1)}\!\sup_{Y'\in\,\CU_{l_1}(\fY_1)}\!\sup_{q_1b_1\in\S_1}\!
\sup_{Z\in\,\CU_k(\q_2\oplus\q_1)}\!\sup_{(q_2,q_1)\in Q_2\times Q_1}\\
&
\frac{\fd_{Q_2}(q'_2)^j\fd_{Q_1}(q'_1)^j\,\|\underline Z_{(q'_2,q'_1)}\widetilde Y_{q_1b_1} 
\widetilde Y'_{q_1b_1} 
\widetilde X_{q_2b_2}\widetilde X'_{q_2b_2} 
\Phi(q_2b_2;q_1b_1;q'_2,q'_1)\|}
{\fd_{\S_2}\!(q_2b_2)^{n_2(k_1,l_1,k,j)}\fd_{Q_2}\!(q_2)^{2k_2}\fd_{\S_1}\!(q_1b_1)^{n_1(k,j)}
\fd_{Q_1}\!(q_1)^{2k_1}
|X|_{k_2}|X'|_{l_2}|Y|_{k_1}|Y'|_{l_1}|Z|_{k}}\;.
\end{align*}
To simplify  the notations, we denote the latter semi-norm by $\|.\|_{k_2,l_2,k_1,l_1,k,j}$.
Then,  for 
$$
(X,X',Y,Y',Z^2,Z^1)\in\CU(\q_2)\times\CU(\fY_2)\times \CU(\q_1)\times\CU(\fY_1)\times
\CU(\q_2)\times\CU(\q_1)\;,
$$
we get within Sweedler's notation:
\begin{align*}
&\underline Z^2_{q'_2}\underline Z^1_{q'_1}\widetilde Y_{q_1b_1} 
\widetilde Y'_{q_1b_1} 
\widetilde X_{q_2b_2}\widetilde X'_{q_2b_2} F^\eta(q_1b_1;q_2b_2;q'_2,q'_1)
=
\sum_{(X)}\sum_{(X')}\sum_{(Y)} \sum_{(X')}\sum_{(Z^2)}\sum_{(Z^1)}\\
&\Big(\!\underline Z^2_{(1)q'_2}\underline Z^1_{(1)q'_1}
 \widetilde Y_{(1)q_1b_1} \widetilde Y'_{(1)q_1b_1} \widetilde X_{(1)q_2b_2}
 \widetilde X'_{(1)q_2b_2}
F(q'_1q_1b_1 q'_2q_2b_2)\!\Big)\!\\
&\qquad\qquad\qquad\times\Big(\!
\underline Z^2_{(2)q'_2}\underline Z^1_{(2)q'_1}
\widetilde Y_{(2)q_1} \widetilde X_{(2)q_2}
\eta\big(q'_2\underline s^e(q_2),q'_1\underline s^e(q_1) \big)\!\Big)\;.
\end{align*}
From the same reasoning as those in the proof of Lemma \ref{symbols} (v),
we deduce that
\begin{align}
\label{NO}
\|F^\eta\|_{k_2,k_1,k,j}\leq&\, C \sup\frac{\fd_{Q_2}(q'_2)^j\fd_{Q_1}(q'_1)^j}{\fd_{\S_2}(q_2b_2)^{n_2(k_1,l_1,k,j)}\,\fd_{Q_2}(q_2)^{2k_2}\,\fd_{\S_1}(q_1b_1)^{n_1(k,j)}\,
\fd_{Q_1}(q_1)^{2k_1}}\nonumber\\
&\times\sup
\frac{\big\|\underline Z^1_{q'_1}\underline Z^2_{q'_2}
 \widetilde Y_{q_1b_1} \widetilde Y'_{q_1b_1} \widetilde X_{q_2b_2}\widetilde X'_{q_2b_2}
F(q'_1q_1b_1 q'_2q_2b_2)\big\|}{|X|_{k_2}|X'|_{l_2}|Y|_{k_1}|Y'|_{l_1}|Z^2|_k|Z^1|_k}
\nonumber\\
&\times\sup
\frac{\big|\underline Z^1_{q'_1}\underline Z^2_{q'_2}
 \widetilde Y_{q_1} \widetilde X_{q_2}
\eta\big(q'_2\underline s^e(q_2),q'_1\underline s^e(q_1) \big)
\big|}{|X|_{k_2}|Y|_{k_1}|Z^2|_k|Z^1|_k}\;,
\end{align}
where the first supremum is over:
$$
(q_2b_2,q_1b_1,q'_2,q'_1)\in \S_2\times\S_1\times Q_2\times Q_1\;,
$$
the second over:
$$
(X,X',Y,Y',Z^2,Z^1)\!\in\!\CU_{k_2}(\q_2)\times\CU_{l_2}(\fY_2)\times \CU_{k_1}(\q_1)\times\CU_{l_1}(\fY_1)\times
\CU_k(\q_2)\times\CU_k(\q_1)\;,
$$
and the third over:
$$
(X,Y,Z^2,Z^1)\in\CU_{k_2}(\q_2)\times \CU_{k_1}(\q_1)\times
\CU_k(\q_2)\times\CU_k(\q_1)\;.
$$
Next, we observe:
\begin{align*}
&\underline Z^1_{q'_1}\underline Z^2_{q'_2}
 \widetilde Y_{q_1b_1} \widetilde Y'_{q_1b_1} \widetilde X_{q_2b_2}\widetilde X'_{q_2b_2}\,
F(q'_1q_1b_1 q'_2q_2b_2)=\\
&
\big(\widetilde{\Ad_{(q'_1q_1b_1 q'_2q_2b_2)^{-1}}(Z^1)}\,\widetilde{\Ad_{(q'_2q_2b_2)^{-1}}(Z^2)}\,
\widetilde {\Ad_{(q'_2q_2b_2)^{-1}}(YY') }\,\widetilde X\,\widetilde X'\,F\big)
(q'_1q_1b_1 q'_2q_2b_2)\;.
\end{align*}
This,  together with Lemma \ref{lemtwo},  entails that the $F$-dependent 
supremum in \eqref{NO}
is, up to a constant,  bounded by:
$$
 \|F\|_{2k+k_1+l_1+k_2+l_2}\,
 \fd_{\S_2\ltimes\S_1}(q'_1q_1b_1 q'_2q_2b_2)^k
 \,\fd_{\S_2\ltimes\S_1}(q'_2q_2b_2)^{k+k_1+l_1}\;,
 $$
 which by sub-multiplicativity of the modular weight, is bounded by:
 $$
  \|F\|_{2k+k_1+l_1+k_2+l_2}\,\fd_{\S_2\ltimes\S_1}(q'_2q_2b_2)^{2k+k_1+l_1}\,
 \fd_{\S_2\ltimes\S_1}(q'_1q_1b_1)^k\;.
$$
Now, by Theorem \ref{TASP}, the pair $\big((\S_2\ltimes\S_1)^2, S_{\rm can}^{\S_1}\oplus 1
+1\oplus  S_{\rm can}^{\S_1}\big)$ is tempered (and admissible and tame), so that by Lemma 
\ref{fd-tempered} the modular weight $\fd_{(\S_2\ltimes\S_1)^2}$ is tempered. 
Then, using the last statement of Lemma
\ref{fd-semi}, together with the methods of Lemmas \ref{lem:mult_inv_temp} 
and \ref{lem:j-algebra-extension-morphism-tempered}, 
we see that $\fd_{\S_2\ltimes\S_1}$ is tempered in (any) adapted coordinates 
(see Definition \ref{adapted}) for $\B=\S_2\ltimes\S_1$. 
This clearly implies that the restriction $\fd_{\S_2\ltimes\S_1}|_{\S_j}$,
$j=1,2$, is also tempered in the adapted coordinates for $\S_j$.
In view of the expressions \eqref{again}, we see that the adapted tempered coordinates and the 
coordinates associated to the one-point pair $(\S_j,\bS_{\S_j})$ (which is tempered by Theorem 
\ref{TP1P}) are related to one another through
a tempered diffeomorphism.
Hence, we deduce that $\fd_{\S_2\ltimes\S_1}|_{\S_j}$ is  tempered in the sense of the
one-point phase function too. Last,
using the explicit expression of the tempered weight $\fd_{\S_j}$ given in Lemma
\ref{fd-explicit}, we deduce that there exist $m_j\in\N$ and $C_j>0$, such that  
$\fd_{\S_2\ltimes\S_1}|_{\S_j}\leq C_j \fd_{\S_j}^{m_j}$, $j=1,2$.
Hence, the $F$-dependent 
supremum in \eqref{NO}
is, up to a constant,  bounded by
$$
  \|F\|_{2k+k_1+l_1+k_2+l_2}\,\fd_{\S_2}(q'_2q_2b_2)^{m_2(2k+k_1+l_1)}\,
 \fd_{\S_1}(q'_1q_1b_1)^{m_1k}\;.
$$

For $\eta$-dependent term in \eqref{NO}, we first note:
$$
\underline Z^1_{q'_1}\underline Z^2_{q'_2}
 \widetilde Y_{q_1} \widetilde X_{q_2}
\eta\big(q'_2\underline s^e(q_2),q'_1\underline s^e(q_1) \big)=
 \widetilde Y_{q_1} \widetilde X_{q_2}
\big(\underline Z^1\underline Z^2\eta\big)
\big(q'_2\underline s^e(q_2),q'_1\underline s^e(q_1) \big)\;,
$$
so that up to a redefinition of $\eta$, we can ignore the right-invariant vector fields.
Next, we observe that with $q=(a,n)$,  $q'=(a',n')$  in the coordinates \eqref{chartQ}, we have:
$$
q'\underline s^e(q)=\big(2a+a',e^{-2a}n'+2n\cosh a\big)\;.
$$
With $H$ the generator of $\a$ and $\{f_j\}_{j=1}^d$ a basis of $\l^\star$, the
associated  left-invariant vector fields on $Q$ read:
$$
\widetilde H=\partial_a-\sum_{j=1}^d n_j\partial_{n_j}\,,\quad \widetilde f_j=\partial_{n_j}\;.
$$
Choosing $\eta=\eta_2\otimes\eta_1$ with $\eta_j\in\CD(Q_j)$, it is enough to treat each variable
separately. So just assume that $\eta\in\CD(Q)$.
 Now, for $N=(N_1,\dots,N_d)\in\N^d$ with $|N|=k$, we have with $\widetilde f^N:=\widetilde f_1^{N_1}\dots \widetilde f_d^{N_d}$, $\partial^N_n:=\partial_{n_1}^{N_1}\dots\partial_{n_d}^{N_d}$ and setting $q_0:=q'\underline s^e(q)\in Q$:
$$
\widetilde f^N_q \eta(q_0)=2^k\cosh^{k}a\,\big(\partial_n^N\eta\big)(q_0)\;.
$$
Since $\cosh a\leq 2\cosh a'/2\cosh a_0/2$,
the latter  and Remark \ref{fd-explicit2} entail that
\begin{align*}
\big|\widetilde f^N_q \eta(q_0)\big|&\leq  C\cosh(a'/2)^{k}\cosh(a_0/2)^{k}
\big|\partial_n^N\eta\big|(q_0)
\\&\leq C\, \fd_Q(q')^k\cosh(a_0/2)^{k}\,\big|\partial_n^N\eta\big|(q_0)\;.
\end{align*}
On the other hand, we have
$$
\widetilde H_q\,\eta(q_0)=2\big(\partial_a\eta\big)(q_0)-2\sum_{j=1}^d(n_je^{-a}+n'_{j}e^{-2a})\big(\partial_{n_j}\eta\big)(q_0)\;.
$$
Since $w_j:=n_je^{-a}+n'_{j}e^{-2a}$ is an eigenvector of $\widetilde H_q$ with eigenvalue $-2$,
we deduce that for $k\in\N$, $\widetilde H_q^k\,\eta(q_0)$ is a linear combinations of the ordinary
derivatives of $\eta$, with coefficients given in  the ring  $\C[w_j]$ of order at most $k$.
Moreover, the rough estimate:
$$
|w|=|(w_1,\dots,w_d)|\leq 4\cosh a_0\cosh a'   (|n_0|+|n'|)\;,
$$
gives by Remark \ref{fd-explicit2}:
$$
\big|\widetilde H_q^k\,\eta(q_0)\big|\leq C\fd_Q(q')^{k}\,\cosh^ka_0|n_0|^k\big|P(\partial_a,\partial_{n_j})\eta\big|(q_0)\;,
$$
for a suitable polynomial $P$. This implies that the $\eta$-dependent term  in \eqref{NO} is,
up to a constant, bounded by:
$$
\fd_{Q_2}(q'_2)^{k_2}\,\fd_{Q_1}(q'_1)^{k_1}\,|\tilde \eta|\big(q'_2\underline s^e(q_2),
q'_1\underline s^e(q_1) \big)\;,
$$
where $\tilde\eta$ belongs to $\CD(Q_2\times Q_1)$ and is obtained from $\eta$ by multiplication
by $\cosh a,|n|$ and by differentiation along all its variables. 
Finally, we deduce (with $m_1,m_2\in\N$ fixed) that
\begin{align*}
&\|F^\eta\|_{k_2,l_2,k_1,l_1,k,j}\leq C\,\|F\|_{2k+k_1+l_1+k_2+l_2}
 \sup_{(q_2b_2,q_1b_1,q'_2,q'_1)\in \S_2\times\S_1\times Q_2\times Q_1} \\
&
 \frac{\fd_{Q_2}\!(q'_2)^{j+k_2}\!\fd_{Q_1}\!(q'_1)^{j+k_1}\!\fd_{\S_2}\!
 (q'_2q_2b_2)^{m_2(2k+k_1+l_1)}\!
 \fd_{\S_1}\!(q'_1q_1b_1)^{m_1k}\!|\tilde \eta|\!\big(\!q'_2\underline s^e(q_2),q'_1\underline s^e(q_1)
 \!  \big)
 }{\fd_{\S_2}(q_2b_2)^{n_2(k_1,l_1,k,j)}\,\fd_{Q_2}(q_2)^{2k_2}\,\fd_{\S_1}(q_1b_1)^{n_1(k,j)}\,
\fd_{Q_1}(q_1)^{2k_1}}
\;.
\end{align*}
Observe then that by Lemma \ref{fd-explicit} and Remark \ref{fd-explicit2}, we have 
$\fd_\S|_Q\leq C\fd_Q^2$.
Then, by the sub-multiplicativity and the invariance under the inversion map of the modular weights, 
we deduce that the fraction above is smaller than (a constant times):
\begin{align*}
 &\frac{ \fd_{Q_2}\big(\underline s^e(q_2)\big)^{j+k_2+2m_2(2k+k_1+l_1)}
 \fd_{Q_1}\big(\underline s^e(q_1)\big)^{j+k_1+2m_1k}
 \fd_{\S_2}( q_2b_2)^{m_2(2k+k_1+l_1)}\,
 }{\fd_{\S_2}(q_2b_2)^{n_2(k_1,l_1,k,j)}\,\fd_{Q_2}(q_2)^{2k_2}\,\fd_{\S_1}(q_1b_1)^{n_1(k,j)}\,
\fd_{Q_1}(q_1)^{2k_1}}\\
&\times  \fd_{\S_1}(q_1b_1)^{m_1k}
\fd_{Q_2}\big(q'_2\underline s^e(q_2)\big)^{j+k+2m_2(2k+k_1+l_1)}
 \fd_{Q_1}\big(q'_1\underline s^e(q_1 )\big)^{j+k_1+2m_1k}\\&\times
|\tilde \eta|\big(q'_2\underline s^e(q_2),q'_1\underline s^e(q_1) \big)\;.
\end{align*} 
Because $\tilde\eta$ is compactly supported, the expression in the last line above
 is smaller than a constant. Also, 
since $\underline s^e(a,n)=(2a,2n\cosh a)$, we deduce, by Remark \ref{fd-explicit2} again, that 
$\fd_Q\circ \underline s^e \leq C\fd_Q^2$. Thus, the expression above is bounded by
(a constant times):
$$
\frac{ \fd_{Q_2}(q_2)^{2j+4m_2(2k+k_1+l_1)}\,
 \fd_{Q_1}(q_1)^{2j+4m_1k}\,
 \fd_{\S_2}( q_2b_2)^{m_2(2k+k_1+l_1)}\,
   \fd_{\S_1}(q_1b_1)^{m_1k}
 }{\fd_{\S_2}(q_2b_2)^{n_2(k_1,l_1,k,j)}\,\fd_{\S_1}(q_1b_1)^{n_1(k,j)}}\;,
 $$
 and one concludes using Lemma \ref{fd-explicit} and Remark \ref{fd-explicit2}, which show that for 
 all $q\in Q$, $b\in\Y$, we have $\fd_Q(q)\leq\fd_\S(qb)$ and, by suitably choosing $n_1(k,j)$ and 
 $n_2(k_1,l_1,k,j)$.
\end{proof}
\begin{rmk}
\label{RM1}
Lemma \eqref{F-eta} admits a straightforward generalization for  symbols valued in
a Fr\'echet algebra.
Namely, if $\CE$ is a Fr\'echet algebra and $\underline\mu$ is a family of tempered weights
on $\B$, then the map $F\mapsto F^\eta$ is continuous  from $\CB^{\underline\mu}(\B,\CE)$ to  
$$
\CB^{\underline\mu_N,\nu_N}\big(\S_N,
\CB^{\underline\mu_{N-1},\nu_{N-1}}
\big(\S_{N-1},\dots\CB^{\underline\mu_1,\nu_1}\big(\S_1,
\CS(Q_N\times\dots\times Q_1,A)\big)\dots\big)\big)\;,
$$ 
where the $\underline\mu_j$'s now depend also of the restriction of $\underline\mu$
to $\S_j$ (and the $\nu_j$'s are unchanged). 
\end{rmk}

We are now ready to prove a non-Abelian (curved) and $C^*$-valued version of the Calder\'on-Vaillancourt 
estimate, the main result of this chapter.
\begin{thm}
\label{C-def}
Let $\B$ be a normal $\bf j$-group, $A$ a $C^*$-algebra faithfully represented
on a separable Hilbert space $\CH$, $F\in\CB(\B,A)$, $\eta
\in\CD(Q_N\times\dots\times Q_1)$ and $\bm\in{\bf\Theta}(\B)$.
Define for $x,y\in\B$, the element of $A$ given by
\begin{align}
\label{FQ}
&\langle\eta_{x},\Omega_{\theta,\bm}(F)\eta_y\rangle_{\CH_\chi}:=
\CF^{\overline \eta}
\Big(\widetilde{\int_{\S_1}\bE_\theta^{\S_1}\, \hat\bm_1}\dots
\Big(\widetilde{\int_{\S_N}\bE_\theta^{\S_N}\, \hat\bm_N}
\Big[g_N\in\S_N\mapsto
 \dots\nonumber \\
 &\qquad\qquad\qquad\qquad \Big[ g_1\in \S_1\mapsto
\big(L^\star_{y^{-1}}F\big)^\eta(g_N,\dots,g_1;.)\Big]\dots\Big]\Big)\dots\Big)
(y^{-1}x)\;,
\end{align}
where $\hat\bm(qb)=\bm(q^{-1})$ and where $\bE_\theta^{\S_j}$ is the one-point phase of
$\S_j$ as defined in \eqref{bS}.  Then we have:
\begin{align*}
\sup_{x\in\B}\int_\B\big\|\langle\eta_{x},\Omega_{\theta,\bm}(F)\eta_y\rangle_{\CH_\chi}\big\|\,{\rm d}_\B(y)<\infty\,,\,\,
\sup_{y\in\B}\int_\B\big\|\langle\eta_{x},\Omega_{\theta,\bm}(F)\eta_y\rangle_{\CH_\chi}\big\|\,{\rm d}_\B(x)<\infty\;.
\end{align*}
Consequently (see Proposition \ref{CV-pre-born}), the operator $\Omega_{\theta,\bm}(F)$ on $\CH_\chi\otimes\CH$ defined by means of  the quadratic form
$$
\Psi,\Phi\in\CH_\chi\otimes\CH\mapsto
\int_{\B\times\B}\big\langle\langle\eta_x,\Psi\rangle_{\CH_\chi},\langle\eta_x,\Omega_{\theta,\bm}(F)\eta_y\rangle_{\CH_\chi}\langle\eta_y,\Psi\rangle_{\CH_\chi}\big\rangle_\CH\,{\rm d}_\B(x)\,{\rm d}_\B(y)\;,
$$
is bounded. Moreover, there exists $k\in\N$ (depending only on $\dim \B$ and on the order
of the polynomial in $\fd_{Q_N}\otimes\dots\otimes\fd_{Q_1}$ that majorizes $|\bm|$)
and $C>0$, such that for all $F\in\CB(\B,A)$ we have
\begin{align*}
\|\Omega_{\theta,\bm}(F)\|\leq C\,\|F\|_{k,\infty}=C\,
\sup_{X\in\CU_k(\b)}\sup_{x\in\B}\big\|\widetilde X\, F(x)\big\|\;.
\end{align*}
\end{thm}
\begin{proof}
To simplify the notation for the matrix element given in \eqref{FQ},
we write
$$
\langle\eta_{x},\Omega_{\theta,\bm}(F)\eta_y\rangle_{\CH_\chi}=
 \CF^{\overline \eta}
\Big(\widetilde{\int_\B\bE_\theta^\B\, \hat\bm}\big[z\mapsto\,\big(L^\star_{y^{-1}}F\big)^\eta(.,z)
\big]\Big)
(y^{-1}x) \;,
$$
where $\bE^\B_\theta$ is given in \eqref{bEB}.
Observe that this notation is coherent with our Fubini type Theorem \ref{Fubini}.
Thus,
\begin{align*}
&\sup_{y\in\B}\int_\B\big\|\langle\eta_{x},\Omega_{\theta,\bm}(F)\eta_y\rangle_{\CH_\chi}\big\|\,{\rm d}_\B(x)\\&\qquad=
\sup_{y\in\B}\int_\B\big\| \CF^{\overline \eta}
\Big(\widetilde{\int_\B\bE^\B_\theta\, \hat\bm}
\big[z\mapsto\,\big(L^\star_{y^{-1}}F\big)^\eta(.,z)\big]\Big)
(y^{-1}x)  \big\|\,{\rm d}_\B(x)\\
&\qquad=
\sup_{y\in\B}\int_\B\big\| \CF^{\overline \eta}
\Big(\widetilde{\int_\B\bE^\B_\theta\, \hat\bm}
\big[z\mapsto\,\big(L^\star_{y^{-1}}F\big)^\eta(.,z)\big]\Big)
(x)  \big\|\,{\rm d}_\B(x)\;.
\end{align*}
The fact that this expression is finite follows then by combining  Propositions \ref{first-type} and
 \ref{CS}
with  Theorem \ref{OIEXT} and Lemma \ref{F-eta} and the fact that 
 $L^\star_{y^{-1}}$ maps $\CB(\B,A)$ to itself isometrically.
The second case is similar since
$$
\langle\eta_{x},\Omega_{\theta,\bm}(F)\eta_y\rangle_{\CH_\chi}^*
=\langle\eta_{y},\Omega_{\theta,\underline\sigma^\star\overline{\bm}}(F^*)\eta_x\rangle_{\CH_\chi}\;.
$$
The final estimation we give is a consequence of Proposition \ref{CV-pre-born} together with the
estimates underlying Lemma \ref{CS}, Theorem \ref{OIEXT} and Lemma \ref{F-eta}.
\end{proof}

\begin{rmk}
\label{RM2}
In view of Remark \ref{RM1}, on may wonder what happens in Theorem \ref{C-def} when 
one choses a symbol  in $\CB^\mu(\B,A)$  instead of 
a symbol in $\CB(\B,A)$. So let $F\in \CB^\mu(\B,A)$, with $\mu$ a tempered weight. Then,
 from the same argument than those of Theorem \ref{C-def} (using 
 by Remark \ref{RM1}, instead of Lemma \ref{F-eta}),  one deduces  that 
\begin{align}
\label{tg}
 \int_\B\big\|\langle\eta_{x},\Omega_{\theta,\bm}(F)\eta_y\rangle_{\CH_\chi}\big\|\,{\rm d}_\B(x)
 <\infty\;.
 \end{align}
However, as the left regular action is no longer isometric on $\CB^\mu(\B,A)$
(see for instance the second item of  Lemma \ref{SmoothFamily}), there is no reason to 
expect that the supremum over $y\in\B$ of the expression given in \eqref{tg} to be finite.
Accordingly (and as expected), there is no chance for the operator $\Omega_{\theta,\bm}(F)$ to be bounded
when $F$ belongs to $\CB^\mu(\B,A)$ with unbounded $\mu$.
\end{rmk}

\section{The deformed $C^*$-norm}
Now, we assume that our $C^*$-algebra $A$ is equipped with a strongly continuous  and isometric action
 $\alpha $ of a normal $\bf j$-group $\B$. 
We stress that the results of this chapter cannot hold true in the more general context of 
tempered actions. This is the main difference between the deformation theory at the level of Fr\'echet
and $C^*$-algebras.
Given an element $a\in A$, we construct as usual the $A$-valued function $\alpha(a)$ on
 $\B$:
 $$
 \alpha(a):=[g\in\B\mapsto\alpha_g(a)\in A]\;.
 $$
 Thus, from Theorem \ref{UDF}, we can deform the Fr\'echet algebra
 structure on the set of smooth vectors $A^\infty$ by means of the deformed product 
 $$
 a\star_{\theta,\bm}^\alpha b:=\big(\alpha(a)\star_{\theta,\bm}\alpha(b)\big)(e)\,,
 \quad a,b\in A^\infty\;.
 $$ 
 We have seen in \eqref{invol} how to modify the original 
 involution at the level of $\CB(\B,A)$. At the level of the Fr\'echet algebra $A^\infty$,
an obvious observation leads to:

 \begin{lem}
 \label{continvol}
 Let $\B$ be a normal $\bf j$-group.
For $\bm\in{\bf\Theta}(\B)$, the following defines a continuous involution of the Fr\'echet
algebra
$(A^\infty,\star_{\theta,\bm}^\alpha)$:
 \begin{equation*}
\ast_{\theta,\bm}:A^\infty\to A^\infty\,,\quad a\mapsto \frac{{\bf m_N}}{\underline s_e^\star \overline\bm_N }
\Big(\tfrac12\arcsinh(\tfrac i\theta E_N^\alpha)\Big)\dots
\frac{{\bf m_1}}{\underline s_e^\star \overline\bm_1 }
\Big(\tfrac12\arcsinh(\tfrac i\theta E_1^\alpha)\Big) a^*\;,
\end{equation*}
where $E_N,\dots,E_1$ are the central elements of the Heisenberg Lie
algebras attached to each elementary factors of $\B$.
\end{lem}
\begin{rmk}
Note that when $\underline s_e^\star \overline\bm_j=\bm_j$, $j=1,\dots,N$, there is no modification
of the involution.
\end{rmk}
The construction of a pre-$C^*$-structure on $(A^\infty,\star_{\theta,\bf m}^\alpha)$
follows then from Theorem \ref{C-def} and from the following immediate result (compare with
Lemma \ref{embed}):
 \begin{lem}
 \label{V}
 Let $(A,\alpha,\B)$ be a $C^*$-algebra endowed with a strongly continuous and isometric
 action of a normal $\bf j$-group.
 Then, we have an isometric equivariant embedding $\alpha:A^\infty\to \CB(\B,A)$.
 \end{lem}
 \begin{proof}
 The equivariance property of $\alpha$ is obvious and implies (with the fact that 
   $\alpha$ is an isometric action of $\B$ on $A$) that for any $k\in\N$:
 \begin{align*}
 \|\alpha(a)\|_{k,\infty}&=\sup_{g\in\B}\sup_{X\in\,\CU_k(\b)}\frac{\|\widetilde X_g\,\alpha_g(a)\|}{|X|_k}
\\&
 =\sup_{g\in\B}\sup_{X\in\,\CU_k(\b)}\frac{\|\alpha_g( X^\alpha\,a)\|}{|X|_k}=
 \sup_{X\in\,\CU_k(\b)}\frac{\|X^\alpha\,a\|}{|X|_k}=
 \| a\|_k\;,
\end{align*}
 and the proof follows.
 \end{proof}

Recall that $\CH$ is any separable Hilbert space carrying a faithful representation of
 our $C^*$-algebra $A$, which therefore, is identified with a $C^*$-subalgebra of
 $\CB(\CH)$. Then by the previous lemma and Theorem \ref{C-def}, we deduce that the map
 $$
 a\in A^\infty\mapsto  \|\Omega_{\theta,\bm}\big(\alpha(a)\big)\|\;,
 $$
 takes finite values. It is also important to observe that the norm above is by construction
 the operator norm on $\CH_\chi\otimes\CH$, that is the spatial $C^*$-norm on 
 $A\otimes\CB(\CH_\chi)$. For future use, we also recall that the spatial $C^*$-norm is the minimal 
 $C^*$-cross-norm. (We invite  the reader unfamiliar with tensor products of $C^*$-algebras
 to consult, for example, Appendix B of \cite{RW}.) More precisely,
 combining the Lemma \ref{V} with Theorem \ref{C-def}, we deduce the following 
 inequality:
 \begin{cor}
 \label{Estim-norm}
 Let $(A,\alpha,\B)$ be a $C^*$-algebra endowed with a strongly continuous action
 of a normal $\bf j$-group and $\bm\in{\bf\Theta}(\B)$.
 Then, there exists $k\in\N$ and $C>0$ such that
  for any  $a\in A^\infty$, we have:
  $$
  \|\Omega_{\theta,\bm}\big(\alpha(a)\big)\|\leq C\|a\|_k:=C\sup_{X\in\,\CU_k(\b)}\Big\{\frac{\|X^\alpha\,a\|}
  {|X|_k}\Big\}\;.
  $$
 \end{cor}

\begin{prop}
\label{Cnorm1}
Let $\bm\in{\bf\Theta}(\B)$. Then, the following defines 
 a $C^*$-norm on the involutive deformed Fr\'echet algebra 
 $(A^\infty,\star^\alpha_{\theta,\bm},\ast_{\theta,\bm})$:
$$
a\in A^\infty\mapsto \|a\|_{\theta,\bm}:=\big\|\Omega_{\theta,\bm}\big(\alpha(a)\big)\big\|\;,
$$
where the 
operator $\Omega_{\theta,\bm}\big(\alpha(a)\big)$ is defined in Theorem \ref{C-def}.
Accordingly, we let $A_{\theta,\bm}$ be the $C^*$-completion of $A^\infty$ that we 
abusively call the $C^*$-deformation of $A$.
\end{prop}
\begin{proof}
By construction we have for all $a,b\in A^\infty$
$$
\Omega_{\theta,\bm}\big(\alpha(a^{\ast_{\theta,\bm}}\star_{\theta,\bm}^\alpha b)\big)=\Omega_{\theta,\bm}\big(\alpha(a)\big)^*
\Omega_{\theta,\bm}\big(\alpha(b)\big)\;,
$$
and the claim follows immediately. 
\end{proof}

\begin{rmk}
We already know that at the level of the deformed 
pre-$C^*$-algebra $(A^\infty,\star_{\theta,\bm}^\alpha)$, 
the action of the group $\B$  is no longer by automorphism. But at the level of the deformed
$C^*$-algebra there is no action of $\B$ at all.
\end{rmk}

In a way very analogous  to Proposition \ref{defofdef}, we can 
show that the $C^*$-deforma\-tion associated with a normal $\bf j$-group coincides with the iterated 
$C^*$-deformations of each of its elementary normal subgroups. To see this, fix
$\B$ be a normal $\bf j$-group with Pyatetskii-Shapiro decomposition $\B=\B'\ltimes\S$
and $A$ a $C^*$-algebra endowed with a strongly continuous and isometric action 
$\alpha$ of $\B$. Of course, 
$\alpha^\S$, the restriction  of $\alpha$ to the subgroup $\S$, is strongly continuous on $A$. 
Let us fix also $\tilde\bm=\bm'\otimes\bm\in{\bf {\bf\Theta}}(\B)$, with $\bm'\in{\bf {\bf\Theta}}(\B')$
and $\bm\in{\bf {\bf\Theta}}(\S)$. Then, we can 
perform the $C^*$-deformation of $A$ by means of the action of $\S$. We call this deformed 
$C^*$-algebra $A_{\theta,\bm}^\S$. Then  $\B'$ acts strongly continuously by $*$-homomorphisms
on $A_{\theta,\bm}^\S$. Indeed, it has been shown in the proof of Proposition \ref{defofdef}
that the subspace of smooth vectors for $\B$ coincides with the subspace
of smooth vectors for $\B'$ within the subspace of smooth vectors for $\S$. In turns,
$A^\infty$, the set of smooth vectors for $\B$ on $A$, is dense in  $A_{\theta,\bm}^\S$.
As the action of $\B'$ is (obviously) strongly continuous and by $\ast$-homomorphisms
(as shown  in the proof of  Proposition \ref{defofdef} too) on $A^\infty$, a density argument yields
the result. Thus, we can perform the $C^*$-deformation of $A_{\theta,\bm}^\S$ by means of the
 action of $\B'$. 
We call this deformed $C^*$-algebra $(A_{\theta,\bm}^\S)_{\theta,\bm'}^{\B'}$. But we could also
perform the $C^*$-deformation of $A$ by means of the action of $\B$ directly. 
We call this deformed $C^*$-algebra $A_{\theta,\tilde\bm}^\B$. Now, the precise result of 
Proposition \ref{defofdef}, is that at  the level of the (common) dense subspace
$A^\infty$, both constructions coincide. Thus it suffices to show that the $C^*$-norms
of  $(A_{\theta,\bm}^\S)_{\theta,\bm'}^{\B'}$ and  $A_{\theta,\tilde\bm}^\B$ coincide  on 
$A^\infty$.
But this easily follows from our construction. Indeed, the $C^*$-norm of 
$(A_{\theta,\bm}^\S)_{\theta,\bm'}^{\B'}$ at the level of $A^\infty$, is by definition the map
$$
a\mapsto\big\|\Omega^{\B'}_{\theta, \bm'}\big(\big[z'\in\B'\mapsto\Omega_{\theta,\bm}^\S\big([
z\in\S\mapsto\alpha_{zz'}(a)]\big)\big]\big)\big\|\;.
$$
But by the construction of  section \ref{QNG}, we precisely have
$$
\Omega^{\B'}_{\theta,\bm'}\big(\big[z'\in\B'\mapsto\Omega_{\theta,\bm}^\S\big([
z\in\S\mapsto\alpha_{zz'}(a)]\big)\big]\big)=\Omega^\B_{\theta, \tilde\bm}\big(\alpha(a)\big)\;.
$$
Thus, we have proved the following:

\begin{prop}
Let $\B$ be a normal $\bf j$-group with Pyatetskii-Shapiro decomposition $\B=\B'\ltimes\S$, where
$\B'$ is a normal $\bf j$-group and  $\S$  is an elementary normal $\bf j$-group. 
Let $A$ be a $C^*$-algebra 
endowed with a strongly continuous isometric action $\alpha$ of 
$\B$. Within the notations displayed above, we have:
$$
A_{\theta,\tilde\bm}^\B=(A_{\theta,\bm}^\S)_{\theta,\bm'}^{\B'}\;.
$$
\end{prop}

In the remaining part of this section we prove that the deformed $C^*$-norm constructed above
coincide  with the $C^*$-norm of bounded and adjointable operators on a $C^*$-module.
This will make clearer the analogies with the construction of Rieffel in \cite{Ri} for the Abelian  
case and it also explains the choice of the spatial tensor product in
Theorem \ref{C-def}.

\begin{dfn}
Let $\bm\in{\bf\Theta}(\B)$. Then, for $f_1,f_2\in\CS^{S_{\rm can}}(\B,A)$, we define the $A$-valued pairing:
\begin{align}
\label{pairing}
\langle f_1,f_2\rangle_{\theta,\bf m}
:=\int_\B \big\langle\eta_x,\Omega_{\theta,\bm}(f^{\ast_{\theta,\bm}}_1\star_{\theta,\bm}f_2)\eta_x\big\rangle\,
{\rm d}_\B(x)\;,
\end{align}
where where $\{\eta_x\}_{x\in\B}\subset\CH_\chi$ is the family of coherent states given in
Definition \ref{CoheState} and the involution $\ast_{\theta,\bm}$ on 
$\CS^{S_{\rm can}}(\B,A)$ is defined by:
$$
\ast_{\theta,\bm}:f\mapsto \frac{{\bf m}_N}{\underline s_e^\star \overline\bm_N }
\Big(\tfrac12\arcsinh(\tfrac i\theta \widetilde E_N)\Big)\circ\dots\circ
\frac{{\bf m_1}}{\underline s_e^\star \overline\bm _1}
\Big(\tfrac12\arcsinh(\tfrac i\theta 
\widetilde E_1)\Big)  f^*\;.
$$
In the last formula, $E_N,\dots,E_1$ denote the central elements in each Heisenberg Lie
algebra attached to each elementary components in $\B$ and $f^*:=[x\in\B\mapsto f(x)^*]\in 
\CS^{S_{\rm can}}(\B,A)$.
\end{dfn}

\begin{prop}
Endowed with the pairing \eqref{pairing} and action 
$$
 \CS^{S_{\rm can}}(\B,A)\times A^\infty\to \CS^{S_{\rm can}}(\B,A)\,,\quad  (f,a)\mapsto \big[g\in\B\mapsto f(g)a\big]\;,
 $$
  the space $\CS^{S_{\rm can}}(\B,A)$ becomes a (right) pre-$C^*$-module for the $C^*$-algebra $A$.
\end{prop}
\begin{proof}
We need first to show that the pairing \eqref{pairing} is well defined. Also, to lighten a little bit
the notations, we assume that the normal $\bf j$-group $\B$ contains only two elementary
factors, i.e$.$ $\B=\S_2\ltimes\S_1$ with $\S_2,\S_1$ elementary normal $\bf j$-groups. 
(The proof for a 
generic normal
$\bf j$-group with an arbitrary number of elementary factors has no essential 
supplementary  difficulties.)
Take first an element $F\in\CB^{\underline\mu}(\B,\CE)$, where $\CE$ is  any Fr\'echet space 
and $\underline\mu$ any family of tempered weights. In this situation (which is slightly more
general than the situation of Theorem  \ref{C-def}), by analogy with 
 equation \eqref{FQ}, it is natural to define:
\begin{align*}
&\langle\eta,\Omega_{\theta,\bm}(F)\eta\rangle_{\CH_\chi}:=\\
&\qquad\CF^{\overline \eta}
\Big(\widetilde{\int_{\S_1}\bE_\theta^{\S_1}\, \hat\bm_1}
\Big(\widetilde{\int_{\S_2}\bE_\theta^{\S_2}\, \hat\bm_2}
\Big[g_2\in\S_2\mapsto
\Big[ g_1\in \S_1\mapsto
F^\eta(g_2,g_1;.)\Big]\Big]\Big)\Big)
(e)\;.
\end{align*}
By a slight adaptation of Lemma \ref{F-eta} (see Remark \ref{RM1}),
we deduce that 
$$
F^\eta\in\CB^{\underline\mu_2,\nu_2}\big(\S_2,\CB^{\underline\mu_1,\nu_1}
\big(\S_1,\CS(Q_2\times Q_1,\CE\big)\big)\;, 
$$
for suitable tempered weights. Now, Theorem \ref{OIEXT} and Proposition \ref{first-type}
entails that $\langle\eta,\Omega_{\theta,\bm}(F)\eta\rangle_{\CH_\chi}\in\CE$.
 This observation being made,
we further remark that  for every  $f\in\CS^{S_{\rm can}}(\B,A)$, the element 
$\dot f\,\in C^\infty(\B,\CS^{S_{\rm can}}(\B,A))$
defined by $\dot{f}(x):=[y\mapsto f(xy)]$ actually lives in 
$\CB^{\underline{\mu}}(\B,\CS^{S_{\rm can}}(\B,A))$ for a suitable family of tempered weights
 $\underline{\mu}$. Applying the preceding reasoning for the Fr\'echet space
 $\CE=\CS^{S_{\rm can}}(\B,A)$, we deduce that the element  
 $\langle\eta\,,\,\Omega_{\theta,\bm}(\dot{f})\eta\rangle$ lives in
  $\CS^{S_{\rm can}}(\B,A)$. Using the $\B$-equivariance of the quantization map 
  $\Omega_{\theta,\bm}$, we then notice that the value at $x\in\B$
of the above element $\langle\eta\,,\,\Omega_{\theta,\bm}(\dot{f})\eta\rangle$ 
equals $\langle\eta_x\,,\,\Omega_{\theta,\bm}(f)\eta_x\rangle$. Which we deduce from that the 
matrix coefficient $[x\mapsto\langle\eta_x\,,\,\Omega_{\theta,\bm}(f)\eta_x\rangle]$ 
belongs to $\CS^{S_{\rm can}}(\B,A)$.
Hence,
$$
\int_\B \big\|\langle\eta_x,\Omega_{\theta,\bm}(f)\eta_x\rangle\big\|\,{\rm d}_\B(x)<\infty\;.
$$
Thus, we conclude
from the stability of $\CS^{S_{\rm can}}(\B,A)$
 under $\star_{\theta,\bf m}$ (see Proposition
\ref{universal-schwartz}),
that the pairing $\langle.,.\rangle_{\theta,\bf m}$ is well defined.

Testing this pairing on the dense subset ${\Span}\{a\vf\,,a\in A\,,\,\vf\in\CS^{S_{\rm can}}(\B)\}$
 of $\CS^{S_{\rm can}}(\B,A)$, we see that 
 $$
 \langle\CS^{S_{\rm can}}(\B,A),\CS^{S_{\rm can}}(\B,A)\rangle_{\theta,\bf m}=A.A\;,
 $$
 which is dense in $A$. Next, we observe that the pairing can be rewritten as
 $$
 \langle f_1,f_2\rangle_{\theta,\bf m}
=\int_\B \big\langle\eta_x,\Omega_{\theta,\bm}(f_1)^{\ast}\,
\Omega_{\theta,\bm}(f_2)\eta_x\big\rangle\,
{\rm d}_\B(x)\;.
$$
This shows  that 
$\langle f_1,f_2\rangle_{\theta,\bf m}^*=\langle f_2,f_1\rangle_{\theta,\bf m}$ and proves positivity and 
non-degeneracy.
Last, it is clear that  $\langle f_1,f_2\rangle_{\theta,\bf m}a= \langle f_1,f_2a\rangle_{\theta,\bf m}$ 
 for all $a\in A$ and all $\CS^{S_{\rm can}}(\B,A)$. 
\end{proof}

\begin{rmk}
It can be shown that the pairing can be rewritten as:
$$
\langle f_1,f_2\rangle_{\theta,\bf m}=\int_\B f_1^{\ast_{\theta,\bm}}\star_{\theta,\bf m}f_2(g)\,
{\rm d}_\B(g)={\rm Tr}\big(\Omega_{\theta,\bm}(f_1)^*\,\Omega_{\theta,\bm}(f_2)\big)\;.
$$
However, this is by far  less convenient expressions, as shown  in the proof of  Theorem 
\ref{CequalsC} below.
\end{rmk}

\begin{dfn}
Let $\bm\in{\bf\Theta}(\B)$.  For  $F\in\CB(\B,A)$, let $L^{\theta,\bf m}(F)$ be the operator on 
$\CS^{S_{\rm can}}(\B,A)$ given by 
$$L^{\theta,\bf m}(F)f=F\star_{\theta,\bf m} f\;.$$
\end{dfn}
By Proposition \ref{universal-schwartz}, the operator $L^{\theta,\bf m}(F)$, $F\in\CB(\B,A)$,  acts 
continuously on $\CS^{S_{\rm can}}(\B,A)$. Moreover, $L^{\theta,\bf m}(F)$ is adjointable, with  adjoint given by 
$L^{\theta,\bf m}(F^{\ast_{\theta,\bm}})$. Indeed, for all $f_1,f_2\in \CS^{S_{\rm can}}(\B,A)$ and $F\in\CB(\B,A)$, 
we have
\begin{align*}
\langle f_1,L^{\theta, \bf m}(F)f_2\rangle_{\theta,\bf m}
&=\int_\B \big\langle\eta_x,\Omega_{\theta,\bm}(f^{\ast_{\theta,\bm}}_1\star_{\theta,\bm}
F\star_{\theta,\bm} f_2)\eta_x\big\rangle\,
{\rm d}_\B(x)\\
&=\int_\B \big\langle\eta_x,\Omega_{\theta,\bm}\big((F^{\ast_{\theta,\bm}}\star_{\theta,\bm }
f_1)^{\ast_{\theta,\bm}}
\star_{\theta,\bm} f_2\big)\eta_x\big\rangle\,
{\rm d}_\B(x)\\&
=\langle L^{\theta,\bf m}(F^{\ast_{\theta,\bm}})f_1,f_2\rangle_{\theta,\bf m}\;.
\end{align*}
Note also that the operators $L^{\theta,\bf m}(F)$ all commute with the right-action of $A$.
But we have more, since in fact  $L^{\theta,\bf m}(F)$, for $F\in\CB(\B,A)$, belongs to the 
$C^*$-algebra of 
$A$-linear adjointable endomorphisms of the pre-$C^*$-module  $\CS^{S_{\rm can}}(\B,A)$. Indeed,
from the operator inequality on $ \CB(\CH_\chi)\otimes A$
\begin{align*}
\Omega_{\theta,\bm}(f^{\ast_{\theta,\bm}}\star_{\theta,\bf m}F^{\ast_{\theta,\bm}} 
\star_{\theta,\bf m}F\star_{\theta,\bf m}f)
&=\Omega_{\theta,\bm}(f)^\ast\big|\Omega_{\theta,\bm}(F)\big|^2\,\Omega_{\theta,\bm}(f)\\&
\leq \|\Omega_{\theta,\bm}(F)\|^2\Omega_{\theta,\bm}(f)^\ast\,\Omega_{\theta,\bm}(f)\\&\quad
=\|\Omega_{\theta,\bm}(F)\|^2\Omega_{\theta,\bm}(f^{\ast_{\theta,\bm}}\star_{\theta,\bf m} f)\;,
\end{align*}
we deduce for $F\in\CB(\B,A)$ and $f\in\CS^{S_{\rm can}}(\B,A)$, the operator inequality on $A$:
$$ 
\langle L^{\theta,\bf m}(F)f,L^{\theta,\bf m}(F)f\rangle_{\theta,\bf m}\leq  \|\Omega_{\theta,\bm}(F)\|^2 \langle f,f\rangle_{\theta,\bf m}\;.
$$
Hence we get
\begin{equation}
\label{CC}
\|L^{\theta,\bf m}(F)\|\leq \|\Omega_{\theta,\bm}(F)\|\;,
\end{equation}
where the norm on the left hand side denotes the norm of  the $C^*$-algebra of 
$A$-linear adjointable endomorphisms of the pre-$C^*$-module  $\CS^{S_{\rm can}}(\B,A)$.
Now, observe the dense embedding of the algebraic tensor product
$$
\CB(\B)\otimes_{\rm alg} A\to\CB(\B,A)\,,\qquad \sum_i \phi_i\otimes a_i\mapsto\big[g\in\B
\mapsto \sum_i \phi_i(g)a_i\in A\big]\;.
$$
Via this embedding, the norm on the right hand side of \eqref{CC} is by construction the restriction
to $\CB(\B)\otimes_{\rm alg} A$ of the minimal
(spatial) $C^*$-norm on $B\otimes_{\rm alg}A$, where 
 $B$ is the $C^*$-completion of $\{\Omega_{\theta,\bm}(F),\,F\in\CB(\B)\}$
in $\CB(\CH_\chi)$. Hence, we deduce that
$$
\|L^{\theta,\bf m}(F)\|\geq \|\Omega_{\theta,\bm}(F)\|\,,\qquad\forall F\in \CB(\B)\otimes_{\rm alg}A \;,
$$
which by density implies that 
\begin{equation*}
\|L^{\theta,\bf m}(F)\|\geq \|\Omega_{\theta,\bm}(F)\|\,,\qquad\forall F\in \CB(\B,A) \;.
\end{equation*}
Thus we have proved the following:
\begin{thm}
\label{CequalsC}
Let $\B$ be a normal $\bf j$-group, $A$ a $C^*$-algebra and $\bm\in{\bf\Theta}(\B)$. Then
$$
\|L^{\theta,\bf m}(F)\|=\|\Omega_{\theta,\bm}(F)\|\,,\qquad \forall F \in \CB(\B,A)\;,
$$
where the norm on the left hand side is the  one  of  the $C^*$-algebra of 
$A$-linear adjointable endomorphisms of the pre-$C^*$-module  $\CS^{S_{\rm can}}(\B,A)$
and the norm on the right hand side is the spatial $C^*$-norm of $ \CB(\CH_\chi)\otimes A$.
\end{thm}

Back to the case where $A$ carries a strongly continuous isometric action $\alpha$, we deduce:

\begin{prop}
\label{same-norms}
Let $(A,\alpha)$ be a $C^*$-algebra endowed with a strongly continuous and isometric action
of  a normal $\bf j$-group $\B$ and let $\bm\in{\bf\Theta}(\B)$. Then, the
  $C^*$-norm  on the involutive Fr\'echet algebra 
  $(A^\infty,\star^\alpha_{\theta,\bm},\ast_{\theta,\bm})$ given by
$$
a\in A^\infty\mapsto \big\|L^{\theta,\bm}\big(\alpha(a)\big)\big\|\;,
$$
coincides with the deformed norm $\|.\|_{\theta,\bm}$ of Proposition \ref{Cnorm1}.
\end{prop}

\section{Functorial properties of the deformation}

In this section, we collect the main functorial properties of the deformation.
We still consider a $C^*$-algebra $A$, endowed with a strongly continuous and isometric action
$\alpha$ of a normal $\bf j$-group $\B$. Given an element $\bm\in{\bf\Theta}(\B)$, we form $A_{\theta,\bm}$, 
the
$C^*$-deformation of $A$. 
We let $ \overline{\CB_{\theta,\bm}}(\B,A)$  and $ \overline{\CS_{\theta,\bm}}(\B,A)$ be the 
$C^*$-completion of the pre-$C^*$-algebras:
 $$
 \big( \CB(\B,A),\star_{\theta,\bm},\ast_{\theta,\bm}\big)\quad\mbox{and}\quad
 \big(\CS^{S_{\rm can}}(\B,A),\star_{\theta,\bm},\ast_{\theta,\bm}\big)\;,
 $$ 
 for the (deformed) $C^*$-norm
$$
 F\mapsto \|\Omega_{\theta,\bm}(F)\|\;.
 $$ 
Firstly, we observe from Proposition \ref{Relation-Weyl2}, the following  isomorphism:
 \begin{lem}
 \label{compact}
  Let $A$ be a $C^*$-algebra  and $\bm\in{\bf\Theta}(\B)$. 
Then we have:
 $$
 \overline{\CS_{\theta,\bm}}(\B,A)\simeq \CK(\CH_\chi)\otimes A\;.
 $$
 \end{lem}

Now, we come to the question of  bounded approximate units for the deformed $C^*$-algebra
$A_{\theta,\bm}$. Since $A$ possesses a bounded approximate unit (as any $C^*$-algebra does),
Proposition \ref{BAU} shows  that the pre-$C^*$-algebra $(A^\infty,\star_{\theta,\bm}^\alpha,\ast_{\theta,\bm})$
possesses a bounded approximate unit as well. Thus, we deduce from Corollary \ref{Estim-norm}:\begin{prop}
 Let $(A,\alpha)$ be a $C^*$-algebra endowed with a strongly continuous action
 of a normal $\bf j$-group $\B$ and $\bm\in{\bf\Theta}(\B)$.  Then $A_{\theta,\bm}$
possesses a  bounded approximate unit $\{e_\lambda\}_{\lambda\in\Lambda}$ 
consisting of elements of $A^\infty$.
\end{prop}

 Next, we observe that   the two-sided
  ideal $(\CS^{S_{\rm can}}(\B,A),\star_{\theta,\bm},\ast_{\theta,\bm})$ is essential
 in  $(\CB(\B,A),\star_{\theta,\bm},\ast_{\theta,\bm})$:
 \begin{prop}
  Let $A$ be a $C^*$-algebra  and $\bm\in{\bf\Theta}(\B)$.
  Then,
  the ideal $\CS^{S_{\rm can}}(\B,A)$ is essential in the pre-$C^*$-algebra $\CB(\B,A)$,
  that is to say
   we have for all $F\in\CB(\B,A)$:
  $$
 \|\Omega_{\theta,\bm}(F)\|=\sup\big\{\|\Omega_{\theta,\bm}(F\star_{\theta,\bm} f)\|\,:\,f\in \CS^{S_{\rm can}}(\B,A)\,,\;
 \|\Omega_{\theta,\bm}(f)\|\leq1\big\}\;.
 $$
 \end{prop}
  \begin{proof}
 This is verbatim the arguments of \cite[Proposition 4.11]{Ri}, combined with the equality
 $ \|\Omega_{\theta,\bm}(F)\|=\|L^{\theta,\bm}(F)\|$ of Proposition \ref{same-norms}, 
 for all $F\in\CB(\B,A)$ (thus  for $f\in\CS^{S_{\rm can}}(\B,A)$ too) and with the existence of bounded
 approximate units of the pre-$C^*$-algebra $(\CS^{S_{\rm can}}(\B,A),\star_{\theta,\bm})$ as shown in 
 Proposition \ref{BAU}.
 \end{proof}

 The proof of the next two results is word for word the one of the corresponding results  
 in the flat situation, given in  \cite[Proposition 4.12 and Proposition 4.15]{Ri}.
 \begin{prop}
  Let $A$ be a $C^*$-algebra, $I$ an essential ideal of $A$  and $\bm\in{\bf\Theta}(\B)$.  
  Then the $C^*$-norm on $\big(\CB(\B,A),\star_{\theta,\bm}\big)$ given in Theorem \ref{C-def} is the 
  same 
  as the $C^*$-norm of Proposition \ref{same-norms} for the restriction of the action of $\CB(\B,A)$ 
  on $\CS^{S_{\rm can}}(\B,I)$. 
  \end{prop}
 \begin{prop}
 Let $A$ be a $C^*$-algebra  and $\bm\in{\bf\Theta}(\B)$. 
  The $C^*$-algebra $ \overline{\CB_{\theta,\bm}}(\B,A)$ is isomorphic to the $C^*$-deformation
 of the algebra of $A$-valued right uniformly continuous and bounded functions
 on $\B$, $C_{ru}(\B,A)$, for the right regular action of $\B$.
 \end{prop}

The following two results treat the question of morphisms and ideals. They can be proved 
exactly as  \cite[Theorem 5.7, Proposition 5.8 and Proposition 5.9]{Ri}, by using our
Propositions \ref{HOM} and \ref{IDEAL}.
\begin{prop}
\label{InjSurj}
Fix $\bm\in{\bf\Theta}(\B)$ and
let $(A,\alpha)$ and $(B,\beta)$ be two $C^*$-algebras endowed with  strongly 
continuous actions of $\B$. Then, if $T:A\to B$ is a continuous homomorphism
which intertwines the actions $\alpha$ and $\beta$,
its restriction $T^\infty:A^\infty\to B^\infty$ extends to a continuous homomorphism
$T_{\theta,\bm}:A_{\theta,\bm}\to B_{\theta,\bm}$.
 If moreover $T$ is injective (respectively surjective) then $T_{\theta,\bm}$
is injective (respectively  surjective) too.
\end{prop}

\begin{prop}
Fix $\bm\in{\bf\Theta}(\B)$ and
let $(A,\alpha)$ be a $C^*$-algebra endowed with a strongly 
continuous  and isometric action of $\B$ and let also $I$ be an $\alpha$-invariant (essential) ideal
of $A$. Then $I_{\theta,\bm}$ is an (essential) ideal of $A_{\theta,\bm}$. 
\end{prop}

\end{document}